# The reduction of $G$-ordinary crystalline representations with $G$-structure

Macarena Peche Irissarry

**Abstract**

Fontaine's $D_{\text{cris}}$ functor allow us to associate an isocrystal to any crystalline representation. For a reductive group $G$, we study the reduction of lattices inside a germ of crystalline representations with $G$-structure $V$ to lattices (which are crystals) with $G$-structure inside $D_{\text{cris}}(V)$. Using Kisin modules theory, we give a description of this reduction in terms of $G$, in the case when the representation $V$ is ($G$-)ordinary. In order to do that, first we need to generalize Fargues' construction of the Harder-Narasimhan filtration for $p$-divisible groups to Kisin modules.

## Contents

# 1 Introduction

## 1.1 The problem

Integral models of Shimura varieties have been constructed by Kisin in [23] and [24], and by Kisin and Pappas in [25]. Given an integral model, we can ask ourselves how to reduce special cycles of the Shimura variety at places of good reduction. This thesis was motivated by that question, in particular, by the study of 0-dimensional special cycles, which thus correspond to CM points of the Shimura variety. If $\gamma$ is an isogeny class of such points, then $\gamma$ may be described as a coset

$$\text{Aut}(\gamma) \backslash \mathcal{X}^p(\gamma) \times \mathcal{X}_p(\gamma).$$

The sets $\mathcal{X}^p(\gamma)$ and $\mathcal{X}_p(\gamma)$ parameter the level structures prime to $p$ and in $p$, respectively, and together they determine the position of a point in the isogeny class. Langlands made a conjecture in [30] about the mod $p$ points of a Shimura variety. Afterwards, the conjecture was made more precise by Kottwitz ([27]) and then by Langlands and Rapoport ([31]). This conjecture has been proved in the case of Shimura varieties of abelian type by Kisin in [21]. The conjecture implies that the points of the Shimura variety in characteristic $p$ may be described analogously as disjoint union of sets of the form

$$\text{Aut}(\overline{\gamma}) \backslash \mathcal{X}^p(\overline{\gamma}) \times \mathcal{X}_p(\overline{\gamma}).$$

While the morphism $\text{Aut}(\gamma) \to \text{Aut}(\overline{\gamma})$ and the bijection $\mathcal{X}^p(\gamma) \simeq \mathcal{X}^p(\overline{\gamma})$ are well known, the map

$$\text{red} \;:\; \mathcal{X}_p(\gamma) \to \mathcal{X}_p(\overline{\gamma})$$

remains quite mysterious. The aim is to describe the latter as concretely as possible. In order to do that, we can use methods from integral $p$-adic Hodge theory and Bruhat-Tits buildings. In particular, we will be working with Kisin modules and various kinds of well-known filtrations in $p$-adic Hodge theory.

Now, let $G/\mathbb{Z}_p$ be the reductive group such that $G_{\mathbb{Q}_p}$ is the group associated to the Shimura variety. For the integral models constructed by Kisin, there is presumably (see the article by Milne [34]) a Tannakian description for the reduction map. An element $x \in \mathcal{X}_p(\gamma)$ corresponds to a $\otimes$-functor

$$x \ : \ \mathrm{Rep}_{\mathbb{Z}_p} G \to \mathrm{Rep}_{\mathbb{Z}_p}^{\mathrm{cr,ab}} \mathrm{Gal}_K$$

where $\mathrm{Rep}_{\mathbb{Z}_p} G$ is the category of $\mathbb{Z}_p$-linear representations of $G$ and $\mathrm{Rep}_{\mathbb{Z}_p}^{\mathrm{cr,ab}} \mathrm{Gal}_K$ is the category of abelian crystalline representations of $\mathrm{Gal}(\overline{\mathbb{Q}}_p/K)$ for some finite extension $K$ of $K_0 = W(\mathbb{F})[\frac{1}{p}]$ and $\mathbb{F} = \overline{\mathbb{F}}_p$. Similarly, an element $y \in \mathcal{X}_p(\overline{\gamma})$ corresponds to a $\otimes$-functor

$$y \ : \ \mathrm{Rep}_{\mathbb{Z}_p} G \to \mathrm{Mod}^\sigma_{W(\mathbb{F})}$$

where $\mathrm{Mod}^\sigma_{W(\mathbb{F})}$ is the category of finite free $W(\mathbb{F})$-modules with a Frobenius isomorphism $\sigma_M : \sigma^* M[\frac{1}{p}] \to M[\frac{1}{p}]$. Here $\sigma$ is the lifting to $W(\mathbb{F})$ of the Frobenius on $\mathbb{F}$. Then, the reduction map that we will consider is given by

$$y = D_{\mathrm{cris}} \circ x$$

where $D_{\mathrm{cris}}$ is Kisin's integral functor described in [5] or [22]. For the Tate module $T_p(X)$ of a $p$-divisible group over $\mathcal{O}_K$, $D_{\mathrm{cris}}(T_p(X))$ is the Dieudonné crystal of the reduction of $X$ to the residue field $\mathbb{F}$. Bruhat-Tits buildings provide yet another convenient description of the source and target sets of the reduction map. The source set $\mathcal{X}_p(\gamma) \simeq G(\mathbb{Q}_p)/G(\mathbb{Z}_p)$ is simply the $G(\mathbb{Q}_p)$-orbit of the hyperspecial point of the extended Bruhat-Tits building $\mathbf{B}^e(G_{\mathbb{Q}_p})$ which corresponds to $G/\mathbb{Z}_p$. The target set is usually more difficult to grasp, but it may still be embedded in the much larger Bruhat-Tits building $\mathbf{B}^e(G_{K_0})$ of $G$ over $K_0$. Nevertheless, it also has a simple description under the ordinary condition on $\gamma$ (that we shall explain later): in this case the target of our reduction map is easy to describe, it is the quotient $M(\mathbb{Q}_p)/M(\mathbb{Z}_p)$ where $M$ is a Levi of a parabolic $P = U \rtimes M$ of $G$ associated to the Newton type of $y$. Moreover, $M(\mathbb{Q}_p)/M(\mathbb{Z}_p) \hookrightarrow G(\mathbb{Q}_p)/G(\mathbb{Z}_p)$ is a fundamental domain for the action of $U(\mathbb{Q}_p)$ on $G(\mathbb{Q}_p)/G(\mathbb{Z}_p)$. Thus, there is a natural retraction

$$G(\mathbb{Q}_p)/G(\mathbb{Z}_p) \twoheadrightarrow M(\mathbb{Q}_p)/M(\mathbb{Z}_p)$$

where the source and target sets respectively correspond to $\mathcal{X}_p(\gamma)$ and $\mathcal{X}_p(\overline{\gamma})$. Our main result, Theorem 5.15thm.5.15, gives us:

**Theorem 1.1.** *For $\gamma$ ordinary, the above map coincides with the reduction map* $\mathrm{red} : \mathcal{X}_p(\gamma) \to \mathcal{X}_p(\overline{\gamma})$.

The proof of the theorem uses the theory of buildings, Wintenberger's work in [43] on abelian crystalline representations and the factorization

$$D_{\mathrm{cris}} \ : \ \mathrm{Rep}_{\mathbb{Z}_p}^{\mathrm{cr,ab}} \mathrm{Gal}_K \xrightarrow{\mathfrak{M}} \mathrm{Mod}^\varphi_{\mathfrak{S},\,\mathrm{fr}} \xrightarrow{\mathrm{mod}\ u} \mathrm{Mod}^\sigma_{W(\mathbb{F})}$$

where $\mathrm{Mod}^\varphi_{\mathfrak{S},\,\mathrm{fr}}$ is the category of Kisin modules, defined by Kisin in [22] together with the functor $\mathfrak{M}$. We thus have to study $G$-Kisin modules, i.e. $\otimes$-functors

$$M \ : \ \mathrm{Rep}_{\mathbb{Z}_p} G \to \mathrm{Mod}^\varphi_{\mathfrak{S},\,\mathrm{fr}}$$

and their mod $p^n$-variants. The key point of the proof is to lift the Newton filtration of $\mathrm{Rep}_{\mathbb{Z}_p} G \to \mathrm{Mod}^\sigma_{W(\mathbb{F})}$ to a filtration of Harder-Narasimhan type on the $G$-Kisin module. Harder-Narasimhan filtrations have already been defined in different contexts such as:

- Fargues developed a Harder-Narasimhan theory for finite flat group schemes and for filtered isocrystals in [15] and [16].

- Fargues also defined the Harder-Narasimhan polygon for $p$-divisible groups in [15, Theorem 2].

- Moonen ([35]) and Shen ([41]) constructed the Hodge-Newton filtration for a $p$-divisible group with additional structures (which corresponds to the case of PEL Shimura varieties), under the ordinarity condition. This filtration lifts the Newton filtration of the isocrystal. Shen's result is more general, since it proves that whenever the Hodge and Newton filtration of the filtered isocrystal coincide in a break point, we can lift the sub-filtered isocrystal to a sub-$p$-divisible group.

## 1.2 The results

Let $\mathfrak{S} = W(\mathbb{F})[[u]]$ for $\mathbb{F}$ a perfect field of characteristic $p > 0$ and $W(\mathbb{F})$ the ring of Witt vectors over $\mathbb{F}$. Let $K$ be a finite extension of $K_0 = \operatorname{Frac} W(\mathbb{F})$, with uniformizer $\pi_K$. Recall that a Kisin module is a finite free $\mathfrak{S}$-module $M$ together with a Frobenius morphism $\varphi_M : \varphi^*M[\frac{1}{E}] \xrightarrow{\sim} M[\frac{1}{E}]$, where $E \in \mathfrak{S}$ is the minimal polynomial of $\pi_K$. These modules were defined by Kisin in [22]. In the same article, he also constructed a functor

$$\mathfrak{M} : \operatorname{Rep}_{\mathbb{Z}_p}^{\mathrm{cr}} \operatorname{Gal}_K \to \operatorname{Mod}_{\mathfrak{S},\,\mathrm{fr}}^{\varphi}$$

and such that $\mathfrak{M}(L)/u\mathfrak{M}(L)[\frac{1}{p}] \simeq D_{\mathrm{cris}}(L \otimes \mathbb{Q}_p)$ for any $L \in \operatorname{Rep}_{\mathbb{Z}_p}^{\mathrm{cr}} \operatorname{Gal}_K$. Some other categories of Kisin modules that will appear are: the categories of isogeny classes of Kisin modules that we denote by $\operatorname{Mod}_{\mathfrak{S}}^{\varphi} \otimes \mathbb{Q}_p$ and $\operatorname{Mod}_{\mathfrak{S}[\frac{1}{p}]}^{\varphi}$, the category of $p^{\infty}$-torsion Kisin modules that we denote by $\operatorname{Mod}_{\mathfrak{S},\,\mathrm{t}}^{\varphi}$ (formed by finitely generated $\mathfrak{S}$-modules killed by a power of $p$ with no $u^{\infty}$-torsion, endowed with a Frobenius) and the category of $p$-torsion Kisin modules denoted by $\operatorname{Mod}_{\mathbb{F}[[u]],\mathrm{fr}}^{\varphi}$ (formed by finite free $\mathbb{F}[[u]]$-modules endowed with a Frobenius). In [19], Genestier and Lafforgue define the category MHP of Hodge-Pink modules and gave $\otimes$-functors $\operatorname{Mod}_{\mathfrak{S},\,\mathrm{fr}}^{\varphi} \to \mathrm{MHP} \to \mathrm{MF}_K^{\sigma}$, where $\mathrm{MF}_K^{\sigma}$ is the category of filtered isocrystals. Then, they proved that there are equivalence of categories between Hodge-Pink modules verifying Griffiths transversality condition and filtered isocrystals, and which sends weakly admissible Hodge-Pink modules (that they define) to weakly admissible filtered isocrystals.

All the following categories admit Harder-Narasimhan filtrations:

$$\operatorname{Rep}_{\mathbb{Q}_p}^{\mathrm{cr}} \operatorname{Gal}_K \xleftarrow{\simeq}{\mathrm{wa}} \mathrm{MF}_K^{\sigma} \xleftarrow{\simeq} \mathrm{MHP}^{\mathrm{wa,Gr}} \hookrightarrow \operatorname{Mod}_{\mathfrak{S}}^{\varphi} \otimes \mathbb{Q}_p \hookrightarrow \operatorname{Mod}_{\mathfrak{S}[\frac{1}{p}]}^{\varphi},$$

where $\mathrm{MHP}^{\mathrm{wa,Gr}}$ denotes the category of weakly admissible Hodge-Pink modules verifying the Griffiths transversality condition. In particular, the Harder-Narasimhan filtration on $\operatorname{Rep}_{\mathbb{Q}_p}^{\mathrm{cr}} \operatorname{Gal}_K$ is given by Fargues in [15], where he also showed that this filtration is compatible with tensor products and with the change of extension $K$. In Proposition 3.11thm.3.11, we show that all the filtrations defined in the categories above are compatible with the functors between the categories, thus the filtrations are all compatible with tensor products, by Fargues' result ([15, Corollaire 6]).

Motivated by the study of special points of Shimura varieties, we want to be able to change the extension $K$ over $K_0$. For that reason, we define the germ of crystalline representations (Definition 4.8defi.4.8) that we can think of as pairs $(V, \rho)$ such that $(V, \rho) \in \operatorname{Rep}_{\mathbb{Q}_p} \operatorname{Gal}_K$ for a sufficiently large finite extension $K$ of $K_0$. We denote by $\operatorname{Rep}_{\mathbb{Q}_p}\{\operatorname{Gal}_{K_0}\}$ the category of germs of crystalline representations. In Proposition 3.8.1 we prove that the Fargues filtration on crystalline representations extends to a filtration on germs of crystalline representation. In Lemma 3.8.2, following Liu's results in [32], we have that the functor $\operatorname{Rep}_{\mathbb{Z}_p} \operatorname{Gal}_K \to \operatorname{Mod}_{W(\mathbb{F})}^{\sigma}$, constructed by Kisin, induces a $\otimes$-functor

$$\operatorname{Rep}_{\mathbb{Z}_p}\{\operatorname{Gal}_{K_0}\} \to \operatorname{Mod}_{W(\mathbb{F})}^{\sigma}.$$

In chapter 5, we give a degree and rank function on objects in $\operatorname{Mod}_{\mathbb{F}[[u]],\mathrm{fr}}^{\varphi}$ by

$$\operatorname{rank}(M, \varphi_M) = \operatorname{rank}_{\mathbb{F}[[u]]} M \quad \text{and} \quad \deg(M, \varphi_M) = \operatorname{Pos}(M, \varphi_M \varphi^*M).$$

These functions define a slope function in the sense of André's formalism and thus Harder-Narasimhan filtration $\mathcal{F}_{\mathrm{F},1}(M)$ and its polygon $\mathbf{t}_{\mathrm{F},1}(M)$ on each object $M \in \operatorname{Mod}_{\mathbb{F}[[u]],\mathrm{fr}}^{\varphi}$. This is a generalization of the results by Fargues in [16], who gave this construction for finite flat group schemes (which are related to $p^{\infty}$-torsion Kisin modules in the same way $p$-divisible groups corresponds to Kisin modules of height 1) and, for this reason, we call Fargues filtration the Harder-Narasimhan filtrations defined on categories of Kisin modules. Levin and Wang Erickson prove in [14] that this filtration is compatible with tensor products.

Afterwards, we work with the category of $p^{\infty}$-torsion Kisin modules, which is quasi abelian and admits a degree and rank function. We follow the strategy given by Fargues in [15] and Shen in [41] to construct a Harder-Narasimhan filtration on $p$-divisible groups under some hypothesis. For a Kisin module $M = (M, \varphi_M) \in \operatorname{Mod}_{\mathfrak{S},\,\mathrm{fr}}^{\varphi}$, we use the polygons defined on the $p^{\infty}$-torsion Kisin modules $M/p^nM$ to define a polygon $\mathbf{t}_{\mathrm{F},\infty}(M)$ on $M$. Even though the category $\operatorname{Mod}_{\mathfrak{S},\,\mathrm{fr}}^{\varphi}$ is not quasi-abelian, we have good degree and rank functions. We define HN-type Kisin modules as the ones for which there exists a filtration by Kisin modules and whose graded pieces are semi-stable Kisin modules. In the rest of the chapter, we prove two important results:

1. Proposition 3.23thm.3.23: A Kisin module $M$ is HN-type if and only if $\mathbf{t}_\infty(M) = \mathbf{t}_{F,1}(M)$,

2. Theorem 3.26thm.3.26: Every Kisin module is isogenous to a HN-type Kisin module.

Both proofs are made in the same spirit as the proofs given by Fargues and Shen for $p$-divisible groups.

In chapter 6, we will be working with objects with $G$-structure, for $G$ a reductive group over $\mathbb{Z}_p$ (related to the reductive group in the Shimura datum). We start then by presenting the filtrations, graduations and space of types for fiber functors

$$\omega_R \;:\; \operatorname{Rep}_\mathcal{O} G \to \operatorname{Bun}_R$$

where $\mathcal{O} \in \{\mathbb{Q}_p, \mathbb{Z}_p, \mathbb{F}_p\}$ and $R$ is a ring. Then, we use the result given by Broshi in [6] to prove the following result, which will be a key tool for us:

**Proposition 1.2.** *Let $R$ be a local strictly henselian and faithfully flat $\mathcal{O}$-algebra. Then, any exact and faithful $\otimes$-functor as above is $\otimes$-isomorphic to the trivial fiber function $\omega_{G,R}$ (the forgetful functor).*

Even if we start working with a trivial (or isomorphic to trivial) crystalline representation with $G$-structure, à priori, the functors that we will apply to it (Fontaine's $D_{\mathrm{cris}}$ functor, Kisin's $\mathfrak{M}$ functor, etc.) will change the functor to a non-trivial one. However, the proposition above will always allow us to reduce to the (isomorphic to the) trivial case. This is very important as the results we will be using from different authors are only stated for the trivial fiber functor. Another important result given in this section is a generalization of Haboush theorem, using Seshadri's results in [40]:

**Proposition 1.3.** *Let $L$ be a field which is an $\mathcal{O}$-algebra. Let $\Gamma \subset \mathbb{R}$ be a subring. Suppose that a fiber functor $\omega_L : \operatorname{Rep}_\mathcal{O} G \to \operatorname{Vect}_L$ admits a factorization through an additive $\otimes$-functor*

$$\mathcal{F} \;:\; \operatorname{Rep}_\mathcal{O} G \to \operatorname{Fil}_L^\Gamma$$

*which is compatible with exterior and symmetric powers. Then $\mathcal{F}$ is exact, thus a filtration on $\omega_L$.*

We use this proposition to prove that the Fargues filtrations and Hodge filtrations defined in previous sections can be generalized to filtrations on the object with $G$-structure.
We define the germs of crystalline representations as faithful $\otimes$-functors

$$V : \operatorname{Rep}_{\mathbb{Q}_p} G \to \operatorname{Rep}_{\mathbb{Q}_p}\{\operatorname{Gal}_{K_0}\}.$$

Using Fontaine's functor, we can associate to $V$ a filtered isocrystal with $G$-structure, that we denote by $D$. Let $\mathbf{C}^\mathbb{Q}(G)$ be a closed Weyl chamber for the split group $G$ over $K_0$, equipped with the dominance order (this is the set where the types of filtrations on objects with $G$-structure live). Then, we define the ordinary germs of crystalline representation with $G$-structure as those verifying $\mathbf{t}_N(D) = \mathbf{t}_H(D)^\#$ in $\mathbf{C}^\mathbb{Q}(G)$, where $\mathbf{t}_N(D)$ and $\mathbf{t}_H(D)$ are the Newton and Hodge types associated to the Newton graduation and Hodge filtration of $D$, and $\#$ is the average of the Galois orbits of a type. In this case, we can already say something about the reduction map: in Corollary 7.3.24, we prove that for an ordinary $V$, there is a factorization of the reduction map $\operatorname{red} : \mathcal{L}(V) \to \mathcal{L}(D)$ by

$$\mathcal{L}(V) \to \mathcal{L}(D, \mathbf{t}_H^t(D)) \hookrightarrow \mathcal{L}(D)$$

where $\mathcal{L}(V)$ and $\mathcal{L}(D)$ are the notations that we will use in the text for the sets corresponding to $\mathcal{X}_p(\gamma)$ and $\mathcal{X}_p(\overline{\gamma})$ in the case of points of Shimura varieties, and $\mathcal{L}(D, \mathbf{t}_H^t(D))$ is the subset of lattices $y$ in $D$ such that $\operatorname{Pos}(y, \sigma_D y) = \mathbf{t}_H(D)$. Note that this factorization is always true (without the ordinary condition) in the case coming from Shimura varieties, as the type $\mathbf{t}_H^t(D)$ is minuscule in that situation. We also prove the following theorem (Theorem 7.3.25):

**Theorem 1.4.** *Suppose $V$ is ordinary, let $x : \operatorname{Rep}_{\mathbb{Z}_p} \to \operatorname{Rep}_{\mathbb{Z}_p}^{\mathrm{cr}} \operatorname{Gal}_K$ for some $K$ such that $V$ is defined over it. Let $M = \mathfrak{M} \circ x$. Then, $M$ is HN-type. Therefore its Fargues filtration $\mathcal{F}_F(M)$ exists.*

In the end of this chapter we define two operators $\Phi_{\acute{e}t}^s$ and $\Phi_{\mathrm{cris}}^s$ on, respectively, lattices inside $V$ and lattices inside $D$ and prove in Proposition 7.3.27 that they commute with the reduction map red in the ordinary case. Then, we obtain the following result (Theorem 4.32thm.4.32):

**Theorem 1.5.** *For $V$ an ordinary germ of crystalline representations with $G$-structure, the map* red *admits a factorization*

$$\mathcal{L}(V) \xrightarrow{\text{red}} \mathcal{L}(D, \mathbf{t}_{\mathrm{H}}^{\iota}(D)).$$
$$\pi \searrow \nearrow$$
$$U(\mathbb{Q}_p) \backslash \mathcal{L}(V)$$

*where $U$ is the unipotent radical of the parabolic subgroup of $G$ stabilizing the Fargues filtration of $V$.*

In the last chapter, we study the particular case when the germ of crystalline representation is abelian, which is the case that interested us originally (since it is the one coming from CM points in Shimura varieties). In [39, 2], Fontaine proves, building on Serre's results in [38], that there is an equivalence of categories

$$V_K^u \; : \; \mathrm{Rep}_{\mathbb{Q}_p} T_K \xrightarrow{\sim} \mathrm{Rep}_{\mathbb{Q}_p}^{\mathrm{cr,ab}} \mathrm{Gal}_K,$$

where $T_K = \varprojlim_{E \subset K} \mathrm{Res}_{E/\mathbb{Q}_p}(\mathbb{G}_{m,E})$. Then, in [43], Wintenberger constructs a $\otimes$-functor

$$D_{\pi_K} \; : \; \mathrm{Rep}_{\mathbb{Q}} T_K \to \mathrm{Mod}_{W(\mathbb{F})}^{\sigma}$$

and proves that when we restrict ourselves to a finitely generated tensor subcategory $\mathcal{V}$ of $\mathrm{Rep}_{\mathbb{Q}_p} T_K$, we have a isomorphism of $\otimes$-functors

$$D_{\pi_K | \mathcal{V}} \simeq (D_{\mathrm{cris}} \circ V_K^u)_{|\mathcal{V}}.$$

We generalize both constructions to germs of crystalline representation with $G$-structure, in Proposition 5.1thm.5.1 and section 5.3.3Wintenberger's functorsubsubsection.5.3.3, and then give an explicit description, as cocharacters, of the Hodge and Newton types in Proposition 5.4thm.5.4. In Proposition 5.6thm.5.6, we prove that in this case, the Fargues filtration defined on a germ of crystalline representation with $G$-structure coincides with the opposed Newton filtration of the $G$-isocrystal $D$ associated to it.

For the rest of the chapter we restrict ourselves to the ordinary case. Using the explicit description of $G$-isocrystals given by Kottwitz in [26], [28] and [29], and by Rapoport and Richartz in [37], we are able to say a little bit more about the reduction map. Putting together the results in section 8.5, we obtain the main theorem (Theorem 5.15thm.5.15):

**Theorem 1.6.** *The reduction map*

$$\mathrm{red} \; : \; \mathcal{L}(V) \to \mathcal{L}(D, \mathbf{t}_{\mathrm{H}}^{\iota}(D))$$

*factors through an $M(\mathbb{Q}_p)$-equivariant bijection*

$$U_{\mathcal{F}_{\mathrm{F}}}(\mathbb{Q}_p) \backslash \mathcal{L}(V) \simeq \mathcal{L}(D, \mathbf{t}_{\mathrm{H}}^{\iota}(D)).$$

#### Acknowledgments

I would like to thank Christophe Cornut for introducing me to this subject, all the helpful conversations and its work about filtration and buildings. I also would like to thank Laurent Fargues for the construction of the Harder-Narasimhan filtration for $p$-divisible groups, as its generalization for Kisin modules is a key point of the proof of the main theorem.

## 2 Preliminaries

Most of the structures described in this section can be found in [2] and [8].

### 2.1 Space of types and the dominance order

Let $\Gamma$ be a nonzero subring of $\mathbb{R}$ and $r \in \mathbb{N}$. We define the space of $\Gamma$-types of length $r$ in three different but equivalent ways :

1. The cone $\Gamma_{\geq}^r = \{(\gamma_1, \ldots, \gamma_r) \in \Gamma^r \mid \gamma_1 \geq \ldots \geq \gamma_r\}$.

2. Consider first the group-ring $\mathbb{Z}[\Gamma] = \left\{ \mathbf{t} = \sum_{\gamma \in \Gamma} n_\gamma \cdot e^\gamma \mid n_\gamma \in \mathbb{Z} \text{ and } \mathbf{t} \text{ has finite support} \right\}$. Inside the group-ring, we have the space of positive elements

$$\mathbb{N}[\Gamma] = \left\{ \mathbf{t} = \sum_{\gamma \in \Gamma} n_\gamma \cdot e^\gamma \mid n_\gamma \in \mathbb{N} \text{ and } \mathbf{t} \text{ has finite support} \right\}$$

and finally, the space of types of length $r$ is given by

$$\mathbb{N}[\Gamma]^r = \left\{ \mathbf{t} = \sum_{\gamma \in \Gamma} n_\gamma \cdot e^\gamma \mid \sum_{\gamma \in \Gamma} n_\gamma = r \text{ and } n_\gamma \in \mathbb{N} \right\},$$

i.e. it is the subspace of positive elements of degree $r$.

3. Concave polygons : Continuous functions $\mathbf{t} : [0, r] \to \mathbb{R}$ with $f(0) = 0$, which are affine on $[i-1, i]$ for every $1 \leq i \leq r$ with slopes $\gamma_i \in \Gamma$ verifying $\gamma_1 \geq \ldots \geq \gamma_r$.

We can also define the entire space of $\Gamma$-types as the union of the spaces of types of length $r$, for all $r \in \mathbb{N}$. In the second case, the entire space of $\Gamma$-types corresponds to $\mathbb{N}[\Gamma]$. We will drop the notation $\Gamma$ when there is no confusion about the indexation of the type. We recall some of the operations that we can consider in the space of types.

- **Degree:** The degree of a type is defined by

  1. The map
  $$\begin{array}{rccc} \deg & : & \Gamma^r_\geq & \to & \Gamma \\ & & (\gamma_1, \ldots, \gamma_r) & \mapsto & \sum_{i=1}^r \gamma_i \end{array}$$

  2. The map
  $$\begin{array}{rccc} \deg & : & \mathbb{N}[\Gamma] & \to & \Gamma \\ & & \sum_{\gamma \in \Gamma} n_\gamma \cdot e^\gamma & \mapsto & \sum_{\gamma \in \Gamma} n_\gamma \cdot \gamma \end{array}$$

  3. For a concave polygon, its degree corresponds to the coordinate of its ending point.

  Degree functions will appear everywhere in this thesis so, in this section, we will show how the degree of a type behave with respect to other operators that will be defined.

- **Dominance order:** We can define a partial order relation called the dominance order over the space of types of length $r$:

  (1) In the cone, the order is given by $\gamma \leq \delta$ if and only if $\sum_{i=1}^k \gamma_i \leq \sum_{i=1}^k \delta_i$ for all $1 \leq k \leq r$ with equality when $k = r$, where $\gamma = (\gamma_1, \ldots, \gamma_r)$ and $\delta = (\delta_1, \ldots, \delta_r)$ are two elements of $\Gamma^r_\geq$.

  (3) Viewed as concave polygons this partial order corresponds to saying that $\delta$ is above $\gamma$ and they both have the same extremities.

  If $\mathbf{t}_1 \leq \mathbf{t}_2$, then $\deg \mathbf{t}_1 = \deg \mathbf{t}_2$.

- **Addition:** The addition of two types of length $r$ is given by:

  (1) The map
  $$\begin{array}{rccc} + & \Gamma^r_\geq \times \Gamma^r_\geq & \to & \Gamma^r_\geq \\ & (\gamma, \delta) & \mapsto & \gamma + \delta = (\gamma_1 + \delta_1, \ldots, \gamma_r + \delta_r) \end{array}$$
  where $\gamma = (\gamma_1, \ldots, \gamma_r)$ and $\delta = (\delta_1, \ldots, \delta_r)$.

  (3) As concave polygons, the addition corresponds to the usual addition of functions.

  From the formula above, it is easy to see that $\deg(\gamma + \delta) = \deg \gamma + \deg \delta$.

- **Norm:** We define the norm of a type as

  $$\begin{array}{rccc} || \cdot || & : & \Gamma^r_\geq & \to & \mathbb{R}_+ \\ & & (\gamma_1, \ldots, \gamma_r) & \mapsto & \sqrt{\gamma_1^2 + \ldots + \gamma_r^2} \end{array}$$

  and it verifies two properties.

  1. We have $||\mathbf{t}_1 + \mathbf{t}_2|| \leq ||\mathbf{t}_1|| + ||\mathbf{t}_2||$ for all $\mathbf{t}_1, \mathbf{t}_2 \in \Gamma^r_\geq$, since $|| \cdot ||$ is the euclidean distance in $\mathbb{R}^r_+ \subset \mathbb{R}^r$.

2. For $\mathbf{t}_1, \mathbf{t}_2 \in \Gamma_{\geq}^r$, the inequality for the dominance order $\mathbf{t}_1 \leq \mathbf{t}_2$ implies an inequality for the respective norms $||\mathbf{t}_1|| \leq ||\mathbf{t}_2||$, moreover $||\mathbf{t}_1|| = ||\mathbf{t}_2||$ implies that $\mathbf{t}_1 = \mathbf{t}_2$. To prove that, let $(\gamma_1, \ldots, \gamma_r) \leq (\delta_1, \ldots, \delta_r)$, thus $\Gamma_i = \sum_{j \leq i} \gamma_j \leq \Delta_i = \sum_{j \leq i} \delta_j$, and the equality holds for $i = r$. Then

$$\begin{aligned}
\sum_{i=1}^r \delta_i^2 - \sum_{i=1}^r \gamma_i^2 &= \sum_{i=1}^r (\delta_i - \gamma_i)(\delta_i + \gamma_i) \\
&= \sum_{i=1}^r (\Delta_i - \Gamma_i)(\delta_i + \gamma_i) - \sum_{i=1}^{r-1} (\Delta_i - \Gamma_i)(\delta_{i+1} + \gamma_{i+1}) \\
&= \sum_{i=1}^{r-1} (\Delta_i - \Gamma_i)(\delta_i - \delta_{i+1} + \gamma_i - \gamma_{i+1}) \\
&\geq 0,
\end{aligned}$$

and the equality holds if and only if $(\Delta_i - \Gamma_i)(\delta_i - \delta_{i+1} + \gamma_i - \gamma_{i+1}) = 0$ for $1 \leq i \leq r-1$, which means that $(\gamma_1, \ldots, \gamma_r) = (\delta_1, \ldots, \delta_r)$.

- **Multiplication by scalars:** The multiplication by scalars is given by

  1. For a type $(\gamma_1, \ldots, \gamma_r)$ and $c > 0$ an element in $\Gamma$, then $c \cdot (\gamma_1, \ldots, \gamma_r) = (c\gamma_1, \ldots, c\gamma_r)$.
  2. For $\sum_{\gamma \in \Gamma} n_\gamma \cdot e^\gamma \in \mathbb{N}[\Gamma]$ and $c > 0$ an element in $\Gamma$, then $c \cdot \sum_{\gamma \in \Gamma} n_\gamma \cdot e^\gamma = \sum_{\gamma \in \Gamma} n_\gamma \cdot e^{c\gamma}$.

  For the degree, we have $\deg(c\gamma) = c \deg \gamma$.

- **Involution:** We define the involution of a type as:

  (1) For a type $\gamma = (\gamma_1, \ldots, \gamma_r) \in \Gamma_{\geq}^r$, it is the map

  $$\begin{array}{rccc}
  \cdot^\iota : & \Gamma_{\geq}^r & \to & \Gamma_{\geq}^r \\
  & \gamma & \mapsto & \gamma^\iota = (-\gamma_r, \ldots, -\gamma_1)
  \end{array}.$$

  (2) For an element of $\mathbb{N}[\Gamma]$, it is the map

  $$\begin{array}{rccc}
  \cdot^\iota : & \mathbb{N}[\Gamma] & \to & \mathbb{N}[\Gamma] \\
  & \sum_{\gamma \in \Gamma} n_\gamma \cdot e^\gamma & \mapsto & \sum_{\gamma \in \Gamma} n_\gamma \cdot e^{-\gamma}
  \end{array}$$

  The degree verifies $\deg(\gamma^\iota) = -\deg \gamma$.

- **Concatenation:** We define the concatenation of two types as

  (1) In the cone, the map

  $$\begin{array}{rccc}
  * : & \Gamma_{\geq}^r \times \Gamma_{\geq}^s & \to & \Gamma_{\geq}^{r+s} \\
  & (\gamma, \delta) & \mapsto & \gamma * \delta = (\epsilon_{\sigma(1)}, \ldots, \epsilon_{\sigma(r+s)})
  \end{array}$$

  where if $\gamma = (\gamma_1, \ldots, \gamma_r)$, $\delta = (\delta_1, \ldots, \delta_s)$, then $(\epsilon_1, \ldots, \epsilon_{r+s}) = (\gamma_1, \ldots, \gamma_r, \delta_1, \ldots, \delta_s)$ and $\sigma \in \mathfrak{S}_{r+s}$ is a permutation such that $\epsilon_{\sigma(1)} \geq \ldots \geq \epsilon_{\sigma(r+s)}$. From this description, we see that

  $$c(\gamma * \delta) = c\gamma * c\delta$$

  for every $\gamma, \delta$ in the space of types and every $c > 0$ in $\Gamma$.

  (2) Inside $\mathbb{N}[\Gamma]$, it is the usual addition.

  (3) As polygons, the concatenation is the smallest concave polygon above the collection of points $\{(i+j, \gamma(i) + \delta(j))\}_{\substack{0 \leq i \leq r \\ 0 \leq j \leq s \\ i, j \in \mathbb{N}}}$ or, equivalently, the polygon defined as the function

  $$\begin{array}{rccc}
  \gamma * \delta : & [0, r+s] & \to & \mathbb{R} \\
  & x & \mapsto & \sup_{\substack{a+b=x \\ 0 \leq a \leq r \\ 0 \leq b \leq s}} \gamma(a) + \delta(b)
  \end{array}$$

  From the first formula, we get $\deg(\gamma * \delta) = \deg(\gamma) + \deg(\delta)$.

- **Tensor product:** The tensor product of two types is given by

(1) In the cone, the map

$$\begin{array}{rccc}\otimes \; : & \Gamma_{\geq}^r \times \Gamma_{\geq}^s & \to & \Gamma_{\geq}^{rs} \\ & (\gamma, \delta) & \mapsto & \gamma \otimes \delta = (\epsilon_{\sigma(1)}, \ldots, \epsilon_{\sigma(rs)})\end{array}$$

where if $\gamma = (\gamma_1, \ldots, \gamma_r)$, $\delta = (\delta_1, \ldots, \delta_s)$, then

$$(\epsilon_1, \ldots, \epsilon_{rs}) = (\gamma_1 + \delta_1, \ldots, \gamma_1 + \delta_s, \ldots, \gamma_r + \delta_1, \ldots, \gamma_r + \delta_s)$$

and $\sigma \in \mathfrak{S}_{rs}$ is a permutation such that $\epsilon_{\sigma(1)} \geq \ldots \geq \epsilon_{\sigma(rs)}$. From this description, we see that $c(\gamma \otimes \delta) = c\gamma \otimes c\delta$, for every $\gamma, \delta$ in the space of types and every $c > 0$ in $\Gamma$.

(2) Inside $\mathbb{N}[\Gamma]$, it is given by the usual multiplication on the group-ring $\mathbb{Z}[\Gamma]$ ($\mathbb{N}[\Gamma]$ is stable under this multiplication).

Again, from the first formula, we have $\deg(\gamma \otimes \delta) = s \deg \gamma + r \deg \delta$, for $\gamma$ and $\delta$ of length $r$ and $s$, respectively.

- **Exterior powers:** For $k \geq 1$, we define the $k$-th exterior power of a type by

$$\begin{array}{rccc}\Lambda^k \; : & \Gamma_{\geq}^r & \to & \Gamma_{\geq}^s \\ & (\gamma_1, \ldots, \gamma_r) & \mapsto & (\underline{\gamma}_{\sigma(1)}, \ldots, \underline{\gamma}_{\sigma(s)})\end{array}$$

where $s = \binom{r}{k}$, $\{\underline{\gamma}_1, \ldots, \underline{\gamma}_s\} = \{\gamma_{i_1} + \ldots + \gamma_{i_k} \mid 1 \leq i_1 < \ldots < i_k \leq r\}$ and $\sigma \in \mathfrak{S}_s$ is a permutation such that $\underline{\gamma}_{\sigma(1)} \geq \ldots \geq \underline{\gamma}_{\sigma(s)}$. If we take $k$ bigger than the length of the type, we will get a zero as its $k$-th exterior power. A straightforward calculation gives us $\deg(\Lambda^k \gamma) = \binom{r-1}{k-1} \deg \gamma$. In particular, for $k = r$, we have $\deg(\Lambda^k \gamma) = \deg \gamma$.

- **Symmetric powers:** For $k \geq 1$, we define the $k$-th symmetric power of a type by

$$\begin{array}{rccc}\mathrm{Sym}^k \; : & \Gamma_{\geq}^r & \to & \Gamma_{\geq}^s \\ & (\gamma_1, \ldots, \gamma_r) & \mapsto & (\underline{\gamma}_{\sigma(1)}, \ldots, \underline{\gamma}_{\sigma(s)})\end{array}$$

where $s = \binom{r+k-1}{k}$, $\{\underline{\gamma}_1, \ldots, \underline{\gamma}_s\} = \{\gamma_{i_1} + \ldots + \gamma_{i_k} \mid 1 \leq i_1 \leq \ldots \leq i_k \leq r\}$ and $\sigma \in \mathfrak{S}_s$ is a permutation such that $\underline{\gamma}_{\sigma(1)} \geq \ldots \geq \underline{\gamma}_{\sigma(s)}$. A straightforward calculation gives us $\deg(\mathrm{Sym}^k \gamma) = r\binom{r+k-1}{k-1} \deg \gamma$.

## 2.2 Filtrations

Let $\Gamma$ a non zero subgroup of $\mathbb{R}$, let $K$ be a field and $V$ a finite dimensional $K$-vector space. A $\Gamma$-filtration $\mathcal{F}$ on $V$ is a collection of $K$-subspaces $(\mathcal{F}^{\geq \gamma} V)_{\gamma \in \Gamma}$ of $V$ which is decreasing, separated, exhaustive and left continuous. For a filtration $\mathcal{F}$ of $V$, we define

$$\mathcal{F}^{>\gamma} V := \bigcup_{\gamma' > \gamma} \mathcal{F}^{\geq \gamma'} V \quad \text{and} \quad \mathrm{Gr}_\mathcal{F}^\gamma V = \mathcal{F}^{\geq \gamma} V / \mathcal{F}^{>\gamma} V.$$

We thus have a short exact sequence of $K$-vector spaces $0 \to \mathcal{F}^{>\gamma} V \to \mathcal{F}^{\geq \gamma} V \to \mathrm{Gr}_\mathcal{F}^\gamma V \to 0$. This yields a $\Gamma$-graded $K$-vector space associated to $\mathcal{F}$,

$$\mathrm{Gr}_\mathcal{F} V := \bigoplus_{\gamma \in \Gamma} \mathrm{Gr}_\mathcal{F}^\gamma V = \mathrm{Gr}_\mathcal{F}^{\gamma_1} V \oplus \ldots \oplus \mathrm{Gr}_\mathcal{F}^{\gamma_s} V$$

where $\{\gamma_1, \ldots, \gamma_s\} = \{\gamma \mid \mathrm{Gr}_\mathcal{F}^\gamma V \neq 0\}$ are called the breaks of $\mathcal{F}$. The type of $\mathcal{F}$ (or, equivalently, the type of $\mathrm{Gr}_\mathcal{F}$) is defined by $\mathbf{t}(\mathcal{F}) = \mathbf{t}(\mathrm{Gr}_\mathcal{F}) = (\gamma_1, \ldots, \gamma_s)$ where each $\gamma_i$ for $1 \leq i \leq s$ is written as many times as the dimension of $\mathrm{Gr}_\mathcal{F}^{\gamma_i} V$. This is an element in the space of types of length $r = \dim_K V$ defined in last subsection and we have

$$\mathbf{t}(\mathcal{F}) = \sum_{\gamma \in \Gamma} \dim_K(\mathrm{Gr}_\mathcal{F}^\gamma V) \cdot e^\gamma \quad \text{in} \quad \mathbb{N}[\Gamma].$$

We define the degree of the filtration as the degree of $\mathbf{t}(\mathcal{F})$. Thus

$$\deg(\mathcal{F}) = \sum_{\gamma \in \Gamma} \dim_K(\mathrm{Gr}_\mathcal{F}^\gamma V) \cdot \gamma.$$

**Example 2.1.** For $\gamma \in \Gamma$, Denote by $V(\gamma)$ the filtration of $V$ such that $\mathrm{Gr}^{\gamma}_{V(\gamma)} = V$ and $\mathrm{Gr}^{\gamma'}_{V(\gamma)} = 0$ if $\gamma' \neq \gamma$. For this filtration, we have

$$\begin{aligned} \mathbf{t}(V(\gamma)) &= \dim_K V \cdot e^{\gamma} \\ \deg(V(\gamma)) &= \dim_K V \cdot \gamma. \end{aligned}$$

**Proposition 2.1.** *Some important properties of the filtrations are the following:*

1. *Let $0 \to W \to V \to V/W \to 0$ be an exact sequence of $K$-vector spaces and let $\mathcal{F}_V$ be a filtration on $V$. We also consider the filtrations on $W$ and $V/W$ induced by $\mathcal{F}_V$ and denote them by $\mathcal{F}_W$ and $\mathcal{F}_{V/W}$ respectively. Then, the types verify*

   $$\mathbf{t}(\mathcal{F}_V) \leq \mathbf{t}(\mathcal{F}_W) * \mathbf{t}(\mathcal{F}_{V/W}).$$

2. *Given a filtration $\mathcal{F}_i$ on a $K$-vector space $V_i$, for $1 \leq i \leq n$, we can define a filtration $\mathcal{F} = \oplus_{i=1}^{n} \mathcal{F}_i$ on $V = \oplus_{i=1}^{n} V_i$, by $\mathcal{F}^{\gamma}(V) = \oplus_{i=1}^{n} \mathcal{F}_i^{\gamma}(V_i)$, and it verifies*

   $$\mathbf{t}(\mathcal{F}) = \mathbf{t}(\mathcal{F}_1) * \ldots * \mathbf{t}(\mathcal{F}_n).$$

3. *Let $\mathcal{F}_1$ and $\mathcal{F}_2$ be two filtrations on a $K$-vector space $V$. Then we have*

   $$\mathbf{t}(\mathcal{F}_2) \leq \mathbf{t}(\mathrm{Gr}_{\mathcal{F}_1}(\mathcal{F}_2))$$

   *where $\mathrm{Gr}_{\mathcal{F}_1}(\mathcal{F}_2)$ is the filtration that $\mathcal{F}_2$ induces on $\mathrm{Gr}_{\mathcal{F}_1} V$. This property is a consequence of the last two properties.*

4. *Given a filtration $\mathcal{F}_i$ on a $K$-vector space $V_i$, for $i = 1, 2$, we can define a filtration $\mathcal{F} = \mathcal{F}_1 \otimes \mathcal{F}_2$ on $V = V_1 \otimes V_2$, by $\mathcal{F}^{\gamma}(V) = \sum_{\substack{\gamma_1 + \gamma_2 = \gamma \\ \gamma_1, \gamma_2 \in \Gamma}} \mathcal{F}_1^{\gamma_1}(V_1) \otimes \mathcal{F}_2^{\gamma_2}(V_2)$. Thus*

   $$\mathrm{Gr}^{\gamma}_{\mathcal{F}}(V) = \bigoplus_{\substack{\gamma_1 + \gamma_2 = \gamma \\ \gamma_1, \gamma_2 \in \Gamma}} \mathrm{Gr}^{\gamma_1}_{\mathcal{F}_1}(V_1) \otimes \mathrm{Gr}^{\gamma_2}_{\mathcal{F}_2}(V_2)$$

   *and*

   $$\mathbf{t}(\mathcal{F}) = \mathbf{t}(\mathcal{F}_1) \otimes \mathbf{t}(\mathcal{F}_2).$$

5. *Given a filtration $\mathcal{F}$ on a $K$-vector space $V$, we can define a filtration $\mathrm{Sym}^k \mathcal{F}$ on $\mathrm{Sym}^k V$ for every $k \geq 1$ by taking the image of $\mathcal{F}$ by $V^{\otimes} \twoheadrightarrow \mathrm{Sym}^k V$, and we have*

   $$\mathbf{t}(\mathrm{Sym}^k \mathcal{F}) = \mathrm{Sym}^k(\mathbf{t}(\mathcal{F})).$$

6. *Given a filtration $\mathcal{F}$ on a $K$-vector space $V$, we can define a filtration $\Lambda^k \mathcal{F}$ on $\Lambda^k V$ for every $k \geq 1$ by taking the image of $\mathcal{F}$ by $V^{\otimes} \twoheadrightarrow \Lambda^k V$, and we have*

   $$\mathbf{t}(\Lambda^k \mathcal{F}) = \Lambda^k(\mathbf{t}(\mathcal{F})).$$

## 2.3 Lattices

Suppose $\mathcal{O}_K$ is a discrete valuation ring with uniformizer $u$, fraction field $K$ and residue field $k$. Let $M_1$ and $M_2$ be two $\mathcal{O}_K$-lattices inside the same finite dimensional $K$-vector space $V$.

Let $\overline{M}_1 = M_1/uM_1$. We can define a $\mathbb{Z}$-filtration $\mathcal{F}(M_1, M_2)$ on the $k$-vector space $\overline{M}_1$ by

$$\mathcal{F}^i(M_1, M_2) = \frac{u^i M_2 \cap M_1 + uM_1}{uM_1}$$

for every $i \in \mathbb{Z}$. The filtration $\mathcal{F}(M_1, M_2)$ allows us to define two operators over the lattices, by

$$\mathrm{Pos}(M_1, M_2) = \mathbf{t}(\mathcal{F}(M_1, M_2)) \quad \text{and} \quad \nu(M_1, M_2) = \deg(\mathcal{F}(M_1, M_2)).$$

We call $\text{Pos}(M_1, M_2)$ the relative position of $M_1$ and $M_2$ and we are going to explain the reason why. There exists a basis $\{e_1, \ldots, e_r\}$ of $V$ adapted to $M_1$ and $M_2$, meaning that

$$M_1 = \bigoplus_{i=1}^r \mathcal{O}_K \cdot e_i, \quad M_2 = \bigoplus_{i=1}^r \mathcal{O}_K \cdot u^{-a_i} e_i$$

for some $a_i \in \mathbb{Z}$ verifying $a_1 > \ldots > a_r$. Then $(a_1, \ldots, a_n)$ does not depend upon the chosen basis. In fact, a direct calculation shows that $\text{Pos}(M_1, M_2) = (a_1, \ldots, a_r)$ and $\nu(M_1, M_2) = \sum_{i=1}^r a_i$.

*Remark* 1. If $M_1 \subset M_2$, then $\text{Pos}(M_1, M_2)$ corresponds to the invariant factors (given by the structure theorem for finitely generated torsion modules over a principal ideal domain) of the quotient $Q = M_2/M_1$ and $\nu(M_1, M_2)$ corresponds to the length of $Q$.

The relative position of two lattices is an element of the space of $\mathbb{Z}$-types. For two $\mathcal{O}_K$-lattices $M_1, M_2$ inside $V$, we can also define the operator

$$\text{d}(M_1, M_2) = ||\text{Pos}(M_1, M_2)||$$

where $||\cdot||$ is the norm of a type.

**Lemma 2.2.** *Let $M_1, M_2, M_3$ be $\mathcal{O}_K$-lattices in a $K$-vector space. Then:*

1. *The relative position verifies the triangular inequality $\text{Pos}(M_1, M_3) \leq \text{Pos}(M_1, M_2) + \text{Pos}(M_2, M_3)$.*

2. *The operator* d *verifies the triangular inequality $\text{d}(M_1, M_3) \leq \text{d}(M_1, M_2) + \text{d}(M_2, M_3)$, thus* d *is a distance.*

*Proof.* The triangular inequality is given in [8, 6.1]. For the second point, we have

$$\begin{aligned}
\text{d}(M_1, M_3) &= ||\text{Pos}(M_1, M_3)|| \\
&\leq ||\text{Pos}(M_1, M_2) + \text{Pos}(M_2, M_3)|| \\
&\leq ||\text{Pos}(M_1, M_2)|| + ||\text{Pos}(M_2, M_3)|| \\
&= \text{d}(M_1, M_2) + \text{d}(M_2, M_3).
\end{aligned}$$

where the first and second inequalities are given, respectively, by the first and second properties of the norm listed in 2.1Space of types and the dominance orderItem.8. □

We study some properties of the filtration associated to two lattices in the next proposition and, afterwards, we will study the properties of the relative position.

**Proposition 2.3.** *The filtration $\mathcal{F}(M_1, M_2)$ verifies the following properties:*

1. *It is compatible with tensor products, i.e. for $M_1, M_2$ (resp. $M_1', M_2'$) two $\mathcal{O}_K$-lattices in $V$ (resp. $V'$), we have*
$$\mathcal{F}(M_1 \otimes M_1', M_2 \otimes M_2') = \mathcal{F}(M_1, M_2) \otimes \mathcal{F}(M_1', M_2').$$

2. *Let $0 \to W \to V \xrightarrow{\pi} V/W \to 0$ be an exact sequence of $K$-vector spaces, $M_1, M_2$ two lattices in $V$, and let $N_i = M_i \cap W$ and $Q_i = \pi(M_i)$ for $i = 1, 2$. Suppose there exists a basis adapted to $M_1$, $M_2$ and $W$. Then, for each $i \in \mathbb{Z}$, there is an exact sequence*
$$0 \to \mathcal{F}^i(N_1, N_2) \to \mathcal{F}^i(M_1, M_2) \to \mathcal{F}^i(Q_1, Q_2) \to 0.$$

3. *Symmetric powers: for $M_1$ and $M_2$ two $\mathcal{O}_K$-lattices in $V$ and for every $k \geq 1$, we have*
$$\mathcal{F}(\text{Sym}^k M_1, \text{Sym}^k M_2) = \text{Sym}^k(\mathcal{F}(M_1, M_2)).$$

4. *Exterior powers: for $M_1$ and $M_2$ two $\mathcal{O}_K$-lattices in $V$ and for every $k \geq 1$, we have*
$$\mathcal{F}(\Lambda^k M_1, \Lambda^k M_2) = \Lambda^k(\mathcal{F}(M_1, M_2)).$$

5. *Direct sums: for $M_1, M_2$ (resp. $M_1', M_2'$) two $\mathcal{O}_K$-lattices in $V$ (resp. $V'$) and for every $i \in \mathbb{Z}$, we have*
$$\mathcal{F}^i(M_1 \oplus M_1', M_2 \oplus M_2') = \mathcal{F}^i(M_1, M_2) \oplus \mathcal{F}^i(M_1', M_2').$$

6. *Graduations:* For a filtration $\mathcal{F}$ and an $\mathcal{O}_K$-lattice $M$ on a $K$-vector space $V$, we define an $\mathcal{O}_K$-lattice in $\mathrm{Gr}_{\mathcal{F}}^{\gamma} V$ by $\mathrm{Gr}_{\mathcal{F}}^{\gamma} M = (M \cap \mathcal{F}^{\geq \gamma})/(M \cap \mathcal{F}^{> \gamma}) \subset \mathrm{Gr}_{\mathcal{F}}^{\gamma} V$. Then, for two $\mathcal{O}_K$-lattices $M_1, M_2$ in $V$ such that there exists an adapted basis for $M_1, M_2$ and $\mathcal{F}$, we have

$$\mathcal{F}(\mathrm{Gr}_{\mathcal{F}} M_1, \mathrm{Gr}_{\mathcal{F}} M_2) = \oplus_{\gamma \in \Gamma} \mathcal{F}(\mathrm{Gr}_{\mathcal{F}}^{\gamma} M_1, \mathrm{Gr}_{\mathcal{F}}^{\gamma} M_2).$$

**Proposition 2.4.** *The relative position verifies the following properties:*

1. *The relative position is antisymmetric, i.e. for $M_1, M_2$ two $\mathcal{O}_K$-lattices, we have*

$$\mathrm{Pos}(M_2, M_1) = \mathrm{Pos}(M_1, M_2)^{\iota}.$$

2. *Let $0 \to W \to V \xrightarrow{\pi} V/W \to 0$ be an exact sequence of $K$-vector spaces, $M_1, M_2$ two lattices in $V$, and let $N_i = M_i \cap W$ and $Q_i = \pi(M_i)$ for $i = 1, 2$. We have*

$$\mathrm{Pos}(N_1, N_2) * \mathrm{Pos}(Q_1, Q_2) \leq \mathrm{Pos}(M_1, M_2),$$

   *with an equality if there exists a basis adapted to $M_1, M_2$ and $W$.*

3. *Let $M_1, M_2$ (respectively, $M_1', M_2'$) be two $\mathcal{O}_K$-lattices inside a $K$-vector space $V$ (respectively, $V'$). Then we have*
$$\mathrm{Pos}(M_1 \oplus M_1', M_2 \oplus M_2') = \mathrm{Pos}(M_1, M_2) * \mathrm{Pos}(M_1', M_2').$$

4. *Let $M_1, M_2$ (respectively, $M_1', M_2'$) be two $\mathcal{O}_K$-lattices in the $K$-vector space $V$ (respectively, $V'$). Then we have*
$$\mathrm{Pos}(M_1 \otimes M_1', M_2 \otimes M_2') = \mathrm{Pos}(M_1, M_2) \otimes \mathrm{Pos}(M_1', M_2').$$

5. *Let $\mathcal{F}$ a filtration on a $K$-vector space $V$ and $M$ a $\mathcal{O}_K$-lattice in $V$. Then, for $M_1, M_2$ two $\mathcal{O}_K$-lattices in $V$, we have*
$$\mathrm{Pos}(\mathrm{Gr}_{\mathcal{F}} M_1, \mathrm{Gr}_{\mathcal{F}} M_2) \leq \mathrm{Pos}(M_1, M_2),$$

   *with equality if there exists a basis adapted to $M_1, M_2$ and $\mathcal{F}$.*

*Proof.* Only point 2 is tricky, and its proof can be found in [8, Proposition 99]. □

**Definition 2.1.** Given a $K$-vector space $V$, a $\mathbb{Z}$-filtration $\mathcal{F}$ on $V$ and an $\mathcal{O}_K$-lattice $M$ in $V$, we can define another $\mathcal{O}_K$-lattice by

$$M + \mathcal{F} = \sum_{i \in \mathbb{Z}} u^{-i} M \cap \mathcal{F}^i V.$$

An easy computation with adapted basis gives us the following proposition:

**Proposition 2.5.** *The addition operator defined above verifies the following properties:*

1. *We have $\mathrm{Pos}(M, M + \mathcal{F}) = \mathbf{t}(\mathcal{F})$.*

2. *It is compatible with tensor products: we have*

$$(M_1 \otimes M_2) + (\mathcal{F}_1 \otimes \mathcal{F}_2) = (M_1 + \mathcal{F}_1) \otimes (M_2 + \mathcal{F}_2)$$

   *for $M_1 \subset V_1$, $M_2 \subset V_2$ $\mathcal{O}_K$-lattices and $\mathcal{F}_1, \mathcal{F}_2$ filtrations on the $K$-vector spaces $V_1$ and $V_2$, respectively.*

## 2.4 The Harder-Narasimhan formalism

In this section, we will recall the Harder-Narasimhan formalism given by André in [2] adding some extra conditions. Let $\mathcal{C}$ be a quasi-abelian category and $\Gamma$ a $\mathbb{Q}$-subspace of $\mathbb{R}$. We denote by $\mathrm{sk}\,\mathcal{C}$ the skeleton of $\mathcal{C}$, i.e. the isomorphism classes of objects of $\mathcal{C}$.

**Definition 2.2.** A *rank function* on $\mathcal{C}$ is a map, $\mathrm{rank} : \mathrm{sk}\,\mathcal{C} \to \mathbb{N}$, that is additive on short exact sequences and takes the value 0 only on the 0 object. A *degree function* on $\mathcal{C}$ with values in $\Gamma$ is a map, $\deg : \mathrm{sk}\,\mathcal{C} \to \Gamma$, taking value 0 (not necessarily only) on the 0 object, that satisfies the following two conditions:

1. It is additive on short exact sequences.

2. For any epi-monic $M_1 \to M_2$, one has $\deg(M_1) \leq \deg(M_2)$.

Given a degree and rank function, we can associate a slope function

$$\mu = \frac{\deg}{\text{rank}} \; : \; \text{sk}\,\mathcal{C}\backslash\{0\} \to \mathbb{Q}.$$

A non-zero object $M$ is called $\mu$-semistable (or just semistable when there is no ambiguity) when for every strict subobject $M' \subset M$, we have $\mu(M') \leq \mu(M)$. We denote by $\mathcal{C}(\gamma)$ the full subcategory of objects of $\mathcal{C}$ which are $\mu$-semistable of slope $\gamma \in \Gamma$, together with the zero object.

Suppose we are given a rank function on $\mathcal{C}$. A $\Gamma$-filtration on $\mathcal{C}$ is a functor $\mathcal{F}^{\geq} : \Gamma^{\text{op}} \times \mathcal{C} \to \mathcal{C}$ which sends any object $(\gamma, M)$ to a strict subobject $\mathcal{F}^{\geq \gamma} M$ of $M$ and verifying:

- It is decreasing, i.e. $\mathcal{F}^{\geq \gamma'} M \subset \mathcal{F}^{\geq \gamma} M$ for $\gamma' \geq \gamma$ for any $M$ in $\mathcal{C}$.

- It is separated, i.e. $\varprojlim \mathcal{F}^{\geq \gamma} M = 0$ for any $M$ in $\mathcal{C}$.

- It is exhaustive, i.e $\varinjlim \mathcal{F}^{\geq \gamma} M = M$ for any $M$ in $\mathcal{C}$.

- It is left continuous, i.e. $\mathcal{F}^{\geq \gamma} M = \varprojlim_{\gamma' < \gamma} \mathcal{F}^{\geq \gamma'} M$ for any $M$ in $\mathcal{C}$.

We can now state the main theorem, which gives us the Harder-Narasimhan filtration.

**Theorem 2.6.** *[2, 1.4.7] Given a slope filtration $\mu$ on $\mathcal{C}$, there exists a unique filtration $\mathcal{F}$ on $\mathcal{C}$ such that for every $M \in \mathcal{C}$, the flag*

$$F(M) \; : \; 0 \hookrightarrow \mathcal{F}^{\geq \gamma_1} M = M_1 \hookrightarrow \ldots \hookrightarrow \mathcal{F}^{\geq \gamma_r} M = M_r = M$$

*attached to $\mathcal{F}^{\geq \cdot}$ is the unique flag on $M$ whose graded pieces $M_i/M_{i-1}$ for $1 \leq i \leq r$ are semistable of decreasing slopes*

$$\gamma_1 > \ldots > \gamma_r.$$

The slope $\Gamma$-filtration $\mathcal{F}^{\geq \cdot}$ of the theorem is called the Harder-Narasimhan filtration of $M$ and the slopes $\gamma_1 > \ldots > \gamma_r$ are called the breaks of the filtration $\mathcal{F}$. Sometimes, we will refer to the flag associated to $\mathcal{F}^{\geq \cdot}$ as the Harder-Narasimhan filtration of $M$.

### 2.4.1 Further assumptions

**Proposition 2.7.** *Let $\mathcal{C}$ be a quasi-abelian category with rank and degree functions such that*

1. *For a mono-epi $M_1 \to M_2$, we have $\text{rank}\,M_1 = \text{rank}\,M_2$.*

2. *The function $\text{rank}$ detects the mono-epi morphisms, i.e. it verifies that a mono (resp. an epi) $M_1 \to M_2$ in $\mathcal{C}$ is an epi (resp. a mono) if and only if $\text{rank}\,M_1 = \text{rank}\,M_2$.*

3. *The functions $\deg$ detects the isomorphisms, i.e. for a mono-epi $M_1 \xrightarrow{f} M_2$ we have $\deg M_1 = \deg M_2$ if and only if $f$ is an isomorphism in $\mathcal{C}$.*

*Then,*

1. *The function $\mu$ also detects the isomorphisms.*

2. *The category $\mathcal{C}(\gamma)$ is abelian, for every $\gamma \in \Gamma$.*

3. *Let $0 \to M_1 \to M_2 \to M_3 \to 0$ be an exact sequence in $\mathcal{C}$. Suppose that $\mathbf{t}_{\text{HN}}(M_2) = \mathbf{t}_{\text{HN}}(M_1) * \mathbf{t}_{\text{HN}}(M_3)$. Then, for every $\gamma \in \Gamma$, we have an exact sequence*

$$0 \to \mathcal{F}^{\geq \gamma}_{\text{HN}} M_1 \to \mathcal{F}^{\geq \gamma}_{\text{HN}} M_2 \to \mathcal{F}^{\geq \gamma}_{\text{HN}} M_3 \to 0.$$

*Proof.* 1. It follows from the condition (1).

2. Let $f : M_1 \to M_2$ be a morphism in $\mathcal{C}(\gamma)$, where $M_1$ and $M_2$ are nonzero objects in $\mathcal{C}(\gamma)$. Then, the inclusion $\ker f \subset M_1$ gives us $\mu(\ker f) \leq \gamma$ and $\operatorname{im} f, \operatorname{coim} f \subset M_2$ gives us $\mu(\operatorname{coim} f) \leq \gamma$ and $\mu(\operatorname{im} f) \leq \gamma$. We have an exact sequence

$$0 \to \ker f \to M_1 \to \operatorname{coim} f \to 0,$$

so $\min\{\mu(\ker f), \mu(\operatorname{coim} f)\} \leq \gamma \leq \max\{\mu(\ker f), \mu(\operatorname{coim} f)\}$, thus we have an equality $\mu(\ker f) = \mu(\operatorname{coim} f) = \gamma$, so $\ker f, \operatorname{coim} f$ are semi-stable. Since $\operatorname{coim} f \hookrightarrow \operatorname{im} f$ is a mono-epi, we have $\gamma = \mu(\operatorname{coim} f) \leq \mu(\operatorname{im} f) \leq \gamma$, so $\mu(\operatorname{im} f) = \gamma$ and then $\operatorname{im} f$ is semi-stable. To finish, the exact sequence

$$0 \to \operatorname{im} f \to M_2 \to \operatorname{coker} f \to 0$$

gives us $\mu(\operatorname{coker} f) = \gamma$, so $\operatorname{coker} f$ is semi-stable.

3. Let $0 \to M_1 \to M_2 \to M_3 \to 0$ be an exact sequence in $\mathcal{C}$ with $\mathbf{t}_{\mathrm{HN}}(M_2) = \mathbf{t}_{\mathrm{HN}}(M_1) * \mathbf{t}_{\mathrm{HN}}(M_3)$. By induction on the rank of $M_2$, it suffices to prove the exactness of

$$0 \to \mathcal{F}_{\mathrm{HN}}^{\geq \gamma} M_1 \to \mathcal{F}_{\mathrm{HN}}^{\geq \gamma} M_2 \to \mathcal{F}_{\mathrm{HN}}^{\geq \gamma} M_3 \to 0.$$

for $\gamma$ the maximal slope in $\mathbf{t}_{\mathrm{HN}}(M_2)$. We have a strict mono $\mathcal{F}^{\geq \gamma} M_1 \hookrightarrow M_1 \hookrightarrow M_2$, so we can consider the quotients

$$0 \to M_1/\mathcal{F}^{\geq \gamma} M_1 \to M_2/\mathcal{F}^{\geq \gamma} M_1 \to M_3 \to 0.$$

Suppose first $\mathcal{F}^{\geq \gamma} M_1 \neq 0$, then $\operatorname{rank}(M_2/\mathcal{F}^{\geq \gamma} M_1) < \operatorname{rank} M_2$. It is easy to see that we have again

$$\mathbf{t}_{\mathrm{HN}}(M_2/\mathcal{F}^{\geq \gamma} M_1) = \mathbf{t}_{\mathrm{HN}}(M_1/\mathcal{F}^{\geq \gamma} M_1) * \mathbf{t}_{\mathrm{HN}}(M_3)$$

so, by induction, we get

$$0 \to \mathcal{F}^{\geq \mu}(M_1/\mathcal{F}^{\geq \gamma} M_1) \to \mathcal{F}^{\geq \mu}(M_2/\mathcal{F}^{\geq \gamma} M_1) \to \mathcal{F}^{\geq \mu} M_3 \to 0$$

for every $\mu \in \Gamma$. In particular, for $\mu = \gamma$, we have $\mathcal{F}^{\geq \gamma}(M_2)/\mathcal{F}^{\geq \gamma}(M_1) \simeq \mathcal{F}^{\geq \gamma}(M_3)$ as we wanted. Suppose then $\mathcal{F}^{\geq \gamma} M_1 = 0$, then since $\mathcal{F}^{\geq \gamma} M_2$ and $\mathcal{F}^{\geq \gamma} M_3$ are semi-stable of slope $\gamma$, the kernel of $\mathcal{F}^{\geq \gamma} M_2 \to \mathcal{F}^{\geq \gamma} M_3$ is also semi-stable of slope $\gamma$ and contained in $M_1$, thus it is zero and the morphism is a mono. Now, $\mathcal{F}^{\geq \gamma} M_2$ and $\mathcal{F}^{\geq \gamma} M_3$ have the same rank by assumption, so by (2), we have that it is a strict epi. Then, $\mathcal{F}^{\geq \gamma} M_2 \to \mathcal{F}^{\geq \gamma} M_3$ is a mono-epi and $\deg(\mathcal{F}^{\geq \gamma} M_2) = \deg(\mathcal{F}^{\geq \gamma} M_3)$. It is then an isomorphism, by (3). □

## 2.5 Modules over $\mathfrak{S}$

Let $\mathbb{F}$ be a perfect field of characteristic $p$ and let $W(\mathbb{F})$ be the ring of Witt vectors over $\mathbb{F}$. The ring $\mathfrak{S} := W(\mathbb{F})[[u]]$ is a complete, regular, local ring of Krull dimension 2, with maximal ideal $\mathfrak{m} = (p, u)$, so it is a unique factorization domain. Moreover, for all $f \in \mathfrak{m}$ nonzero, the ring $\mathfrak{S}[\frac{1}{f}]$ is a unique factorization domain of Krull dimension 1, therefore a principal ideal domain.

### 2.5.1 Classification of finitely generated $\mathfrak{S}$-modules

Let $M$ be a finitely generated $\mathfrak{S}$-module. Then there is a unique exact sequence

$$0 \to M_{\mathrm{tors}} \to M \to M_{\mathrm{fr}} \to \overline{M} \to 0$$

where $M_{\mathrm{tors}}$ is a torsion module, $M_{\mathrm{fr}}$ is a free module and $\overline{M}$ is a finite length module. This dévissage of $M$ is given by Bhatt-Morrow-Scholze in [4].

**Definition 2.3.** A finitely generated $\mathfrak{S}$-module is called pseudo-null if it has finite length. A morphism between two $\mathfrak{S}$-modules is called a pseudo-isomorphism if its kernel and cokernel are pseudo-null.

The following theorem is the equivalent of the structure theorem for finitely generated modules over a PID. It is a consequence of [36, 5.1.10] and [36, 5.3.7].

**Theorem 2.8** (Structure theorem of finitely generated $\mathfrak{S}$-modules). *Let $M$ be a finitely generated $\mathfrak{S}$-module. Then there exists a pseudo-isomorphism*

$$M \sim \mathfrak{S}^r \oplus \bigoplus_{i=0}^{s} \mathfrak{S}/(p^{n_i}) \oplus \bigoplus_{j=0}^{t} \mathfrak{S}/(P_j(u)^{m_j})$$

*for some $r \geq 0$, $s, t \geq 0$, $n_i \geq 0$ for $0 \leq i \leq s$, and $m_j \geq 0$ and $P_j(u)$ distinguished irreducible $\mathfrak{S}$-polynomials for $0 \leq j \leq t$. These quantities are unique up to reordering.*

The structure theorem gives us two important invariants:

$$\mu_{\mathrm{IW}}(M) = \sum_{i=0}^{r} n_i, \quad \lambda_{\mathrm{IW}}(M) = \sum_{j=0}^{t} m_j \deg(f_j).$$

*Remark* 2. After localization by $(p)$ in the structure theorem, we can see that $\mu_{\mathrm{IW}}(M) = \mathrm{length}_{\mathfrak{S}_{(p)}}(M_{(p)})$. As a consequence, we have that $\mu_{\mathrm{IW}}$ is additive on short exact sequences of torsion $\mathfrak{S}$-modules.

### 2.5.2 Cohomological properties

Let $M$ be a finitely generated $\mathfrak{S}$-module. There is an exact sequence $0 \to M[\mathfrak{m}^\infty] \to M_{\mathrm{tors}} \to M_{\mathrm{t}} \to 0$, where $M[\mathfrak{m}^\infty]$ is pseudo-null.

**Definition 2.4.** We define the depth of $M$ as

$$\mathrm{depth}(M) := \min\{i \mid \mathrm{Ext}^i(\mathbb{F}, M) \neq 0\}$$

or, equivalently, the maximal length of an $M$-regular $\mathfrak{m}$-sequence, i.e. a sequence of elements $x_1, \ldots, x_r \in \mathfrak{m}$ such that $x_1$ is not a zero divisor of $M$ and $x_i$ is not a zero divisor in $M/(x_1, \ldots, x_{i-1})M$ for $2 \leq i \leq r$.

The projective dimension of $M$, that we denote $\mathrm{pd}(M)$, is the minimal length of a projective resolution of $M$. In the case of modules over noetherian local rings, there is an interesting relation between the projective dimension and the Tor functor. For a finitely generated $\mathfrak{S}$-module, we have

$$\mathrm{pd}(M) = \max\{i \mid \mathrm{Tor}^i(\mathbb{F}, M) \neq 0\}$$

and the Auslander-Buchsbaum theorem ([3]) tells us that

$$\mathrm{pd}(M) + \mathrm{depth}(M) = 2.$$

As a consequence, for finite free and finitely generated torsion modules over $\mathfrak{S}$, we have

- For a finitely generated $\mathfrak{S}$-module $M$, the following conditions are equivalent: $M$ is projective, $M$ flat, $M$ is free, $M$ is reflexive and $\mathrm{Tor}^1(\mathbb{F}, M) = 0$.

- For a finitely generated torsion $\mathfrak{S}$-module $M$, the following conditions are equivalent: $\mathrm{pd}(M) = 1$, $\mathrm{depth}(M) = 1$ and $\mathrm{Ext}^0(\mathbb{F}, M) \simeq \mathrm{Tor}^2(\mathbb{F}, M) = M[\mathfrak{m}] = 0$.

**Proposition 2.9.** *Let $0 \to K \to M \to Q \to 0$ with $M$ a finite free $\mathfrak{S}$-module and $\mathrm{pd}(Q) \leq 1$. Then $K$ is a finite free $\mathfrak{S}$-module.*

*Proof.* Since $M$ is a finite free $\mathfrak{S}$-module, we know that $\mathrm{Tor}^i(\mathbb{F}, M) = 0$ for $i = 1, 2$, and thus $\mathrm{Tor}^1(\mathbb{F}, K) = \mathrm{Tor}^2(\mathbb{F}, Q) = 0$ since $\mathrm{pd}(Q) \leq 1$. □

Another way to check if $M$ is free is given by the following proposition.

**Proposition 2.10.** *Let $M$ be a finitely generated $\mathfrak{S}$-module. Then $M$ is free if and only if*

1. *The reduction $M/p^n M$ has no $u$-torsion for every $n \geq 1$.*

2. *For every $n \geq 1$, we have $\mu_{\mathrm{IW}}(M/p^n M) = n\, \mu_{\mathrm{IW}}(M/pM)$.*

*Proof.* First, we remark that the conditions are necessary. Conversely, by hypothesis we know that the $\mathbb{F}[[u]]$-module $M/pM$ is free. Let $m$ be its rank. By Nakayama's lemma, we know that there is a system of $m$ generators of $M$, giving us the exact sequence

$$0 \to \ker f \to \mathfrak{S}^m \xrightarrow{f} M \to 0.$$

Considering reductions modulo $p^n$ for each $n \geq 1$ we get

$$0 \to M[p^n] \to \ker f / p^n \ker f \to (\mathfrak{S}/p^n\mathfrak{S})^m \to M/p^n M \to 0$$

that we can cut into two exact sequences:

$$0 \to M[p^n] \to \ker f / p^n \ker f \to K_n \to 0$$
$$0 \to K_n \to (\mathfrak{S}/p^n\mathfrak{S})^m \to M/p^n M \to 0.$$

Now, the second hypothesis tells us that

$$\mu_{\text{IW}}((\mathfrak{S}/p^n\mathfrak{S})^m) - \mu_{\text{IW}}(M/p^n M) = mn - n\mu_{\text{IW}}(M/pM) = 0.$$

So by the additivity of $\mu_{\text{IW}}$ on short exact sequences of torsion $\mathfrak{S}$-modules seen in Remark 2rem.2 we get $\mu_{\text{IW}}(K_n) = 0$. Being also killed by $p^n$, $K_n$ has to be pseudo-null. Since $(\mathfrak{S}/p^n\mathfrak{S})^m$ contains no nonzero pseudo-null submodule, actually $K_n = 0$. Therefore $(\mathfrak{S}/p^n\mathfrak{S})^m \simeq M/p^n M$. Since $\mathfrak{S}$ and $M$ are $p$-adically complete, it follows that $\mathfrak{S}^m \simeq M$. □

### 2.5.3 Some categories of $\mathfrak{S}$-modules

Denote $\text{Mod}_\mathfrak{S}$ the category of finitely generated $\mathfrak{S}$-modules and let $\mathcal{A}$ be the category of finitely generated $\mathfrak{S}$-modules killed by a power of $p$. It is a full subcategory of $\text{Mod}_\mathfrak{S}$ which is abelian and stable by subobjects, quotients and extensions. Inside $\mathcal{A}$, we consider two subcategories:

- The full subcategory $\mathcal{T}$ of objects with finite length, i.e. the $\mathfrak{S}$-modules such that $M[\mathfrak{m}^\infty] = M$ (i.e. $M[u^\infty] = M$, since the $\mathfrak{m}^\infty$-torsion and the $u^\infty$-torsion coincide in this category).

- The full subcategory $\text{Mod}_{\mathfrak{S},\text{t}}$ of objects with projective dimension 1. This is a full subcategory of $\text{Mod}_\mathfrak{S}$, it contains $\mathfrak{S}/p\mathfrak{S} = \mathbb{F}[[u]]$ and it is stable by subobjects and extensions, but not stable by quotients.

For any module $M$ in $\mathcal{A}$, consider the short exact sequence $0 \to M[u^\infty] \to M \to M' \to 0$, where $M[u^\infty]$ is an object in $\mathcal{T}$ and $M' = M/M[u^\infty]$ is an object in $\text{Mod}_{\mathfrak{S},\text{t}}$. There are no nonzero morphisms from $\mathcal{T}$ to $\text{Mod}_{\mathfrak{S},\text{t}}$ since the objects in $\text{Mod}_{\mathfrak{S},\text{t}}$ have no $u^\infty$-torsion. This proves that $\text{Mod}_{\mathfrak{S},\text{t}}$ is a quasi-abelian category.

There is a rank function in the sense of the Harder-Narasimhan formalism on the abelian category $\mathcal{A}$, given by the Iwasawa invariant $\mu_{\text{IW}}$ and we can reinterpret it in the following way:

**Proposition 2.11.** *Let $M$ be an object in $\mathcal{A}$, then*

$$\mu_{\text{IW}}(M) = \text{length}_{\mathfrak{S}_{(p)}}(M_{(p)}) = \text{length}_{\mathfrak{S}[\frac{1}{u}]}(M[\tfrac{1}{u}])$$

*and for $M$ an object of the subcategory $\text{Mod}_{\mathfrak{S},\text{t}}$, we have*

$$\mu_{\text{IW}}(M) = \text{length}_{\mathfrak{S}_{(p)}}(M_{(p)}) = \text{length}_{\mathfrak{S}[\frac{1}{u}]}(M[\tfrac{1}{u}]) = \text{length}_{W(\mathbb{F})}(M/uM).$$

## 3 The Fargues filtration on Kisin modules

Let $K_0 = \text{Frac}\, W(\mathbb{F})$ and let $K$ be a finite totally ramified extension of $K_0$, with ring of integers $\mathcal{O}_K$, uniformizer $\pi$ and residue field $\mathbb{F}$. Let $E(u) \in \mathfrak{S}$ be the minimal polynomial of $\pi$, which is Eisenstein. Then, $E$ is an irreducible distinguished polynomial, hence a prime element in $\mathfrak{S}$, and $\mathfrak{S}/E\mathfrak{S} \simeq \mathcal{O}_K$ by $u \mapsto \pi$.

Consider the category $\text{Mod}_\mathfrak{S}^\varphi$, whose objects are pairs $(M, \varphi_M)$ with $M$ a module in $\text{Mod}_\mathfrak{S}$ and $\varphi_M$ is a Frobenius isomorphism $\varphi_M : \varphi^* M[\tfrac{1}{E}] \xrightarrow{\sim} M[\tfrac{1}{E}]$, where $\varphi$ is the Frobenius over $\mathfrak{S}$. When there is no confusion, we will denote only by $M$ the objects in $\text{Mod}_\mathfrak{S}^\varphi$. A morphism $f : (M, \varphi_M) \to (N, \varphi_N)$ in $\text{Mod}_\mathfrak{S}^\varphi$ is a morphism of $\mathfrak{S}$-modules $f : M \to N$ compatible with the Frobenius defined on $M$ and $N$. The category $\text{Mod}_\mathfrak{S}^\varphi$ is abelian and there is a unique dévissage

$$
\begin{array}{c}
0 \\
\downarrow \\
(M[\mathfrak{m}^\infty], \varphi_{\mathfrak{m}^\infty}) \\
\downarrow \\
0 \longrightarrow (M_{\text{tors}}, \varphi_{\text{tors}}) \longrightarrow (M, \varphi_M) \longrightarrow (M_{\text{fr}}, \varphi_{\text{fr}}) \longrightarrow (\overline{M}, \varphi_{\overline{M}}) \longrightarrow 0 \\
\downarrow \\
(M_{\text{t}}, \varphi_{\text{t}}) \\
\downarrow \\
0
\end{array}
$$

where $M_{\text{fr}}$ is a free $\mathfrak{S}$-module, $M_{\text{tors}}$ is killed by a power of $p$ (i.e. it is an object in the category $\mathcal{A}$ from last section), $M_{\text{t}}$ is an object in $\text{Mod}_{\mathfrak{S},\text{t}}$, and $\overline{M}$ and $M[\mathfrak{m}^\infty]$ are $\mathfrak{S}$-modules of finite length. The proof that $M_{\text{tors}}$ is killed by a power of $p$ is given by Bhatt, Morrow and Scholze in [4]. This decomposition above allows us to define the following categories:

1. The category $\text{Mod}^\varphi_{\mathfrak{S},\,\text{fr}}$: It is the additive full subcategory of $\text{Mod}^\varphi_{\mathfrak{S}}$ whose objects are pairs $(M, \varphi_M)$ with $M$ a free module. We call its objects Kisin modules.

2. The category $\text{Mod}^\varphi_{\mathfrak{S},\text{tors}}$: It is the abelian full subcategory of $\text{Mod}^\varphi_{\mathfrak{S}}$ whose objects are pairs $(M, \varphi_M)$ with $M$ a module in $\mathcal{A}$. By point 2 in the remark above, we see that the objects in $\text{Mod}^\varphi_{\mathfrak{S},\,\text{t}}$ do not depend on $K$, $\pi$ or $E$.

Using the exact sequence

$$0 \to (M[\mathfrak{m}^\infty], \varphi_{\mathfrak{m}^\infty}) \to (M_{\text{tors}}, \varphi_{\text{tors}}) \to (M_{\text{t}}, \varphi_{\text{t}}) \to 0$$

above, we can define a torsion theory on $\text{Mod}^\varphi_{\mathfrak{S},\text{tors}}$, using the following two categories:

(3) The category $\text{Mod}^\varphi_{\mathfrak{S},\,\text{t}}$: It is the full subcategory of $\text{Mod}^\varphi_{\mathfrak{S},\text{tors}}$ whose objects are pairs $(M, \varphi_M)$ with $M$ a module in $\text{Mod}_{\mathfrak{S},\text{t}}$. We call its object $p^\infty$-torsion Kisin modules. This will play the part of the torsion-free subcategory. It is then a quasi-abelian category.

(4) The category $\text{Mod}^\varphi_{\mathfrak{S},\text{fl}}$: It is the abelian full subcategory of $\text{Mod}^\varphi_{\mathfrak{S}}$ whose objects are pairs $(M, \varphi_M)$ with $M$ a finite length module, i.e. $M[\mathfrak{m}^\infty] = M$. We have seen that $\varphi_{\mathfrak{m}^\infty} = 0$ in Remark 5, so we have an equivalence of categories $\text{Mod}^\varphi_{\mathfrak{S},\text{fl}} \simeq \text{Mod}_{\mathfrak{S},\text{fl}}$. This will play the part of the torsion subcategory.

One last category will be studied, the $p$-torsion category:

(5) The category $\text{Mod}^\varphi_{\mathbb{F}[[u]],\text{fr}}$: It is the additive full subcategory of $\text{Mod}^\varphi_{\mathfrak{S},\,\text{t}}$ whose objects are killed by $p$. This corresponds to finite free $\mathbb{F}[[u]]$-modules (since they do not have $u$-torsion, as $M[E] = M[u]$ for every $M \in \text{Mod}_\mathfrak{S}$ killed by $p$ and objects in $\text{Mod}^\varphi_{\mathfrak{S},\,\text{t}}$ have no $E$-torsion) endowed with a Frobenius $\varphi_M : \varphi^* M[\frac{1}{u}] \xrightarrow{\sim} M[\frac{1}{u}]$. We call the objects in this category $p$-torsion Kisin modules. This is a quasi-abelian category, since it is fully embedded in $\text{Mod}^\varphi_{\mathfrak{S},\,\text{t}}$ and it is stable by kernels and cokernels.

**Definition 3.1.** Let $M$ be an object in $\text{Mod}^\varphi_\mathfrak{S}$ such that $M[\mathfrak{m}^\infty] = 0$. Then $M \subset M[\frac{1}{E}]$ and $\varphi^* M \subset \varphi^* M[\frac{1}{E}]$. We say that $M$ is effective when $\varphi_M(\varphi^* M) \subset M$. Thus $\varphi_M : \varphi^* M \to M$ is an injective morphism with cokernel killed by a power of $E$.

*Remark* 3. It is easy to see that if $M$ is effective, then a subobject $N$ of $M$ is effective if $M/N$ has no $\mathfrak{m}^\infty$-torsion and a quotient of $M$ is effective if it has no $\mathfrak{m}^\infty$-torsion.

### 3.1 Isogenies classes of Kisin modules

Let $\text{Mod}^\varphi_\mathfrak{S} \otimes \mathbb{Q}_p$ be the category whose objects are the same objects as in $\text{Mod}^\varphi_\mathfrak{S}$ and whose morphisms are morphisms in $\text{Mod}^\varphi_\mathfrak{S}$ tensored by $\mathbb{Q}_p$. It is the isogeny category of Kisin modules. For $(M, \varphi_M)$ in $\text{Mod}^\varphi_\mathfrak{S}$, we denote by $(M, \varphi_M) \otimes \mathbb{Q}_p$ the corresponding object in $\text{Mod}^\varphi_\mathfrak{S} \otimes \mathbb{Q}_p$ (or $M \otimes \mathbb{Q}_p$ as an abuse of notation when there is no confusion). The adjoint functors

$$
\begin{array}{ccc}
\text{Mod}^\varphi_\mathfrak{S} & \to & \text{Mod}^\varphi_{\mathfrak{S},\,\text{fr}} \\
M & \mapsto & M_{\text{fr}}
\end{array}
\quad \text{and} \quad
\begin{array}{ccc}
\text{Mod}^\varphi_{\mathfrak{S},\,\text{fr}} & \to & \text{Mod}^\varphi_\mathfrak{S} \\
M & \mapsto & M
\end{array}
$$

induce an equivalence of categories $\text{Mod}^\varphi_\mathfrak{S} \otimes \mathbb{Q}_p \simeq \text{Mod}^\varphi_{\mathfrak{S},\,\text{fr}} \otimes \mathbb{Q}_p$. Both categories $\text{Mod}^\varphi_\mathfrak{S} \otimes \mathbb{Q}_p$ and $\text{Mod}^\varphi_{\mathfrak{S},\,\text{fr}} \otimes \mathbb{Q}_p$ are abelian, since $\text{Mod}^\varphi_\mathfrak{S}$ is abelian.

**Definition 3.2.** Let $M_1, M_2 \in \mathrm{Mod}_{\mathfrak{S}}^{\varphi}$. An isogeny between $M_1$ and $M_2$ is a morphism $f : M_1 \to M_2$ which becomes an isomorphism in $\mathrm{Mod}_{\mathfrak{S}}^{\varphi} \otimes \mathbb{Q}_p$. Equivalently, it is a morphism whose kernel and cokernel are objects in $\mathrm{Mod}_{\mathfrak{S},\mathrm{tors}}^{\varphi}$. We say that two objects are isogenous if there is an isogeny between them, i.e. if they become isomorphic in the category $\mathrm{Mod}_{\mathfrak{S}}^{\varphi} \otimes \mathbb{Q}_p$.

There is a fully faithfully functor, studied in [19],

$$\begin{array}{rcl} \mathrm{Mod}_{\mathfrak{S}}^{\varphi} \otimes \mathbb{Q}_p & \to & \mathrm{Mod}_{\mathfrak{S}[\frac{1}{p}]}^{\varphi} \\ (M \otimes \mathbb{Q}_p, \varphi_M) & \mapsto & (M[\frac{1}{p}], \varphi_M \otimes 1) \end{array}$$

where $\mathrm{Mod}_{\mathfrak{S}[\frac{1}{p}]}^{\varphi}$ is the abelian category whose objects are finitely generated $\mathfrak{S}[\frac{1}{p}]$-modules $N$ together with a Frobenius $\varphi_N : \varphi^* N[\frac{1}{E}] \xrightarrow{\sim} N[\frac{1}{E}]$, and whose morphisms are the morphisms between modules compatible with the Frobenius. When an object in $\mathrm{Mod}_{\mathfrak{S}[\frac{1}{p}]}^{\varphi}$ comes from an object $M$ in $\mathrm{Mod}_{\mathfrak{S},\mathrm{fr}}^{\varphi}$, we will denote it by $M[\frac{1}{p}]$. By proceeding as in [4, Proposition 4.3], we can prove that all the modules in $\mathrm{Mod}_{\mathfrak{S}[\frac{1}{p}]}^{\varphi}$ are free modules.

**Lemma 3.1.** *The essential image of the functor* $\mathrm{Mod}_{\mathfrak{S}}^{\varphi} \otimes \mathbb{Q}_p \to \mathrm{Mod}_{\mathfrak{S}[\frac{1}{p}]}^{\varphi}$ *given above is stable by subobjects and quotients. In particular, it is an abelian subcategory of* $\mathrm{Mod}_{\mathfrak{S}[\frac{1}{p}]}^{\varphi}$.

*Proof.* It suffices to show that for any submodule $N$ of $M[\frac{1}{p}]$ in $\mathrm{Mod}_{\mathfrak{S}[\frac{1}{p}]}^{\varphi}$, where $M \in \mathrm{Mod}_{\mathfrak{S},\mathrm{fr}}^{\varphi}$, there is a module $N'$ in $\mathrm{Mod}_{\mathfrak{S}}^{\varphi}$ such that $N'[\frac{1}{p}] = N$. We can take $N' = N \cap M$ and it is easy to show that this is a module in $\mathrm{Mod}_{\mathfrak{S}}^{\varphi}$, by the compatibily between $\varphi_M$, $\varphi_M[\frac{1}{p}]$ and $\varphi_N$. Then, each quotient is also included in the essential image, since the localization functor is exact, thus the essential image forms an abelian category. The lemma follows since the essential image of $\mathrm{Mod}_{\mathfrak{S},\mathrm{fr}}^{\varphi} \to \mathrm{Mod}_{\mathfrak{S}[\frac{1}{p}]}^{\varphi}$ equals the essential image of $\mathrm{Mod}_{\mathfrak{S}}^{\varphi} \otimes \mathbb{Q}_p \to \mathrm{Mod}_{\mathfrak{S}[\frac{1}{p}]}^{\varphi}$. $\square$

Now, we want to show that there are Harder-Narasimhan filtrations in the two categories, and compare them. For $\mathrm{Mod}_{\mathfrak{S}[\frac{1}{p}]}^{\varphi}$, we put

$$\begin{array}{rclcrcl} \mathrm{rank} & : & \mathrm{Mod}_{\mathfrak{S}[\frac{1}{p}]}^{\varphi} & \to & \mathbb{Z}, & & \\ & & N & \mapsto & \mathrm{rank}_{\mathfrak{S}[\frac{1}{p}]} N & & \end{array} \qquad \begin{array}{rclcrcl} \deg & : & \mathrm{Mod}_{\mathfrak{S}[\frac{1}{p}]}^{\varphi} & \to & \mathbb{Z} \\ & & N & \mapsto & \nu(N \otimes \hat{\mathfrak{S}}, \varphi_N \varphi^* N \otimes \hat{\mathfrak{S}}). \end{array}$$

An easy calculation shows that the rank and degree functions verify the conditions of a Harder-Narasimhan rank and degree function. Together with the fact that $\mathrm{Mod}_{\mathfrak{S}[\frac{1}{p}]}^{\varphi}$ is an abelian category, we can define a Harder-Narasimhan filtration on objects of $\mathrm{Mod}_{\mathfrak{S}[\frac{1}{p}]}^{\varphi}$ that we will denote by $\mathcal{F}_{\mathrm{F},\mathfrak{S}[\frac{1}{p}]}$.

For $\mathrm{Mod}_{\mathfrak{S}}^{\varphi} \otimes \mathbb{Q}_p$, let

$$\begin{array}{rclcrcl} \mathrm{rank} & : & \mathrm{Mod}_{\mathfrak{S}}^{\varphi} \otimes \mathbb{Q}_p & \to & \mathbb{Z}, & & \\ & & M \otimes \mathbb{Q}_p & \mapsto & \mathrm{rank}_{\mathfrak{S}}(M) & & \end{array} \qquad \begin{array}{rclcrcl} \deg & : & \mathrm{Mod}_{\mathfrak{S}}^{\varphi} \otimes \mathbb{Q}_p & \to & \mathbb{Z} \\ & & M \otimes \mathbb{Q}_p & \mapsto & \nu(M \otimes \hat{\mathfrak{S}}, \varphi_M \varphi^* M \otimes \hat{\mathfrak{S}}). \end{array}$$

*Remark* 4. The degree and rank function on $M \otimes \mathbb{Q}_p$ and $M[\frac{1}{p}]$ coincide.

Since the functor $\mathrm{Mod}_{\mathfrak{S}}^{\varphi} \otimes \mathbb{Q}_p \to \mathrm{Mod}_{\mathfrak{S}[\frac{1}{p}]}^{\varphi}$ is exact, and the degree and rank functions coincide in the two categories, it follows that they verify the properties of a Harder-Narasimhan degree and rank function. Together with the fact that $\mathrm{Mod}_{\mathfrak{S}}^{\varphi} \otimes \mathbb{Q}_p$ is an abelian category, we obtain a Harder-Narasimhan filtration on $\mathrm{Mod}_{\mathfrak{S}}^{\varphi} \otimes \mathbb{Q}_p$ that we denote by $\mathcal{F}_{\mathrm{F},\circ}$ and whose polygon we will denote by $\mathbf{t}_{\mathrm{F},\circ}$. Next proposition proves that the functor $\mathrm{Mod}_{\mathfrak{S}}^{\varphi} \otimes \mathbb{Q}_p \to \mathrm{Mod}_{\mathfrak{S}[\frac{1}{p}]}^{\varphi}$ is compatible with the Harder-Narasimhan filtrations given in each category. As a consequence of Lemma 3.1thm.3.1, we get:

**Proposition 3.2.** *Let* $M \in \mathrm{Mod}_{\mathfrak{S}}^{\varphi} \otimes \mathbb{Q}_p$, *then for every* $\gamma \in \mathbb{Q}$, *we have*

$$\mathcal{F}_{\mathrm{F},\mathfrak{S}[\frac{1}{p}]}^{\geq \gamma}(M[\tfrac{1}{p}]) = (\mathcal{F}_{\mathrm{F},\circ}^{\geq \gamma}(M))[\tfrac{1}{p}].$$

## 3.2 Hodge-Pink modules

We recall in this section the definitions and some results about Hodge-Pink modules that can be found in the article of Genestier and Lafforgue [19]. Let $K$ be a totally ramified finite extension of $K_0$ with uniformizer $\pi_K$, let $E \in \mathfrak{S}$ be the minimal polynomial of $\pi_K$ and denote by $\hat{\mathfrak{S}}$ the completion of $\mathfrak{S}[\frac{1}{p}]$ with respect to the ideal generated by $E$.

**Definition 3.3.** A Hodge-Pink module is a 3-tuple $(D, \varphi_D, V_D)$ where:
- $D$ is a finite dimensional $K_0$-vector space,
- $\sigma_D : \sigma^* D \to D$ is an isomorphism of $K_0$-vector spaces,
- $V_D$ is a Hodge-Pink structure on $D$: a free $\hat{\mathfrak{S}}$-module $V_D$ which is a lattice in $\varphi^* D \otimes_{K_0} \hat{\mathfrak{S}}[\frac{1}{E}]$.

We denote by MHP the category of Hodge-Pink modules. A morphism of Hodge-Pink modules $D \to D'$ is a morphism $f$ compatible with $\sigma_D$ and $\sigma_{D'}$ such that $f(V_D) \subset V_{D'}$. This is a quasi-abelian category.

*Remark* 5. If $(D, \sigma_D, V_D)$ is a Hodge-Pink module and $(D', \sigma_{D'})$ a subobject of the isocrystal $(D, \sigma_D)$, we can endow $D'$ with a Hodge-Pink structure by setting $V_{D'} = V_D \cap (\sigma^* D' \otimes_{K_0} \hat{\mathfrak{S}}[\frac{1}{E}])$.

We denote by $D$ the Hodge-Pink module $(D, \sigma_D, V_D)$ when there is no confusion. We define the Newton and Hodge types of $D$ by

$$\mathbf{t}_\mathrm{N}(D, \sigma_D, V_D) = \mathbf{t}_\mathrm{N}(D, \sigma_D) \quad \text{and} \quad \mathbf{t}_\mathrm{H}(D, \sigma_D, V_D) = \mathrm{Pos}(\sigma^* D \otimes_{K_0} \hat{\mathfrak{S}}, V_D)$$

where $\mathbf{t}_\mathrm{N}(D, \sigma_D)$ denotes the type associated to the Newton polygon of the isocrystal $(D, \sigma_D)$ given by the Dieudonné-Manin decomposition. We denote their degrees by

$$\mathrm{t}_\mathrm{N}(D) = \deg(\mathbf{t}_\mathrm{N}(D, \sigma_D, V_D)) \quad \text{and} \quad \mathrm{t}_\mathrm{H}(D) = \deg(\mathbf{t}_\mathrm{H}(D, \sigma_D, V_D)) = \nu(\sigma^* D \otimes_{K_0} \hat{\mathfrak{S}}, V_D).$$

The usual rank function and the degree function given by

$$\deg(D) = \mathrm{t}_\mathrm{H}(D) - \mathrm{t}_\mathrm{N}(D)$$

define a Harder-Narasimhan theory on MHP.

**Definition 3.4.** We say that a Hodge-Pink module $D$ is weakly admissible when it is semi-stable of slope 0 for the Harder-Narasimhan theory defined above. We denote by MHP$^\mathrm{wa}$ the full subcategory of weakly admissible Hodge-Pink modules. It is an abelian category.

We say that a Hodge-Pink module $D$ verifies the Griffiths transversality condition when

$$1 \otimes u \frac{d}{du}(V_D) \subset E^{-1} V_D.$$

We denote by MHP$^\mathrm{Gr}$ the full subcategory of Hodge-Pink modules verifying the Griffiths transversality condition, and MHP$^\mathrm{wa,Gr}$ the full subcategory of weakly admissible Hodge-Pink modules verifying the Griffiths transversality condition.

**Proposition 3.3.** *The category* MHP$^\mathrm{wa,Gr}$ *is a full subcategory of* MHP$^\mathrm{wa}$ *stable by subobjects.*

*Proof.* Let $D = (D, \sigma_D, V_D)$ be an object in MHP$^\mathrm{Gr,wa}$, $D' = (D', \sigma_{D'}, V_{D'})$ be a weakly admissible subobject of $D$ and $D'' = (D', \sigma_{D'}, V'_{D'})$ the image of $D'$ in the category MHP, i.e.

$$V'_{D'} = V_D \cap \sigma^* D' \otimes \hat{\mathfrak{S}}[\tfrac{1}{E}].$$

We have

$$\mathrm{t}_\mathrm{H}(D') \leq \mathrm{t}_\mathrm{H}(D'') \leq \mathrm{t}_\mathrm{N}(D'') = \mathrm{t}_\mathrm{N}(D') = \mathrm{t}_\mathrm{H}(D')$$

where the first inequality is due to the definition of $\mathrm{t}_\mathrm{H}$, the second inequality is given by the weakly admissibility of $D$ and the last equality is given by the weakly admissibility of $D'$. Thus $\mathrm{t}_\mathrm{H}(D') = \mathrm{t}_\mathrm{H}(D'')$, so $V_{D'} = V'_{D'}$. Now, by hypothesis we have $u\frac{d}{du}(V_D) \subset E^{-1} V_D$, so

$$\begin{aligned}
u\tfrac{d}{du}(V_{D'}) &= u\tfrac{d}{du}\left(V_D \cap \sigma^* D' \otimes_{K_0} \hat{\mathfrak{S}}[\tfrac{1}{E}]\right) \\
&\subset u\tfrac{d}{du}(V_D) \cap u\tfrac{d}{du}\left(\sigma^* D' \otimes_{K_0} \hat{\mathfrak{S}}[\tfrac{1}{E}]\right) \\
&\subset E^{-1}\left(V_D \cap \sigma^* D' \otimes_{K_0} \hat{\mathfrak{S}}[\tfrac{1}{E}]\right) \\
&= E^{-1} V_{D'},
\end{aligned}$$

therefore $D'$ verifies the Griffiths transversality condition. □

The following theorem gives us the relation between Hodge-Pink modules and Kisin modules. It can be found in [19, Theorem 0.4].

**Theorem 3.4.** *There is an equivalence of $\otimes$-categories*

$$\mathrm{Mod}^{\varphi}_{\mathfrak{S},\,\mathrm{fr}} \otimes \mathbb{Q}_p \simeq \mathrm{MHP}^{\mathrm{wa}}\,.$$

We denote by $\mathrm{Mod}^{\varphi,\log}_{\mathfrak{S},\mathrm{fr}}$ the full subcategory of $\mathrm{Mod}^{\varphi}_{\mathfrak{S},\,\mathrm{fr}}$ whose image in $\mathrm{MHP}^{\mathrm{wa}}$ lies in $\mathrm{MHP}^{\mathrm{wa},\mathrm{Gr}}$.

The following proposition gives us the relation between the Hodge types of the Kisin module and the Hodge-Pink module.

**Proposition 3.5.** *We have $\mathbf{t}_{\mathrm{H}}(M) = \mathbf{t}_{\mathrm{H}}(D)^{\iota}$.*

*Proof.* We have

$$\mathrm{Pos}\left(\sigma^* D \otimes_{K_0} \hat{\mathfrak{S}}, V_D\right) = \mathrm{Pos}\left(\xi(D \otimes_{K_0} \hat{\mathfrak{S}}), M \otimes_{\mathfrak{S}} \hat{\mathfrak{S}}\right) = \mathrm{Pos}\left(\varphi_M \otimes 1\left(\varphi^* M \otimes_{\mathfrak{S}} \hat{\mathfrak{S}}\right), M \otimes_{\mathfrak{S}} \hat{\mathfrak{S}}\right)$$

where $\xi$ is the isomorphism given in [19, Lemma 3.5] or [22, 1.2.6], the first equality is given by the construction of the equivalence of categories $\mathrm{Mod}^{\varphi}_{\mathfrak{S},\,\mathrm{fr}} \otimes \mathbb{Q}_p \simeq \mathrm{MHP}^{\mathrm{wa}}$ and the second equality is given in [22, Lemma 1.2.6]. □

### 3.2.1 Filtered isocrystals

**Definition 3.5.** A filtered isocrystal is a 3-tuple $(D, \sigma_D, \mathcal{F}_{\mathrm{H}})$ where:

- The pair $(D, \sigma_D)$ is an isocrystal,
- $\mathcal{F}_{\mathrm{H}}$ is a $\mathbb{Z}$-filtration on $D_K = D \otimes_{K_0} K$ by $K$-subspaces.

We denote by $\mathrm{MF}^{\sigma}_K$ the category of filtered isocrystals (with the filtration defined on $K$), whose objects are filtered isocrystals and the morphisms are morphisms between isocrystals $f : D \to D'$ such that $f(\mathcal{F}^{\geq i}_{\mathrm{H}} D_K) \subset \mathcal{F}^{\geq i}_{\mathrm{H}} D'_K$ for every $i \in \mathbb{Z}$. This is a quasi-abelian category. If $D'$ is a subobject of the underlying isocrystal of $D = (D, \sigma_D, \mathcal{F}_{\mathrm{H}})$, we can endow $D'$ with a Hodge filtration given by $\mathcal{F}^{\geq i}_{\mathrm{H}} D'_K = \mathcal{F}^{\geq i}_{\mathrm{H}} D_K \cap (D' \otimes_{K_0} K)$.

We can define the Newton and Hodge type by

$$\mathbf{t}_{\mathrm{N}}(D, \sigma_D, \mathcal{F}_{\mathrm{H}}) = \mathbf{t}_{\mathrm{N}}(D, \sigma_D) \quad \text{and} \quad \mathbf{t}_{\mathrm{H}}(D, \varphi_D, \mathcal{F}_{\mathrm{H}}) = \mathbf{t}(\mathcal{F}_{\mathrm{H}}).$$

We denote their degrees by

$$\mathrm{t}_{\mathrm{N}}(D) = \deg(\mathbf{t}_{\mathrm{N}}(D, \sigma_D, \mathcal{F}_{\mathrm{H}})) \quad \text{and} \quad \mathrm{t}_{\mathrm{H}}(D) = \deg(\mathcal{F}_{\mathrm{H}}) = \sum_{i \in \mathbb{Z}} i \cdot \dim_K(\mathrm{Gr}^{i}_{\mathcal{F}_{\mathrm{H}}}(D_K)).$$

It is easy to check (analogously to the Hodge-Pink modules case) that the usual rank function and the degree function given by $\deg(D) = \mathrm{t}_{\mathrm{H}}(D) - \mathrm{t}_{\mathrm{N}}(D)$ define a Harder-Narasimhan theory on $\mathrm{MF}^{\sigma}_K$.

**Definition 3.6.** We say that a filtered isocrystal $D$ is weakly admissible when it is semi-stable of slope 0 for the Harder-Narasimhan theory defined above. We denote by $^{\mathrm{wa}}\mathrm{MF}^{\sigma}_K$ the full subcategory of weakly admissible filtered isocrystals. It is an abelian category.

The category $^{\mathrm{wa}}\mathrm{MF}^{\sigma}_K$ is abelian and in [15], we see that there is a good Harder-Narasimhan formalism by taking for every object $D = (D, \varphi_D, \mathcal{F}_{\mathrm{H}})$

$$\deg D = -\mathrm{t}_{\mathrm{H}}(D) \quad \text{and} \quad \mathrm{rank}\, D = \dim_{K_0} D.$$

We denote by $\mathcal{F}_{\mathrm{F,wa}}$ the Harder-Narasimhan filtration given by these degree and rank functions, and $\mathbf{t}_{\mathrm{F,wa}}$ will be its polygon.

There is a functor

$$\begin{array}{ccc} \mathrm{MHP} & \to & \mathrm{MF}^{\sigma}_K \\ D_{\mathrm{HP}} = (D, \sigma_D, V_D) & \mapsto & D_{\mathrm{Fil}} = (D, \sigma_D, \mathcal{F}_{\mathrm{H}}) \end{array}$$

where $\mathcal{F}_{\mathrm{H}} = \mathcal{F}(D \otimes_{K_0} \hat{\mathfrak{S}}, \sigma_D \otimes 1(V_D))$, since both arguments are $\hat{\mathfrak{S}}$-lattices in $D \otimes \hat{\mathfrak{S}}[\frac{1}{E}]$, so we obtain a filtration on the reduction of $D \otimes_{K_0} \hat{\mathfrak{S}}$, which is $D \otimes_{K_0} K$. It is clear that this is a $\otimes$-functor. The following theorem is given by Genestier and Lafforgue in [19, Lemma 1.3] and [19, Lemma 1.4].

**Theorem 3.6.** *This construction yields an equivalence of $\otimes$-categories*
$$\mathrm{MHP}^{\mathrm{Gr}} \xrightarrow{\sim} \mathrm{MF}_K^{\sigma}$$
*inducing an equivalence of $\otimes$-categories ${}^{\mathrm{wa}}\mathrm{MF}_K^{\sigma} \xrightarrow{\sim} \mathrm{MHP}^{\mathrm{wa},\mathrm{Gr}}$.*

The next proposition gives us the comparison between the Hodge types defined on the Hodge-Pink module and the filtered isocrystal.

**Proposition 3.7.** *Let $D_{\mathrm{HP}}$ be a Hodge-Pink module whose image by the functor above is $D_{\mathrm{Fil}}$. We have*
$$\mathbf{t}_{\mathrm{H}}(D_{\mathrm{HP}}) = \mathbf{t}_{\mathrm{H}}(D_{\mathrm{Fil}}).$$

*Proof.* We have $\mathbf{t}_{\mathrm{H}}(D_{\mathrm{HP}}) = \mathrm{Pos}\left(\sigma^*D \otimes_{K_0} \hat{\mathfrak{S}}, V_D\right) = \mathrm{Pos}\left(D \otimes_{K_0} \hat{\mathfrak{S}}, \sigma_D \otimes 1(V_D)\right) = \mathbf{t}_{\mathrm{H}}(D_{\mathrm{Fil}})$. $\square$

To summarize, we get a diagram

$$\begin{array}{ccccc}
\mathrm{Mod}_{\mathfrak{S}}^{\varphi} \otimes \mathbb{Q}_p & \xrightarrow{\simeq} & \mathrm{MHP}^{\mathrm{wa}} & \hookrightarrow & \mathrm{MHP} \\
\uparrow & & \uparrow & & \uparrow \\
\mathrm{Mod}_{\mathfrak{S},\mathrm{fr}}^{\varphi,\mathrm{log}} \otimes \mathbb{Q}_p & \xrightarrow{\simeq} & \mathrm{MHP}^{\mathrm{wa},\mathrm{Gr}} & \hookrightarrow & \mathrm{MHP}^{\mathrm{Gr}} \\
& \underset{\tilde{\Theta}}{\overset{\alpha}{\leftrightarrows}} & \downarrow{\simeq} & & \downarrow{\simeq} \\
& & {}^{\mathrm{wa}}\mathrm{MF}_K^{\sigma} & \hookrightarrow & \mathrm{MF}_K^{\sigma}
\end{array}$$

where the functor
$$\alpha \;:\; \mathrm{Mod}_{\mathfrak{S},\mathrm{fr}}^{\varphi,\mathrm{log}} \otimes \mathbb{Q}_p \to {}^{\mathrm{wa}}\mathrm{MF}_K^{\sigma}$$
is described by associating to $(M, \varphi_M) \otimes \mathbb{Q}_p$ the filtered isocrystal $(M(0), \varphi_{M(0)}, \mathcal{F}_{\mathrm{H}})$ where

- The underlying module is $M(0) = (M/uM)[\frac{1}{p}]$,
- The Frobenius $\varphi_{M(0)}$ is the reduction of $\varphi_M$ to $M(0)$,
- The filtration on $M(0)_K = M(0) \otimes_{K_0} K$ is given, via the isomorphism $M(0) \otimes_{K_0} K \xrightarrow{\sim} \varphi_M \varphi^* M \otimes_{\hat{\mathfrak{S}}} K$, by the filtration
$$\mathcal{F}(\varphi_M \varphi^* M \otimes \hat{\mathfrak{S}}, M \otimes \hat{\mathfrak{S}}) = \frac{\varphi_M \varphi^* M \cap E(u)^i M + E(u) \varphi_M \varphi^* M}{E(u) \varphi_M \varphi^* M} \subset \varphi_M \varphi^* M \otimes_{\hat{\mathfrak{S}}} K,$$

and
$$\tilde{\Theta} \;:\; {}^{\mathrm{wa}}\mathrm{MF}_K^{\sigma} \to \mathrm{Mod}_{\mathfrak{S},\mathrm{fr}}^{\varphi,\mathrm{log}} \otimes \mathbb{Q}_p$$

is its inverse, the functor defined in [22]. From the diagram above and the results for Hodge-Pink modules, we get:

**Proposition 3.8.** *The category $\mathrm{Mod}_{\mathfrak{S},\mathrm{fr}}^{\varphi,\mathrm{log}} \otimes \mathbb{Q}_p$ is a full subcategory of $\mathrm{Mod}_{\mathfrak{S}}^{\varphi} \otimes \mathbb{Q}_p$ stable by subobjects.*

As a consequence, we get the compatibility for the filtrations, as we did in 3.2thm.3.2.

**Proposition 3.9.** *Let $D \in {}^{\mathrm{wa}}\mathrm{MF}_K^{\sigma}$, then for every $\gamma \in \mathbb{Q}$, we have*
$$\mathcal{F}_{\mathrm{F,o}}^{\geq \gamma}(\widetilde{\Theta}(D)) = \widetilde{\Theta}(\mathcal{F}_{\mathrm{F,wa}}^{\geq \gamma}(D)).$$

### 3.2.2 Fontaine's functors

Using the same notations from last section, let $K$ be a finite totally ramified extension of $K_0$ contained in $\overline{K}_0$, $\pi_K$ an uniformizer of $K$, $K_\infty = \cup_{n \geq 1} K(\sqrt[p^n]{\pi_K})$, $\mathrm{Gal}_K = \mathrm{Gal}(\overline{K}_0/K)$ and $\mathrm{Gal}_{K_\infty} = \mathrm{Gal}(\overline{K}_0/K_\infty)$. We denote by $\mathrm{Rep}_{\mathbb{Q}_p}^{\mathrm{cr}} \mathrm{Gal}_K$ the category of crystalline $\mathrm{Gal}_K$-representations, and by $\mathrm{Rep}_{\mathbb{Z}_p}^{\mathrm{cr}} \mathrm{Gal}_K$ the category of $\mathrm{Gal}_K$-stable $\mathbb{Z}_p$-lattices $L$ such that $L \otimes \mathbb{Q}_p \in \mathrm{Rep}_{\mathbb{Q}_p}^{\mathrm{cr}} \mathrm{Gal}_K$.

For a crystalline representation $V$, Fontaine functor $D_{\mathrm{cris}}(V) = (B_{\mathrm{cris}} \otimes_{\mathbb{Q}_p} V)^{\mathrm{Gal}_K}$, defined in [17, 2.3.3], where $D_{\mathrm{cris}}(V)$ is an object in $\mathrm{MF}_K^{\varphi}$, and it is compatible with the change of the extension $K$ of $K_0$. We get an equivalence of $\otimes$-categories
$$D_{\mathrm{cris}} : \mathrm{Rep}_{\mathbb{Q}_p}^{\mathrm{cr}} \mathrm{Gal}_K \xrightarrow{\simeq} {}^{\mathrm{wa}}\mathrm{MF}_K^{\sigma},$$

whose inverse is denoted by $V_{\text{cris}}$. Then, we also get a functor

$$\mathfrak{N} = \tilde{\Theta} \circ D_{\text{cris}} \ : \ \text{Rep}_{\mathbb{Q}_p}^{\text{cr}} \text{Gal}_K \to \text{Mod}_{\mathfrak{S}}^{\varphi} \otimes \mathbb{Q}_p.$$

By construction, there is a canonical isomorphism of $\otimes$-functors

$$\alpha \circ \mathfrak{N} \simeq D_{\text{cris}} \ : \ \text{Rep}_{\mathbb{Q}_p} \text{Gal}_K \to^{\text{wa}} \text{MF}_K^{\sigma}.$$

Denote by $\mathcal{O}_{\mathcal{E}}$ the $p$-adic completion of $\mathfrak{S}_{(p)}$, a discrete valuation ring with fraction field $\mathcal{E} = \text{Frac}\,\mathcal{O}_{\mathcal{E}}$ and residue field $\mathbb{F}((u))$. Let $\text{Mod}_{\mathcal{E}}^{\varphi,\text{ét}}$ be the category whose objects are finite dimensional $\mathcal{E}$-vector spaces $M$ together with a Frobenius isomorphism $\varphi_M : \varphi^* M \xrightarrow{\sim} M$, and the morphisms are morphisms between $\mathcal{E}$-vector spaces compatible with the Frobenius. Also, let $\text{Mod}_{\mathcal{O}_{\mathcal{E}}}^{\varphi,\text{ét}}$ be the category whose objects are finite free $\mathcal{O}_{\mathcal{E}}$-modules together with a Frobenius isomorphism as above. Then, Fontaine gives in [17] two equivalences of $\otimes$-categories

$$\text{Mod}_{\mathcal{O}_{\mathcal{E}}}^{\varphi,\text{ét}} \xrightarrow{\sim} \text{Rep}_{\mathbb{Z}_p} \text{Gal}_{K_{\infty}} \quad \text{and} \quad \text{Mod}_{\mathcal{E}}^{\varphi,\text{ét}} \xrightarrow{\sim} \text{Rep}_{\mathbb{Q}_p} \text{Gal}_{K_{\infty}}$$

which are exact. These two isomorphism induce a diagram

$$\begin{array}{ccccc}
\text{Mod}_{\mathfrak{S},\text{fr}}^{\varphi,\log} & \xrightarrow{-\otimes \mathcal{O}_{\mathcal{E}}} & \text{Mod}_{\mathcal{O}_{\mathcal{E}}}^{\varphi,\text{ét}} & \xrightarrow{\simeq} & \text{Rep}_{\mathbb{Z}_p} \text{Gal}_{K_{\infty}} \\
\downarrow^{-\otimes \mathbb{Q}_p} & & & & \downarrow^{-\otimes \mathbb{Q}_p} \\
\text{Mod}_{\mathfrak{S},\text{fr}}^{\varphi,\log} \otimes \mathbb{Q}_p & \xrightarrow{-\otimes \mathcal{E}} & \text{Mod}_{\mathcal{E}}^{\varphi,\text{ét}} & \xrightarrow{\simeq} & \text{Rep}_{\mathbb{Q}_p} \text{Gal}_{K_{\infty}} \\
\| & & & & \uparrow \\
\text{Mod}_{\mathfrak{S},\text{fr}}^{\varphi,\log} \otimes \mathbb{Q}_p & \xrightarrow{\alpha} & ^{\text{wa}}\text{MF}_K^{\sigma} & \xrightarrow{\simeq} & \text{Rep}_{\mathbb{Q}_p}^{\text{cr}} \text{Gal}_K
\end{array}$$

where the upper diagram is strictly commutative by construction and Kisin proves in [24, Proposition 2.1.5] that the lower diagram is also commutative, i.e. the two functors $\text{Mod}_{\mathfrak{S},\text{fr}}^{\varphi,\log} \otimes \mathbb{Q}_p \to \text{Rep}_{\mathbb{Q}_p} \text{Gal}_{K_{\infty}}$ are canonically isomorphic.

### 3.2.3 Kisin's functor

Using the Fontaine's functors, we obtain an exact and fully faithful $\otimes$-functor $\text{Mod}_{\mathfrak{S},\text{fr}}^{\varphi,\log} \to \text{Rep}_{\mathbb{Z}_p}^{\text{cr}}(\text{Gal}_K, \text{Gal}_{K_{\infty}})$, which is an equivalence of categories by [22, 2.1.15], where $\text{Rep}_{\mathbb{Z}_p}^{\text{cr}}(\text{Gal}_K, \text{Gal}_{K_{\infty}})$ denotes the category whose objects are the couples $(T, V)$ such that $T \in \text{Rep}_{\mathbb{Z}_p} \text{Gal}_{K_{\infty}}$ is a lattice in $V \in \text{Rep}_{\mathbb{Q}_p}^{\text{cr}} \text{Gal}_K$. We denote by $\mathfrak{M}$ its inverse and also the fully faithful $\otimes$-functor

$$\mathfrak{M} \ : \ \text{Rep}_{\mathbb{Z}_p}^{\text{cr}} \text{Gal}_K \to \text{Mod}_{\mathfrak{S},\,\text{fr}}^{\varphi}$$

which is obtained by pre and post composing this inverse with the fully faithful embeddings $\text{Rep}_{\mathbb{Z}_p}^{\text{cr}} \text{Gal}_K \hookrightarrow \text{Rep}_{\mathbb{Z}_p}^{\text{cr}}(\text{Gal}_K, \text{Gal}_{K_{\infty}})$ and $\text{Mod}_{\mathfrak{S},\text{fr}}^{\varphi,\log} \hookrightarrow \text{Mod}_{\mathfrak{S},\text{fr}}^{\varphi}$. The functor $\mathfrak{M}$ is compatible with the formation of symmetric and exterior powers, by [24, Theorem 1.2.1]. By construction,

- We have $\mathfrak{M} \otimes \mathbb{Q}_p \simeq \mathfrak{N}(- \otimes \mathbb{Q}_p)$. In particular, $\mathfrak{M} \otimes \mathbb{Q}_p$ and $\mathfrak{M}[\frac{1}{p}]$ are exact.

- The functor $\mathfrak{M} \otimes \mathcal{O}_{\mathcal{E}}$ is isomorphic to Fontaine's functor $\text{Rep}_{\mathbb{Z}_p}^{\text{cr}} \text{Gal}_K \to \text{Rep}_{\mathbb{Z}_p} \text{Gal}_{K_{\infty}} \to \text{Mod}_{\mathcal{O}_{\mathcal{E}}}^{\varphi,\text{ét}}$. In particular, it is exact and so is $\mathfrak{M}_{(p)} = \mathfrak{M} \otimes_{\mathfrak{S}} \mathfrak{S}_{(p)}$.

*Remark* 6. It follows that for every $\mathfrak{p} \in \text{Spec}\,\mathfrak{S} \setminus \{\mathfrak{m}\}$, the localized functor $\mathfrak{M}_{\mathfrak{p}}$ is exact. For an exact sequence $0 \to L_1 \to L_2 \to L_3 \to 0$ in $\text{Rep}_{\mathbb{Z}_p}^{\text{cr}} \text{Gal}_K$, the sequence

$$0 \to \mathfrak{M}(L_1) \to \mathfrak{M}(L_2) \to \mathfrak{M}(L_3)$$

is exact and the cokernel of $\mathfrak{M}(L_2) \to \mathfrak{M}(L_3)$ has finite length.

In the sequel, a $\mathbb{Z}$-filtration on an $A$-module $M$ is a decreasing collection $\mathcal{F} = (\mathcal{F}^{\geq i})_{i \in \mathbb{Z}}$ of direct summands of $M$ such that $\mathcal{F}^{\geq i} = M$ for $i \ll 0$ and $\mathcal{F}^{\geq i} = 0$ for $i \gg 0$. Any such filtration as a (non-unique) splitting: a $\mathbb{Z}$-graduation $M = \oplus_{i \in \mathbb{Z}} \mathcal{G}_i$ such that $\mathcal{F}^{\geq i} = \oplus_{j \geq i} \mathcal{G}_j$ for all $i \in \mathbb{Z}$. For $M \in \text{Mod}_{\mathfrak{S},\text{fr}}$ and a $\mathbb{Z}$-filtration $\mathcal{F} = (\mathcal{F}^{\geq i})_{i \in \mathbb{Z}}$ on $N = M[\frac{1}{p}]$, we set

$$M + \mathcal{F} = \sum_{i \in \mathbb{Z}} p^{-i} M \cap \mathcal{F}^{\geq i}.$$

Since $M + \mathcal{F} = \sum_{i=-r}^{r} p^{-i}M \cap \mathcal{F}^{\geq i}$ for $r \gg 0$, this is a finitely generated $\mathfrak{S}$-submodule of $N$, and so is therefore also
$$M +_{\text{fr}} \mathcal{F} = (M + \mathcal{F})_{\text{fr}}.$$
The following conditions on $(M, \mathcal{F})$ are equivalent:

1. We have $\mathcal{F} \cap M = (\mathcal{F}^{\geq i} \cap M)_{i \in \mathbb{Z}}$ is a $\mathbb{Z}$-filtration on $M$,

2. There is a $\mathbb{Z}$-filtration $\mathcal{F}' = (\mathcal{F}'^{\geq i})_{i \in \mathbb{Z}}$ on $M$ such that $\mathcal{F}'[\frac{1}{p}] = \mathcal{F}$.

Indeed $(1) \Rightarrow (2)$ with $\mathcal{F}' = \mathcal{F} \cap M$ and $(2) \Rightarrow (1)$ since then $\mathcal{F} \cap M = \mathcal{F}'$. They are also equivalent to:

3. For each $i \in \mathbb{Z}$, $M/\mathcal{F}^{\geq i} \cap M$ is free over $\mathfrak{S}$.

Then
$$M + \mathcal{F} = M +_{\text{fr}} \mathcal{F}.$$
Indeed, let $M = \oplus \mathcal{G}_i$ be a splitting of $\mathcal{F} \cap M$. Then each $\mathcal{G}_i$ is a free $\mathfrak{S}$-module and
$$\begin{aligned}
p^{-i}M \cap \mathcal{F}^{\geq i} &= \left(\oplus_{j \in \mathbb{Z}} p^{-i}\mathcal{G}_j\right) \cap \left(\oplus_{j \geq i} \mathcal{G}_j[\tfrac{1}{p}]\right) \\
&= \oplus_{j \geq i} p^{-i}\mathcal{G}_j \\
\text{thus}\quad M + \mathcal{F} &= \sum_{i \in \mathbb{Z}} \left(\oplus_{j \geq i} p^{-i}\mathcal{G}_j\right) \\
&= \oplus_{j \in \mathbb{Z}} \left(\sum_{j \geq i} p^{-i}\mathcal{G}_j\right) \\
&= \oplus_{j \in \mathbb{Z}} p^{-j}\mathcal{G}_j
\end{aligned}$$
is already a free $\mathfrak{S}$-submodule of $N$.

If now $M \in \text{Mod}^{\varphi}_{\mathfrak{S},\,\text{fr}}$ and $\mathcal{F}$ is a $\varphi_N$-stable $\mathbb{Z}$-filtration on $N = M[\frac{1}{p}] \in \text{Mod}^{\varphi}_{\mathfrak{S}[\frac{1}{p}]}$, then $M + \mathcal{F}$ and $(M + \mathcal{F})_{\text{fr}}$ are $\varphi_N$-stable $\mathfrak{S}$-submodules of $N$. We may thus view them as objects in $\text{Mod}^{\varphi}_{\mathfrak{S}}$ and $\text{Mod}^{\varphi}_{\mathfrak{S},\,\text{fr}}$, respectively.

**Proposition 3.10.** *Let $L \in \text{Rep}^{\text{cr}}_{\mathbb{Z}_p} \text{Gal}_K$, let $V = L \otimes \mathbb{Q}_p$, and let $\mathcal{F}$ be a $\text{Gal}_K$-stable $\mathbb{Z}$-filtration on $V$. Then we have*
$$\mathfrak{M}(L + \mathcal{F}) = \mathfrak{M}(L) +_{\text{fr}} \mathfrak{N}(\mathcal{F}).$$

*Proof.* Plainly, $L + \mathcal{F}$ is a $\text{Gal}_K$-stable sub-lattice of $V$. Applying Kisin's functor $\mathfrak{M}$ and $\mathfrak{N}$ to
$$\begin{array}{ccccccccc}
0 & \to & p^{-i}L \cap \mathcal{F}^{\geq i} & \to & p^{-i}L & \to & p^{-i}L/p^{-i}L \cap \mathcal{F}^{\geq i} & \to & 0 \\
& & \cap & & \cap & & \cap & & \\
0 & \to & \mathcal{F}^{\geq i} & \to & V & \to & V/\mathcal{F}^{\geq i} & \to & 0
\end{array}$$
we obtain a commutative diagram with exact rows
$$\begin{array}{ccccccccc}
0 & \to & \mathfrak{M}(p^{-i}L \cap \mathcal{F}^{\geq i}) & \to & \mathfrak{M}(p^{-i}L) & \to & \mathfrak{M}(p^{-i}L/p^{-i}L \cap \mathcal{F}^{\geq i}) & & \\
& & \cap & & \cap & & \cap & & \\
0 & \to & \mathfrak{N}(\mathcal{F}^{\geq i}) & \to & \mathfrak{N}(V) & \to & \mathfrak{N}(V/\mathcal{F}^{\geq i}) & \to & 0
\end{array}$$
Since also $\mathfrak{M}\left(p^{-i}L\right) = p^{-i}\mathfrak{M}(L)$ in $\mathfrak{N}(V)$, we obtain
$$\begin{aligned}
\mathfrak{M}\left(p^{-i}L \cap \mathcal{F}^{\geq i}\right) &= \mathfrak{M}\left(p^{-i}L\right) \cap \mathfrak{N}\left(\mathcal{F}^{\geq i}\right) \\
&= p^{-i}\mathfrak{M}(L) \cap \mathfrak{N}\left(\mathcal{F}^{\geq i}\right)
\end{aligned}$$
inside $\mathfrak{N}(V)$. Fix $r \geq 0$ such that $\mathcal{F}^{\geq -r} = V$ and $\mathcal{F}^{\geq r} = 0$. Applying Kisin's (additive) functors $\mathfrak{M}$ and $\mathfrak{N}$ to the diagram
$$\begin{array}{ccccc}
\oplus_{i=-r}^{r} p^{-i}L \cap \mathcal{F}^{\geq i} & \longrightarrow & L + \mathcal{F} & \longrightarrow & 0 \\
\cap & & \cap & & \\
\oplus_{i=-r}^{r} \mathcal{F}^{\geq i} & \longrightarrow & V & \longrightarrow & 0
\end{array}$$

we thus obtain a commutative diagram

$$\begin{array}{ccc} \oplus_{i=-r}^{r} p^{-i}\mathfrak{M}(L) \cap \mathfrak{N}(\mathcal{F}^{\geq i}) & \longrightarrow & \mathfrak{M}(L+\mathcal{F}) \\ \cap & & \cap \\ \oplus_{i=-r}^{r} \mathfrak{N}(\mathcal{F}^{\geq i}) & \longrightarrow & \mathfrak{N}(V) & \longrightarrow & 0 \end{array}$$

The image of the top map is the $\mathfrak{S}$-submodule $\mathfrak{M}(L) + \mathfrak{N}(\mathcal{F})$ of $\mathfrak{N}(V)$ and its cokernel has finite length, thus indeed

$$(\mathfrak{M}(L) +_{\text{fr}} \mathfrak{N}(\mathcal{F})) = (\mathfrak{M}(L) + \mathfrak{N}(\mathcal{F}))_{\text{fr}} = \mathfrak{M}(L+\mathcal{F}).$$

$\square$

### 3.2.4 Harder-Narasimhan filtrations

Denote by $C = \hat{\overline{K}}$, the completion of an algebraic closure $\overline{K}$ of $K$. In [15], Fargues constructs a Harder-Narasimhan filtration $\mathcal{F}_{\text{F,cr}}$ on $\text{Rep}_{\mathbb{Q}_p}^{\text{HT}} \text{Gal}_K$, the abelian category of Hodge-Tate $\mathbb{Q}_p$-representations of $\text{Gal}_K$, by considering the dimension as rank function and the degree function defined by $\deg(V) = d$ where $\wedge^{\text{rank } V} V_C \simeq C(d)$, for $V_C = V \otimes_{\mathbb{Q}_p} C$. Then, he construct a Harder-Narasimhan filtration on $\text{Rep}_{\mathbb{Q}_p}^{\text{cr}} \text{Gal}_K$ via the fully faithful functor $\text{Rep}_{\mathbb{Q}_p}^{\text{cr}} \text{Gal}_K \to \text{Rep}_{\mathbb{Q}_p}^{\text{HT}} \text{Gal}_K$, since every subrepresentation of a crystalline representation is crystalline.

Consider the category $\text{VectFil}_{C/\mathbb{Q}_p}$ of $\mathbb{Q}_p$-vector spaces endowed with a filtration on $V_C$. This is a quasi-abelian category admitting a Harder-Narasimhan formalism for the rank function defined as the dimension of the $\mathbb{Q}_p$-vector space and the degree function defined by $\deg(V, \text{Fil}^\bullet V_C) = \sum_{i \in \mathbb{Z}} i \cdot \dim_C(\text{Gr}^i_{\text{Fil}^\bullet} V_C)$. Moreover, the Harder-Narasimhan filtration obtained is compatible with tensor products. Fargues also shows that there is an exact and faithful $\otimes$-functor $\mathcal{G} : \text{Rep}_{\mathbb{Q}_p}^{\text{HT}} \text{Gal}_K \to \text{VectFil}_{C/\mathbb{Q}_p}$ such that $\deg V = \deg(\mathcal{G}(V))$ and $\text{rank } V = \text{rank}(\mathcal{G}(V))$.

To summarize, we have the following relations between all the categories presented before:

$$\begin{array}{ccccccc} \text{MHP}^{\text{wa,Gr}} & \xrightarrow{\simeq} & \text{Mod}_{\mathfrak{S}}^{\varphi,\log} \otimes \mathbb{Q}_p & \hookrightarrow & \text{Mod}_{\mathfrak{S}}^{\varphi} \otimes \mathbb{Q}_p & \hookrightarrow & \text{Mod}_{\mathfrak{S}[\frac{1}{p}]}^{\varphi} \\ \downarrow \simeq & & & & & & \\ {}^{\text{wa}}\text{MF}_K^{\sigma} & \xrightarrow{\simeq} & \text{Rep}_{\mathbb{Q}_p}^{\text{cr}} \text{Gal}_K & \hookrightarrow & \text{Rep}_{\mathbb{Q}_p}^{\text{HT}} \text{Gal}_K & \hookrightarrow & \text{VectFil}_{C/\mathbb{Q}_p} \end{array}$$

where all the functors are exact and $\otimes$-functors. Also, there is a Harder-Narasimhan filtration defined in every category above (the filtration in $\text{MHP}^{\text{wa,Gr}}$ is defined by the equivalence of categories with $\text{Mod}_{\mathfrak{S}}^{\varphi,\log} \otimes \mathbb{Q}_p$).

By the results in previous sections and the results given by Fargues in [15], we have:

**Proposition 3.11.** *All the Harder-Narasimhan filtrations defined on the categories above are compatible.*

### 3.2.5 Germs of crystalline representations

Assume $\mathbb{F}$ algebraically closed.

**Definition 3.7.** For a $\mathbb{Q}_p$-vector space $V$, denote by $\text{Cr}(V, K)$ the set of all morphisms $\rho : \text{Gal}_K \to \text{GL}(V)$ such that $(V, \rho)$ is a crystalline representation of $\text{Gal}_K$, and set

$$\text{Cr}(V) := \varinjlim_{K_0 \subset K \subset \overline{K}_0} \text{Cr}(V, K)$$

where the transition maps $\text{Cr}(V, K) \to \text{Cr}(V, K')$ are induced by the inclusion $\text{Gal}_{K'} \hookrightarrow \text{Gal}_K$, for finite (totally ramified) extensions $K_0 \subset K \subset K' \subset \overline{K}_0$. For an element $\rho \in \text{Cr}(V)$, we denote by $\rho_K$ the morphism $\rho_K : \text{Gal}_K \to \text{GL}(V)$, for $K$ large enough.

We define a germ of crystalline representations of $\text{Gal}_{K_0}$ as an object $(V, \rho)$ where $V$ is a finite dimensional $\mathbb{Q}_p$-vector space and $\rho \in \text{Cr}(V)$. A morphism $f : (V, \rho) \to (V', \rho')$ of germs of representations is a $\mathbb{Q}_p$-linear morphism $f : V \to V'$ verifying that for every $K$ and $K'$ such that $\rho = \rho_K$ and $\rho' = \rho_{K'}$, there exists a finite extension $K_0 \subset K'' \subset \overline{K}_0$ containing $K$ and $K'$ such that $f \circ \rho_{K|K''} = \rho_{K'|K''} \circ f$, where $\rho_{K|K''}$ and $\rho_{K'|K''}$ are the restrictions of $\rho_K$ and $\rho_{K'}$ to $K''$. We define the germs of integral crystalline representations of $\text{Gal}_{K_0}$ analogously.

The category of germs of crystalline representations and the category of germs of integral crystalline representations are denoted by $\text{Rep}_{\mathbb{Q}_p}^{\text{cr}}\{\text{Gal}_{K_0}\}$ and $\text{Rep}_{\mathbb{Z}_p}^{\text{cr}}\{\text{Gal}_{K_0}\}$, respectively.

As a consequence of Proposition 13 in [15], we have next proposition which gives us the necessary compatibility between them in order to get a Harder-Narasimhan filtration on $\text{Rep}_{\mathbb{Q}_p}^{\text{cr}}\{\text{Gal}_K\}$.

**Proposition 3.12.** *The Harder-Narasimhan filtration defined on objects in $\text{Rep}_{\mathbb{Q}_p}^{\text{cr}} \text{Gal}_K$ is compatible with base change. In particular, it defines a Harder-Narasimhan filtration on objects in $\text{Rep}_{\mathbb{Q}_p}^{\text{cr}}\{\text{Gal}_{K_0}\}$.*

For a germ of crystalline representations $(V, \rho)$, there always exists a large enough finite extension $K_0 \subset K \subset \overline{K}_0$ such that $(V, \rho_K) \in \text{Rep}_{\mathbb{Q}_p} \text{Gal}_K$. Then, we can consider Fontaine's $D_{\text{cris}}$ functor to get a functor
$$\text{Rep}_{\mathbb{Q}_p}^{\text{cr}}\{\text{Gal}_{K_0}\} \to \text{Vect}_{K_0}^\sigma$$
(which is just the composition of $D_{\text{cris}}$ with the functor forgetting the filtration) and by Fontaine's construction, this functor does not depend on $K$.

We would like to define the $D_{\text{cris}}$ functor on germs of integral crystalline representations using Kisin's construction but, in order to do that, we have to prove that this construction is independent of the choice of $K$. For a finite extension $K_0 \subset K \subset \overline{K}_0$ with uniformizer $\pi_K$. Let $\eta$ be the $\otimes$-isomorphism between Kisin's $\otimes$-functor
$$D_{\text{cris}}^K \; : \; \text{Rep}_{\mathbb{Q}_p}^{\text{cr}} \text{Gal}_K \xrightarrow{\mathfrak{M}} \text{Mod}_{\mathfrak{S},\text{fr}}^{\varphi,\log} \otimes \mathbb{Q}_p \xrightarrow{\alpha} \text{MF}_K^\sigma$$
and Fontaine's $D_{\text{cris}}$ functor. Define $\otimes$-functors

$$\begin{array}{rcl} D_{\text{cris}}^K \; : \; \text{Rep}_{\mathbb{Z}_p}^{\text{cr}} \text{Gal}_K & \to & \text{Mod}_{W(\mathbb{F})}^\sigma, \\ L & \mapsto & \mathfrak{M}(L)/u\mathfrak{M}(L) \end{array} \quad \text{and} \quad \begin{array}{rcl} D_{\text{cris}} \; : \; \text{Rep}_{\mathbb{Z}_p}^{\text{cr}} \text{Gal}_K & \to & \text{Mod}_{W(\mathbb{F})}^\sigma \\ L & \mapsto & \eta_{L \otimes \mathbb{Q}_p}(D_{\text{cris}}^K(L)) \end{array}$$

so that $D_{\text{cris}}^K(L)$ is a lattice in $D_{\text{cris}}^K(L \otimes \mathbb{Q}_p)$ and $D_{\text{cris}}(L)$ is a lattice in $D_{\text{cris}}(L \otimes \mathbb{Q}_p)$.

**Lemma 3.13.** *This construction induces a $\otimes$-functor $D_{\text{cris}} \; : \; \text{Rep}_{\mathbb{Z}_p}^{\text{cr}}\{\text{Gal}_{K_0}\} \to \text{Mod}_{W(\mathbb{F})}^\sigma$.*

*Proof.* It follows from the more general result given by Liu in [32, Proposition 2.2.4] for semi-stable representations, knowing that in the case of a crystalline representation we have $D_{\text{st}}(V) = D_{\text{cris}}(V)$. $\square$

## 3.3 The Fargues filtration on $\text{Mod}_{\mathbb{F}[[u]],\text{fr}}^\varphi$

We have already seen that $\text{Mod}_{\mathbb{F}[[u]],\text{fr}}^\varphi$ is an abelian category. In order to define a Harder-Narasimhan filtration on $\text{Mod}_{\mathbb{F}[[u]],\text{fr}}^\varphi$, we define the rank as
$$\text{rank}(M) = \text{rank}_{\mathbb{F}[[u]]}(M).$$

*Remark* 7. We see that another description of the rank function is given by $\text{rank}(M) = \mu_{\text{IW}}(M)$, for $M \in \text{Mod}_{\mathbb{F}[[u]],\text{fr}}^\varphi$ viewed as an $\mathfrak{S}$-module.

We define the degree function by
$$\deg(M) = \nu(M, \varphi_M(\varphi^*M)),$$
for $\nu$ the operator defined in 2.3Latticessubsection.2.3. When $M$ is effective, we have $\deg M = -\text{length}_{\mathbb{F}[[u]]} Q$ for $Q = M/\varphi_M(\varphi^*M)$ the cokernel of $\varphi_M$ in the category of $\mathbb{F}[[u]]$-modules. Also, define
$$\mu(M) = \frac{\deg M}{\text{rank } M}.$$

It is easy to see that rank and deg are a rank and degree function in the sense of the Harder-Narasimhan formalism. They thus define a Harder-Narasimhan filtration on every object of $\text{Mod}_{\mathbb{F}[[u]],\text{fr}}^\varphi$. This filtration was already announced by Carl Wang-Erickson and Brandon Levin in [14] and was originally inspired by Fargues' theory [16] of Harder-Narasimhan filtrations for finite flat group schemes. This is why, from now on, we will refer to our filtration as the the Fargues filtrations. We will denote by $\mathcal{F}_{\text{F},1}$ the Fargues filtration and by $\mathbf{t}_{\text{F},1}$, the polygon associated to it.

Levin and Wang Erickson prove in [14] that this filtration is compatible with tensor products. Using the tensor compatibility and a calculation with slopes, we obtain that it is also compatible with symmetric and exterior powers.

We define the Hodge polygon of a $p$-torsion Kisin module $M$ by

$$\mathbf{t}_{\mathrm{H}}(M) = \mathrm{Pos}(M, \varphi_M(\varphi^* M)).$$

**Proposition 3.14.** *We have $\mathbf{t}_{\mathrm{F}}(M) \leq \mathbf{t}_{\mathrm{H}}(M)$.*

*Proof.* Since we have $\mathbf{t}_{\mathrm{F}}(M) = \mathbf{t}_{\mathrm{F}}(\mathrm{Gr}_{\mathcal{F}_{\mathrm{F}}} M) = \mathbf{t}_{\mathrm{F}}(\mathrm{Gr}_{\mathcal{F}_{\mathrm{F}}}^{\gamma_1} M) * \ldots * \mathbf{t}_{\mathrm{F}}(\mathrm{Gr}_{\mathcal{F}_{\mathrm{F}}}^{\gamma_r} M)$ where $\{\gamma_1, \ldots, \gamma_r\}$ are the breaks in the Fargues filtration of $M$, it suffices to prove it for $M$ semi-stable the statement is true. For $M$ semi-stable, it is verified because the Fargues polygon of a semi-stable module has only one slope and both polygons have the same terminal points. □

## 3.4 The Fargues filtration on $\mathrm{Mod}^{\varphi}_{\mathfrak{S},\mathrm{t}}$

We have already seen that the category $\mathrm{Mod}^{\varphi}_{\mathfrak{S},\mathrm{t}}$ is a quasi-abelian category. To give a mono-epi in this category, is the same as to give a monomorphism $M_1 \hookrightarrow M_2$ such that $M_1\left[\frac{1}{u}\right] \xrightarrow{\sim} M_2\left[\frac{1}{u}\right]$, which is the same as to say that we have an exact sequence in the category $\mathcal{A}$ defined in 2.5.3Some categories of $\mathfrak{S}$-modulessubsubsection.2.5.3 of torsion $\mathfrak{S}$-modules $0 \to M_1 \to M_2 \to Q \to 0$ with $Q$ an object in $\mathcal{T}$, i.e. a finite length $\mathfrak{S}$-module.

We define the $i$-twist of a $p^\infty$-torsion Kisin module as

$$(M(i), \varphi_{M(i)}) = (M, u^i \varphi_M),$$

and we have that for all $p^\infty$-torsion Kisin modules, there exists $i \geq 0$ such that $M(i)$ is effective.

Now we need to define a slope function. First, we start by defining the rank function as

$$\mathrm{rank}(M) = \mathrm{length}_{W(\mathbb{F})}(M/uM).$$

We cannot use relative position in this context to define the degree function since we do not have a DVR anymore, but recall that the degree function of a $p$-torsion Kisin module corresponded to the length of a quotient when the $p$-torsion Kisin module was effective. Since the twist turns a module effective for a sufficiently large integer, we define the degree function by

$$\deg(M) = -\mathrm{length}_{\mathfrak{S}}(Q(i)) + i \, \mathrm{rank}\, M,$$

for $i$ large enough and $Q(i) = M(i)/(\varphi_{M(i)}\varphi^* M(i))$, i.e. the cokernel of $\varphi_{M(i)}$ viewed as a morphism in the category of $\mathfrak{S}$-modules. We can easily check that deg is independent of the choice of a large enough $i \in \mathbb{N}$ and that both degree and rank functions are a degree and rank function in the sense of Harder-Narasimhan. Moreover, for $p$-torsion Kisin modules, it coincides with the degree function defined in last section. We will denote by $\mathcal{F}_{\mathrm{F},\mathrm{t}}(M)$ the Fargues filtration of a $p^\infty$-torsion Kisin module $M$ and by $\mathbf{t}_{\mathrm{F},\mathrm{t}}(M)$, the type associated to $\mathcal{F}_{\mathrm{F},\mathrm{t}}(M)$.

## 3.5 The Fargues filtration on $\mathrm{Mod}^{\varphi}_{\mathfrak{S},\mathrm{fr}}$

### 3.5.1 The category $\mathrm{Mod}^{\varphi}_{\mathfrak{S},\mathrm{fr}}$

Some operators we defined for objects in $\mathrm{Mod}^{\varphi}_{\mathbb{F}[[u]],\mathrm{fr}}$ can be generalized to objects in $\mathrm{Mod}^{\varphi}_{\mathfrak{S},\mathrm{fr}}$:

1. Tensor products : For two Kisin modules $(M_1, \varphi_{M_1})$ and $(M_2, \varphi_{M_2})$, define

   $$(M_1 \otimes M_2, \varphi_{M_1 \otimes M_2}) = (M_1 \otimes_{\mathfrak{S}} M_2, \varphi_{M_1} \otimes \varphi_{M_2}).$$

   This definition works since we have $\varphi^*(M_1 \otimes_{\mathfrak{S}} M_2)[\frac{1}{E}] = (\varphi^* M_1[\frac{1}{E}]) \otimes_{\mathfrak{S}[\frac{1}{E}]} (\varphi^* M_2[\frac{1}{E}])$ and $(M_1 \otimes_{\mathfrak{S}} M_2)[\frac{1}{E}] = (M_1[\frac{1}{E}]) \otimes_{\mathfrak{S}[\frac{1}{E}]} (M_2[\frac{1}{E}])$. The identity object for the tensor product is the Kisin module $\mathbf{1} = (\mathfrak{S}, \varphi)$ where $\varphi$ is the identity on $\varphi^*\mathfrak{S} \simeq \mathfrak{S}$.

2. Duality : For a Kisin module $(M, \varphi_M)$, we define the dual of $(M, \varphi_M)$ as $(M^\vee, \varphi_{M^\vee})$ where $M^\vee = \mathrm{Hom}_{\mathfrak{S}}(M, \mathfrak{S})$ and

$$\begin{array}{rccc} \varphi_{M^\vee} : & (\varphi^* M^\vee)[\frac{1}{E}] & \to & M^\vee[\frac{1}{E}] \\ & f & \mapsto & \varphi \circ f \circ \varphi_M^{-1} \end{array}$$

since $\varphi^* \mathrm{Hom}_{\mathfrak{S}}(M, \mathfrak{S})[\frac{1}{E}] \simeq \mathrm{Hom}_{\varphi^* \mathfrak{S}}(\varphi^* M, \varphi^* \mathfrak{S})[\frac{1}{E}] \simeq \mathrm{Hom}_{\varphi^* \mathfrak{S}[\frac{1}{E}]}(\varphi^* M[\frac{1}{E}], \varphi^* \mathfrak{S}[\frac{1}{E}])$.

3. Internal homomorphisms: Since we have defined the tensor product and duality, we can define an internal Hom by $\underline{\mathrm{Hom}}((M_1, \varphi_{M_1}), (M_2, \varphi_{M_2})) = (M_1, \varphi_{M_1})^\vee \otimes (M_2, \varphi_{M_2})$ for all pair of objects $(M_1, \varphi_1), (M_2, \varphi_{M_2})$ in $\mathrm{Mod}^\varphi_{\mathfrak{S}, \mathrm{fr}}$.

4. Twist: We define the $i$-twist of a Kisin module $(M, \varphi_M)$ as $(M(i), \varphi_{M(i)}) = (M, E^i \varphi_M)$, and we have that for all Kisin modules, there exists $i \geq 0$ such that $M(i)$ is effective. We have $M(i) = M \otimes \mathbf{1}(i)$ for $i \in \mathbb{Z}$.

Let $\mathrm{Mod}^\sigma_{W(\mathbb{F})}$ denote the category of crystals, whose objects are finite free $W(\mathbb{F})$-modules $D$ together with a morphism $\sigma_D : D[\frac{1}{u}] \xrightarrow{\sim} D[\frac{1}{u}]$ with cokernel killed by a power of $u$. We can always associate a crystal to a Kisin module, by sending $(M, \varphi_M)$ to the crystal $(M/uM, \varphi_M \bmod u)$.

For a Kisin module $M \in \mathrm{Mod}^\varphi_{\mathfrak{S}, \mathrm{fr}}$, we define a rank function by

$$\mathrm{rank}(M) = \mathrm{rank}_{\mathfrak{S}} M.$$

Also, we can define two Hodge filtrations and their associated types on $M$ by

$$\begin{aligned} \mathcal{F}_{\mathrm{H}}(M) &= \mathcal{F}(M_{(E)}, \varphi_M \varphi^* M_{(E)}) \quad \text{and} \quad \mathcal{F}_{\mathrm{H},u}(M) = \mathcal{F}(M/uM, \sigma_M \sigma^*(M/uM)), \\ \mathbf{t}_{\mathrm{H}}(M) &= \mathrm{Pos}(M_{(E)}, \varphi_M \varphi^* M_{(E)}) \quad \text{and} \quad \mathbf{t}_{\mathrm{H},u}(M) = \mathrm{Pos}(M/uM, \sigma_M \sigma^*(M/uM)) \end{aligned}$$

in $\mathbb{Z}^{\mathrm{rank}\, M}_{\geq}$, where $\sigma$ is the Frobenius on $\mathfrak{S}/u\mathfrak{S} = W(\mathbb{F})$ and $\sigma_M$ is the Frobenius on $M/uM$. Then, we can also define a degree function given by

$$\deg(M) = \deg(\mathbf{t}_{\mathrm{H}}(M)) = \nu(M_{(E)}, \varphi_M \varphi^* M_{(E)}).$$

Suppose $M$ is effective, then and $\deg(M) = -\mathrm{length}_{\mathfrak{S}_{(E)}}(Q_{(E)})$. If $M$ is an object in $\mathrm{Mod}^\varphi_{\mathfrak{S}, \mathrm{fr}}$, then $M/p^n M$ is an object in $\mathrm{Mod}^\varphi_{\mathfrak{S}, \mathrm{t}}$ for every $n \geq 1$, the Frobenius being the reduction of the Frobenius on $M$. In particular $M/pM$ is an object in $\mathrm{Mod}^\varphi_{\mathbb{F}[[u]], \mathrm{fr}}$.

**Proposition 3.15.** *Let* $M \in \mathrm{Mod}_{\mathfrak{S}, \mathrm{fr}}$ *and* $\overline{M} = M/pM$. *Then*

$$\begin{aligned} \mathbf{t}_{\mathrm{H}}(\overline{M}) &\leq e \cdot \mathbf{t}_{\mathrm{H}}(M) \\ \mathbf{t}_{\mathrm{H},u}(M) &\leq \mathbf{t}_{\mathrm{H}}(M). \end{aligned}$$

*Proof.* We can reduce to the effective case. Let $M$ be an effective Kisin module and $Q = \mathrm{coker}\, \varphi_M$ and $m = \mathrm{rank}\, M$. We fix a pseudo-isomorphism $Q \sim Q' = \bigoplus_{i=1}^{m} \mathfrak{S}/(E(u)^{n_i})$, which becomes an isomorphism after localizing by $(E)$, with $n_1 \geq \ldots \geq n_m \geq 0$ (completing by 0 the invariant factors of $Q$). For $N = n_1$, consider the filtration

$$0 = Q'_0 \subsetneq Q'_{\frac{1}{N}} \subsetneq \ldots \subsetneq Q'_{\frac{N}{N}} = Q'$$

with $Q'_{\frac{i}{N}} = Q'[E^i]$. Let $Q_{\frac{i}{N}} = Q \cap Q'_{\frac{i}{N}} = Q[E^i]$, so we have a filtration

$$0 = Q_0 \subsetneq Q_{\frac{1}{N}} \subsetneq \ldots \subsetneq Q_{\frac{N}{N}} = Q.$$

Also, an exact sequence

$$0 \to Q/Q_{\frac{i}{N}} \hookrightarrow Q'/Q'[E^i] \twoheadrightarrow Q'/(Q + Q'[E^i]) \to 0$$

and an isomorphism

$$Q'/Q'[E^i] \simeq \bigoplus_{n_j > i} \mathfrak{S}/E^{n_j - i} \mathfrak{S}$$

for $1 \leq i \leq N$, so $\operatorname{Tor}^2\left(Q/Q_{\frac{i}{N}}, \mathbb{F}\right) = \left(Q/Q_{\frac{i}{N}}\right)[\mathfrak{m}] = 0$. Let $\varphi^* M \subset M_{\frac{i}{N}} \subset M$ with

$$M_{\frac{i}{N}}/\varphi^* M = Q_{\frac{i}{N}} \subset Q = M/\varphi^* M.$$

Since $M/M_{\frac{i}{N}} = Q/Q_{\frac{i}{N}}$, we have $\operatorname{Tor}^2\left(M/M_{\frac{i}{N}}, \mathbb{F}\right) = 0$, $\operatorname{Tor}^1(M_{\frac{i}{N}}, \mathbb{F}) = 0$ and $M_{\frac{i}{N}}$ is free. So we have a filtration

$$\varphi^* M \subsetneq M_{\frac{1}{N}} \subsetneq \ldots \subsetneq M_{\frac{N}{N}} = M$$

with quotients

$$M_{\frac{i}{N}}/M_{\frac{i-1}{N}} = Q_{\frac{i}{N}}/Q_{\frac{i-1}{N}} = Q[E^i]/Q[E^{i-1}] = R_i,$$

which are contained in

$$R_i' = Q'[E^i]/Q'[E^{i-1}] \simeq \bigoplus_{n_j \geq i} \mathfrak{S}/E\mathfrak{S}$$

and with cokernel of finite length. Since $\mathfrak{S}/E\mathfrak{S}$ is a DVR, we also have $R_i \simeq \bigoplus_{n_j \geq i} \mathfrak{S}/E\mathfrak{S}$.

Let $f = p$ or $u$, and let quotient by $f$. We obtain a sequence of free $\overline{\mathfrak{S}} = \mathfrak{S}/f\mathfrak{S}$-modules of rank $h$,

$$\overline{\varphi^* M} \subsetneq \overline{M}_{\frac{1}{N}} \subsetneq \ldots \subsetneq \overline{M}_{\frac{N}{N}} = \overline{M},$$

with torsion $\overline{\mathfrak{S}}$-modules as quotients $\overline{R}_i = \bigoplus_{n_j \geq i} \overline{\mathfrak{S}}/E\overline{\mathfrak{S}}$. The triangle inequality for relative positions tells us

$$\operatorname{Pos}(\overline{\varphi^* M}, \overline{M}) \leq \sum_{i=1}^{N} \operatorname{Pos}(\overline{M}_{\frac{i-1}{N}}, \overline{M}_{\frac{i}{N}})$$

and we know that $\operatorname{Pos}(\overline{M}_{\frac{i-1}{N}}, \overline{M}_{\frac{i}{N}}) = (e, \ldots, e, 0, \ldots, 0) \in \mathbb{Z}_{\geq}^m$, where $e = \overline{v}(E)$ is the valuation of $E$ in $\overline{\mathfrak{S}}$ (i.e. $e = \deg f$ for $f = p$ and $e = 1$ for $f = u$) and the multiplicity of $e$ is $\#\{j \,:\, n_j \geq i\}$. So we get

$$\operatorname{Pos}(\overline{\varphi^* M}, \overline{M}) \leq e \cdot (n_1, n_2, \ldots, n_m) = e \cdot \operatorname{Pos}(\varphi^* M_{(E)}, M_{(E)}).$$

$\square$

For $M \in \operatorname{Mod}_{\mathfrak{S},\,\mathrm{fr}}^{\varphi}$, it is easy to see that we have

$$\begin{array}{rcl} \deg(M/p^n M) & = & en \deg(M) \\ \operatorname{rank}(M/p^n M) & = & n \operatorname{rank}(M), \end{array}$$

for $n \geq 1$, where $e = \deg E$. Also, for the isogeny class, we have $\deg(M \otimes \mathbb{Q}_p) = \deg(M)$ and $\operatorname{rank}(M \otimes \mathbb{Q}_p) = \operatorname{rank}(M)$. This way, we get an invariant

$$\mu_M = \frac{\deg M}{\operatorname{rank} M}$$

which only depends on the isogeny class of $M$. We want to get rid of the coefficient $e$. For that, it suffices to define a new degree function on $\operatorname{Mod}_{\mathfrak{S},\,\mathrm{t}}^{\varphi}$ by putting

$$\deg_{\mathrm{new}}(M) = \tfrac{1}{e} \deg_{\mathrm{old}}(M)$$

for every $M \in \operatorname{Mod}_{\mathfrak{S},\,\mathrm{t}}^{\varphi}$. The new degree function verifies the same properties as the old one and gives us a new slope function, Fargues filtration and Fargues polygon. From now on, we will only work with the new degree function, filtration and polygon and drop it from the notation. Thus, for all $n \in \mathbb{N}$,

$$\begin{array}{rcl} n \deg(M) & = & \deg(M/p^n M) \\ n \operatorname{rank}(M) & = & \operatorname{rank}(M/p^n M) \\ \mu_M & = & \mu(M/p^n M). \end{array}$$

### 3.5.2 Minimal slopes

**Definition 3.8.** For $M \in \text{Mod}^{\varphi}_{\mathfrak{S},\,\text{t}}$, let $M^{\min}$ be the last quotient of the Harder-Narasimhan filtration of $M$. We define the minimal slope of $M$ as

$$\mu_{\min}(M) := \mu(M^{\min}) = \max\{\lambda \mid \mathcal{F}^{\geq \lambda} M = M\}.$$

It is easy to see that then $M$ is semi-stable if and only if $\mu(M) = \mu_{\min}(M)$.

**Proposition 3.16.** *Let $M \in \text{Mod}^{\varphi}_{\mathfrak{S},\,\text{t}}$. Then for every nonzero strict quotient $M \twoheadrightarrow Q$ in $\text{Mod}^{\varphi}_{\mathfrak{S},\,\text{t}}$, we have*

$$\mu(Q) \geq \mu_{\min}(Q) \geq \mu_{\min}(M)$$

*and the equality $\mu(Q) = \mu_{\min}(M)$ holds if and only if $Q$ is a semi-stable quotient of $M^{\min}$.*

*Proof.* 1. The filtration $\mathcal{F}_{\text{F},\text{t}}$ is functorial, so have morphisms $\mathcal{F}^{\geq \lambda}_{\text{F},\text{t}} M \to \mathcal{F}^{\geq \lambda}_{\text{F},\text{t}} Q$ for every $\lambda \in \mathbb{R}$. In particular, for $\lambda = \mu_{\min}(M)$, we get

$$M = \mathcal{F}^{\geq \mu_{\min}(M)}_{\text{F},\text{t}} M \longrightarrow \mathcal{F}^{\geq \mu_{\min}(M)}_{\text{F},\text{t}} Q$$
$$\searrow \qquad \qquad \swarrow$$
$$Q$$

so $\mathcal{F}^{\geq \mu_{\min}(M)}_{\text{F},\text{t}} Q = Q$. Since, by definition, we have $\mu_{\min}(Q) = \max\{\lambda \in \mathbb{R} \mid \mathcal{F}^{\geq \lambda}_{\text{F},\text{t}} Q = Q\}$ it follows that $\mu_{\min}(M) \leq \mu_{\min}(Q)$. Also, we know that $\mu_{\min}(Q) \leq \mu(Q)$ by concavity of $\mathbf{t}_{\text{F},\text{t}}(Q)$.

2. Suppose that $\mu(Q) = \mu_{\min}(M)$, thus Then $Q$ is semi-stable and $M \twoheadrightarrow Q$ maps $\mathcal{F}^{> \mu_{\min}(M)}_{\text{F},\text{t}} M$ to $\mathcal{F}^{> \mu_{\min}(M)}_{\text{F},\text{t}} Q = 0$, i.e the morphsim $M \twoheadrightarrow Q$ factors through $M \twoheadrightarrow M^{\min} \twoheadrightarrow Q$ and $Q$ is a semi-stable quotient of $M^{\min}$. $\square$

For $M \in \text{Mod}^{\varphi}_{\mathfrak{S},\,\text{fr}}$, we set $\mu_{\min}(M) := \mu_{\min}(M/pM)$.

**Lemma 3.17.** *Let $M \in \text{Mod}^{\varphi}_{\mathfrak{S},\,\text{fr}}$ nonzero and $M_n = M/p^n M$. Then:*

1. *For every $n \geq 1$, we have $\mu_{\min}(M) = \mu_{\min}(M_n) = \mu(M^{\min}_n)$.*
2. *For $1 \leq m \leq n$, we have $p^m M^{\min}_n = \text{im}(M^{\min}_{n-m} \to M^{\min}_n)$.*
3. *For $1 \leq m \leq n$, we have $M^{\min}_n / p^m M^{\min}_n = M^{\min}_m$.*
4. *The system $(N)_{n \geq 1}$, for $N_n = \ker(M_n \to M^{\min}_n)$, is Mittag-Leffler surjective.*

*Proof.* 1. Let $M \in \text{Mod}^{\varphi}_{\mathfrak{S},\,\text{fr}}$ and consider the exact sequence

$$0 \to M^{\min}_n[p^m] \to M^{\min}_n \xrightarrow{\times p^m} M^{\min}_n \to M^{\min}_n / p^m M^{\min}_n \to 0$$

for all $1 \leq m \leq n$. By the properties of our Harder-Narasimhan formalism, seen in 2.7thm.2.7, we know that $M^{\min}_n / p^m M^{\min}_n$ is semi-stable and

$$\mu(M^{\min}_n) = \mu(M^{\min}_n / p^m M^{\min}_n).$$

Since $\mu_{\min}(M/p^n M) = \mu(M^{\min}_n)$ (and $\mu_{\min}(M/p^m M) = \mu(M^{\min}_m)$ respectively), and we have the quotients

$$M/p^n M \twoheadrightarrow M/p^m M \twoheadrightarrow M^{\min}_m$$

and

$$M/p^m M \simeq M_n / p^m M_n \twoheadrightarrow M^{\min}_n / p^m M^{\min}_n$$

we get the following inequalities

$$\mu(M^{\min}_n) \leq \mu(M^{\min}_m) \leq \mu(M^{\min}_n / p^m M^{\min}_n)$$

which must be equalities in order to have $\mu(M^{\min}_n / p^m M^{\min}_n) = \mu(M^{\min}_n)$. In particular $\mu(M^{\min}_n) = \mu(M^{\min}_1) = \mu_{\min}(M)$, for every $n \geq 1$.

2. We have

$$\begin{array}{rcl}
\mathrm{im}(M_{n-m}^{\min} \to M_n^{\min}) & = & \mathrm{im}(M_{n-m} \twoheadrightarrow M_{n-m}^{\min} \to M_n^{\min}) \\
& = & \mathrm{im}(M_{n-m} \hookrightarrow M_n \twoheadrightarrow M_n^{\min}) \\
& = & \mathrm{im}(p^m M_n \hookrightarrow M_n \twoheadrightarrow M_n^{\min}) \\
& = & p^m M_n^{\min}.
\end{array}$$

3. If we apply last proposition to the surjection $M_m \simeq M_n/p^m M_n \twoheadrightarrow M_n^{\min}/p^m M_n^{\min}$, we get

$$M_m^{\min} \twoheadrightarrow M_n^{\min}/p^m M_n^{\min}.$$

On the other hand, applying last proposition to the surjection $M_n \twoheadrightarrow M_m \twoheadrightarrow M_m^{\min}$ we obtain that $M_n^{\min} \twoheadrightarrow M_m^{\min}$ and it is obvious that it factors through $M_n^{\min}/p^m M_n^{\min}$, giving us a surjection

$$M_n^{\min}/p^m M_n^{\min} \twoheadrightarrow M_m^{\min},$$

which ends the proof.

4. Denote by $N_n^m = \ker(N_n \to N_m)$, $Q_n^m = \ker(M_n^{\min} \to M_m^{\min})$, for $n \geq m \geq 1$ and let $C = \mathrm{coker}(N_n \to N_m)$. Then, we have $N_{n-m} \subset N_n^m = M_{n-m} \cap N_n$ and

$$0 \to N_n^m/N_{n-m} \to M_{n-m}^{\min} \to Q_n^m \to C \to 0.$$

Since both $M_{n-m}^{\min}$ and $M_{n-m}^{\min}$ are semi-stable of slope $\mu_{\min}(M)$, then $Q_n^m$ is semi-stable of slope $\mu_{\min}(M)$. Suppose $C \neq 0$, then $C$ is semi-stable of slope $\mu_{\min}(M)$ and also a quotient of $N_m$, so

$$\mu_{\min}(M) = \mu(C) \geq \mu(N_m) > \mu(M_m^{\min}) = \mu_{\min}(M)$$

which is a contradiction, therefore $C = 0$ and $(N_n)_{n \geq 1}$ is Mittag-Leffler surjective. $\square$

### 3.5.3 A polygon on $\mathrm{Mod}_{\mathfrak{S},\,\mathrm{fr}}^{\varphi}$

Even though we have rank and degree functions defined on Kisin modules, $\mathrm{Mod}_{\mathfrak{S},\,\mathrm{fr}}^{\varphi}$ is not a quasi-abelian category, so we cannot use André's formalism in order to endow this category with a Harder-Narasimhan filtration. However, we can consider the projective limit of Fargues polygons on torsion Kisin modules to try to construct a polygon on a Kisin module, as follows: for $M$ a Kisin module, we can renormalize the polygons $\mathbf{t}_{\mathrm{F},\mathrm{t}}(M/p^n M)$ by considering the concave functions $\mathbf{t}_{\mathrm{F},n}(M)$ defined as

$$\mathbf{t}_{\mathrm{F},n}(M)(x) = \tfrac{1}{n}\mathbf{t}_{\mathrm{F},\mathrm{t}}(M/p^n M)(nx)$$

for every $x \in [0, \mathrm{rank}\, M]$, for every $n \geq 1$. Such functions are concave polygons with end point at $(\mathrm{rank}\, M, \deg M)$ but whose break points may not have integer abscissas, so we cannot consider them as types, but we can extend all the properties for types to these concave functions. In particular, the dominance order still makes sense for them. Moreover, the next lemma shows that they form a decreasing system.

**Lemma 3.18.** *For every Kisin module $M$, we have $\mathbf{t}_{F,nm}(M) \leq \mathbf{t}_{F,n}(M)$ for every $m, n \geq 1$.*

*Proof.* By induction using the exact sequences

$$0 \to M/p^n M \to M/p^{mn} M \to M/p^{m(n-1)} M \to 0$$

we get

$$\mathbf{t}_{\mathrm{F},\mathrm{t}}(M/p^{mn} M) \leq \mathbf{t}_{\mathrm{F},\mathrm{t}}(M/p^n M)^{*m}.$$

Now, the polygon $\mathbf{t}_{\mathrm{F},\mathrm{t}}(M/p^n M)^{*m}$ is just an homotethy of $\mathbf{t}_{\mathrm{F},\mathrm{t}}(M/p^n M)$ and we have

$$\mathbf{t}_{\mathrm{F},\mathrm{t}}(M/p^n M)^{*m}(x) = m\mathbf{t}_{\mathrm{F},\mathrm{t}}(M/p^n M)(\tfrac{x}{m})$$

for $0 \leq x \leq mn\, \mathrm{rank}\, M$. Thus

$$\tfrac{1}{mn}\mathbf{t}_{\mathrm{F},\mathrm{t}}(M/p^{mn} M)(mnx) \leq \tfrac{1}{n}\mathbf{t}_{\mathrm{F},\mathrm{t}}(M/p^n M)(nx)$$

for $0 \leq x \leq \mathrm{rank}\, M$ and we obtain $\mathbf{t}_{\mathrm{F},mn}(M) \leq \mathbf{t}_{\mathrm{F},n}(M)$. $\square$

Then, by Ascoli's theorem, we obtain:

**Theorem 3.19.** *Let $M \in \text{Mod}_{\mathfrak{S},\,\text{fr}}^{\varphi}$. Let $r = \text{rank}\, M$ and $d = \deg M$. The sequence of functions $\mathbf{t}_{\text{F},n}(M) : [0,r] \to \mathbb{R}$ converges uniformly for the divisibility order to a continuous concave function*

$$\mathbf{t}_{\text{F},\infty}(M) \;:\; [0,r] \to \mathbb{R}$$

*which is equal to $\inf_{n \geq 1} \mathbf{t}_{\text{F},n}(M)$ and, moreover, verifying $\mathbf{t}_{\text{F},\infty}(M)(0) = 0$ and $\mathbf{t}_{\text{F},\infty}(M)(r) = d$. We call this function the Fargues polygon of $M$.*

The Fargues polygon of a Kisin module only depends on the isogeny class of the module, as we can see in the next proposition.

**Proposition 3.20.** *Let $M_1$, $M_2 \in \text{Mod}_{\mathfrak{S},\,\text{fr}}^{\varphi}$ such that $M_1 \otimes \mathbb{Q}_p \simeq M_2 \otimes \mathbb{Q}_p$ in $\text{Mod}_{\mathfrak{S}}^{\varphi} \otimes \mathbb{Q}_p$. Then $\mathbf{t}_{\text{F},\infty}(M_1) = \mathbf{t}_{\text{F},\infty}(M_2)$.*

*Proof.* The proof is analogous to the one in [15, Proposition 3]. By reflexivity of the isogeny relation, it suffices to prove one inequality. Fix an isogeny $M_1 \to M_2$ giving rise to an exact sequence

$$0 \to M_1 \to M_2 \to Q \to 0$$

in $\text{Mod}_{\mathfrak{S}}^{\varphi}$ with $p^n Q = 0$ for $n$ large enough. By multiplication by $p^n$, we obtain an exact sequence

$$0 \to Q \to M_1/p^n M_1 \to M_2/p^n M_2 \to Q \to 0$$

in $\text{Mod}_{\mathfrak{S},\,\text{t}}^{\varphi}$. This yields two exact sequences

$$0 \to Q \to M_1/p^n M_1 \to N_n \to 0$$
$$0 \to N_n \to M_2/p^n M_2 \to Q \to 0$$

in $\text{Mod}_{\mathfrak{S},\,\text{t}}^{\varphi}$, with $nr = r_n + q$, for $r = \text{rank}\, M$, $r_n = \text{rank}\, N_n$ and $q = \text{rank}\, Q$. Then $\frac{r_n}{n} = r - \frac{q}{n}$. For $n$ large enough, let $x \in [0, \frac{r_n}{n}] = [0, r - \frac{q}{n}]$, then the second exact sequence gives us

$$\mathbf{t}_{\text{F},\text{t}}(N_n)(nx) \leq \mathbf{t}_{\text{F},\text{t}}(M_2/p^n M_2)(nx)$$

and the first exact sequence implies

$$\begin{aligned}\mathbf{t}_{\text{F},\text{t}}(M_1/p^n M_1)(nx) &\leq (\mathbf{t}_{\text{F},\text{t}}(N_n) * \mathbf{t}_{\text{F},\text{t}}(Q))(nx) \\ &= \mathbf{t}_{\text{F},\text{t}}(N_n)(nx - \delta_n) + \mathbf{t}_{\text{F},\text{t}}(Q)(\delta_n)\end{aligned}$$

for some $0 \leq \delta_n \leq nx, q$, by definition of the concatenation. By concavity of $\mathbf{t}_{\text{F},\text{t}}(N_n)$, we have

$$\mathbf{t}_{\text{F},\text{t}}(N_n)(nx - \delta_n) \leq \mathbf{t}_{\text{F},\text{t}}(N_n)(nx) - \mu_{\min}(N_n) \cdot \delta_n.$$

Since $N_n$ is a quotient of $M_1/p^n M_1$ by proposition 3.16thm.3.16, we obtain

$$\mu_{\min}(N_n) \geq \mu_{\min}(M_1/p^n M_1) = \mu_{\min}(M_1),$$

thus

$$\begin{aligned}\mathbf{t}_{\text{F},\text{t}}(M_1/p^n M_1)(nx) &\leq \mathbf{t}_{\text{F},\text{t}}(N_n)(nx) + \mathbf{t}_{\text{F},\text{t}}(Q)(\delta_n) - \mu_{\min}(M_1) \cdot \delta_n \\ &\leq \mathbf{t}_{\text{F},\text{t}}(N_n)(nx) + C\end{aligned}$$

for $C = \max\{\mathbf{t}_{\text{F},\text{t}}(Q)(\delta) - \mu_{\min}(M_1) \cdot \delta \mid \delta \in [0,q]\}$, which is independent of $n$. Combining all the inequalities above, we have

$$\mathbf{t}_{\text{F},\text{t}}(M_1/p^n M_1)(nx) \leq \mathbf{t}_{\text{F},\text{t}}(N_n)(nx)) + C \leq \mathbf{t}_{\text{F},\text{t}}(M_2/p^n M_2)(nx) + C,$$

so

$$\mathbf{t}_{\text{F},n}(M_1)(x) \leq \mathbf{t}_{\text{F},n}(M_2)(x) + \frac{C}{n}$$

for $n$ large enough and $x \in [0, r - \frac{q}{n}]$, therefore

$$\mathbf{t}_{\text{F},\infty}(M_1)(x) \leq \mathbf{t}_{\text{F},\infty}(M_2)(x)$$

for every $x \in [0, r)$ and we already know the equality for $x = r$. □

### 3.5.4 Semi-stability and type HN

There is not a Harder-Narasimhan formalism, but we can define semi-stability on objects in $\text{Mod}^{\varphi}_{\mathfrak{S}, \text{fr}}$.

**Proposition 3.21.** *Let $M \in \text{Mod}^{\varphi}_{\mathfrak{S}, \text{fr}}$. The following conditions are equivalent:*

1. *The module $M/pM$ is semi-stable (as a $p$-torsion Kisin module),*
2. *For every $n \geq 1$, the module $M/p^n M$ is semi-stable (as a $p^{\infty}$-torsion Kisin module),*
3. *For every $p^{\infty}$-torsion Kisin module $Q$ which is a quotient of $M$ we have $\mu_M \leq \mu(Q)$,*

*Proof.* It is easy to see that (1) implies (2), since the category $\mathcal{C}(\mu)$ of semi-stable objects of slope $\mu$ is stable by extension. For (2) implies (3), we have that for any $p^{\infty}$-torsion quotient $Q$ of $M$, there exists a $n \geq 1$, such that we have a factorization $M \twoheadrightarrow M/p^n M \twoheadrightarrow Q$, thus $\mu_M = \mu(M/p^n M) \leq \mu(Q)$. For (3) implies (1), we have that any quotient $Q$ of $M/pM$ is a $p^{\infty}$-torsion quotient of $M$ therefore $\mu(M/pM) = \mu_M \leq \mu(Q)$ and $M/pM$ is semi-stable. □

**Definition 3.9.** We say that a Kisin module $M$ is semi-stable if it verifies one of the conditions above.

*Remark 8.* If a Kisin module $M$ is semi-stable, then $\mathbf{t}_{F,\infty}(M)$ is isoclinic of slope $\mu_M$.

In the next proposition, we compare the semi-stability on $\text{Mod}^{\varphi}_{\mathfrak{S}, \text{fr}}$ with the semi-stability on $\text{Mod}^{\varphi}_{\mathfrak{S}} \otimes \mathbb{Q}_p$.

**Proposition 3.22.** *Let $M, M' \in \text{Mod}^{\varphi}_{\mathfrak{S}, \text{fr}}$ isogenous and $M'$ semi-stable in $\text{Mod}^{\varphi}_{\mathfrak{S}, \text{fr}}$. Then for every subobject $N \neq 0$ of $M \otimes \mathbb{Q}_p$, we have $\mu_N \leq \mu_M$, i.e. $M \otimes \mathbb{Q}_p$ is semi-stable in $\text{Mod}^{\varphi}_{\mathfrak{S}} \otimes \mathbb{Q}_p$.*

*Proof.* As we have already discussed, an object in $\text{Mod}^{\varphi}_{\mathfrak{S}} \otimes \mathbb{Q}_p$ is semi-stable of slope $\lambda$ if and only if its image in $\text{Mod}^{\varphi}_{\mathfrak{S}[\frac{1}{p}]}$ is semi-stable of slope $\lambda$. We use this characterization. Let $M \in \text{Mod}^{\varphi}_{\mathfrak{S}, \text{fr}}$ and $N = M[\frac{1}{p}]$ its image in $\text{Mod}^{\varphi}_{\mathfrak{S}[\frac{1}{p}]}$. An exact sequence $0 \to N_1 \to N \to N_2 \to 0$ in $\text{Mod}^{\varphi}_{\mathfrak{S}[\frac{1}{p}]}$ induces an exact sequence $0 \to M_1 \to M \to M_2 \to 0$ in $\text{Mod}^{\varphi}_{\mathfrak{S}}$, where $M_1 = M \cap N_1$ and $M_2$ is the image of $M$ in $N_2$. The cokernel of $M_2 \hookrightarrow M_{2,\text{fr}}$ is killed by $p^m$ for some $m \geq 1$, so multiplication by $p^m$ gives us an isogeny $f: M_{2,\text{fr}} \hookrightarrow M_2$. We have thus a diagram

$$
\begin{array}{ccccccccc}
0 & \longrightarrow & M_1 & \longrightarrow & M & \longrightarrow & M_2 & \longrightarrow & 0 \\
& & \downarrow i & & \downarrow i & & \| & & \\
0 & \longrightarrow & M_{1,\text{fr}} & \longrightarrow & M' & \longrightarrow & M_2 & \longrightarrow & 0 \\
& & \| & & f \uparrow & & f \uparrow & & \\
0 & \longrightarrow & M_{1,\text{fr}} & \longrightarrow & M'' & \longrightarrow & M_{2,\text{fr}} & \longrightarrow & 0
\end{array}
$$

by push out of $i: M_1 \to M_{1,\text{fr}}$ and pull-back by $f: M_{2,\text{fr}} \to M_2$. Since $M''$ is an extension of two free modules, it is itself free and, for every $n \geq 1$, we have an exact sequence

$$0 \to M_{1,\text{fr}}/p^n M_{1,\text{fr}} \to M''/p^n M'' \to M_{2,\text{fr}}/p^n M_{2,\text{fr}} \to 0.$$

Suppose $M$ is semi-stable of slope $\lambda$. Then, $M/pM$ is semi-stable of slope $\lambda$, so $\mathbf{t}_{F,1}(M)$ is isoclinic of slope $\lambda$, which implies that $\mathbf{t}_{F,\infty}(M)$ is also isoclinic of slope $\lambda$, since $\mathbf{t}_{F,\infty}(M) \leq \mathbf{t}_{F,1}(M)$. Now, $\mathbf{t}_{F,\mathfrak{S}[\frac{1}{p}]}(N)$ has the same starting and ending point, so

$$\mathbf{t}_{F,\infty}(M) \leq \mathbf{t}_{F,\mathfrak{S}[\frac{1}{p}]}(N).$$

Let $\lambda_{\max}$ be the maximal slope of $\mathbf{t}_{F,\mathfrak{S}[\frac{1}{p}]}(N)$. Then, $\mathbf{t}_{F,\mathfrak{S}[\frac{1}{p}]}(N)$ is isoclinic of slope $\lambda$ if and only if $\lambda_{\max} = \lambda$, if and only if $\lambda_{\max} \leq \lambda$, since we already have the other inequality. Let $N_1$ be the first submodule in the Fargues filtration of $N$, such that $\mu(N_1) = \lambda_{\max}$ and $r_1 = \text{rank } N_1$. We have that

$$r_1 \lambda_{\max} = \deg N_1 = \deg M_{1,\text{fr}} = \mathbf{t}_{F,n}(M_{1,\text{fr}})(r_1) \leq \mathbf{t}_{F,n}(M'')(r_1)$$

so $r_1 \lambda_{\max} \leq \mathbf{t}_{F,\infty}(M'')(r_1)$. In 3.20thm.3.20, we have seen that $\mathbf{t}_{F,\infty}$ is invariant by isogeny, so we have $\mathbf{t}_{F,\infty}(M'') = \mathbf{t}_{F,\infty}(M'_{1,\text{fr}}) = \mathbf{t}_{F,\infty}(M)$ which implies that

$$r_1 \lambda_{\max} \leq \mathbf{t}_{F,\infty}(M'')(r_1) = \mathbf{t}_{F,\infty}(M)(r_1) = r_1 \lambda,$$

hence $\lambda_{\max} \leq \lambda$ and $\mathbf{t}_{F,\mathfrak{S}[\frac{1}{p}]}(N) = \mathbf{t}_{F,\circ}(M \otimes \mathbb{Q}_p) = \mathbf{t}_{F,\infty}(M)$ is isoclinic of slope $\lambda$. □

We cannot use André's formalism to give a Harder-Narasimhan filtration on Kisin modules, but some of them may have a Harder-Narasimhan filtration for the semi-stable definition given above.

**Proposition 3.23.** *Let $M \in \mathrm{Mod}^{\varphi}_{\mathfrak{S},\,\mathrm{fr}}$. The following statements are equivalent:*

1. *The Kisin module $M$ admits a Harder-Narasimhan filtration by Kisin submodules or, equivalently, an increasing flag*
$$0 = M_0 \subset M_1 \subset \ldots \subset M_r = M$$
*of Kisin submodules such that, for $1 \leq i \leq r$, the modules $M_i/M_{i-1}$ are free and semi-stable of slopes $\mu_i$ verifying $\mu_1 > \ldots > \mu_r$.*

2. *We have $\mathbf{t}_{\mathrm{F},\infty}(M) = \mathbf{t}_{\mathrm{F},1}(M)$.*

*Moreover, if the conditions above are verified, we have:*

(a) *The flag is uniquely determined by $M_i = \varprojlim_n \mathcal{F}^{\mu_i}_{\mathrm{F},n}(M)$ for $1 \leq i \leq r$ and we have $\mathbf{t}_{\mathrm{F},\infty}(M) = (\mu_1, \ldots, \mu_r)$ (with multiplicity).*

(b) *The reduction of the flag modulo $p^n$ gives us the Fargues filtration of $M/p^n M$.*

(c) *The flag is compatible with the one obtained for the isogeny categories, i.e. for $1 \leq i \leq r$*
$$M_i \otimes \mathbb{Q}_p = \mathcal{F}^{\mu_i}_{\mathrm{F},\circ}(M \otimes \mathbb{Q}_p) \quad \text{and} \quad M_i[\tfrac{1}{p}] = \mathcal{F}^{\mu_i}_{\mathrm{F},\mathfrak{S}[\frac{1}{p}]}(M[\tfrac{1}{p}]) \ .$$

*Proof.* It is easy to see that (1) implies (2) since the semi-stability condition is preserved by reducing modulo $p^n$, so we obtain a flag for $M_n$ whose quotients are semi-stable, and the Fargues filtration is unique. Moreover, the slopes are also preserved by reduction and by the compatibility of the rank and slopes, we get that the polygons modulo $p^n$ are constant and we have $\mathbf{t}_{\mathrm{F},n}(M) = (\mu_1, \ldots, \mu_r)$, thus $\mathbf{t}_{F,\infty}(M) = \lim \mathbf{t}_{\mathrm{F},n}(M) = (\mu_1, \ldots, \mu_r)$.

Conversely, suppose that $\mathbf{t}_{\mathrm{F},\infty}(M) = \mathbf{t}_{\mathrm{F},1}(M)$. We consider the sequences
$$0 \to M_n \to M_{m+n} \to M_m \to 0,$$
for $n, m \geq 1$. Our hypothesis implies that all the polygons $\mathbf{t}_{\mathrm{F},n}(M)$ coincide for $n \geq 1$, so we have
$$\mathbf{t}_{\mathrm{F},\mathrm{t}}(M/p^{n+m}M) = \mathbf{t}_{\mathrm{F},\mathrm{t}}(M/p^n M) * \mathbf{t}_{\mathrm{F},\mathrm{t}}(M/p^m M)$$
for $n, m \geq 1$. Then, we can apply the point (3) in 2.7thm.2.7 to get a diagram

$$\begin{array}{ccccccccc}
& & 0 & & 0 & & 0 & & \\
& & \downarrow & & \downarrow & & \downarrow & & \\
0 & \longrightarrow & \mathcal{F}^{>\gamma}_{\mathrm{F},\mathrm{t}} M_n & \longrightarrow & \mathcal{F}^{>\gamma}_{\mathrm{F},\mathrm{t}} M_{m+n} & \longrightarrow & \mathcal{F}^{>\gamma}_{\mathrm{F},\mathrm{t}} M_m & \longrightarrow & 0 \\
& & \downarrow & & \downarrow & & \downarrow & & \\
0 & \longrightarrow & \mathcal{F}^{\geq\gamma}_{\mathrm{F},\mathrm{t}} M_n & \longrightarrow & \mathcal{F}^{\geq\gamma}_{\mathrm{F},\mathrm{t}} M_{m+n} & \longrightarrow & \mathcal{F}^{\geq\gamma}_{\mathrm{F},\mathrm{t}} M_m & \longrightarrow & 0 \\
& & \downarrow & & \downarrow & & \downarrow & & \\
0 & \to & \mathrm{Gr}^{\geq\gamma}_{\mathcal{F}_{\mathrm{F},\mathrm{t}}} M_n & \to & \mathrm{Gr}^{\gamma}_{\mathcal{F}_{\mathrm{F},\mathrm{t}}} M_{m+n} & \to & \mathrm{Gr}^{\gamma}_{\mathcal{F}_{\mathrm{F},\mathrm{t}}} M_m & \to & 0 \\
& & \downarrow & & \downarrow & & \downarrow & & \\
& & 0 & & 0 & & 0 & &
\end{array}$$

for all $\gamma \in \mathbb{Q}$. In particular $(\mathcal{F}^{>\gamma}_{\mathrm{F},\mathrm{t}} M_n)_{n \geq 1}$, $(\mathcal{F}^{\geq\gamma}_{\mathrm{F},\mathrm{t}} M_n)_{n \geq 1}$ and $(\mathrm{Gr}^{\gamma}_{\mathcal{F}_{\mathrm{F},\mathrm{t}}} M_n)_{n \geq 1}$ are Mittag-Leffler surjective. Call $\mathcal{F}^{>\gamma}_{\mathrm{F}} M$, $\mathcal{F}^{\geq\gamma}_{\mathrm{F}} M$ and $\mathrm{Gr}^{\gamma}_{\mathcal{F}_{\mathrm{F}}} M$ their projective limit. Since $(\mathcal{F}^{>\gamma}_{\mathrm{F},\mathrm{t}} M)_{n \geq 1}$ is Mittag-Leffler surjective, we obtain an exact sequence
$$0 \to \mathcal{F}^{>\gamma}_{\mathrm{F}} M \to \mathcal{F}^{\geq\gamma}_{\mathrm{F}} M \to \mathrm{Gr}^{\gamma}_{\mathcal{F}_{\mathrm{F}}} M \to 0.$$

On the other hand, an exact sequence
$$0 \to X_n \xrightarrow{\times p^m} X_{m+n} \to X_m \to 0$$
for all $n, m \geq 1$, induces an exact sequence
$$0 \to X \xrightarrow{\times p^m} X \to X_m \to 0$$
where $X = \varprojlim_n X_n \simeq \varprojlim_n X_{n+m}$. We apply this result to

- The family $(X_n)_{n\geq 1} = (\mathcal{F}_{\mathrm{F,t}}^{>\gamma} M_n)_{n\geq 1}$ to obtain that $\mathcal{F}_{\mathrm{F}}^{>\gamma} M/p^m \mathcal{F}_{\mathrm{F}}^{>\gamma} M \simeq \mathcal{F}_{\mathrm{F,t}}^{>\gamma} M_m$ for every $m \geq 1$, so $\mathcal{F}_{\mathrm{F}}^{>\gamma} M$ is free of rank $\max\{i \mid \mu_i > \gamma\}$.

- The family $(X_n)_{n\geq 1} = (\mathcal{F}_{\mathrm{F,t}}^{\geq\gamma} M_n)_{n\geq 1}$ to obtain that $\mathcal{F}_{\mathrm{F}}^{\geq\gamma} M/p^m \mathcal{F}_{\mathrm{F}}^{\geq\gamma} M \simeq \mathcal{F}_{\mathrm{F,t}}^{\geq\gamma} M_m$ for every $m \geq 1$, so $\mathcal{F}_{\mathrm{F}}^{\geq\gamma} M$ is free of rank $\max\{i \mid \mu_i \geq \gamma\}$.

- The family $(X_n)_{n\geq 1} = (\mathrm{Gr}_{\mathcal{F}_{\mathrm{F,t}}}^{\gamma} M_n)_{n\geq 1}$ to obtain that $\mathrm{Gr}_{\mathcal{F}_{\mathrm{F}}}^{\gamma} M/p^m \mathrm{Gr}_{\mathcal{F}_{\mathrm{F}}}^{\gamma} M \simeq \mathrm{Gr}_{\mathcal{F}_{\mathrm{F,t}}}^{\gamma} M_m$ for every $m \geq 1$, so $\mathrm{Gr}_{\mathcal{F}_{\mathrm{F}}}^{\gamma} M$ is free of rank $\max\{i \mid \mu_i = \gamma\}$.

In particular, $\mathrm{Gr}_{\mathcal{F}_{\mathrm{F}}} M$ is semi-stable of slope $\gamma$, for $\gamma \in \{\mu_1, \ldots, \mu_r\}$, and the zero object for $\gamma$ otherwise. Since $\mathcal{F}_{\mathrm{F}}^{\gamma} M = M$ for $\gamma \leq \mu_r$, this proves (1).

Assertions ($a$) and ($b$) have already been established in the proof of (1) if and only if (2), and (3) follows from Proposition 3.22thm.3.22. $\square$

**Definition 3.10.** Let $M \in \mathrm{Mod}_{\mathfrak{S},\,\mathrm{fr}}^{\varphi}$. We say that $M$ is HN-type if it verifies one of the conditions above.

### 3.5.5 The algorithm

The aim of this subsection is to prove that every Kisin module is isogenous to a HN-type Kisin module.

We would like to define the last submodule in the Harder Narasimhan filtration of $M \in \mathrm{Mod}_{\mathfrak{S},\,\mathrm{fr}}^{\varphi}$ as

$$\theta M = \ker(M \to M^{\min}) = \varprojlim_i (\ker(M_i \twoheadrightarrow M_i^{\min}))$$

for $M^{\min} = \varprojlim_i M_i^{\min}$ and $M_i^{\min}$ defined as above. The problem is that we cannot ensure that $M/\theta M$ is free. In the next proposition, we find a Kisin module $M'$ isogenous to $M$ and such that $M'/\theta M$ is also a Kisin module.

**Proposition 3.24.** *Let $M \in \mathrm{Mod}_{\mathfrak{S},\,\mathrm{fr}}^{\varphi}$. Then:*

1. *The inclusion $\theta M \to M$ factors (not canonically) through*

$$\theta M \to M' \to M$$

   *where $\theta M, M' \in \mathrm{Mod}_{\mathfrak{S},\,\mathrm{fr}}^{\varphi}$, $M'/\theta M \in \mathrm{Mod}_{\mathfrak{S},\,\mathrm{fr}}^{\varphi}$ is either zero or semi-stable of slope $\mu_{\min}(M)$ and $M' \to M$ is an isogeny with cokernel in $\mathrm{Mod}_{\mathfrak{S},\,\mathrm{t}}^{\varphi}$ which is semi-stable of slope $\mu_{\min}(M)$.*

2. *If $\theta M \neq 0$, then $\mu_{\min}(\theta M) > \mu_{\min}(M)$.*

3. *If $M$ is effective, then $\theta M$ is effective.*

4. *For every $r \geq 1$, we have $\theta(M(r)) = (\theta M)(r)$.*

*Proof.* 1. Let $M \in \mathrm{Mod}_{\mathfrak{S},\mathrm{fr}}^{\varphi}$. We have the following commutative diagram

$$\begin{array}{ccccccccc}
0 & \to & \ker \pi_{k-1} & \to & M_k^{\min} & \overset{\pi_{k-1}}{\to} & M_{k-1}^{\min} & \to & 0 \\
& & \downarrow [p] & & \downarrow [p] & & \downarrow [p] & & \\
0 & \to & \ker \pi_k & \to & M_{k+1}^{\min} & \overset{\pi_k}{\to} & M_k^{\min} & \to & 0
\end{array}$$

where the morphism $\pi_k$ is the reduction modulo $p^k$ and $[p]$ is multiplication by $p$, coming from

$$M/p^k M \xrightarrow{p} pM/p^{k+1}M \hookrightarrow M/p^{k+1}M.$$

Multiplication by $p$ restricted to the kernels is surjective. Indeed, let $x \in \ker \pi_k = p^k M_{k+1}^{\min}$, so there exist an $y \in M_{k+1}^{\min}$ such that $x = p^k y$. Consider $y' = \pi_k(y)$, and $x' = p^{k-1} y'$ which is in $p^{k-1} M_k^{\min} = \ker \pi_{k-1}$. Since $([p] \circ \pi_k)(y) = py$, we get that $[p](x') = x$.

This way, for $a_k = \mathrm{rank}(\ker \pi_k)$, we get a decreasing family of integers $(a_k)_{k\geq 1}$ that will stabilize at some point. Let $k_0 \geq 1$ such that $a_{k+k_0} = a_{k_0}$ for all $k \geq 0$. This means that $\ker \pi_{k+k_0} \simeq \ker \pi_{k_0}$ for all $k \geq 0$, a fact that we will use later. In the case that $a_{k_0} = 0$ we get that $M_k^{\min} \simeq M_{k_0}^{\min}$ for all $k \geq k_0$ and thus $M^{\min} = M_{k_0}^{\min}$.

Let $K_i = \ker(M_{i+k_0}^{\min} \twoheadrightarrow M_{k_0}^{\min})$ for such a $k_0 \geq 1$. Thus,

$$K_i = p^{k_0} M_{i+k_0}^{\min} = \operatorname{im}(M_i^{\min} \to M_{i+k_0}^{\min})$$

by Lemma 3.17thm.3.17. We have the following commutative diagram

$$\begin{array}{ccccccccc} 0 & \to & K_{i+1} & \to & M_{i+1+k_0}^{\min} & \to & M_{k_0}^{\min} & \to & 0 \\ & & \downarrow \pi_{i+k_0} & & \downarrow \pi_{i+k_0} & & \| & & \\ 0 & \to & K_i & \to & M_{i+k_0}^{\min} & \to & M_{k_0}^{\min} & \to & 0 \end{array}$$

so $(K_i)_{i \geq 1}$ is Mittag-Leffler surjective and we can prove by induction that $\operatorname{rank}(K_i) = i \operatorname{rank}(\ker \pi_{k_0})$. We would like to prove that $K = \varprojlim K_i$ is a Kisin module using Proposition 2.10thm.2.10. It only remains to show that $K/p^n K \simeq K_n$ for $n \geq 1$, since $K_n$ has no $u$-torsion and the ranks are compatible.

First, we remark that for a fixed $n \geq 1$ and for all $i \geq n$, we have

$$\ker(K_i \twoheadrightarrow K_n) = \ker(M_{i+k_0}^{\min} \twoheadrightarrow M_{n+k_0}^{\min}) = p^{n+k_0} M_{i+k_0}^{\min} = p^n K_i$$

since $K_i = p^{k_0} M_{i+k_0}^{\min}$ for every $i \geq 1$, so $K_n \simeq K_i/p^n K_i$. Moreover, using the properties of the last quotient studied in 3.17thm.3.17, we have

$$K_i = \ker(M_{i+k_0}^{\min} \to M_{k_0}) = p^{k_0} M_{i+k_0}^{\min} \simeq M_i^{\min} \simeq p^{n+k_0} M_{i+n+k_0}^{\min} \simeq p^n K_{i+n}$$

for every $i, n \geq 1$, so there are exact sequences

$$0 \to K_n \xrightarrow{\times p^i} K_{i+n} \to K_i \to 0,$$

for every $i, n \geq 1$, which gives us

$$\varprojlim_i K_{i+n} \simeq \varprojlim_i K_i$$

and

$$0 \to K \xrightarrow{\times p^n} K \to K_n \to 0,$$

thus $K/p^n K = K_n$ for all $n \geq 1$, hence $K$ is free. It is semi-stable of slope $\mu_{\min}(M)$ since so is each $K_n = K/p^n K$ being the kernel of $M_{n+k_0}^{\min} \to M_{k_0}^{\min}$, both semi-stable of slope $\mu_{\min}(M)$.

For $i \geq k_0$, let $M_i' = \ker(M_i \twoheadrightarrow M_{k_0}^{\min})$ and $(\theta M)_i = \ker(M_i \to M_i^{\min})$. There is a diagram

$$\begin{array}{ccccccccc} & & & & 0 & & 0 & & \\ & & & & \downarrow & & \downarrow & & \\ 0 & \to & (\theta M)_i & \to & M_i' & \to & K_{i-k_0} & \to & 0 \\ & & \| & & \downarrow & & \downarrow & & \\ 0 & \to & (\theta M)_i & \to & M_i & \to & M_i^{\min} & \to & 0 \\ & & & & \downarrow & & \downarrow & & \\ & & & & M_{k_0}^{\min} & = & M_{k_0}^{\min} & & \\ & & & & \downarrow & & \downarrow & & \\ & & & & 0 & & 0 & & \end{array}$$

where $(M_i)_{i \geq k_0}$, $(M_i')_{i \geq k_0}$ and $(K_{i-k_0})_{i \geq k_0}$ are Mittag-Leffler surjective. By Lemma 3.17thm.3.17 that $((\theta M)_i)_{i \geq 1}$ is also Mittag-Leffler surjective. Taking limits in the above diagram, we obtain

$$
\begin{array}{ccccccc}
 & & 0 & & 0 & & \\
 & & \downarrow & & \downarrow & & \\
0 & \to & \theta M & \longrightarrow & M' & \longrightarrow & K & \longrightarrow & 0 \\
 & & \| & & \downarrow & & \downarrow & & \\
0 & \to & \theta M & \longrightarrow & M & \longrightarrow & M^{\min} & \to & 0 \\
 & & & & \downarrow & & \downarrow & & \\
 & & & & M^{\min}_{k_0} & = & M^{\min}_{k_0} & & \\
 & & & & \downarrow & & \downarrow & & \\
 & & & & 0 & & 0 & &
\end{array}
$$

Using Proposition 2.9thm.2.9 on the isogeny $0 \to M' \to M \to M^{\min}_{k_0} \to 0$, we have that $M' \in \mathrm{Mod}^{\varphi}_{\mathfrak{S},\,\mathrm{fr}}$ and the same proposition for $0 \to \theta M \to M' \to K \to 0$ tells us that $\theta M \in \mathrm{Mod}^{\varphi}_{\mathfrak{S},\,\mathrm{fr}}$, since we have already proved that $K \in \mathrm{Mod}^{\varphi}_{\mathfrak{S},\,\mathrm{fr}}$. Moreover, $M'/\theta M \simeq K \in \mathrm{Mod}^{\varphi}_{\mathfrak{S},\,\mathrm{fr}}$ and it is semi-stable of slope $\mu_{\min}(M)$.

2. If $Q$ is a nonzero quotient of $\theta M$ in $\mathrm{Mod}^{\varphi}_{\mathfrak{S},\,\mathrm{t}}$, there exists $i \geq 1$ such that $Q$ is a quotient of $(\theta M)_i$, therefore $\mu(Q) \geq \mu(\theta M_i) > \mu_{\min}(M_i) \geq \mu_{\min}(M)$. Thus, if $\theta M \neq 0$,

$$\mu_{\min}(\theta M) = \mu_{\min}(\theta M / p \theta M) > \mu_{\min}(M).$$

3. By remark 3rem.3, $M'$ is effective since $M/M' \simeq M^{\min}_{k_0}$ has no $u^\infty$-torsion and so $\theta M$ is effective since $M'/\theta M \simeq K$ which also has no $u^\infty$-torsion.

4. We have
$$
\begin{array}{rcl}
\theta(M(r)) & = & \ker((M, E^r \varphi_M) \to \varprojlim_n ((M, E^r \varphi_M)_n)^{\min}) \\
 & = & \ker((M, E^r \varphi_M) \to \varprojlim_n (M_n, E^r \varphi_{M_n})^{\min}) \\
 & = & \ker((M, E^r \varphi_M) \to \varprojlim_n (M_n^{\min}, E^r \varphi_{M_n^{\min}})) \\
 & = & (\ker(M \to \varprojlim_n M_n^{\min}), E^r \varphi_M) \\
 & = & (\theta M)(r)
\end{array}
$$
where the third equality is verified since twisting by $E^r$ does not change the Harder-Narasimhan filtration on $M_n$, it just change the slopes. $\square$

We get a family of Kisin modules $(M^{(i)})_{i \geq 0}$, where $M^{(i)} = \theta^{(i)} M$, together with inclusions

$$\iota_i \;:\; M^{(i+1)} \hookrightarrow M^{(i)}$$

verifying the properties of last proposition. This algorithm stops when $M^{(i)} = 0$ for some $i \geq 0$.

**Proposition 3.25.** *Let $M \in \mathrm{Mod}^{\varphi}_{\mathfrak{S},\,\mathrm{fr}}$. The algorithm above stops in a finite number of steps.*

*Proof.* We can reduce to the case when $M$ is an effective Kisin module for all $i \geq 1$. Suppose that the algorithm does not stop, then, by point (2) in last proposition $\mu_i = \mu_{\min}(\theta^{(i)} M)$ is strictly increasing for $i \geq 1$. However, for $e = \deg E$, we have

$$\mu_{\min}(\theta^{(i)} M) = \frac{\deg(\theta^{(i)} M / p \theta^{(i)} M)^{\min}}{\mathrm{rank}(\theta^{(i)} M / p \theta^{(i)} M)^{\min}} \in \frac{-\frac{1}{e}\mathbb{N}}{\mathrm{rank}(\theta^{(i)} M)} \subset -\frac{1}{e \cdot \mathrm{rank}\, M!}\mathbb{N},$$

since $(\theta^{(i)} M / p \theta^{(i)} M)^{\min}$ is effective by hypothesis. This is a discrete set, so we have a contradiction. $\square$

**Theorem 3.26.** *Every Kisin module is isogenuous to a HN-type Kisin module.*

*Proof.* Given three families $(M^{(i)})_{0 \leq i \leq r}$, $(M'^{(i)})_{0 \leq i \leq r-1}$ and $(Q^{(i)})_{0 \leq i \leq r-1}$ of Kisin modules verifying

- $M^{(0)} = M$,
- There is an isogeny between $M'^{(i)}$ and $M^{(i)}$ for all $0 \leq i \leq r-1$,
- For $0 \leq i \leq r-1$ we have exact sequences $0 \to M^{(i+1)} \to M'^{(i)} \to Q^{(i)} \to 0$,
- For $0 \leq i \leq r-1$, the Kisin module $Q^{(i)}$ is semi-stable of slope $\mu_i$ and we have $\mu_0 < \ldots < \mu_{r-1}$.

Then we can prove by induction on $r$ that $M$ is isogenuous to a HN-type Kisin module. Indeed, suppose that $f : M^{(2)} \to N$ is an isogeny with $N$ a HN-type Kisin module. So we have a diagram

$$\begin{array}{ccccccccc} 0 & \to & M^{(2)} & \to & M'^{(1)} & \to & Q^{(1)} & \to & 0 \\ & & \downarrow f & & \downarrow & & \downarrow \simeq & & \\ 0 & \longrightarrow & N & \longrightarrow & N' & \longrightarrow & Q^{(1)} & \to & 0 \end{array}$$

where $N'$ is the pushout of $M^{(2)} \to M'^{(1)}$ and $M^{(2)} \to N$. The isogeny $f$ gives us an isogeny between $M'^{(1)}$ and $N'$ and so an isogeny between $M = M^{(1)}$ and $N'$. If we put $\mathrm{Fil}_r N' = N'$, $\mathrm{Fil}_{r-1} N' = N$ and $\mathrm{Fil}_i N' = \iota(\mathrm{Fil}_i N)$ for all $1 \leq i \leq r-1$, where $\mathrm{Fil}_\bullet N$ is the HN flag defined on $N$, we obtain a HN-type Kisin module isogenuous to $M$.

It remains to show that such three families of modules exist for any Kisin module $M$, but this is equivalent to show that the algorithm given before stops, since that way we can take $M^{(i)} = \theta^{(i)} M$, $M'^{(i)}$ corresponds to the $M'$ given by 3.24thm.3.24 for $M^{(i)}$ and $Q^{(i)} = M'^{(i)}/\theta^{(i+1)} M$. By proposition 3.25thm.3.25, the algorithm stops. □

As a consequence, we have

**Corollary 3.27.** *Notations as above, we have* $\mathbf{t}_{\mathrm{F},\infty}(M) = \mathbf{t}_{\mathrm{F},\circ}(M)$ *in* $\mathbb{Q}_{\geq}^r$.

The following lemma gives us all the inequalities and equalities between the polygons presented before.

**Lemma 3.28.** *Let $M$ be a Kisin module. Then, we have inequalities of types*

$$\mathbf{t}_{\mathrm{F},\circ}(M) = \mathbf{t}_{\mathrm{F},\infty}(M) \leq \mathbf{t}_{\mathrm{F},1}(M) \leq \mathbf{t}_{\mathrm{H}}(M/pM) \leq \mathbf{t}_{\mathrm{H}}(M).$$

*Proof.* Last proposition together with Proposition 3.23thm.3.23 gives us the first equality. The inequality $\mathbf{t}_{\mathrm{F},\infty} \leq \mathbf{t}_{\mathrm{F},n}$ is given by Proposition 3.18thm.3.18 and Proposition 3.19thm.3.19. The inequality $\mathbf{t}_{\mathrm{F},1}(M) \leq \mathbf{t}_{\mathrm{H}}(M/pM)$ is given in Proposition Proposition 3.14thm.3.14 and the last inequality is given by Proposition 3.15thm.3.15. □

# 4 Crystalline representations with $G$-structure

## 4.1 The partially ordered commutative monoid $\mathbb{C}^\Gamma(G)$

For a reductive group $G$ over a base scheme $S$ and a totally ordered commutative group $\Gamma = (\Gamma, +, \leq) \neq 0$, Cornut constructs a sequence of $S$-schemes

$$\mathbb{G}^\Gamma(G) \xrightarrow{\mathrm{Fil}} \mathbb{F}^\Gamma(G) \xrightarrow{t} \mathbb{C}^\Gamma(G)$$

in [8, 2]. The construction is compatible with base change on $S$ and covariantly functorial in $G$ and $\Gamma$. On $\mathbb{G}^\Gamma(G)$ there is an involution $\iota$, compatible with the order and addition, inducing an involution on $\mathbb{C}^\Gamma(G)$. Moreover, $\mathbb{C}^\Gamma(G)$ is an étale, partially ordered commutative monoid over $S$. The partial order is the weak dominance order of [8, 2.2.12]. It is compatible with the functoriality in $G$ and $\Gamma$, but the monoid structure is only compatible with the functoriality in $\Gamma$. The functors represented by these schemes are related to $\Gamma$-graduations and $\Gamma$-filtrations on various fiber functors, as explained in [8, 3].

Suppose that $S$ is a connected normal scheme, $\Gamma$ is (uniquely) divisible and $T$ is a connected $S$-scheme. Then $\mathbb{C}^\Gamma(G)(S) \to \mathbb{C}^\Gamma(G)(T)$ is a monomorphism, with a canonical additive retraction $\sharp : \mathbb{C}^\Gamma(G)(T) \to \mathbb{C}^\Gamma(G)(S)$, functorial in the connected $S$-scheme $T$, by [8, 3.11.8]. If $T = s$ is a geometric point of $S$, corresponding to the fundamental group $\pi(S, s)$, then $\pi(S, s)$ acts on $\mathbb{C}^\Gamma(G)(s)$ with finite orbits, $\mathbb{C}^\Gamma(G)(S)$ is the fixed point set of this action and the retraction $\sharp : \mathbb{C}^\Gamma(G)(s) \to \mathbb{C}^\Gamma(G)(S)$ maps $x$ to the average of its orbit in the (uniquely) divisible monoid $\mathbb{C}^\Gamma(G)(s)$.

If $\mathrm{Spec}(R)$ is an affine $S$-scheme, we set

$$\begin{aligned} \mathbf{G}^\Gamma(G_R) &= \mathbb{G}^\Gamma(G)(R) \\ \mathbf{F}^\Gamma(G_R) &= \mathbb{F}^\Gamma(G)(R) \\ \mathbf{C}^\Gamma(G_R) &= \text{image of } \mathbb{F}^\Gamma(G)(R) \to \mathbb{C}^\Gamma(G)(R) \end{aligned}$$

Suppose that $R$ is local, as in [8, 4.1]. Then $\mathbf{C}^\Gamma(G_R)$ is a (partially ordered) commutative submonoid of $\mathbb{C}^\Gamma(G)(R)$. Moreover,
$$G_R \text{ is quasi-split} \iff \mathbf{C}^\Gamma(G_R) = \mathbb{C}^\Gamma(G)(R).$$

This follows from [8, 4.1.8] and [1, XXVI 3.8]. If $R$ is a valuation ring with fraction field $K$, then $\mathbf{C}^\Gamma(G_R) = \mathbf{C}^\Gamma(G_K)$ by [8, 4.1.18]. If $R$ is Henselian with residue field $k$, then $\mathbf{C}^\Gamma(G_R) = \mathbf{C}^\Gamma(G_k)$ by [8, 4.1.17], and $G_R$ is quasi-split if and only if $G_k$ is quasi-split. This is for instance the case if $k$ is finite (by Lang's theorem) or algebraically closed. If $G_R$ is split, then $\mathbb{C}^\Gamma(G_R)$ is the constant partially ordered commutative monoid $\mathbf{C}^\Gamma(G_R)_{\mathrm{Spec}(R)}$ ([8, 2.2.11]). Then for any morphism of local $S$-algebras $R \to R'$, the base change map $\mathbf{C}^\Gamma(G_R) \to \mathbf{C}^\Gamma(G_{R'})$ is an isomorphism. If $G = GL_n$, then $\mathbf{C}^\Gamma(G_R) = \Gamma_\geq^n$.

In the sequel, we will have a reductive group $G$ over $\mathbb{Z}_p$. Thus $G$ is quasi-split over $\mathbb{Z}_p$ and split over a finite unramified extension $\mathbb{Z}_{p^N}$ of $\mathbb{Z}_p$. In particular, applying $\mathbf{C}^\Gamma(G_-)$ to the diagram

$$\begin{array}{ccc}
 & W(\overline{\mathbb{F}}_p) \longrightarrow \overline{\mathbb{F}}_p & \\
 & \downarrow \qquad \downarrow & \\
K_0 \longleftarrow & W(\mathbb{F}) \longrightarrow \mathbb{F} & \\
\downarrow & \downarrow \qquad \downarrow & \\
 & \mathfrak{S} \longrightarrow \mathbb{F}[[u]] & \\
 & \downarrow & \\
K \longleftarrow & \mathfrak{S}_{(E)} &
\end{array}$$

all the arrows become isomorphisms. We will drop the index relative to the ring and note it simply by $\mathbf{C}^\Gamma(G)$ (this is an abuse of notation, but we will denote by $\mathbf{C}^\Gamma(G_{\mathbb{Z}_p})$ when we want to make clear that we are considering the group defined over $\mathbb{Z}_p$). The Frobenius $\sigma$ acts on it, with fixed point set $\mathbf{C}^\Gamma(G_{\mathbb{Z}_p}) = \mathbf{C}^\Gamma(G_{\mathbb{F}_p})$ and $\sigma^N \equiv \mathrm{Id}$ on $\mathbf{C}^\Gamma(G)$. The partial order on $\mathbf{C}^\Gamma(G)$ then has the following Tannakian caracterisation, which can be found in [8, 3.11.8], see also [11, 9.4.2].

**Proposition 4.1.** *Suppose $\Gamma$ is divisible and $G$ is defined over $\mathbb{Z}_p$. For every connected $\mathbb{Z}_p$-scheme $\mathrm{Spec}\,R$ and $\mathbf{t}_1, \mathbf{t}_2 \in \mathbf{C}^\Gamma(G_R)$, consider the following conditions:*

1. *$\mathbf{t}_1 \leq \mathbf{t}_2$ in $\mathbf{C}^\Gamma(G_R)$,*

2. *For every $\tau \in \mathrm{Rep}_R G$, $\mathbf{t}_1(\tau) \leq \mathbf{t}_2(\tau)$ in $\mathbf{C}^\Gamma(G_R)(\mathrm{GL}(V(\tau)) = \Gamma_\geq^{\mathrm{rank}_R(\tau)}$,*

3. *$\mathbf{t}_1^\# \leq \mathbf{t}_2^\#$ in $\mathbf{C}^\Gamma(G_{\mathbb{Z}_p})$,*

4. *For every $\tau \in \mathrm{Rep}_{\mathbb{Z}_p} G$, $\mathbf{t}_1(\tau_R) \leq \mathbf{t}_2(\tau_R)$ in $\mathbf{C}^\Gamma(G)(\mathrm{GL}(V(\tau_R)) = \Gamma_\geq^{\mathrm{rank}_{\mathbb{Z}_p}(\tau)}$.*

*Then, we have (1) $\iff$ (2), (4) $\implies$ (3) and we have (3) $\iff$ (4) if $\mathbf{t}_1 = \mathbf{t}_1^\#$.*

## 4.2 Fiber functors

In this section, let $G$ be a reductive group defined over $\mathcal{O}$, for $\mathcal{O} \in \{\mathbb{Q}_p, \mathbb{Z}_p, \mathbb{F}_p\}$. For a scheme $X$, we define the category $\mathrm{Bun}_X$ as the category of finite locally free sheaves of $\mathcal{O}_X$-modules. In particular, for $X = \mathrm{Spec}\,R$, it corresponds to the category of finite projective $R$-modules and we denote this category by $\mathrm{Bun}_R$. We define the category $\mathrm{Rep}_\mathcal{O} G$ as the category of algebraic representations of $G$ on finite free $\mathcal{O}$-modules. We say that a functor
$$\omega_X : \mathrm{Rep}_\mathcal{O} G \to \mathrm{Bun}_X$$
is a fiber functor when $\omega_X$ is an exact $\otimes$-functor, for a scheme over $\mathrm{Spec}\,\mathcal{O}$. We denote by $\mathrm{Aut}^\otimes(\omega_X)$ the group of $\otimes$-automorphisms of $\omega_X$. We denote by $\mathbf{G}^\Gamma(\omega_R)$ and $\mathbf{F}^\Gamma(\omega_R)$ the set of $\Gamma$-graduations and $\Gamma$-filtrations on $\omega_R$, respectively, as defined in [8, 3.2].

Let $\mathrm{Gr}^\Gamma \mathrm{Bun}_R$ is the category of finite projective $R$-modules endowed with a $\Gamma$-graduation and $\mathrm{Fil}^\Gamma \mathrm{Bun}_R$ is the category of finite projective $R$-modules endowed with a $\Gamma$-filtration by direct summands, as defined in [8, 3.3.1]. There is an exact $\otimes$-functor $\mathrm{Fil} : \mathrm{Gr}^\Gamma \mathrm{Bun}_R \to \mathrm{Fil}^\Gamma \mathrm{Bun}_R$ sending a graduation $\mathcal{G}$ to the filtration $\mathrm{Fil}(\mathcal{G}) = \mathcal{F}$ defined by $\mathcal{F}^{\geq \gamma} = \oplus_{\eta \geq \gamma} \mathcal{G}^\eta$. It induces an $\mathrm{Aut}^\otimes(\omega_R)$-equivariant map

$$\mathrm{Fil} : \mathbf{G}^\Gamma(\omega_R) \to \mathbf{F}^\Gamma(\omega_R).$$

We define the types on $\omega_R$ by
$$\mathbf{C}^\Gamma(\omega_R) := \mathrm{Aut}^\otimes(\omega_R)\backslash \mathbf{F}^\Gamma(\omega_R).$$

The three sets $\mathbf{G}^\Gamma(\omega_R)$, $\mathbf{F}^\Gamma(\omega_R)$ and $\mathbf{C}^\Gamma(\omega_R)$ are compatible with change of the $\mathcal{O}$-algebra $R$, so we can define three presheaves $\mathbf{G}^\Gamma(\omega)$, $\mathbf{F}^\Gamma(\omega)$ and $\mathbf{C}^\Gamma(\omega)$ together with maps

$$\mathbf{G}^\Gamma(\omega_R) \xrightarrow{\mathrm{Fil}} \mathbf{F}^\Gamma(\omega_R) \xrightarrow{\mathbf{t}} \mathbf{C}^\Gamma(\omega_R).$$

Now, consider the standard fiber functor

$$\omega_{G,R} \ : \ \mathrm{Rep}_\mathcal{O} G \to \mathrm{Bun}_R$$

given by the composition of the forgetting functor $\omega_G : \mathrm{Rep}_\mathcal{O} G \to \mathrm{Bun}_\mathcal{O}$ with the tensorization by $R$. This fiber functor is called the trivial fiber functor. The next proposition gives a criterion for the existence of a $\otimes$-isomorphism $\omega_{G,R} \simeq \omega_R$.

**Proposition 4.2.** *Let $R$ be a local strictly henselian and faithfully flat $\mathcal{O}$-algebra. Then, any exact and faithful $\otimes$-functor $\omega_R : \mathrm{Rep}_\mathcal{O} G \to \mathrm{Bun}_R$ is $\otimes$-isomorphic to the standard fiber functor $\omega_{G,R}$.*

*Proof.* We can use the results in [6], where Broshi proved that there is an equivalence of categories between exact and faithful $\otimes$-functors taking values in $\mathrm{Bun}_R$ and $G$-torsors over $\mathrm{Spec}(R)$. Now, $G$-torsors over $\mathrm{Spec}(R)$ are classified by the étale cohomology group $H^1(\mathrm{Spec}(R), G)$. We have that

$$H^1(\mathrm{Spec}(R), G) \simeq H^1(\mathbb{F}, G_\mathbb{F}) \simeq 0$$

under our hypothesis, the first isomorphism using [1, Exposé XXIV, Proposition 8.1], and the second one s given in [42, Theorem 1.9]. Thus, all $G$-torsors over $R$ are trivial, and $\omega_R \simeq \omega_{G,R}$ as wanted, by the equivalence of categories given by Broshi. $\square$

The choice of a $\otimes$-isomorphism $\omega_R \simeq \omega_{G,R}$ (when it exists) yields isomorphisms

$$\begin{array}{rcl}
\mathbf{G}^\Gamma(\omega_R) & \simeq & \mathbf{G}^\Gamma(\omega_{G,R}) \ =: \ \mathbf{G}^\Gamma(G_R) \\
\mathbf{F}^\Gamma(\omega_R) & \simeq & \mathbf{F}^\Gamma(\omega_{G,R}) \ =: \ \mathbf{F}^\Gamma(G_R) \\
\mathbf{C}^\Gamma(\omega_R) & \simeq & \mathbf{C}^\Gamma(\omega_{G,R}) \ =: \ \mathbf{C}^\Gamma(G_R)
\end{array}$$

so, in this case, we can use the results in [8]. The first two isomorphisms are not canonical, since they depend on the chosen isomorphism $\omega \simeq \omega_{G,R}$. Let $\eta, \eta'$ be two different isomorphisms $\omega \simeq \omega_{G,R}$, then we have $\eta' = \kappa \eta$ for some $\kappa \in \mathrm{Aut}^\otimes(\omega_{G,R}) = G(R)$. Then, the induced isomorphisms

$$\eta, \eta' \ : \ \mathbf{C}^\Gamma(\omega_R) = \mathrm{Aut}^\otimes(\omega_R)\backslash \mathbf{G}^\Gamma(\omega_R) \simeq G(R)\backslash \mathbf{G}^\Gamma(G_R) = \mathbf{C}^\Gamma(G_R)$$

are equal, so $\mathbf{C}^\Gamma(\omega)$ is canonically isomorphic to $\mathbf{C}^\Gamma(G_R)$.

The definition of a filtration on a fiber functor implies that the filtration is compatible with exterior and symmetric powers. The next proposition shows that the converse is also sometimes true.

**Proposition 4.3.** *Let $L$ be a field which is an $\mathcal{O}$-algebra. Suppose that a fiber functor $\omega_L : \mathrm{Rep}_\mathcal{O} G \to \mathrm{Vect}_L$ admits a factorization through an additive $\otimes$-functor*

$$\mathcal{F} \ : \ \mathrm{Rep}_\mathcal{O} G \to \mathrm{Fil}^\Gamma_L$$

*which is compatible with $\Lambda$ and $\mathrm{Sym}$. Then $\mathcal{F}$ is exact, thus $\mathcal{F}$ is a filtration of $\omega_L$.*

*Proof.* We say that $\mathcal{F}$ is exact if it verifies the following equivalent conditions. For every exact sequence

$$0 \to \tau_1 \to \tau_2 \xrightarrow{\pi} \tau_3 \to 0$$

with $\tau_1, \tau_2, \tau_3 \in \mathrm{Rep}_\mathcal{O} G$:

- The sequence $0 \to \mathcal{F}^\gamma(\tau_1) \to \mathcal{F}^\gamma(\tau_2) \to \mathcal{F}^\gamma(\tau_3) \to 0$ is exact, for every $\gamma \in \Gamma$.

- We have $\mathcal{F}^\gamma(\tau_1) = \mathcal{F}^\gamma(\tau_2) \cap \omega_L(\tau_1)$, and $\mathcal{F}^\gamma(\tau_3) = \pi(\mathcal{F}^\gamma(\tau_2))$, for every $\gamma \in \Gamma$.

- The functor $\mathcal{F}$ transforms strict monomorphisms (respectively strict epimorphisms) into strict monomorphisms (respectively, strict epimorphisms).

We will use the last characterization of an exact functor. By [13, II, Proposition 1.9], an additive ⊗-functor between two rigid categories is also compatible with duality, i.e. $\mathcal{F}(\tau)^\vee = \mathcal{F}(\tau^\vee)$. This means that we only need to check that $\mathcal{F}$ transforms strict monomorphisms into strict monomorphisms.

Let $\tau_1 \hookrightarrow \tau_2$ be a strict monomorphism. Set $\mathcal{G}^\gamma(\tau_1) = \mathcal{F}^\gamma(\tau_2) \cap \omega_L(\tau_1)$. Since $\mathcal{F}$ is a functor, we have a commutative diagram in $\mathrm{Fil}_L^\Gamma$

$$\mathcal{F}(\tau_1) \xrightarrow{\iota} \mathcal{F}(\tau_2)$$
$$\searrow^{i_1} \quad \nearrow^{i_2}$$
$$\mathcal{G}(\tau_1)$$

where $i_1$ is a mono-epi and $i_2$ is a strict mono. We want to show that $\mathcal{F}(\tau_1) = \mathcal{G}(\tau_1)$, i.e. that $i_1$ is an isomorphism. First, we prove that we can reduce to the case where $\mathrm{rank}_\mathcal{O}\,\tau_1 = 1$. Indeed, let $d = \mathrm{rank}_\mathcal{O}\,\tau_1$ and consider $\Lambda^d \tau_1 \hookrightarrow \Lambda^d \tau_2$, applying $\mathcal{F}$, and since $\mathcal{F}$ is compatible with exterior powers, we get a diagram

$$\Lambda^d \mathcal{F}(\tau_1) \xrightarrow{\Lambda^d \iota} \Lambda^d \mathcal{F}(\tau_2)$$
$$\searrow^{\Lambda^d i_1} \quad \nearrow^{\Lambda^d i_2}$$
$$\Lambda^d \mathcal{G}(\tau_1)$$

where $\Lambda^d i_1$ is again a mono-epi by the properties of the exterior power, and by a direct calculation considering a splitting, we can check that $\Lambda^d i_2$ is a strict monomorphism. Since

$$\mu(\Lambda^d \mathcal{F}(\tau_1)) = d\mu(\mathcal{F}(\tau_1)) \quad \text{and} \quad \mu(\Lambda^d \mathcal{G}(\tau_1)) = d\mu(\mathcal{G}(\tau_1))$$

we obtain that $i_1$ is an isomorphism if and only if $\Lambda^d i_1$ is an isomorphism.

Now, we have $\tau_1 \hookrightarrow \tau_2$ where $\tau_1$ is a character of $G$. Since $\mathcal{F}$ is compatible with tensor products and duality, twisting the initial sequence by the inverse of that character, we may assume that $\tau_1 = \mathbf{1}$.

Let thus $\mathbf{1} \hookrightarrow \tau$ be a strict monomorphism. The generalization of Haboush's theorem given by Seshadri in [40, Theorem 1] tells us that there is an $r \geq 1$ such that $\mathbf{1} = \mathrm{Sym}^r \mathbf{1} \hookrightarrow \mathrm{Sym}^r(\tau)$ is split. This implies that $\mathcal{F}(\mathbf{1}) \to \mathcal{F}(\mathrm{Sym}^r(\tau)) = \mathrm{Sym}^r(\mathcal{F}(\tau))$ is a strict monomorphism. Since $\mathrm{Sym}^r i_2$ is also a strict mono, the mono-epi

$$\mathrm{Sym}^r i_1 \ : \ \mathrm{Sym}^r \mathcal{F}(\mathbf{1}) \to \mathrm{Sym}^r \mathcal{G}(\mathbf{1})$$

is an isomorphism. Then, we have $\mathcal{F}(\mathbf{1}) \simeq \mathcal{G}(\mathbf{1})$, since

$$\mu(\mathrm{Sym}^r \mathcal{F}(\tau_1)) = r\mu(\mathrm{Sym}^r(\tau_1)) \quad \text{and} \quad \mu(\mathrm{Sym}^r \mathcal{G}(\tau_1)) = r\mu(\mathcal{G}(\tau_1)).$$

□

*Remark* 9. If $\mathcal{O} = \mathbb{Q}_p$, exact sequences in $\mathrm{Rep}_\mathcal{O}\,G$ are split and the exactness of $\mathcal{F}$ follows from its additivity.

### 4.2.1 Lattices

Let $\mathcal{O}_L$ be an $\mathcal{O}$-algebra which is a discrete valuation ring, $L$ its fraction field, $\pi_L$ a uniformizer and $l$ its residue field. For $G$ a reductive group over $\mathcal{O}$, let

$$V \ : \ \mathrm{Rep}_\mathcal{O}\,G \to \mathrm{Bun}_L$$

be a ⊗-functor. We denote by $\mathcal{L}'(V)$ the set of ⊗-functors

$$x \ : \ \mathrm{Rep}_\mathcal{O}\,G \to \mathrm{Bun}_{\mathcal{O}_L}$$

such that

1. For every $\tau \in \mathrm{Rep}_\mathcal{O}\,G$, its image $x(\tau)$ is an $\mathcal{O}_L$-lattice in $V(\tau)$,

2. For every $\tau_1 \xrightarrow{f} \tau_2$, the morphism $x(f) : x(\tau_1) \to x(\tau_2)$ is induced by $V(f) : V(\tau_1) \to V(\tau_2)$.

Then, there is a canonical $\otimes$-isomorphism $x_L \simeq V$, where $x_L = x \otimes L$.

Conversely, to give an element of $\mathcal{L}'(V)$ amounts to give an isomorphism class of pairs $(x', \eta)$ where $x' : \text{Rep}_\mathcal{O} G \to \text{Bun}_{\mathcal{O}_L}$ is a $\otimes$-functor and $\eta$ is an isomorphism of $\otimes$-functors $x'_L \simeq V$: we associate to $(x', \eta)$ the functor $x(\tau) = \eta_\tau(x'(\tau))$.

We denote by $\mathcal{L}^{\text{ex}}(V) \subset \mathcal{L}'(V)$ the subset of lattices $x$ which are exact. Then $\mathcal{L}^{\text{ex}}(V) \neq \emptyset$ implies that $V$ is exact. We suppose $V$ exact henceforth, i.e. $V : \text{Rep}_\mathcal{O} G \to \text{Bun}_L$ is a fiber functor. Denote by $\mathcal{L}(V) \subset \mathcal{L}^{\text{ex}}(V)$ the subset of lattices $x$ which are isomorphic to the standard fiber functor $\omega_{G, \mathcal{O}_L}$ (they are exact since the latter is exact). Then, if $\mathcal{L}(V) \neq \emptyset$, we have $V \simeq \omega_{G,L}$ and, conversely, if $V \simeq \omega_{G,L}$, then we have $\omega_{G, \mathcal{O}_L} \in \mathcal{L}(\omega_{G,L}) \simeq \mathcal{L}(V) \neq \emptyset$. We suppose $V \simeq \omega_{G,L}$ henceforth.

*Remark* 10. In the case when $\mathcal{O}_L$ is a local strictly henselian faithfully flat $\mathcal{O}$-algebra and $V : \text{Rep}_\mathcal{O} G \to \text{Bun}_L$ is faithful, we get $\mathcal{L}(V) = \mathcal{L}^{\text{ex}}(V)$, using Proposition 4.2thm.4.2, since $V$ faithful implies $x$ faithful for every $x \in \mathcal{L}'(V)$.

Let $\eta : V_1 \to V_2$ be a $\otimes$-isomorphism. Then $\eta$ induces a bijection $\mathcal{L}^*(V_1) \xrightarrow{\sim} \mathcal{L}^*(V_2)$, by $x \mapsto \eta x$ where $\eta x(\tau) = \eta_\tau(x(\tau))$ (so $\eta$ is also a $\otimes$-isomorphism $x \simeq \eta x$). In particular, the group $\text{Aut}^\otimes(V)$ acts on $\mathcal{L}(V) \subset \mathcal{L}^{\text{ex}}(V) \subset \mathcal{L}'(V)$. The action is transitive on $\mathcal{L}(V)$, since two elements $x, y \in \mathcal{L}(V)$ are isomorphic as functors, through a $\otimes$-isomorphism $\eta : x \to y$ inducing a $\otimes$-isomorphism $\eta_L : V \to V$, i.e. an element $g \in \text{Aut}^\otimes(V)$, which sends $x$ to $y$ by definition. Thus for any $x \in \mathcal{L}(V)$ with stabilizer $\text{Aut}^\otimes(x)$ in $\text{Aut}^\otimes(V)$, the map $g \mapsto g \cdot x$ yields an $\text{Aut}^\otimes(V)$-equivariant bijection

$$\text{Aut}^\otimes(V) / \text{Aut}^\otimes(x) \simeq \mathcal{L}(V).$$

*Remark* 11.  1. For $V = \omega_{G,L}$ and $x = \omega_{G, \mathcal{O}_L}$, we obtain $G(L)/G(\mathcal{O}_L) \simeq \mathcal{L}(V)$, since $\text{Aut}^\otimes(x) = G(\mathcal{O}_L)$ inside $\text{Aut}^\otimes(V) = G(L)$.

2. For $\tilde{G} = G_{\mathcal{O}_L}$, $V = \omega_{G,L}$ and $\tilde{V} = \omega_{\tilde{G}, L}$, the map $\mathcal{L}(\tilde{V}) \to \mathcal{L}(V)$ defined by $\tilde{x} \mapsto x$, where

$$x \ : \ \text{Rep}_\mathcal{O} G \xrightarrow{\otimes \mathcal{O}_L} \text{Rep}_{\mathcal{O}_L} \tilde{G} \xrightarrow{\tilde{x}} \text{Bun}_{\mathcal{O}_L}$$

is $G(L) = \text{Aut}^\otimes(\tilde{V}) = \text{Aut}^\otimes(V)$-equivariant, thus a bijection since $\mathcal{L}(\tilde{V}) \simeq G(L)/G(\mathcal{O}_L) \simeq \mathcal{L}(V)$.

3. For $\mathcal{O} = \mathcal{O}_L$ Henselian and $V = \omega_{G,L}$, there is a $G(L)$-equivariant embedding

$$\mathcal{L}(V) \hookrightarrow \mathbf{B}^e(\omega_G, L)$$

mapping $x$ to the gauge norms $\alpha_x \in \mathbf{B}^e(\omega_G, L)$ defined by $\alpha_x(\tau)(v) = \inf\{|\lambda| \mid \lambda \in L, \ v \in \lambda x(\tau)\}$ for $\tau \in \text{Rep}_\mathcal{O} G$ and $v \in \omega_{G,L}(\tau)$ and $\mathbf{B}^e(\omega_G, L)$ is the space of $L$-norms on $\omega_G$ defined in [8, 6.4].

**Definition 4.1.** Let $\mathbf{F}^\mathbb{Z}(V) = \mathbf{F}^\mathbb{Z}(\omega(V))$. We define an addition operator between lattices and filtrations

$$\begin{array}{rcl} + \ : \ \mathcal{L}(V) \times \mathbf{F}^\mathbb{Z}(V) & \to & \mathcal{L}(V) \\ (x, \mathcal{F}) & \mapsto & x + \mathcal{F} \end{array}$$

where, for every $\tau \in \text{Rep}_{\mathbb{Z}_p} G$, we have $(x + \mathcal{F})(\tau) = x(\tau) + \mathcal{F}(\tau \otimes \mathbb{Q}_p) = \sum_{i \in \mathbb{Z}} \pi_L^{-1} x(\tau) \cap \mathcal{F}^{\geq i}(\tau \otimes \mathbb{Q}_p)$.

**Proposition 4.4.** *The operator defined above is well-defined.*

*Proof.* We may assume that $V$ is trivial. Fix $x \in \mathcal{L}(V)$ and $\mathcal{F} \in \mathbf{F}^\mathbb{Z}(V)$. We know that $x + \mathcal{F}$ is functorial and, by Proposition 2.5thm.2.5, compatible with tensor products. We have then an element $x + \mathcal{F} \in \mathcal{L}'(V)$, so it suffices to see that it is isomorphic to the trivial lattice. Up to multiplication with an element of $G(L)$, we may assume that $x = \omega_{G, \mathcal{O}_L}$. Since $\mathcal{F} \in \mathbf{F}^\mathbb{Z}(G_L) = \mathbf{F}^\mathbb{Z}(G_{\mathcal{O}_L})$, there exists a cocharacter $\chi : \mathbb{G}_{m, \mathcal{O}_L} \to G_{\mathcal{O}_L}$ splitting $\mathcal{F}$. For every $\tau \in \text{Rep}_{\mathbb{Z}_p} G$, let $x(\tau) = \oplus_{i \in \mathbb{Z}} x(\tau)_i$ be the weight decomposition of $(\tau \otimes \mathcal{O}_L) \circ \chi : \mathbb{G}_{m, \mathcal{O}_L} \to \text{GL}_{\mathcal{O}_L}(x(\tau))$. Then $\mathcal{F}^{\geq j}(\tau \otimes \mathbb{Q}_p) = \oplus_{i \geq j} V(\tau)_j$ where $V(\tau)_j = x(\tau)_j \otimes L$, so

$$\begin{array}{rcl} x(\tau) + \mathcal{F}^{\geq j}(\tau \otimes \mathbb{Q}_p) & = & \sum_j \pi_L^{-j} x(\tau) \cap \mathcal{F}^{\geq j}(\tau \otimes \mathbb{Q}_p) \\ & = & \sum_j \oplus_{i \geq j} \pi_L^{-j} x(\tau)_i \\ & = & \oplus_i \sum_{j \geq i} \pi_L^{-j} x(\tau)_i \\ & = & \oplus_i \pi_L^{-i} x(\tau)_i \\ & = & \tau(\chi(\pi_L^{-1})) \cdot x(\tau) \\ & = & (\chi(\pi_L^{-1}) \cdot x)(\tau), \end{array}$$

and $x + \mathcal{F} = \chi(\pi_L^{-1}) \cdot x$, which indeed belongs to $\mathcal{L}(V)$. $\square$

### 4.2.2 A variant: From $\operatorname{Rep}_{\mathbb{Q}_p} G$ to $\operatorname{Rep}_{\mathbb{Z}_p} G$

For the rest of the chapter, we will use a slightly different framework. If $\mathcal{O} = \mathbb{Z}_p$, $G$ a reductive group over $\mathbb{Z}_p$ and $V : \operatorname{Rep}_{\mathbb{Q}_p} G_{\mathbb{Q}_p} \to \operatorname{Bun}_L$ is a fiber functor, then we denote

$$\mathcal{L}(V) = \mathcal{L}(V')$$

where

$$V' \ : \ \operatorname{Rep}_{\mathbb{Z}_p} G \to \operatorname{Bun}_L$$

is the fiber functor induced by precomposition of $V$ with $\operatorname{Rep}_{\mathbb{Z}_p} G \to \operatorname{Rep}_{\mathbb{Q}_p} G$. We have

$$\begin{aligned}\mathcal{L}(V) \neq \emptyset \ &\Leftrightarrow \ V \text{ is isomorphic to the trivial fiber functor } \omega_{G_{\mathbb{Q}_p},L} \\ &\Leftrightarrow \ V' \text{ is isomorphic to the trivial fiber functor } \omega_{G,L}.\end{aligned}$$

Indeed, we already know that $\mathcal{L}(V) \neq \emptyset$ if and only if $\mathcal{L}(V') \neq \emptyset$ if and only if $V' \simeq \omega_{G,L}$. Suppose $V \simeq \omega_{G_{\mathbb{Q}_p},L}$, then it is obvious that $V \simeq \omega_{G,L}$. Conversely, suppose we have an element $a : V' \simeq \omega_{G,L}$ in $\operatorname{Iso}^\otimes(V', \omega_{G,L})$. By [12, Theorem 1.12 and Remark 1.13], there exists a finite Galois extension $L'/L$ such that there exists an isomorphism $b : V_{L'} \simeq \omega_{G_{\mathbb{Q}_p}, L'}$. Consider the diagram

$$\begin{array}{ccc} \operatorname{Iso}^\otimes(V, \omega_{G_{\mathbb{Q}_p},L}) & \longrightarrow & \operatorname{Iso}^\otimes(V', \omega_{G,L}) \ni a \\ \downarrow & & \downarrow \\ b \in \operatorname{Iso}^\otimes(V_{L'}, \omega_{G_{\mathbb{Q}_p},L'}) & \longrightarrow & \operatorname{Iso}^\otimes(V'_{L'}, \omega_{G,L'}). \end{array}$$

Changing $b$ by an element $g \in G_{\mathbb{Q}_p}(L') = G(L') = \operatorname{Aut}^\otimes(\omega_{G,L'})$, we may assume that $b$ and $a$ have the same image in $\operatorname{Iso}^\otimes(V'_{L'}, \omega_{G,L'})$. Then $b$ is fixed by the action of $\operatorname{Gal}(L'/L)$ on $\operatorname{Iso}^\otimes(V_{L'}, \omega_{G_{\mathbb{Q}_p},L'})$, thus $b \in \operatorname{Iso}^\otimes(V, \omega_{G_{\mathbb{Q}_p},L})$ and indeed $V \simeq \omega_{G_{\mathbb{Q}_p},L}$.

### 4.2.3 Vectorial distance

For $x, y \in \mathcal{L}'(V)$, we denote by

$$\mathcal{F}(x,y) \ : \operatorname{Rep}_\mathcal{O} G \to \operatorname{Fil}^\mathbb{Z} \operatorname{Bun}_l$$

the functor given by $\mathcal{F}(x,y)(\tau) = \mathcal{F}(x(\tau), y(\tau))$, with underlying fiber functor $\overline{x} = x \otimes l$. It is clear that $\overline{x}$ is a $\otimes$-functor, and so is $\mathcal{F}(x,y)$, by Proposition 2.3thm.2.3.

**Proposition 4.5.** *Suppose $x, y \in \mathcal{L}^{\operatorname{ex}}(V)$, then the $\otimes$-functor $\mathcal{F}(x,y)$ is exact, i.e.*

$$\mathcal{F}(x,y) \in \mathbf{F}^\mathbb{Z}(\overline{x}).$$

*Proof.* Let $x, y \in \mathcal{L}^{\operatorname{ex}}(V)$. Then, $\overline{x}$ is exact and, by Proposition 4.3thm.4.3, it suffices to show that $\mathcal{F}(x,y)$ is compatible with exterior and symmetric powers. We check it for exterior powers. We have

$$\Lambda^r \mathcal{F}(x(\tau), y(\tau)) = \mathcal{F}(\Lambda^r x(\tau), \Lambda^r y(\tau)) = \mathcal{F}(x(\Lambda^r \tau), y(\Lambda^r \tau))$$

for every $\tau \in \operatorname{Rep}_\mathcal{O} G$ and every $r \geq 0$, where the first equality is given by Proposition 2.3thm.2.3 and the second equality follows from the exactness of $x$ and $y$. Thus $\mathcal{F}(x,y)$ is exact. □

Using this filtration, we define the vectorial distance (also called relative position) of $x, y \in \mathcal{L}(V)$, as

$$\operatorname{Pos}(x,y) = \mathbf{t}(\mathcal{F}(x,y)) \in \mathbf{C}^\mathbb{Z}(\overline{x}) = \mathbf{C}^\mathbb{Z}(G_l).$$

(where the equality $\mathbf{C}^\mathbb{Z}(\overline{x}) = \mathbf{C}^\mathbb{Z}(G_l)$ is due to the fact that $\overline{x}$ is isomorphic to the trivial fiber functor).

**Example 4.1.** Let $V = \omega_{G,L}$, $x = \omega_{G,\mathcal{O}_L}$ and $y = \mu(\pi_L) \cdot \omega_{G,\mathcal{O}_L}$, where $\mu \ : \ \mathbb{G}_{m,\mathcal{O}_L} \to G_{\mathcal{O}_L}$. We also denote by $\mu$ the type of the graduation given by $\mu$. Then, we have

$$\operatorname{Pos}(x,y) = \mu^\iota \quad \text{in} \quad \mathbf{C}^\mathbb{Z}(G_l).$$

*Remark* 12. If $\mathcal{O}_L$ is Henselian, with the proper normalization of the multiplicative valuation on $L$, there is a commutative diagram

$$\begin{array}{ccccc}
\mathcal{L}(\omega_{G,L}) & \times & \mathcal{L}(\omega_{G,L}) & \xrightarrow{\text{Pos}} & \mathbf{C}^{\mathbb{Z}}(G_l) \\
\downarrow & & \downarrow & & \downarrow \\
\mathbf{B}^e(\omega_{\tilde{G}},L) & \times & \mathbf{B}^e(\omega_{\tilde{G}},L) & \xrightarrow{\mathbf{d}} & \mathbf{C}^{\mathbb{R}}(G_L)
\end{array}$$

where $\tilde{G} = G_{\mathcal{O}_L}$, $\tilde{V} = \omega_{\tilde{G},L}$ and $\mathcal{L}(V) \leftarrow \mathcal{L}(\tilde{V}) \hookrightarrow \mathbf{B}^e(\omega_{\tilde{G}},L)$ are as above and $\mathbf{C}^{\mathbb{Z}}(G_l) \simeq \mathbf{C}^{\mathbb{Z}}(G_{\mathcal{O}_L}) \simeq \mathbf{C}^{\mathbb{Z}}(G_L) \subset \mathbf{C}^{\mathbb{R}}(G_L)$, by [8, 4.1.17] and [8, 4.1.18]. The vectorial distance $\mathbf{d}$ is introduced in [8, 5.2.8], using [8, 6.2] and [8, 6.4]. We explain the diagram a bit more in detail. For the calculation of relative positions, we may assume that $\mathcal{O} = \mathcal{O}_L$ and $V = \omega_{G,L}$. Both maps $\mathcal{L}(\omega_{G,L})^2 \to \mathbf{C}^{\mathbb{R}}(G_L)$ are $G(L)$-equivariant, so we may fix the first component of $(x,y) \in \mathcal{L}(\omega_{G,L})^2$ to $x = \omega_{G,\mathcal{O}_L}$, which maps to $\alpha_{G,L}^{\circ}$ in $\mathbf{B}^{\circ}(\omega_G, L)$. By [8, 6.4.8], it suffices to check the commutativity of

$$\begin{array}{cccc}
\text{Pos}(\omega_{G,\mathcal{O}_L}, -): & \mathcal{L}(\omega_{G,L}) & \to & \mathbf{F}^{\mathbb{Z}}(G_l) \\
& \downarrow & & \downarrow \\
\text{loc}: & \mathbf{B}^e(\omega_G, L) & \to & \mathbf{F}^{\mathbb{R}}(G_l)
\end{array}$$

The two maps $\mathcal{L}(\omega_{G,L}) \to \mathbf{F}^{\mathbb{R}}(G_l)$ on the diagram send $y \in \mathcal{L}(\omega_{G,L})$ to the filtration defined by

$$\mathcal{F}_1^i(\tau) = \frac{\omega_{G,\mathcal{O}_L}(\tau) \cap \pi_L^i y(\tau) + \pi_L \omega_{G,\mathcal{O}_L}(\tau)}{\pi_L \omega_{G,\mathcal{O}_L}} \subset \omega_l(\tau)$$

$$\mathcal{F}_2^{\gamma}(\tau) = \frac{\omega_{G,\mathcal{O}_L}(\tau) \cap \overline{B}(\alpha_y(\tau),\gamma) + \pi_L \omega_{G,\mathcal{O}_L}(\tau)}{\pi_L \omega_{G,\mathcal{O}_L}} \subset \omega_l(\tau)$$

for every $i \in \mathbb{Z}$, every $\gamma \in \mathbb{R}$ and every $\tau \in \text{Rep}_{\mathbb{Z}_p} G$, where $\overline{B}(\alpha_y(\tau),\gamma) = \{v \mid \alpha_y(\tau)(v) \leq \exp(-\gamma)\}$. If $|L^{\times}| = \exp(\mathbb{Z})$, i.e. $|\pi_L| = \exp(-1)$, then $\mathcal{F}_2 \in \mathbf{F}^{\mathbb{Z}}(G_l)$ and

$$\overline{B}(\alpha_y(\tau),i) = \{v \mid \alpha_y(\tau)(v) \leq \exp(-i)\} = \pi_L^i y(\tau),$$

therefore $\mathcal{F}_1 = \mathcal{F}_2$.

**Proposition 4.6.** *The relative position defined above verifies*

1. *The triangular inequality: for every $x,y,z \in \mathcal{L}(V)$, we have $\text{Pos}(x,z) \leq \text{Pos}(x,y) + \text{Pos}(y,z)$ in the partially ordered commutative monoid $\mathbf{C}^{\mathbb{Z}}(G_l)$.*

2. *Compatibility with the involution: For every $x,y \in \mathcal{L}(V)$, we have $\text{Pos}(y,x) = \text{Pos}(x,y)^{\iota}$.*

*Proof.* We may assume that $V = \omega_{G,L}$. We can now use all the results in [8] as follows: as we have seen above, the functor sending a lattice to its Gauge norm embeds $\mathcal{L}(V)$ in $\mathbf{B}(\omega_{\tilde{G}}, L)$ and in the last remark we have seen that for a suitable normalizaton of the valuation on $L$, the operator Pos corresponds to the vectorial distance $\mathbf{d}$ defined on $\mathbf{B}(\omega_{\tilde{G}}, L)$ via 6.4.10 and 6.2 in [8], which verifies itself the triangular inequality. The compatibility with the involution is an obvious property of the vectorial distance $\mathbf{d}$ given in [8, 5.2.8]. □

Another property that will be useful is given by the next proposition:

**Proposition 4.7.** *Let $\varphi: L \to L'$ be a finite extension with ramification index $e$ and set $V' = V \otimes L'$. Then we have*
$$\mathcal{F}(x', y') = e \cdot \mathcal{F}(x,y) \otimes l' \quad \text{in} \quad \mathbf{F}^{\mathbb{Z}}(\overline{x}') = \mathbf{F}^{\mathbb{Z}}(\overline{x}_{l'})$$

*for every $x,y \in \mathcal{L}(V)$, $x' = x \otimes \mathcal{O}_{L'}$, $y' = y \otimes \mathcal{O}_{L'}$ in $\mathcal{L}(V')$ and $l'$ the residue field of $\mathcal{O}_{L'}$. As a consequence, we have*
$$\text{Pos}(x', y') = e \cdot \varphi \text{Pos}(x,y) \quad \text{in} \quad \mathbf{C}^{\mathbb{Z}}(G_{l'}).$$

*Proof.* It suffices to prove that $\mathcal{F}(x'(\tau), y'(\tau)) = e \cdot \mathcal{F}(x(\tau), y(\tau)) \otimes l'$ for every $\tau \in \text{Rep}_{\mathcal{O}} G$. Let $x(\tau) = \mathcal{O}_L e_1 \oplus \ldots \oplus \mathcal{O}_L e_r$, $y(\tau) = \mathcal{O}_L \pi_L^{-n_1} e_1 \oplus \ldots \oplus \mathcal{O}_L \pi_L^{-n_r} e_r$, where $\pi_L$ and $\pi_{L'}$ are the uniformizers of $\mathcal{O}_L$ and $\mathcal{O}_{L'}$, respectively. Then

$$\begin{aligned}
x(\tau) \cap \pi_L^n y(\tau) &= \oplus_{i=1}^r \mathcal{O}_L \pi_L^{\max\{0, n-n_i\}} e_i \\
x(\tau) \cap \pi_L^n y(\tau) + \pi_L x(\tau) &= \oplus_{i=1}^r \mathcal{O}_L \pi_L^{\min\{1, \max\{0, n-n_i\}\}} e_i \\
\mathcal{F}^n(x(\tau), y(\tau)) &= \oplus_{n_i \leq n} l \overline{e}_i.
\end{aligned}$$

Now, we have $x'(\tau) = \mathcal{O}_{L'}e_1 \oplus \ldots \oplus \mathcal{O}_{L'}e_r$, $y'(\tau) = \mathcal{O}_{L'}\pi_{L'}^{-en_1}e_1 \oplus \ldots \oplus \mathcal{O}_{L'}\pi_{L'}^{-en_r}e_r$, so

$$\mathcal{F}^n(x'(\tau), y'(\tau)) = \oplus_{en_i \leq n} l'\overline{e}_i = \oplus_{n_i \leq \frac{n}{e}} l'\overline{e}_i = \mathcal{F}^{\frac{n}{e}}(x(\tau), y(\tau)) \otimes l' = e \cdot \mathcal{F}^n(x(\tau), y(\tau)) \otimes l'.$$

$$\begin{array}{rcl}
\mathcal{F}^n(x'(\tau), y'(\tau)) & = & \oplus_{en_i \leq n} l'\overline{e}_i \\
& = & \oplus_{n_i \leq \frac{n}{e}} l'\overline{e}_i \\
& = & \mathcal{F}^{\frac{n}{e}}(x(\tau), y(\tau)) \otimes l' \\
& = & e \cdot \mathcal{F}^n(x(\tau), y(\tau)) \otimes l'.
\end{array}$$

$\square$

*Remark* 13. We can apply the last proposition in the case where $L' = L$ and $\varphi$ is a morphism which may be ramified, with ramification index $e$. Suppose there is an isomorphism $\varphi_V : \varphi^*V \to V$ of $\otimes$-functors, where $\varphi^*V = V \otimes_{L,\varphi} L$. As an abuse of notation, we denote by $\varphi_V : \mathcal{L}(V) \xrightarrow{\varphi} \mathcal{L}(\varphi^*V) \xrightarrow{\mathcal{L}(\varphi_V)} \mathcal{L}(V)$ the induced morphism on lattices. Then we have

$$\mathrm{Pos}(\varphi_V(x), \varphi_V(y)) = e\varphi \mathrm{Pos}(x, y) \quad \text{in} \quad \mathbf{C}^{\mathbb{Z}}(G_l),$$

where $l$ is the residue field of $L$, for every $x, y \in \mathcal{L}(V)$.

## 4.3 Crystalline representations with $G$-structure

Let $\mathbb{F}$ be an algebraically closed field of characteristic $p > 0$, $W(\mathbb{F})$ the ring of Witt vectors over $\mathbb{F}$, $K_0 = \mathrm{Frac}\, W(\mathbb{F})$ and we fix an algebraic closure of $K_0$ denoted by $\overline{K}_0$. We denote $\mathrm{Gal}_K = \mathrm{Gal}(\overline{K}_0/K)$, for every (totally ramified) extension $K_0 \subset K \subset \overline{K}_0$. Let $G$ be a reductive group over $\mathbb{Z}_p$.

The aim of this section is to define and study the diagram

$$\begin{array}{ccc}
\mathcal{L}(V, K) \to \mathcal{L}(N) & \to & \mathcal{L}(D', \leq \mathbf{t}_{\mathrm{H}}) \\
\downarrow & & \downarrow \simeq \\
\mathcal{L}(V) & \longrightarrow & \mathcal{L}(D, \leq \mathbf{t}_{\mathrm{H}})
\end{array}$$

for $V, N, D$ and $D'$ some $\otimes$-functors on $\mathrm{Rep}_{\mathbb{Q}_p} G$ that will be defined later.

### 4.3.1 Isocrystals with $G$-structure

An isocrystal with $G$-structure (or $G$-isocrystal) is an exact and faithful $\otimes$-functor

$$D \;:\; \mathrm{Rep}_{\mathbb{Q}_p} G \to \mathrm{Mod}_{K_0}^{\sigma}.$$

We denote by $\omega(D)$ its underlying fiber functor $\omega(D) : \mathrm{Rep}_{\mathbb{Q}_p} G \to \mathrm{Bun}_{K_0}$. We say that $D$ is trivial when $\omega(D) = \omega_{G,K_0}$.

**Lemma 4.8.** *There is a canonical correspondance between trivial $G$-isocrystals and elements of $G(K_0)$.*

*Proof.* A trivial $G$-isocrystal is a $\otimes$-isomorphism $\sigma_D : \sigma^*\omega_{G,K_0} \simeq \omega_{G,K_0}$. Since $\sigma^*\omega_{G,K_0}$ is canonically isomorphic to $\omega_{G,K_0}$, a trivial $G$-isocrystal corresponds to an element in $\mathrm{Aut}^{\otimes}(\omega_{G,K_0}) = G(K_0)$. $\square$

Composing $D$ with the Newton slope graduation $\mathrm{Mod}_{K_0}^{\sigma} \to \mathrm{Gr}_{K_0}^{\mathbb{Q}}$, we obtain a $\mathbb{Q}$-graduation

$$\mathcal{G}_{\mathrm{N}} \;:\; \mathrm{Rep}_{\mathbb{Q}_p} G \to \mathrm{Gr}_{K_0}^{\mathbb{Q}}.$$

It lives in $\mathbf{G}^{\mathbb{Q}}(\omega(D))$, and so does the opposed graduation, defined by $\mathcal{G}_{\mathrm{N}}^{\iota\gamma} = \mathcal{G}_{\mathrm{N}}^{-\gamma}$, for every $\gamma \in \mathbb{Q}$. We get two filtrations $\mathcal{F}_{\mathrm{N}} = \mathrm{Fil}(\mathcal{G}_{\mathrm{N}})$ and $\mathcal{F}_{\mathrm{N}}^{\iota} = \mathrm{Fil}(\mathcal{G}_{\mathrm{N}}^{\iota})$, living in $\mathbf{F}^{\mathbb{Q}}(\omega(D))$, called the Newton filtrations. We denote by $\mathbf{t}_{\mathrm{N}}(D)$ and $\mathbf{t}_{\mathrm{N}}^{\iota}(D)$ the types of these filtrations, which are elements in $\mathbf{C}^{\mathbb{Q}}(\omega(D))$.

**Proposition 4.9.** *Suppose $D$ is isomorphic to a trivial $G$-isocrystal, then*

$$\mathbf{t}_{\mathrm{N}}(D) = \mathbf{t}_{\mathrm{N}}(D)^{\#} \quad \text{and} \quad \mathbf{t}_{\mathrm{N}}^{\iota}(D) = \mathbf{t}_{\mathrm{N}}^{\iota}(D)^{\#} \quad \text{in} \quad \mathbf{C}^{\mathbb{Q}}(G_{K_0}).$$

*Proof.* This is well known, see [26, 4.4]. The Frobenius of $D$ induces an isomorphism of fiber functors $\sigma_D : \sigma^*\omega(D) \to \omega(D)$. It gives rises to a commutative diagram

$$\begin{array}{ccc}
\mathbf{F}^{\mathbb{Q}}(\omega(D)) & \xrightarrow{\sigma} \mathbf{F}^{\mathbb{Q}}(\sigma^*\omega(D)) & \xrightarrow{\sigma_D} \mathbf{F}^{\mathbb{Q}}(\omega(D)) \\
\downarrow t & \downarrow t & \downarrow t \\
\mathbf{C}^{\mathbb{Q}}(G_{K_0}) & \xrightarrow{\sigma} \mathbf{C}^{\mathbb{Q}}(G_K) & \xrightarrow{\mathrm{Id}} \mathbf{C}^{\mathbb{Q}}(G_K)
\end{array}$$

Since $\mathcal{G}_{\mathrm{N}}$ is a graduation by sub-isocrystals, the top map fixes $\mathcal{F}_{\mathrm{N}}$ and $\mathcal{F}_{\mathrm{N}}^{\iota}$. Thus $\sigma \mathbf{t}_{\mathrm{N}}(D) = \mathbf{t}_{\mathrm{N}}(D)$ and $\sigma \mathbf{t}_{\mathrm{N}}^{\iota}(D) = \mathbf{t}_{\mathrm{N}}^{\iota}(D)$. □

Set $\mathcal{L}(D) = \mathcal{L}(\omega(D))$. Thus $\mathcal{L}(D) \neq \emptyset$ if and only if $D$ is isomorphic to a trivial $G$-isocrystal, by section 4.2.2A variant: From $\mathrm{Rep}_{\mathbb{Q}_p} G$ to $\mathrm{Rep}_{\mathbb{Z}_p} G$subsubsection.4.2.2, which we assume from now on. As we have seen in Remark 13rem.13, the Frobenius on $D$ induces a bijection

$$\mathcal{L}(D) = \mathcal{L}(\omega(D)) \xrightarrow{\sigma} \mathcal{L}(\sigma^*\omega(D)) \xrightarrow{\sigma_D} \mathcal{L}(\omega(D)) = \mathcal{L}(D)$$

which we simply denote by $\sigma_D$. By Remark 13rem.13, we have

$$\mathrm{Pos}(\sigma_D x, \sigma_D y) = \sigma \mathrm{Pos}(x, y)$$

in $\mathbf{C}^{\mathbb{Z}}(G_{\mathbb{F}})$, for every $x, y \in \mathcal{L}(D)$. For any $y \in \mathcal{L}(D)$, we define the Hodge filtration by

$$\mathcal{F}_{\mathrm{H}}(y) = \mathcal{F}(y, \sigma_D y).$$

Proposition 4.5thm.4.5 shows that it lives in $\mathbf{F}^{\mathbb{Z}}(G_{\mathbb{F}})$. We denote by

$$\mathbf{t}_{\mathrm{H}}(y) = \mathrm{Pos}(y, \sigma_D y)$$

its type in $\mathbf{C}^{\mathbb{Z}}(G_{\mathbb{F}}) = \mathbf{C}^{\mathbb{Z}}(G_{K_0})$, which verifies Mazur's inequality given by the next proposition.

**Proposition 4.10.** *The types above verify the inequality $\mathbf{t}_{\mathrm{N}}^{\iota}(D) \leq \mathbf{t}_{\mathrm{H}}(y)^{\#}$ in $\mathbf{C}^{\mathbb{Q}}(G_{K_0})$ for any $y \in \mathcal{L}(D)$.*

*Proof.* This is well-known, see [37, Theorem 4.2], [29, 4.1 and 4.10] or [10, 4.2]. The inequality $\mathbf{t}_{\mathrm{N}}^{\iota}(D)(\tau) \leq \mathbf{t}_{\mathrm{H}}(y)(\tau)$ for every $\tau \in \mathrm{Rep}_{\mathbb{Z}_p} G$ and every $y \in \mathcal{L}(D)$ is a well-known inequality, stated by Mazur in [33] and whose proof can be found in [20, Theorem 1.4.1]. Thus, by proposition 4.1thm.4.1 and since $\mathbf{t}_{\mathrm{N}}(D) = \mathbf{t}_{\mathrm{N}}(D)^{\#}$, we get $\mathbf{t}_{\mathrm{N}}^{\iota}(D) \leq \mathbf{t}_{\mathrm{H}}(y)^{\#}$. □

For any $\mu \in \mathbf{C}^{\mathbb{Z}}(G_{K_0})$, we define two subsets of lattices in $\mathcal{L}(D)$ by

$$\mathcal{L}(D, \leq \mu) = \{y \in \mathcal{L}(D) \mid \mathbf{t}_{\mathrm{H}}(y) \leq \mu\} \quad \text{and} \quad \mathcal{L}(D, \mu) = \{y \in \mathcal{L}(D) \mid \mathbf{t}_{\mathrm{H}}(y) = \mu\}.$$

These objects have already been studied by various authors, in particular we can give an improved version of Mazur's inequality. First, Kottwitz in [28], and Rapoport and Richartz in [37], classify trivial isocrystals with $G$-structure, i.e. they describe the set $B(G)$ of $\sigma$-conjugacy classes in $G(K_0)$. The key ingredient for the classification is a map

$$\begin{array}{rcl}
B(G) & \to & \mathbf{C}^{\mathbb{Q}}(G_{K_0}) \times \pi_1(G)_{\Gamma} \\
{[b]} & \mapsto & (\nu_G[b], \kappa_G[b])
\end{array}$$

which is shows to be injective in [28, 4.13], see also [37]. The first component is our Newton type, constructed by Kottwitz in [26, 3]: if $D_b$ is the $G$-isocrystal corresponding to $b \in G(K_0)$, then $\nu_G[b] = \mathbf{t}_{\mathrm{N}}(D_b)$. The second component is the Kottwitz map, defined in [26, Proposition 5.6] or [37, 1.13], with values in the $\Gamma$-coinvariants of the algebraic fundamental group $\pi_1(G)$, there $\Gamma = \mathrm{Gal}_{\mathbb{Q}_p}$. In our quasi-split case, $\mathbf{C}^{\mathbb{Q}}(G)^{\Gamma} = \mathbf{C}^{\mathbb{Q}}(G_{\mathbb{Q}_p})$ and for any Borel pair $(T, B)$ in $G_{\mathbb{Q}_p}$, there are canonical $\Gamma$-equivariant isomorphisms

$$\mathbf{C}^{\mathbb{Z}}(G) \simeq X_*(T)^{B-\mathrm{dom}}, \quad \mathbf{C}^{\mathbb{Q}}(G) \simeq (X_*(T) \otimes \mathbb{Q})^{B-\mathrm{dom}} \quad \text{and} \quad \pi_1(G) \simeq X_*(T)/\mathbb{Z} \cdot R_G^{\vee}$$

where $X_*(T)^{B-\mathrm{dom}}$ and $(X_*(T) \otimes \mathbb{Q})^{B-\mathrm{dom}}$ are the cones of $B$-dominant elements in $X_*(T)$ and $X_*(T) \otimes \mathbb{Q}$, respectively, while $\mathbb{Z} \cdot R_G^{\vee}$ is the subgroup of $X_*(T)$ spanned by the coroots $R_G^{\vee}$ of $T$ in $G$. Moreover, the maps

$$[-]_G : \mathbf{C}^{\mathbb{Z}}(G) \twoheadrightarrow \pi_1(G) \quad \text{and} \quad [-]_{G,\Gamma} : \mathbf{C}^{\mathbb{Z}}(G) \twoheadrightarrow \pi_1(G)_{\Gamma}$$

induced by these isomorphisms, given in [37, 4.1] do not depend upon the chosen Borel pair.

Let $(T, B)$ be a Borel pair in $G = G_{\mathbb{Z}_p}$. For $\mu \in X_*(T)^{B-\mathrm{dom}} \simeq \mathbf{C}^{\mathbb{Z}}(G_{K_0})$ and $b \in G(K_0)$, set

$$X_\mu^G(b) = \{x \in G(K_0)/G(W(\mathbb{F})) \mid x^{-1} b \sigma(x) \in G(W(\mathbb{F}))\mu(p)G(W(\mathbb{F}))\}.$$

The following improved version of Mazur's inequality was established in [37, Theorem 4.2], see also Theorem 4.1 and section 4.4 of [29]:

$$X_\mu^G(b) \neq \emptyset \implies \begin{cases} \nu_G[b] \leq \mu^\# & \text{in } \mathbf{C}^{\mathbb{Q}}(G) \quad \text{and} \\ \kappa_G[b] = [\mu]_{G,\Gamma} & \text{in } \pi_1(G)_\Gamma. \end{cases}$$

The last equality is also explained in [10, 4.2]. The converse implication was established by Gashi in [18], building on a strategy proposed by Kottwitz in [29]. Thus:

$$X_\mu^G(b) \neq \emptyset \iff \begin{cases} \nu_G[b] \leq \mu^\# & \text{in } \mathbf{C}^{\mathbb{Q}}(G) \quad \text{and} \\ \kappa_G[b] = [\mu]_{G,\Gamma} & \text{in } \pi_1(G)_\Gamma. \end{cases}$$

On the other hand, the $G(K_0)$-equivariant bijection

$$\begin{array}{rcl} G(K_0)/G(W(\mathbb{F})) = \mathrm{Aut}^\otimes(\omega_{G,K_0})/\mathrm{Aut}^\otimes(\omega_{G,W(\mathbb{F})}) & \to & \mathcal{L}(\omega_{G,K_0}) = \mathcal{L}(D_b) \\ g & \mapsto & g \cdot \omega_{G,W(\mathbb{F})} \end{array}$$

induces a bijection

$$X_\mu^G(b) \simeq \mathcal{L}(D_b, \mu^\iota),$$

as we can check using Example 4.1ex.4.1. Thus,

$$\mathcal{L}(D, \mu) \neq \emptyset \iff \begin{cases} \mathbf{t}_{\mathrm{N}}^\iota(D) \leq \mu^\# & \text{in } \mathbf{C}^{\mathbb{Q}}(G_{K_0}) \quad \text{and} \\ \kappa(D) = -[\mu]_{G,\Gamma} & \text{in } \pi_1(G)_\Gamma \end{cases}$$

for every $G$-isocrystal $D$.

**Definition 4.2.** For $\mu \in \mathbf{C}^{\mathbb{Z}}(G_{K_0})$, we say that $D$ is $\mu$-ordinary when $\mathcal{L}(D, \mu) \neq \emptyset$ and $\mathbf{t}_{\mathrm{N}}^\iota(D) = \mu^\#$.

**Proposition 4.11.** *Up to isomorphism, there is a unique $D$ which is $\mu$-ordinary.*

*Proof.* We will prove the existence later. For the unicity, if $D$ is $\mu$-ordinary, $\mathcal{L}(D) \neq \emptyset$ so we may assume that $D$ is trivial. From the discussion above, we obtain that then there is at most one $[b] \in B(G)$ with $X_\mu(b) \neq \emptyset$ and $\nu_G[b] = \mu^\#$, the one associated to the element $(\mu^\#, -[\mu]_{G,\Gamma}) \in \mathbf{C}^{\mathbb{Q}}(G_{K_0}) \times \pi_1(G)_\Gamma$. □

### 4.3.2 Torsion Kisin modules with $G$-structure

**Definition 4.3.** An exact $\otimes$-functor

$$M_{\mathrm{tors}} : \mathrm{Rep}_{\mathbb{Z}_p} G \to \mathrm{Mod}_{\mathbb{F}[[u]],\mathrm{fr}}^\varphi$$

is called a torsion Kisin module with $G$-structure. We denote by $\omega(M_{\mathrm{tors}})$ its underlying fiber functor and we say that $M_{\mathrm{tors}}$ is trivial when $\omega(M_{\mathrm{tors}}) = \omega_{G,\mathbb{F}[[u]]}$.

**Proposition 4.12.** *Every torsion Kisin module with $G$-structure is isomorphic a trivial one.*

*Proof.* There is a correspondence between torsion Kisin modules with $G$-structure and $\otimes$-fiber functors $\mathrm{Rep}_{\mathbb{F}_p} G_{\mathbb{F}_p} \to \mathrm{Mod}_{\mathbb{F}[[u]],\mathrm{fr}}^\varphi$ which are always isomorphic to the trivial fiber functor, by Proposition 4.2thm.4.2. □

**Definition 4.4.** An isogeny class of torsion Kisin modules with $G$-structure is an exact $\otimes$-functor

$$\overline{X} : \mathrm{Rep}_{\mathbb{Z}_p} G \to \mathrm{Mod}_{\mathbb{F}((u))}^\varphi.$$

We denote by $\omega(\overline{X})$ its underlying fiber functor and we say what $\overline{X}$ is trivial when $\omega(\overline{X}) = \omega_{G,\mathbb{F}((u))}$.

Set $\mathcal{L}(\overline{X}) = \mathcal{L}(\omega(\overline{X}))$. Thus $\mathcal{L}(\overline{X}) \neq \emptyset$ if and only if $\overline{X}$ is isomorphic to a trivial isogeny class of torsion Kisin module with $G$-structure, which we assume from now on. Any $z \in \mathcal{L}(\overline{X})$ has a unique factorization through a $\otimes$-functor $\operatorname{Rep}_{\mathbb{Z}_p} G \to \operatorname{Mod}^\varphi_{\mathbb{F}[[u]],\mathrm{fr}}$, i.e. yields a torsion Kisin module with $G$-structure.

As explained in Remark 13rem.13, the Frobenius morphism induces a map

$$\varphi_{\overline{X}} \; ; \; \mathcal{L}(\overline{X}) = \mathcal{L}(\omega(\overline{X})) \xrightarrow{\varphi^*} \mathcal{L}(\varphi^*\omega(\overline{X})) \xrightarrow{\varphi_{\overline{X}}} \mathcal{L}(\omega(\overline{X})) = \mathcal{L}(\overline{X}).$$

For every $z \in \mathcal{L}(\overline{X})$ and $n \geq 1$, we set

$$\mathcal{F}_{\mathrm{H}}(z) = \mathcal{F}(z, \varphi_{\overline{X}}(z)) \quad \text{in} \quad \mathbf{F}^{\mathbb{Z}}(\omega(z \otimes \mathbb{F})),$$

by Proposition 4.5thm.4.5. Its type lives in $\mathbf{C}^{\mathbb{Z}}(G_{\mathbb{F}}) = \mathbf{C}^{\mathbb{Z}}(G_{K_0})$.

The filtration $\mathcal{F}_{\mathrm{F},1}$ over $\operatorname{Mod}^\varphi_{\mathbb{F}[[u]],\mathrm{fr}}$ defines a $\otimes$-functor

$$\mathcal{F}_{\mathrm{F},1}(z) = \mathcal{F}_{\mathrm{F},1} \circ z \; : \; \begin{array}{rcl} \operatorname{Rep}_{\mathbb{Z}_p} G & \to & \operatorname{Fil}^{\mathbb{Q}}_{\mathbb{F}[[u]]} \\ \tau & \mapsto & \mathcal{F}_{\mathrm{F},1}(z(\tau)) \end{array}$$

and the next proposition shows that this is a filtration in $\mathbf{F}^{\mathbb{Q}}(\omega(z))$.

**Proposition 4.13.** *The functor*

$$\mathcal{F}_{\mathrm{F},1}(z) = \mathcal{F}_{\mathrm{F},1} \circ z \; : \; \operatorname{Rep}_{\mathbb{Z}_p} G \to \operatorname{Fil}^{\mathbb{Q}}_{\mathbb{F}[[u]]},$$

*factoring the fiber functor $\omega(z)$, is an exact $\otimes$-functor.*

*Proof.* We have to show that $\mathcal{F}_{\mathrm{F},1}(z)$ is exact. Since it is a filtration by strict subobjects, it is sufficient to show that $\mathcal{F}_{\mathrm{F},1}(z) \otimes \mathbb{F}((u)) : \operatorname{Rep}_{\mathbb{Z}_p} G \to \operatorname{Fil}^{\mathbb{Q}}_{\mathbb{F}((u))}$ is exact. It follows from Proposition 4.3thm.4.3. $\square$

**Definition 4.5.** If $z \in \mathcal{L}(\overline{X})$, we define the Fargues type $\mathbf{t}_{\mathrm{F},1}(z)$ as the type associated to the filtration $\mathcal{F}_{\mathrm{F},1}$ of $z$. It lives in $\mathbf{C}^{\mathbb{Q}}(G_{\mathbb{F}[[u]]}) = \mathbf{C}^{\mathbb{Q}}(G)$.

**Proposition 4.14.** *The types defined above on $z \in \mathcal{L}(\overline{X})$ verify the following inequalities:*

1. $\mathbf{t}_{\mathrm{F},1}(z) = \mathbf{t}_{\mathrm{F},1}(z)^{\#}$

2. $\mathbf{t}_{\mathrm{F},1}(z) \leq \mathbf{t}_{\mathrm{H}}(z)^{\#}$ *for every $n \geq 1$.*

*Proof.*  1. The commutative diagram

$$\begin{array}{ccccc} \mathbf{F}^{\mathbb{Q}}(\omega(\overline{X})) & \xrightarrow{\varphi} & \mathbf{F}^{\mathbb{Q}}(\varphi^*\omega(\overline{X})) & \xrightarrow{\varphi_{\overline{X}}} & \mathbf{F}^{\mathbb{Q}}(\omega(\overline{X})) \\ \downarrow t & & \downarrow t & & \downarrow t \\ \mathbf{C}^{\mathbb{Q}}(G_{\mathbb{F}((u))}) & \xrightarrow{\sigma} & \mathbf{C}^{\mathbb{Q}}(G_{\mathbb{F}((u))}) & \xrightarrow{\mathrm{Id}} & \mathbf{C}^{\mathbb{Q}}(G_{\mathbb{F}((u))}) \end{array}$$

gives us

$$\mathbf{t}_{\mathrm{F},1}(z) = \mathbf{t}\left(\mathcal{F}_{\mathrm{F},1}(z) \otimes \mathbb{F}((u))\right) = \mathbf{t}(\varphi_{\overline{X}}\varphi^*(\mathcal{F}_{\mathrm{F},1}(z) \otimes \mathbb{F}((u)))) = \sigma \cdot \mathbf{t}\left(\mathcal{F}_{\mathrm{F},1}(z) \otimes \mathbb{F}((u))\right) = \sigma \cdot \mathbf{t}_{\mathrm{F},1}(z)$$

since $\mathcal{F}_{\mathrm{F},1}(z) \otimes \mathbb{F}((u))$ is stable under the upper row map.

2. We have $\mathbf{t}_{\mathrm{F},1}(z)(\tau) \leq \mathbf{t}_{\mathrm{H}}(z)(\tau)$ for every $\tau \in \operatorname{Rep}_{\mathbb{Z}_p} G$, by Proposition 3.14thm.3.14, so by proposition 4.1thm.4.1 and the previous point, we have $\mathbf{t}_{\mathrm{F},1}(z) \leq \mathbf{t}_{\mathrm{H}}(z)^{\#}$. $\square$

### 4.3.3 Kisin modules with $G$-structure

We fix a finite extension $K$ of $K_0$, an uniformizer $\pi_K$ of $K$, $E \in \mathfrak{S}$ the minimal polynomial of $\pi_K$. Let $e = [K : K_0]$. Let $\hat{\mathfrak{S}}$ be the completion of $\mathfrak{S}_{(E)}$.

**Definition 4.6.** An exact and faithful $\otimes$-functor

$$M \ : \ \mathrm{Rep}_{\mathbb{Z}_p} G \to \mathrm{Mod}^{\varphi}_{\mathfrak{S},\,\mathrm{fr}}$$

is called a Kisin module with $G$-structure. Let $\omega(M)$ be its underlying fiber functor. We say that $M$ is trivial if $\omega(N) = \omega_{G,\mathfrak{S}}$. An exact and faithful $\otimes$-functor

$$N \ : \ \mathrm{Rep}_{\mathbb{Q}_p} G \to \mathrm{Mod}^{\varphi}_{\mathfrak{S}[\frac{1}{p}]}$$

is called an isogeny class of Kisin modules with $G$-structure. We denote by $\omega(N)$ its underlying fiber functor. We say that $N$ is trivial if $\omega(M) = \omega_{G,\mathfrak{S}[\frac{1}{p}]}$.

*Remark* 14. To any Kisin module with $G$-structure $M$, we can associate a torsion Kisin module with $G$-structure $\overline{M}$, obtained by composition

$$\overline{M} \ : \ \mathrm{Rep}_{\mathbb{Z}_p} G \xrightarrow{M} \mathrm{Mod}^{\varphi}_{\mathfrak{S},\,\mathrm{fr}} \xrightarrow{\mathrm{mod}\ p} \mathrm{Mod}^{\varphi}_{\mathbb{F}[[u]],\mathrm{fr}}.$$

We fix an isogeny class of Kisin modules with $G$-structure $N$. We define the set of lattices inside $N$, and denote it by $\mathcal{L}(N)$, as the set of Kisin modules with $G$-structure $M$ such that $M[\frac{1}{p}]$ corresponds to the $\otimes$-functor $N'$ given by precomposition of $N$ with $\mathrm{Rep}_{\mathbb{Z}_p} G \to \mathrm{Rep}_{\mathbb{Q}_p} G$.

**Proposition 4.15.** *1. Every Kisin module with $G$-structure is isomorphic a the trivial one.*

*2. If $\mathcal{L}(N) \neq \emptyset$, then $N$ is isomorphic to a trivial isogeny class of Kisin modules with $G$-structure*

*Proof.* 1. This follows from Proposition 4.2thm.4.2.

2. Let $M \in \mathcal{L}(N)$. Then $\omega(M) \simeq \omega_{G,\mathfrak{S}}$, thus $\omega(N') = \omega(M) \otimes \mathfrak{S}[\frac{1}{p}] \simeq \omega_{G,\mathfrak{S}[\frac{1}{p}]}$. This also implies that $\omega(N) \simeq \omega_{G_{\mathbb{Q}_p},\mathfrak{S}[\frac{1}{p}]}$, as in 4.2.2A variant: From $\mathrm{Rep}_{\mathbb{Q}_p} G$ to $\mathrm{Rep}_{\mathbb{Z}_p} G$subsubsection.4.2.2. □

We assume that $\mathcal{L}(N) \neq \emptyset$ from now on. We can construct a functor

$$\mathcal{F}_{\mathrm{H}}(N) = \mathcal{F}(N \otimes \hat{\mathfrak{S}}, \varphi_N \varphi^* N \otimes \hat{\mathfrak{S}}) \ : \ \mathrm{Rep}_{\mathbb{Q}_p} G \to \mathrm{Fil}(\omega(N)_K).$$

Then, Proposition 4.5thm.4.5 tells us that it is a filtration in $\mathbf{F}^{\mathbb{Z}}(\omega_K(N))$. Its type it is denoted by $\mathbf{t}_{\mathrm{H}}(N)$ and it lives in $\mathbf{C}^{\mathbb{Z}}(G_K) = \mathbf{C}^{\mathbb{Z}}(G_{K_0})$.

We can associate an isocrystal with $G$-structure to $N$, by setting

$$D' = N/uN : \mathrm{Rep}_{\mathbb{Q}_p} G \to \mathrm{Mod}^{\sigma}_{K_0},$$

which is a faithful functor since its faithfulness only depends on the fiber functor $\omega(D')$ which is faithful as it is isomorphic to the trivial fiber functor $\omega_{G,K_0}$ by assumption. We thus get a map $\mathcal{L}(N) \to \mathcal{L}(D')$ by sending $M$ to the exact $\otimes$-functor

$$y' = M/uM \ : \ \mathrm{Rep}_{\mathbb{Z}_p} G \to \mathrm{Mod}^{\sigma}_{W(\mathbb{F})}.$$

**Proposition 4.16.** *Let $M \in \mathcal{L}(N)$ and $\overline{M}$ the torsion Kisin module with $G$-structure associated to $M$ and $y'$ the crystal with $G$-structure associated to $M$. We have*

*1. The inequality $\mathbf{t}_{\mathrm{H}}(y') \leq \mathbf{t}_{\mathrm{H}}(N)$, i.e. the map $\mathcal{L}(N) \to \mathcal{L}(D')$ factors through $\mathcal{L}(D', \leq \mathbf{t}_{\mathrm{H}}(N))$.*

*2. The inequality $\mathbf{t}_{\mathrm{H},1}(\overline{M}) \leq e \cdot \mathbf{t}_{\mathrm{H}}(N)$.*

*Proof.* We treat the two points in one. We can suppose that $M$ is trivial, so the Frobenius on $M$ is given by an element $g \in G(\mathfrak{S}[\frac{1}{E}]) = \mathrm{Aut}^\otimes(\omega_{\mathfrak{S}[\frac{1}{E}]})$. Then

$$\begin{array}{rcl} \mathbf{t}_{\mathrm{H}}(y') & = & \mathrm{Pos}(\omega_{G,W(\mathbb{F})}, g(0) \cdot \omega_{G,W(\mathbb{F})}) \\ \mathbf{t}_{\mathrm{H},1}(\overline{M}) & = & \mathrm{Pos}(\omega_{G,\mathbb{F}[[u]]}, \overline{g} \cdot \omega_{G,\mathbb{F}[[u]]}) \\ \mathbf{t}_{\mathrm{H}}(N) & = & \mathrm{Pos}(\omega_{G,\hat{\mathfrak{S}}}, g \cdot \omega_{G,\hat{\mathfrak{S}}}) \end{array}$$

in, respectively, $\mathbf{C}^{\mathbb{Z}}(G_\mathbb{F}) = \mathbf{C}^{\mathbb{Z}}(G)$, $\mathbf{C}^{\mathbb{Z}}(G_\mathbb{F}) = \mathbf{C}^{\mathbb{Z}}(G)$ and $\mathbf{C}^{\mathbb{Z}}(G_K) = \mathbf{C}^{\mathbb{Z}}(G)$. Thus, by 4.1thm.4.1, it suffices to show that for every $\tau \in \mathrm{Rep}_{\mathfrak{S}} G_{\mathfrak{S}}$ (or, respectively, $\mathrm{Rep}_{W(\mathbb{F})} G_{W(\mathbb{F})}$ and $\mathrm{Rep}_{\mathbb{F}[[u]]} G_{\mathbb{F}[[u]]}$), we have $\mathrm{Pos}(X/uX, Y/uY) \leq \mathrm{Pos}(X \otimes \tilde{\mathfrak{S}}, Y \otimes \tilde{\mathfrak{S}})$ (resp. $\mathrm{Pos}(X/pX, Y/pY) \leq e \cdot \mathrm{Pos}(X \otimes \tilde{\mathfrak{S}}, Y \otimes \tilde{\mathfrak{S}})$), where $X = M(\tau)$ and $Y = \tau(g)M(\tau)$. This follows from the proof of 3.15thm.3.15. $\square$

**Definition 4.7.** A Kisin module with $G$-structure is said to be HN-type if for every $\mathbb{Z}_p$-representation $\tau$ of $G$, the Kisin module $M(\tau)$ is HN-type.

Under the type HN-type hypothesis, we can define a Fargues type on $M$.

**Proposition 4.17.** *Suppose $M$ is HN-type, then the functor $\mathcal{F}_\mathrm{F}(M)$ sending a $\mathbb{Z}_p$-representation $\tau$ of $G$ to the Fargues filtration of the HN-type Kisin module $M(\tau)$ is in $\mathbf{F}^{\mathbb{Q}}(\omega(M))$.*

*Proof.* It is a filtration by direct summands, so we need to check the exactness and compatibility with tensor products. For exactness, let

$$0 \to \tau_1 \to \tau_2 \to \tau_3 \to 0$$

be an exact sequence of representations of $G$. Then, we have a commutative diagram where the second row is exact

$$\begin{array}{ccccccc} 0 & \to & \mathcal{F}^\gamma_\mathrm{F}(\tau_1) & \to & \mathcal{F}^\gamma_\mathrm{F}(\tau_2) & \xrightarrow{f} & \mathcal{F}^\gamma_\mathrm{F}(\tau_3) \\ & & \downarrow & & \downarrow & & \downarrow \\ 0 & \to & M(\tau_1) & \to & M(\tau_2) & \to & M(\tau_3) & \to & 0 \end{array}$$

Reducing the first sequence modulo $p$, we obtain the exact sequence

$$0 \to \mathcal{F}^\gamma_{\mathrm{F},1}(\tau_1) \to \mathcal{F}^\gamma_{\mathrm{F},1}(\tau_2) \to \mathcal{F}^\gamma_{\mathrm{F},1}(\tau_3) \to 0$$

so by Nakayama's lemma, the map $\mathcal{F}^\gamma_\mathrm{F}(\tau_2) \to \mathcal{F}^\gamma_\mathrm{F}(\tau_3)$ is surjective. Then, $\ker f = \ker(\mathcal{F}^\gamma_\mathrm{F}(\tau_2) \to \mathcal{F}^\gamma_\mathrm{F}(\tau_3))$ is a free $\mathfrak{S}$-module, and we also have $\mathcal{F}^\gamma_\mathrm{F}(\tau_1) \hookrightarrow \ker f$, which becomes an isomorphism modulo $p$, so again by Nakayama's lemma, we have $\mathcal{F}^\gamma_\mathrm{F}(\tau_1) \simeq \ker f$.

For the compatibility with tensor products, let $\mathcal{G}(\tau_1 \otimes \tau_2) = \mathcal{F}_\mathrm{F}(\tau_1) \otimes \mathcal{F}_\mathrm{F}(\tau_2)$. Then, we have to show that the Kisin module

$$X = \mathrm{Gr}^\gamma_\mathcal{G}(\tau_1 \otimes \tau_2) = \bigoplus_{\gamma_1 + \gamma_2 = \gamma} \mathrm{Gr}^{\gamma_1}_\mathrm{F}(M(\tau_1)) \otimes \mathrm{Gr}^{\gamma_2}_\mathrm{F}(M(\tau_2))$$

is semi-stable of slope $\gamma$. We can reduce modulo $p$ to obtain

$$\begin{array}{rcl} \overline{X} = \mathrm{Gr}^\gamma_{\overline{\mathcal{G}}}(\tau_1 \otimes \tau_2) & = & \bigoplus_{\gamma_1 + \gamma_2 = \gamma} \mathrm{Gr}^{\gamma_1}_\mathrm{F}(\overline{M}(\tau_1)) \otimes \mathrm{Gr}^{\gamma_2}_\mathrm{F}(\overline{M}(\tau_2)) \\ & = & \mathrm{Gr}^\gamma_{\mathcal{F}_{\mathrm{F},1}}(\tau_1 \otimes \tau_2) \end{array}$$

Then $\overline{X}$ is semi-stable of slope $\gamma = \gamma_1 + \gamma_2$ and so is therefore also $X$. $\square$

For $M$ a HN-type Kisin module, we can define the Fargues type $\mathbf{t}_\mathrm{F}(M)$ associated to the filtration $\mathcal{F}_\mathrm{F}(M)$. It is an element in $\mathbf{C}^{\mathbb{Q}}(G_\mathfrak{S})$.

### 4.3.4 Germs of crystalline representations with $G$-structure

**Definition 4.8.** A germ of crystalline representations with $G$-structure is a faithful $\otimes$-functor

$$V : \mathrm{Rep}_{\mathbb{Q}_p} G \to \mathrm{Rep}^{\mathrm{cr}}_{\mathbb{Q}_p}\{\mathrm{Gal}_{K_0}\}.$$

We denote by $\omega(V)$ its underlying fiber functor. We say that $V$ is trivial when $\omega(V) \simeq \omega_{G,\mathbb{Q}_p}$.

Let $V$ be a germ of crystalline representations with $G$-structure. We can define a Fargues filtration on $V$ by

$$\begin{array}{rccc}\mathcal{F}_{\mathrm{F}}(V) & : & \mathrm{Rep}_{\mathbb{Q}_p} G & \to & \mathrm{Fil}^{\mathbb{Q}}_{\mathbb{Q}_p} \\ & & \tau & \mapsto & \mathcal{F}_{\mathrm{F,cr}}(V(\tau))\end{array}$$

for $\mathcal{F}_{\mathrm{F,cr}}$ the filtration defined on crystalline representation in section 2. This functor a $\otimes$-functor by [15, Corollary 6] and exact since it is defined over $\mathbb{Q}_p$ and every exact sequence in $\mathrm{Rep}_{\mathbb{Q}_p} G$ is split, so we have $\mathcal{F}_{\mathrm{F}}(V) \in \mathbf{F}^{\mathbb{Q}}(\omega(V))$. We denote its type by $\mathbf{t}_{\mathrm{F}}(V)$ and it lives in $\mathbf{C}^{\mathbb{Q}}(G_{\mathbb{Q}_p})$.

*Remark* 15. The type $\mathbf{t}_{\mathrm{F}}(V)$ is invariant under $\#$, since it is defined in $\mathbf{C}^{\mathbb{Q}}(G_{\mathbb{Q}_p})$ and $\sigma$ acts trivially on this set.

The functor $V$ induces a filtered isocrystal with $G$-structure by taking

$$D = D_{\mathrm{cr}} \circ V \; : \; \mathrm{Rep}_{\mathbb{Q}_p} G \to^{\mathrm{wa}} \mathrm{MF}^{\sigma}_{\overline{K}_0}.$$

We (still) denote by $D$ the underlying $G$-isocrystal, by $\omega(D) : \mathrm{Rep}_{\mathbb{Q}_p} G \to \mathrm{Bun}_{K_0}$ the underlying fiber functor, and by

$$\begin{array}{rccc}\mathcal{F}_{\mathrm{H}}(D) & : & \mathrm{Rep}_{\mathbb{Q}_p} G & \to & \mathrm{Fil}^{\mathbb{Z}}_{\overline{K}_0}\end{array}$$

the Hodge filtration of the filtered isocrystal, which is an exact $\otimes$-functor by construction, i.e. an element of $\mathbf{F}^{\mathbb{Z}}(\omega(D) \otimes_{K_0} \overline{K}_0)$. Since $\overline{K}_0$ is algebraically closed, $\omega(D) \otimes_{K_0} \overline{K}_0 \simeq \omega_{G,\overline{K}_0}$ and the type $\mathbf{t}_{\mathrm{H}}(D)$ of $\mathcal{F}_{\mathrm{H}}(D)$ lives in $\mathbf{C}^{\mathbb{Z}}(G_{\overline{K}_0}) = \mathbf{C}^{\mathbb{Z}}(G_{K_0})$. We can also define a Fargues filtration on weakly admissible filtered isocrystals with $G$-structure by

$$\begin{array}{rccc}\mathcal{F}_{\mathrm{F}}(D) & : & \mathrm{Rep}_{\mathbb{Q}_p} G & \to & \mathrm{Fil}^{\mathbb{Q}}_{K_0} \\ & & \tau & \mapsto & \mathcal{F}_{\mathrm{F,wa}}(D(\tau))\end{array}$$

where $\mathcal{F}_{\mathrm{F,wa}}$ is the Fargues filtration on weakly admissible filtered isocrystals. This is an exact $\otimes$-functor by [9, Theorem 12], i.e. an element in $\mathbf{F}^{\mathbb{Q}}(\omega(D))$. If $\omega(D)$ is isomorphic to the trivial fiber functor, we obtain a type $\mathbf{t}_{\mathrm{F}}(D) \in \mathbf{C}^{\mathbb{Q}}(G_{K_0})$.

**Proposition 4.18.** *Suppose that $D$ is isomorphic to a trivial $G$-isocrystal. We have*

1. *The equalities $\mathbf{t}_{\mathrm{F}}(D) = \mathbf{t}_{\mathrm{F}}(D)^{\#}$ and $\mathbf{t}_{\mathrm{F}}(V) = \mathbf{t}_{\mathrm{F}}(D)$.*

2. *The inequality $\mathbf{t}_{\mathrm{F}}(V) \leq \mathbf{t}^{\iota}_{\mathrm{N}}(D)$.*

3. *If $\mathbf{t}_{\mathrm{N}}(D) = \mathbf{t}_{\mathrm{H}}(D)^{\#}$, then $\mathcal{F}^{\iota}_{\mathrm{N}}(D) = \mathcal{F}_{\mathrm{F}}(D)$ and $\mathbf{t}^{\iota}_{\mathrm{N}}(D) = \mathbf{t}_{\mathrm{F}}(D)$.*

*Proof.* 1. There is a commutative diagram

$$\begin{array}{ccccc}\mathbf{F}^{\mathbb{Q}}(\omega(D)) & \xrightarrow{\sigma} & \mathbf{F}^{\mathbb{Q}}(\sigma^*\omega(D)) & \xrightarrow{\sigma_D} & \mathbf{F}^{\mathbb{Q}}(\omega(D)) \\ \downarrow t & & \downarrow t & & \downarrow t \\ \mathbf{C}^{\mathbb{Q}}(G_{K_0}) & \xrightarrow{\sigma} & \mathbf{C}^{\mathbb{Q}}(G_{K_0}) & \xrightarrow{\mathrm{Id}} & \mathbf{C}^{\mathbb{Q}}(G_{K_0})\end{array}$$

and since $\mathcal{F}_{\mathrm{F}}(D)$ is stable under the upper row map, we get

$$\mathbf{t}_{\mathrm{F}}(D) = \mathbf{t}(\mathcal{F}_{\mathrm{F}}(D)) = \mathbf{t}(\sigma_D \sigma^*_{K_0} \mathcal{F}_{\mathrm{F}}(D)) = \sigma \cdot \mathbf{t}(\mathcal{F}_{\mathrm{F}}(D)) = \sigma \cdot \mathbf{t}_{\mathrm{F}}(D).$$

For the other equality, we have $\mathbf{t}_{\mathrm{F}}(V)(\tau) = \mathbf{t}_{\mathrm{F}}(D)(\tau)$ for every $\tau \in \mathrm{Rep}_{\mathbb{Z}_p} G$, by Proposition 3.11thm.3.11, so by proposition 4.1thm.4.1, we have $\mathbf{t}_{\mathrm{F}}(V) = \mathbf{t}_{\mathrm{F}}(V)^{\#} = \mathbf{t}_{\mathrm{F}}(D)^{\#} = \mathbf{t}_{\mathrm{F}}(D)$.

2. The inequality $\mathbf{t}_{\mathrm{F}}(V)(\tau) \leq \mathbf{t}^{\iota}_{\mathrm{N}}(D)(\tau)$ is given by Fargues in [15, Theorem 6], for every $\tau \in \mathrm{Rep}_{\mathbb{Q}_p} G$, thus $\mathbf{t}_{\mathrm{F}}(V) = \mathbf{t}_{\mathrm{F}}(V)^{\#} \leq \mathbf{t}^{\iota}_{\mathrm{N}}(D)^{\#} = \mathbf{t}^{\iota}_{\mathrm{N}}(D)$ by Proposition 4.1thm.4.1 and Remark 15rem.15.

3. We may assume $D$ is trivial and so, we can use the results in [9]. We drop $D$ from the notation. We have to show that $\mathbf{t}_N = \mathbf{t}_H^\#$ implies $\mathcal{F}_N^\iota = \mathcal{F}_F$. In [9], we see that if we equip $\mathbf{F}^\mathbb{R}(G_{K_0})$ with a $\varphi$-invariant CAT(0)-distance as defined in [8, Corollary 88], then

$$\begin{array}{rcl}\mathbf{F}^\mathbb{R}(G_{K_0},\varphi) & = & \{\mathcal{F} \in \mathbf{F}^\mathbb{R}(G_{K_0}) \mid \mathcal{F}^\gamma(\tau) \in \mathrm{Mod}_{K_0}^\sigma \text{ for every } \gamma \in \mathbb{R} \text{ and every } \tau \in \mathrm{Rep}_{\mathbb{Q}_p} G\} \\ \mathbf{F}^\mathbb{R}(G_{K_0},\varphi,\mathcal{F}_H) & = & \{\mathcal{F} \in \mathbf{F}^\mathbb{R}(G_{K_0}) \mid \mathcal{F}^\gamma(\tau) \in^{\mathrm{wa}} \mathrm{MF}_{\overline{K_0}}^\sigma \text{ for every } \gamma \in \mathbb{R} \text{ and every } \tau \in \mathrm{Rep}_{\mathbb{Q}_p} G\}.\end{array}$$

are closed convex subsets of $\mathbf{F}^\mathbb{R}(G_{K_0})$ and there is a convex projection

$$\mathbf{F}^\mathbb{R}(G_{K_0},\varphi) \to \mathbf{F}^\mathbb{R}(G_{K_0},\varphi,\mathcal{F}_H).$$

In [9, Proposition 13], we see that $\mathcal{F}_F$ is the convex projection of $\mathcal{F}_N^\iota$, so it suffices to prove that $\mathcal{F}_N^\iota \in \mathbf{F}^\mathbb{R}(G_{K_0},\varphi,\mathcal{F}_H)$. In [9, Lemma 11], we see that

$$\begin{array}{rcl}\mathbf{F}^\mathbb{R}(G_{K_0},\varphi,\mathcal{F}_H) & = & \{\Xi \in \mathbf{F}^\mathbb{R}(G_{K_0},\varphi) \mid \langle \mathcal{F}_H, \Xi\rangle = \langle \mathcal{F}_N, \Xi\rangle\} \\ & = & \{\Xi \in \mathbf{F}^\mathbb{R}(G_{K_0},\varphi) \mid \langle \mathcal{F}_H, \Xi\rangle \geq \langle \mathcal{F}_N, \Xi\rangle\}\end{array}$$

so we need to prove that, under the assumption $\mathbf{t}_N(D) = \mathbf{t}_H(D)^\#$, we have $\langle \mathcal{F}_H, \mathcal{F}_N^\iota\rangle \geq \langle \mathcal{F}_N, \mathcal{F}_N^\iota\rangle$. The operator

$$\langle x, y\rangle^{\mathrm{tr}} = \inf\{\langle X, Y\rangle \mid \mathbf{t}(X) = x, \mathbf{t}(Y) = y\ X, Y \in \mathbf{F}^\mathbb{R}(G_{K_0})\}$$

defined in [8, 4.2.5] verifies that for two opposed filtrations $\mathcal{F}$ and $\mathcal{F}^\iota$, we have $\langle \mathbf{t}(\mathcal{F}), \mathbf{t}(\mathcal{F}^\iota)\rangle^{\mathrm{tr}} = \langle \mathcal{F}, \mathcal{F}^\iota\rangle$. Thus $\langle \mathbf{t}_N, \mathbf{t}_N^\iota\rangle^{\mathrm{tr}} = \langle \mathcal{F}_N, \mathcal{F}_N^\iota\rangle$ and $\langle \mathcal{F}_H, \mathcal{F}_N^\iota\rangle \geq \langle \mathbf{t}_H, \mathbf{t}_N^\iota\rangle^{\mathrm{tr}}$. Therefore,

$$\langle \mathcal{F}_H, \mathcal{F}_N^\iota\rangle \geq \langle \mathbf{t}_H, \mathbf{t}_N^\iota\rangle^{\mathrm{tr}} = \langle \mathbf{t}_H^\#, \mathbf{t}_N^\iota\rangle^{\mathrm{tr}} = \langle \mathbf{t}_N, \mathbf{t}_N^\iota\rangle^{\mathrm{tr}} = \langle \mathcal{F}_N, \mathcal{F}_N^\iota\rangle,$$

where the first equality is true because the operator $\langle \cdot, \cdot\rangle^{\mathrm{tr}}$ is invariant under the action of Galois and additive by [8, 4.2.7], and the second equality is given by our hypothesis. Thus, $\mathcal{F}_N^\iota = \mathcal{F}_F$. □

*Remark* 16. 1. We will show in Proposition 17rem.17, that if $V$ is isomorphic to a trivial germ of crystalline representations with $G$-structure, then $D$ is isomorphic to a trivial isocrystal with $G$-structure.

2. The filtration $\mathcal{F}_F(D)$ is the image of $\mathcal{F}_F(V)$ by $D_{\mathrm{cris}}$. However, since $D_{\mathrm{cris}}$ is not induced by an isomorphism between fiber functors, we cannot deduce directly that $\mathbf{t}_F(D) = \mathbf{t}_F(V)$.

Set $\mathcal{L}(V) = \mathcal{L}(\omega(V))$. Thus $\mathcal{L}(V) \neq \emptyset$ if and only if $V$ is isomorphic to a trivial germ of crystalline representations with $G$-structure, by Remark 4.2.2A variant: From $\mathrm{Rep}_{\mathbb{Q}_p} G$ to $\mathrm{Rep}_{\mathbb{Z}_p} G$subsubsection.4.2.2, which we assume from now on. Any $x$ in $\mathcal{L}(V)$ has a unique factorization through an exact and faithful ⊗-functor

$$\mathrm{Rep}_{\mathbb{Z}_p} G \to \mathrm{Rep}_{\mathbb{Z}_p}^{\mathrm{cr}}\{\mathrm{Gal}_{K_0}\}$$

that we also denote by $x$. We call such a ⊗-functor a germ of integral crystalline representations with $G$-structure.

**Definition 4.9.** 1. A finite extension $K \subset \overline{K}_0$ of $K_0$ is a field of definition of $V$ if $V$ factors through the full subcategory $\mathrm{Rep}_{\mathbb{Q}_p}^{\mathrm{cr}} \mathrm{Gal}_K$ of $\mathrm{Rep}_{\mathbb{Q}_p}^{\mathrm{cr}}\{\mathrm{Gal}_{K_0}\}$.

2. A finite extension $K \subset \overline{K}_0$ of $K_0$ is a field of definition of $x \in \mathcal{L}(V)$ if $x$ factors through the full subcategory $\mathrm{Rep}_{\mathbb{Z}_p}^{\mathrm{cr}} \mathrm{Gal}_K$ of $\mathrm{Rep}_{\mathbb{Z}_p}^{\mathrm{cr}}\{\mathrm{Gal}_{K_0}\}$. We denote by $\mathcal{L}(V, K)$ the set of lattices in $\mathcal{L}(V)$ having $K$ as field of definition.

**Lemma 4.19.** *There exists a field of definition of $V$.*

*Proof.* Let $\tau$ be a ⊗-generator of $\mathrm{Rep}_{\mathbb{Q}_p} G$. Then $V(\tau) \in \mathrm{Rep}_{\mathbb{Q}_p}^{\mathrm{cr}} \mathrm{Gal}_K$ for some large enough $K$, in which case also $V(\tau') \in \mathrm{Rep}_{\mathbb{Q}_p}^{\mathrm{cr}} \mathrm{Gal}_K$ for every $\tau' \in \mathrm{Rep}_{\mathbb{Q}_p} G$. Then $V$ factors through the full-subcategory $\mathrm{Rep}_{\mathbb{Q}_p}^{\mathrm{cr}} \mathrm{Gal}_K$ of $\mathrm{Rep}_{\mathbb{Q}_p}^{\mathrm{cr}}\{\mathrm{Gal}_{K_0}\}$. □

Let $K$ be a field of definition of $V$. Then $\mathrm{Gal}_K$ acts on $\mathcal{L}(V)$ by

$$(g \cdot x)(\tau) = g \cdot x(\tau)$$

for $g \in \mathrm{Gal}_K$, $x \in \mathcal{L}(V)$ and $\tau \in \mathrm{Rep}_{\mathbb{Z}_p} G$. Plainly, $\mathcal{L}(V, K') = \mathcal{L}(V)^{\mathrm{Gal}_{K'}}$ for every finite extension $K'$ of $K$.

**Lemma 4.20.** *We have $\mathcal{L}(V) = \cup \mathcal{L}(V, K')$, for $K' \subset \overline{K}_0$ running through the finite extensions of $K$. In particular, every $x \in \mathcal{L}(V)$ has a field of definition.*

*Proof.* We may assume that $V$ is trivial. Then $V : \text{Rep}_{\mathbb{Q}_p} G \to \text{Rep}_{\mathbb{Q}_p}^{\text{cr}} \text{Gal}_K$ is induced by a continuous morphism $\text{Gal}_K \to \text{Aut}^{\otimes}(V) = G(\mathbb{Q}_p)$. Since any $x \in \mathcal{L}(V) \simeq G(\mathbb{Q}_p)/G(\mathbb{Z}_p)$ has an open stabilizer in $G(\mathbb{Q}_p)$, it also has an open stabilizer in $\text{Gal}_K$, which proves the lemma. $\square$

Fix $K$, a field of definition of $V$ such that $\mathcal{L}(V, K) \neq \emptyset$, with uniformizer $\pi_K$ and $E_K$ the minimal polynomial of $\pi_K$. Then, $V$ induces an isogeny class of Kisin modules with $G$-structure $N$, defined by the exact and faithful $\otimes$-functor

$$N = \mathfrak{N} \circ V \;:\; \text{Rep}_{\mathbb{Q}_p} G \to \text{Mod}_{\mathfrak{S}[\frac{1}{p}]}^{\varphi, E_K}.$$

Let

$$D' = N/uN : \text{Rep}_{\mathbb{Q}_p} G \to^{\text{wa}} \text{MF}_K^{\sigma}$$

the weakly admissible filtered $G$-isocrystal associated to $N$. We also denote by $D'$ the underlying $G$-isocrystal.

*Note* 1. Until now, we have avoided the superscript $E_K$ on categories of Kisin modules, to ease the notations, as we were working with a fixed field $K$. Since we will be changing $K$, depending on the field of definition of germs of integral crystalline representations with $G$-structure, we use now the rigorous notation, including $E_K$.

Given an integral crystalline representation with $G$-structure $x \in \mathcal{L}(V, K)$, we can associate a Kisin module with $G$-structure $M$, with isogeny class $N$, by composition with the Kisin functor $\mathfrak{M}$. Next proposition tells us that $M$ is an exact functor (we already know that it is a faithful $\otimes$-functor) giving us a map

$$\begin{array}{rcl} \mathcal{L}(V, K) & \to & \mathcal{L}(N) \\ x & \mapsto & \mathfrak{M} \circ x. \end{array}$$

**Proposition 4.21.** *Let $x : \text{Rep}_{\mathbb{Z}_p} G \to \text{Rep}_{\mathbb{Z}_p}^{\text{cr}} \text{Gal}_K$ be an integral crystalline representation with $G$-structure. Then the functor $M = \mathfrak{M}(x)$ is exact.*

*Proof.* We view $\omega(M)$ as a faithful $\otimes$-functor $\omega : \text{Rep}_{\mathbb{Z}_p} G \to \text{Bun}_{\mathfrak{S}}$. Let $\omega_U : \text{Rep}_{\mathbb{Z}_p} G \to \text{Bun}_U$ be the restriction of $\omega(M)$ to the open set $U = \text{Spec}\,\mathfrak{S} \setminus \{\mathfrak{m}\} \subset \text{Spec}\,\mathfrak{S}$. Then $\omega_U$ is a faithful $\otimes$-functor which is exact by the exactness of $x$ and the properties of the Kisin functor $\mathfrak{M}$. By Broshi [6], such functors are classified by $H^1(U, G)$. By [7, Theorem 6.13], we have

$$H^1(U, G) = H^1(\mathfrak{S}, G)$$

and we have seen in the proof of 4.2thm.4.2 that $H^1(\mathfrak{S}, G) = 0$. Thus, $\omega_U \simeq \omega_{G,U}$. But then, we have

$$\omega = \Gamma(U, -) \circ \omega_U \simeq \Gamma(U, -) \circ \omega_{G,U} \simeq \omega_{G,\mathfrak{S}},$$

so $\omega \simeq \omega_{G,\mathfrak{S}}$ is exact. $\square$

*Remark* 17. In particular, $\mathcal{L}(N) \neq \emptyset$ since $\mathcal{L}(V, K) \neq \emptyset$. Thus $\omega(N) \simeq \omega_{G_{\mathbb{Q}_p}, \mathfrak{S}[\frac{1}{p}]}$ by Proposition 4.15thm.4.15. Therefore,

$$\omega(D') = \omega(N) \otimes K_0 \simeq \omega_{G_{\mathbb{Q}_p}, \mathfrak{S}[\frac{1}{p}]} \otimes K_0 = \omega_{G_{\mathbb{Q}_p}, K_0}$$

so $D'$ is isomorphic to a trivial $G$-isocrystal. By Kisin's construction, there is an isomorphism $\eta : D \to D'$, so $D$ is isomorphic to a trivial $G$-isocrystal, as claimed in last Remark.

**Proposition 4.22.** *We have $\mathbf{t}_{\text{H}}(N) = \mathbf{t}_{\text{H}}(D')^{\iota} = \mathbf{t}_{\text{H}}(D)^{\iota}$.*

*Proof.* The isomorphism $D \simeq D'$ gives $\mathbf{t}_{\text{H}}(D) = \mathbf{t}_{\text{H}}(D')$. On the other hand, the $\otimes$-isomorphism

$$\omega(D') \otimes_{K_0} K \to \omega(\varphi_N \varphi^* N) \otimes_{\mathfrak{S}} K$$

maps $\mathcal{F}_H(D')$ to $\mathcal{F}_H^\iota(N) = \mathcal{F}(\omega(\varphi_N\varphi^*N) \otimes \hat{\mathfrak{S}}, \omega(N) \otimes \hat{\mathfrak{S}})$ as we have seen in section 3.2.1Filtered isocrystalssubsubsection.3.2.1. Therefore,

$$\begin{aligned}\mathbf{t}_H(D') &= \mathbf{t}(\mathcal{F}_H(D'))\\ &= \mathbf{t}(\mathcal{F}_H^\iota(N))\\ &= \operatorname{Pos}\left(\omega(\varphi_N\varphi^*N) \otimes \hat{\mathfrak{S}}, \omega(N) \otimes \hat{\mathfrak{S}}\right)\\ &= \operatorname{Pos}\left(\omega(N) \otimes \hat{\mathfrak{S}}, \omega(\varphi_N\varphi^*N) \otimes \hat{\mathfrak{S}}\right)^\iota\\ &= \mathbf{t}_H^\iota(N)\end{aligned}$$

□

The next proposition gives us the relation between the Fargues type defined on $V$ and the types defined on the torsion Kisin module with $G$-structure $\overline{M}$, for $x \in \mathcal{L}(V,K)$ and $M = \mathfrak{M} \circ x \in \mathcal{L}(N)$.

**Proposition 4.23.** *We have*

1. *The inequality $e \cdot \mathbf{t}_F(V) \leq \mathbf{t}_H(\overline{M})^\#$.*

2. *If $\overline{M}$ is aligned, then $e \cdot \mathbf{t}_F(V) \leq \mathbf{t}_{F,1}(\overline{M})$. If, moreover, $e \cdot \mathbf{t}_F(V) = \mathbf{t}_{F,1}(\overline{M})$, then, $M$ is HN-type.*

*Proof.* 1. In section 4 we have seen that for all $\tau \in \operatorname{Rep}_{\mathbb{Z}_p} G$, there are inequalities

$$e \cdot \mathbf{t}_F(V)(\tau) = e \cdot \mathbf{t}_F(N(\tau)) \leq \mathbf{t}_{F,1}(\overline{M})(\tau) \leq \mathbf{t}_H(\overline{M})(\tau)$$

(the coefficient $e$ appears since we have not renormalized the types defined on $\overline{M}$ as we did in section 5), since proposition 3.15thm.3.15 gives us $\mathbf{t}_{F,1}(\overline{M}(\tau)) \leq \mathbf{t}_H(\overline{M}(\tau))$ for all $\tau$. Using proposition 4.1thm.4.1, we get

$$e \cdot \mathbf{t}_F(V) = e \cdot \mathbf{t}_F(V)^\# \leq \mathbf{t}_H(\overline{M})^\#,$$

the first equality given by Remark 15rem.15.

2. For every $\tau \in \operatorname{Rep}_{\mathbb{Z}_p} G$, we have

$$e \cdot \mathbf{t}_F(V)(\tau) = e \cdot \mathbf{t}_F(N(\tau)) \leq \mathbf{t}_{F,1}(\overline{M})(\tau),$$

thus $e \cdot \mathbf{t}_F(V) \leq \mathbf{t}_{F,1}(\overline{M})$ by Proposition 4.1thm.4.1, since both types are invariant by $\#$ by Remark 15rem.15 and Proposition 4.14thm.4.14, respectively. The equality implies that

$$e \cdot \mathbf{t}_F(V)(\tau) = e \cdot \mathbf{t}_F(N)(\tau) = \mathbf{t}_{F,1}(\overline{M})(\tau)$$

for every $\tau \in \operatorname{Rep}_{\mathbb{Z}_p} G$, so $M(\tau)$ is HN-type by Proposition 3.23thm.3.23, and $M$ is HN-type. □

For every finite extension $K_0 \subset K \subset \overline{K}_0$, we can construct maps

$$\begin{array}{rccc} D_{\operatorname{cris}}^K : & \mathcal{L}(V,K) & \to & \mathcal{L}(D) \\ & x & \mapsto & y \end{array}$$

where $y = \eta(y')$ for $y' \in \mathcal{L}(D')$ the image of $M = \mathfrak{M} \circ x$ for the map $\mathcal{L}(N) \to \mathcal{L}(D')$ constructed in last section, and $\eta : \mathcal{L}(D') \simeq \mathcal{L}(D)$ induced by the isomorphism between Fontaine's $D_{\operatorname{cris}}$ and Kisin's functor $\mathfrak{M} \otimes K_0$. This construction extends to a map

$$\operatorname{red} : \mathcal{L}(V) \to \mathcal{L}(D)$$

thanks to Liu's results in 3.13thm.3.13 which gives us compatibility of $D_{\operatorname{cris}}$ with base change for $K'$ a totally ramified extension of $K$, and Proposition 4.20thm.4.20 which gives us $\mathcal{L}(V) = \cup \mathcal{L}(V,K)$.

We can still say a bit more about the image of this functor.

**Proposition 4.24.** *There is a factorization*

$$\begin{array}{ccc} \mathcal{L}(V) & \xrightarrow{\operatorname{red}} & \mathcal{L}(D) \\ & \searrow^{\operatorname{red}} & \nearrow \\ & \mathcal{L}(D, \leq \mathbf{t}_H^\iota(D)) & \end{array}$$

*where the maps are equivariant under the action of* $\mathrm{Aut}^\otimes(V) \to \mathrm{Aut}^\otimes(D)$.

*Proof.* We have $\mathcal{L}(D) \simeq \mathcal{L}(D')$ and $\mathcal{L}(D, \leq \mathbf{t}_\mathrm{H}^\iota(D)) \simeq \mathcal{L}(D', \mathbf{t}_\mathrm{H}^\iota(D'))$, by Proposition 4.22thm.4.22. We conclude by the analogous result for $D'$ given in Proposition 4.16thm.4.16, since $\mathbf{t}_\mathrm{H}(N) = \mathbf{t}_\mathrm{H}(D)^\iota$ by Proposition 4.22thm.4.22. □

### 4.3.5 The ordinary case

Fix a germ of crystalline representations with $G$-structure $V$ inducing a filtered $G$-isocrystal $D = D_\mathrm{cris} \circ V$.

**Proposition 4.25.** *Suppose $V$ isomorphic to a trivial germ of crystalline representations with $G$-structure, then $D$ is isomorphic to a trivial filtered $G$-isocrystal and we have $\mathbf{t}_\mathrm{N}(D) \leq \mathbf{t}_\mathrm{H}(D)^\#$ in $\mathbf{C}^\mathbb{Q}(G_{K_0})$.*

*Proof.* We have $V$ isomorphic to a trivial germ of crystalline representations with $G$-structure if and only if $\mathcal{L}(V) \neq \emptyset$, which implies that $\mathcal{L}(D) \neq \emptyset$, so $D$ is isomorphic to a trivial filtered $G$-isocrystal. The inequality is obtained assembling Proposition 4.10thm.4.10, Proposition 4.16thm.4.16 and Proposition 4.22thm.4.22. □

**Definition 4.10.** We say that $V$ is ordinary when

$$\mathbf{t}_\mathrm{N}(D) = \mathbf{t}_\mathrm{H}(D)^\# \quad \text{in} \quad \mathbf{C}^\mathbb{Q}(G_{K_0}).$$

Under the ordinarity condition, the more general proposition in last subsection becomes:

**Corollary 4.26.** *Suppose that $V$ is ordinary, then there is a factorization of $\mathrm{red} : \mathcal{L}(V) \to \mathcal{L}(D)$, by*

$$\mathcal{L}(V) \xrightarrow{\mathrm{red}} \mathcal{L}(D, \mathbf{t}_\mathrm{H}^\iota(D)) \hookrightarrow \mathcal{L}(D).$$

*In particular, $D$ is $\mathbf{t}_\mathrm{H}^\iota(D)$-ordinary.*

*Proof.* Let $x \in \mathcal{L}(V)$ with image $\mathrm{red}(x) = y \in \mathcal{L}(D)$. From $\mathbf{t}_\mathrm{N}^\iota(D) \leq \mathbf{t}_\mathrm{H}(y)^\# \leq \mathbf{t}_\mathrm{H}^\iota(D)^\#$ and the ordinary hypothesis $\mathbf{t}_\mathrm{N}^\iota(D) = \mathbf{t}_\mathrm{H}^\iota(D)^\#$, we get the equality $\mathbf{t}_\mathrm{H}(y)^\# = \mathbf{t}_\mathrm{H}^\iota(D)^\#$. Moreover, we have seen that $\mathbf{t}_\mathrm{H}(y) \leq \mathbf{t}_\mathrm{H}^\iota(D)$. Thus actually, $\mathbf{t}_\mathrm{H}(y) = \mathbf{t}_\mathrm{H}^\iota(D)$. □

**Theorem 4.27.** *Suppose $V$ is ordinary, with field of definition $K$. Given $x \in \mathcal{L}(V, K)$ and $N$ as before, we have that for every $M \in \mathcal{L}(N)$ (for instance $M = \mathfrak{M} \circ x$), $M$ is HN-type. Therefore $\mathcal{F}_\mathrm{F}(M)$ exists.*

*Proof.* We use the fact that $\overline{M}$ is aligned, which means that $\mathcal{F}_{\mathrm{F},1}$ is an exact $\otimes$-functor, so a filtration in $\mathbf{F}^\mathbb{Q}(G_{\mathbb{F}[[u]]})$ and we can consider its type $\mathbf{t}_{\mathrm{F},1}$ as an element in $\mathbf{C}^\mathbb{Q}(G_{\mathbb{F}[[u]]})$. Then, we get inequalities

$$e \cdot \mathbf{t}_\mathrm{F}(V) \leq \mathbf{t}_{\mathrm{F},1}(\overline{M}) \leq \mathbf{t}_\mathrm{H}(\overline{M})^\# \leq e \cdot \mathbf{t}_\mathrm{H}(N)^\#$$

by Proposition 4.23thm.4.23, Proposition 4.14thm.4.14 and Proposition 4.16thm.4.16. Again, by ordinarity, since $\mathbf{t}_\mathrm{F}(V) = \mathbf{t}_\mathrm{F}(D) = \mathbf{t}_\mathrm{H}^\iota(D)^\# = \mathbf{t}_\mathrm{H}(N)^\#$, we get $e \cdot \mathbf{t}_\mathrm{F}(V) = \mathbf{t}_{\mathrm{F},1}(\overline{M})$ and, by Proposition 4.23thm.4.23, $M$ is HN-type. □

Next proposition allows us to associate a $\mathbb{Z}$-filtration to $\mathcal{F}_\mathrm{F}(V)$, which will be used to define some operators:

**Lemma 4.28.** *There exists an integer $s > 0$ such that $s\mathcal{F}_\mathrm{F}(V) \in \mathbf{F}^\mathbb{Z}(\omega(V))$.*

*Proof.* We may suppose that $V$ is trivial. Then $\mathcal{F}_\mathrm{F}(V)$ is (non-canonically) split by a morphism $\mathcal{G} : \mathbb{D}(\mathbb{Q}) \to G_{\mathbb{Q}_p}$. The kernel of this morphism is $\mathbb{D}(\mathbb{Q}/N)$, where $N$ is a finitely generated subgroup of $\mathbb{Q}$, i.e. there exists an integer $s > 0$ such that $sN \subset \mathbb{Z}$. Thus $s\mathcal{G} : \mathbb{D}(\mathbb{Q}) \to G_{\mathbb{Q}_p}$ factors through

$$s\mathcal{G} : \mathbb{D}(\mathbb{Q}) \to \mathbb{D}(\mathbb{Z}) \to G_{\mathbb{Q}_p}$$

and $\mathrm{Fil}(s\mathcal{G}) = s\,\mathrm{Fil}(\mathcal{G}) = s\mathcal{F}_\mathrm{F}(V)$ is a $\mathbb{Z}$-filtration. □

We fix $s > 0$ such that $s\mathcal{F}_F(V) \in \mathbf{F}^{\mathbb{Z}}(\omega(V))$. We define the étale and crystalline operators by

$$\Phi^s_{\text{ét}} : \begin{array}{ccc} \mathcal{L}(V) & \to & \mathcal{L}(V) \\ x & \mapsto & \Phi^s_{\text{ét}}(x) = x + s\mathcal{F}_F(V) \end{array}$$

and

$$\Phi^s_{\text{cris}} : \begin{array}{ccc} \mathcal{L}(D) & \to & \mathcal{L}(D) \\ y & \mapsto & \Phi^s_{\text{cris}}(y) = y + s\mathcal{F}^\iota_N(D), \end{array}$$

where the addition of a lattice and a filtration is given in Definition 4.1defi.4.1. The following proposition gives us the compatibility between the two operators.

**Proposition 4.29.** *Suppose that $V$ is ordinary and let $x \in \mathcal{L}(V)$. Then, $\text{red}(\Phi^s_{\text{ét}}(x)) = \Phi^s_{\text{cris}}(\text{red}(x))$.*

*Proof.* Fix a field of definition $K$ for $x$, so that $x \in \mathcal{L}(V,K)$. Then also $x' = \Phi^s_{et}(x) \in \mathcal{L}(V,K)$ since $\mathcal{F}_F(V)$ is a $\text{Gal}_K$-stable filtration on $V$. Fix also a uniformizer $\pi_K$ of $K$, giving rise to $N = \mathfrak{N}(V)$, $M = \mathfrak{M}(x)$ and $M' = \mathfrak{M}(x')$ in $\mathcal{L}(N)$, as well as their images $y$ and $y'$ in $\mathcal{L}(D')$, where $D'$ is the (weakly admissible filtered) $G$-isocrystal obtained from $N$ by reduction modulo $u$. We have to show that

$$y' = \Phi^s_{\text{cris}}(y) \quad \text{in} \quad \mathcal{L}(D')$$

where $\Phi^s_{\text{cris}}(y) = y + s\mathcal{F}^\iota_N(D')$: the isomorphism of $G$-isocrystals $\eta : D' \to D$ transports this equality in $\mathcal{L}(D')$ to the desired equality in $\mathcal{L}(D)$.

Since $V$ is ordinary, $M$ is HN-type and $\mathcal{F}_F(M)$ exists. By compatibility of the various Fargues filtrations,

$$\mathcal{F}_F(N) := \mathcal{F}_F(M)[\tfrac{1}{p}] = \mathfrak{N}(\mathcal{F}_F(V))$$

and this reduces to $\mathcal{F}_F(D')$, which equals $\mathcal{F}^\iota_N(D')$ by Proposition 4.18thm.4.18. In particular, multiplying any of these $\mathbb{Q}$-filtrations by $s$ yields $\mathbb{Z}$-filtrations. Let $\omega(M) = \oplus_{a \in \mathbb{Z}} \omega_a(M)$ be a splitting of $\mathcal{F}_F(M)$. Note that such a splitting exists by Proposition 4.15thm.4.15 and [8, 3.11.3]. It gives rise to splittings

$$\omega(N) = \oplus_{a \in \mathbb{Z}} \omega_a(N) \quad \text{and} \quad \omega(D') = \oplus_{a \in \mathbb{Z}} \omega_a(D')$$

of $\mathcal{F}_F(N)$ and $\mathcal{F}_F(D') = \mathcal{F}^\iota_N(D')$. Consider now some $\tau \in \text{Rep}_{\mathbb{Z}_p} G$. By Proposition 3.10thm.3.10 and the computations preceeding it,

$$\begin{aligned} M'(\tau) &= M(\tau) +_{\text{fr}} s\mathcal{F}_F(N)(\tau) \\ &= M(\tau) + s\mathcal{F}_F(N)(\tau) \end{aligned}$$

with underlying $\mathfrak{S}$-module $\omega(M')(\tau) = \oplus_{a \in \mathbb{Z}} p^{-a} \omega_a(M)(\tau)$. This reduces modulo $u$ to the $W(\mathbb{F})$-lattice $y'(\tau) = \oplus_{a \in \mathbb{Z}} p^{-a} \omega_a(y)(\tau)$ in $D'(\tau)$, where $y = \oplus_{a \in \mathbb{Z}} \omega_a(y)$ is the reduction of $\omega(M) = \oplus_{a \in \mathbb{Z}} \omega_a(M)$, so that also $\omega_a(y)[\tfrac{1}{p}] = \omega_a(D')$ for every $a \in \mathbb{Z}$. On the other hand,

$$\begin{aligned} \Phi^s_{\text{cris}}(y)(\tau) &= y(\tau) + s\mathcal{F}^\iota_N(D')(\tau) \\ &= \sum_{i \in \mathbb{Z}} p^{-i} y(\tau) \cap \left(s\mathcal{F}^\iota_N(D')\right)^i(\tau) \\ &= \sum_{i \in \mathbb{Z}} \left(\oplus_a p^{-i} \omega_a(y)(\tau)\right) \cap \left(\oplus_{a \geq i} \omega_a(D')(\tau)\right) \\ &= \sum_{i \in \mathbb{Z}} \left(\oplus_{a \geq i} p^{-i} \omega_a(y)(\tau)\right) \\ &= \oplus_{a \in \mathbb{Z}} \sum_{i \leq a} p^{-i} \omega_a(y)(\tau) \\ &= \oplus_{a \in \mathbb{Z}} p^{-a} \omega_a(y)(\tau). \end{aligned}$$

Thus $\Phi^s_{\text{cris}}(y)(\tau) = y'(\tau)$ for every $\tau \in \text{Rep}_{\mathbb{Z}_p}(G)$ and $\Phi^s_{\text{cris}}(y) = y'$. $\square$

**Lemma 4.30.** *Let $x, y \in \mathcal{L}(V)$ and let*

$$\pi : \mathcal{L}(V) \to U_{\mathcal{F}_F}(\mathbb{Q}_p) \backslash \mathcal{L}(V)$$

*be the usual projection, for $\mathcal{F}_F$ the Fargues filtration on $V$ and $U_{\mathcal{F}_F}$ the unipotent radical of the parabolic subgroup of $G$ associated to $\mathcal{F}_F$. Then, we have $\pi(x) = \pi(y)$ if and only if*

$$(\Phi^s_{\text{ét}})^{(n)}(x) = (\Phi^s_{\text{ét}})^{(n)}(y)$$

*for $n$ large enough.*

*Proof.* We may assume $V$ is trivial. We use the results in [8] as follows: In 6.2.7, we see that $\mathbf{B}^e(\omega_G)$ is a tight building, so verifies the axiom $ST$. It also verifies $L(s)$, thus it verifies $UN^+$, by Lemma 114. Then, the lemma holds by 5.6.2. □

**Lemma 4.31.** *The operator $\Phi_{\mathrm{cris}}^s$ is bijective on the set $\mathcal{L}(D, \mathbf{t}_{\mathrm{H}}^\iota(D))$.*

*Proof.* As we will see (Remark 18rem.18), the lemma is true in the abelian case. We can construct an abelian germ of crystalline representations $V'$ such that $D' = D_{\mathrm{cris}}(V')$ is $\mathbf{t}_{\mathrm{H}}^\iota(D)$-ordinary: fix a Borel pair $(T_G, B_G)$ in $G$ and lift $\mathbf{t}_{\mathrm{H}}(D)$ to a $B_G$-dominant cocharacter $\tilde{\mathbf{t}}_{\mathrm{H}}(D) : \mathbb{G}_{m,K_0} \to T_{G,K_0}$. Let $V'$ be the abelian germ of crystalline representations associated, by Proposition 5.7thm.5.7, to $\tilde{\mathbf{t}}_{\mathrm{H}}(D)$. By Proposition 5.9thm.5.9, $V'$ is ordinary and, by Proposition 4.26thm.4.26, $D'$ is $\mathbf{t}_{\mathrm{H}}^\iota(D)$-ordinary. Now, by Proposition 4.11thm.4.11, there is a unique $\mathbf{t}_{\mathrm{H}}^\iota(D)$-ordinary $G$-isocrystal up to isomorphism, we get $D' \simeq D$ and the lemma follows. □

As a consequence of the three previous results, we obtain:

**Theorem 4.32.** *Notations as above. For $V$ an ordinary germ of crystalline representations with $G$-structure, the map* red *admits a factorization*

$$\mathcal{L}(V) \xrightarrow{\mathrm{red}} \mathcal{L}(D, \mathbf{t}_{\mathrm{H}}^\iota(D)).$$
$$\searrow_\pi \quad \nearrow$$
$$U_{\mathcal{F}_{\mathrm{F}}}(\mathbb{Q}_p) \backslash \mathcal{L}(V)$$

*Proof.* Suppose we have $\pi(x) = \pi(y)$ for $x, y \in \mathcal{L}(V)$. By last lemma, we have $(\Phi_{\mathrm{\acute{e}t}}^s)^{(n)}(x) = (\Phi_{\mathrm{\acute{e}t}}^s)^{(n)}(y)$ for $n$ large enough. We apply the corollary above to obtain

$$(\Phi_{\mathrm{cris}}^s)^{(n)}(\mathrm{red}(x)) = \mathrm{red}((\Phi_{\mathrm{\acute{e}t}}^s)^{(n)}(x)) = \mathrm{red}((\Phi_{\mathrm{\acute{e}t}}^s)^{(n)}(y)) = (\Phi_{\mathrm{cris}}^s)^{(n)}(\mathrm{red}(y))$$

for $n$ large enough, and the theorem follows by Lemma 4.31thm.4.31. □

# 5 A particular case: the abelian case

The aim of this section is to describe more precisely the results obtained in last section in the case where the germ of crystalline representations with $G$-structure is abelian. We use the same notations from last section. Let $\mathbb{F}$ be an algebraically closed field of characteristic $p > 0$, $K_0 = \mathrm{Frac}\, W(\mathbb{F})$, we fix an algebraic closure $\overline{K}_0$ of $K_0$ and an embedding $\iota_0 : \overline{\mathbb{Q}}_p \subset \overline{K}_0$. We shall use the notation $E$ for a finite extension $\mathbb{Q}_p \subset E \subset \overline{\mathbb{Q}}_p$, with ring of integers $\mathcal{O}_E$ and uniformizer $\pi_E$, and the notation $K$ for a finite extension $K_0 \subset K \subset \overline{K}_0$, with ring of integers $\mathcal{O}_K$ and uniformizer $\pi_K$.

## 5.1 Preliminaries

### 5.1.1 The pro-tori $T_K$ and $T$

Denote by $T_E = \mathrm{Res}_{E/\mathbb{Q}_p}(\mathbb{G}_{m,E})$ the torus over $\mathbb{Q}_p$ associated to the group $E^\times$, i.e. for a $\mathbb{Q}_p$-algebra $R$, we define $T_E(R)$ as the group of units of $E \otimes_{\mathbb{Q}_p} R$. For $E \subset E' \subset K$, the norm $N_{E'/E}$ induces a morphism $T_{E'} \to T_E$. We define two pro-tori by

$$T_K = \varprojlim_{E \subset K} T_E \quad \text{and} \quad T = \varprojlim_{E \subset \overline{K}_0} T_E$$

where the norms $N_{E'/E}$ are the transition maps. For every $E \subset K$, there is a continuous morphism

$$\chi_{K,E} : \mathrm{Gal}_K \to \mathrm{Gal}(E^{\mathrm{ab}}/E^{\mathrm{ur}}) \simeq \mathcal{O}_E^\times \subset E^\times = T_E(\mathbb{Q}_p)$$

and since these morphisms are compatible with the norm maps, they define a continuous morphism

$$\chi_K : \mathrm{Gal}_K \to T_K(\mathbb{Q}_p) = \varprojlim_{E \subset K} T_E(\mathbb{Q}_p).$$

The maps $T \to T_K \to T_E$ give us maps between the characters groups $X^*(T_E) \to X^*(T_K) \to X^*(T)$ and the cocharacters groups $X_*(T) \to X_*(T_K) \to X_*(T_E)$. For a finite extension $E$ of $\mathbb{Q}_p$, we have

$$X^*(T_E) \simeq \mathcal{C}(\mathrm{Hom}(E, \overline{\mathbb{Q}}_p), \mathbb{Z}),$$

the set of functions $f : \mathrm{Hom}(E, \overline{\mathbb{Q}}_p) \to \mathbb{Z}$. The character $\chi_f : T_{E|\overline{\mathbb{Q}}_p} \to \mathbb{G}_{m,\overline{\mathbb{Q}}_p}$ corresponding to $f$ induces

$$\begin{array}{ccc} T_E(\overline{\mathbb{Q}}_p) & \xrightarrow{\chi_f} & \mathbb{G}_m(\overline{\mathbb{Q}}_p) \\ \| & & \| \\ (E \otimes \overline{\mathbb{Q}}_p)^\times & \longrightarrow & \overline{\mathbb{Q}}_p^\times \\ x & \mapsto & \prod_{\iota : E \hookrightarrow \overline{\mathbb{Q}}_p} (\iota \otimes \mathrm{id})(x)^{f(\iota)}. \end{array}$$

The action of $\mathrm{Gal}_{\mathbb{Q}_p}$ on $X^*(T_E)$ is given by

$$(\sigma \cdot f)(\iota) = f(\sigma^{-1} \circ \iota), \quad \text{for } \sigma \in \mathrm{Gal}_{\mathbb{Q}_p},\ \iota : E \hookrightarrow \overline{\mathbb{Q}}_p.$$

Since $T = \varprojlim T_E$, for $E$ running through all finite extensions of $\mathbb{Q}_p$ contained in $\overline{K}_0$,

$$X^*(T) = \varinjlim X^*(T_E) \simeq \mathcal{C}(\mathrm{Gal}_{\mathbb{Q}_p}, \mathbb{Z}),$$

the set of locally constant functions $f : \mathrm{Gal}_{\mathbb{Q}_p} \to \mathbb{Z}$ with the Galois action

$$(\sigma \cdot f)(\sigma') = f(\sigma^{-1}\sigma') \quad \text{for } \sigma, \sigma' \in \mathrm{Gal}_{\mathbb{Q}_p}\ .$$

Similarly, for $T_K = \varprojlim_{E \subset K} T_E$,

$$X^*(T_K) = \varinjlim_{E \subset K} X^*(T_E) \simeq \mathcal{C}(\mathrm{Gal}_{\mathbb{Q}_p}\,/\,\mathrm{Gal}_{K \cap \overline{\mathbb{Q}}_p}, \mathbb{Z}),$$

the set of locally constant, right $\mathrm{Gal}_{K \cap \overline{\mathbb{Q}}_p}$-invariant functions $f : \mathrm{Gal}_{\mathbb{Q}_p} \to \mathbb{Z}$.

### 5.1.2 The Hodge cocharacters

We define the Hodge cocharacter

$$\mu : \mathbb{G}_{m, \overline{\mathbb{Q}}_p} \to T_{|\overline{\mathbb{Q}}_p}$$

as the cocharacter corresponding to the evaluation at 1 in the dual group

$$X_*(T) = \mathrm{Hom}_{\mathrm{Gr}}(X^*(T), \mathbb{Z}) = \mathrm{Hom}_{\mathrm{Gr}}(\mathcal{C}(\mathrm{Gal}_{\mathbb{Q}_p}, \mathbb{Z}), \mathbb{Z}),$$

i.e. $\langle \cdot, \mu \rangle = \mathrm{ev}_1$, where $\mathrm{ev}_1$ sends $f$ to $f(1)$. Thus, for any $x \in \overline{\mathbb{Q}}_p^\times = \mathbb{G}_m(\overline{\mathbb{Q}}_p)$ and $f \in \mathcal{C}(\mathrm{Gal}_{\mathbb{Q}_p}, \mathbb{Z})$,

$$\chi_f \circ \mu(x) = x^{\langle \chi_f, \mu \rangle} = x^{f(1)}.$$

We denote by $\mu_K \in X^*(T_K)$ and $\mu_E \in X^*(T_E)$ the corresponding Hodge cocharacters of $T_K$ and $T_E$, which are respectively defined over $K \cap \overline{\mathbb{Q}}_p$ and $E$:

$$\mu_K : \mathbb{G}_{m, K \cap \overline{\mathbb{Q}}_p} \to T_{K\,|\,K \cap \overline{\mathbb{Q}}_p} \quad \text{and } \mu_E : \mathbb{G}_{m, E} \to T_{E\,|\,E}.$$

For every $E$-algebra $R$, the canonical decomposition $E \otimes E \simeq E \times E'$ where $E' = \ker m$ for

$$\begin{array}{ccc} E \otimes E & \xrightarrow{m} & E \\ x \otimes y & \mapsto & xy, \end{array}$$

yields a canonical decomposition $E \otimes R \simeq R \times R'$, and

$$\begin{array}{ccc} \mathbb{G}_{m,E}(R) & \xrightarrow{\mu_E} & T_E(R) \\ \| & & \| \\ R^\times & \longrightarrow & R^\times \times R'^\times \\ x & \longmapsto & (x, 1). \end{array}$$

In the dual group $X_*(T_E)$ of $X^*(T_E) = \mathcal{C}(\mathrm{Hom}(E, \overline{\mathbb{Q}}_p), \mathbb{Z})$, the cocharacter $\mu_E$ corresponds to evaluation at the given embedding $E \hookrightarrow \overline{\mathbb{Q}}_p$.

### 5.1.3 The Newton cocharacter

Similarly, integration against the Haar measure $\mu_{\text{Haar}}$ of $\text{Gal}_{\mathbb{Q}_p}$ defines a $\text{Gal}_{\mathbb{Q}_p}$-invariant morphism

$$i_{\text{Haar}} \;:\; \begin{array}{rcl} \mathcal{C}(\text{Gal}_{\mathbb{Q}_p}, \mathbb{Z}) & \to & \mathbb{Q} \\ f & \mapsto & \int_{\text{Gal}_{\mathbb{Q}_p}} f d\mu_{\text{Haar}} \end{array}$$

which yields the Newton cocharacter, defined over $\mathbb{Q}_p$, $\nu : \mathbb{D}(\mathbb{Q}) \to T$. We denote by

$$\nu_K : \mathbb{D}(\mathbb{Q}) \to T_K \quad \text{and} \quad \nu_E : \mathbb{D}(\mathbb{Q}) \to T_E$$

the corresponding Newton cocharacters for $T_K$ and $T_E$. Inside

$$\text{Hom}\left(\mathbb{D}(\mathbb{Q})_{|\overline{\mathbb{Q}}_p}, T_{E|\overline{\mathbb{Q}}_p}\right) = \text{Hom}\left(\mathcal{C}(\text{Hom}(E, \overline{\mathbb{Q}}_p), \mathbb{Q})\right)$$

the Newton cocharacter $\nu_E$ corresponds to the morphism

$$\begin{array}{rcl} \mathcal{C}(\text{Hom}(E, \overline{\mathbb{Q}}_p), \mathbb{Z}) & \to & \mathbb{Q} \\ f & \to & \frac{1}{[E:\mathbb{Q}_p]} \sum_{\iota: E \hookrightarrow \overline{\mathbb{Q}}_p} f(\iota). \end{array}$$

## 5.2 The Fontaine-Serre functor

**Definition 5.1.** We say that a crystalline representation $(V, \rho) \in \text{Rep}_{\mathbb{Q}_p}^{\text{cr}} \text{Gal}_K$ is abelian when its image $\rho(\text{Gal}_K) \subset \text{GL}(V)$ is abelian. We denote by $\text{Rep}_{\mathbb{Q}_p}^{\text{cr,ab}} \text{Gal}_K$ the category of abelian crystalline representations of $\text{Gal}_K$.

Denote by $\text{CrAb}(V, K) \subset \text{Cr}(V, K)$ the set of all morphisms $\rho : \text{Gal}_K \to \text{GL}(V)$ such that $(V, \rho)$ is an abelian crystalline representation of $\text{Gal}_K$, and set

$$\text{CrAb}(V) := \varinjlim_{\substack{K_0 \subseteq K \subset \overline{K}_0 \\ f}} \text{CrAb}(V, K) \subset \text{Cr}(V).$$

We say that a germ of crystalline representations $(V, \rho)$ is abelian when $\rho \in \text{CrAb}(V)$. We denote by $\text{Rep}_{\mathbb{Q}_p}^{\text{cr,ab}}\{\text{Gal}_{K_0}\}$ the category of germs of abelian crystalline representations, a full Tannakian subcategory of $\text{Rep}_{\mathbb{Q}_p}^{\text{cr}}\{\text{Gal}_{K_0}\}$.

In [39, 2], Fontaine proves, building on Serre's results in [38], that there is an equivalence of categories

$$V_K^u \;:\; \text{Rep}_{\mathbb{Q}_p} T_K \xrightarrow{\sim} \text{Rep}_{\mathbb{Q}_p}^{\text{cr,ab}} \text{Gal}_K$$

given by the composition of $\chi_K$ with a representation $T_K(\mathbb{Q}_p) \to \text{GL}(V)$.

**Proposition 5.1.** *The equivalences of neutralized Tannakian categories*

$$V_K^u \;:\; \text{Rep}_{\mathbb{Q}_p} T_K \xrightarrow{\sim} \text{Rep}_{\mathbb{Q}_p}^{\text{cr,ab}} \text{Gal}_K,$$

*for every finite extension $K_0 \subset K \subset \overline{K}_0$, extends to a diagram*

$$\begin{array}{ccc} \text{Rep}_{\mathbb{Q}_p} T_K & \hookrightarrow & \text{Rep}_{\mathbb{Q}_p} T \\ \downarrow V_K^u & & \downarrow V^u \\ \text{Rep}_{\mathbb{Q}_p}^{\text{cr,ab}} \text{Gal}_K & \hookrightarrow & \text{Rep}_{\mathbb{Q}_p}^{\text{cr,ab}}\{\text{Gal}_{K_0}\} \end{array}$$

*where $V^u$ is an equivalence of neutralized Tannakian categories.*

*Proof.* Let $(V, \rho) \in \text{Rep}_{\mathbb{Q}_p} T$. For a large enough finite extension $K$ of $K_0$, we have $(V, \rho) \in \text{Rep}_{\mathbb{Q}_p} T_K$, thus $(V, \rho \circ \chi_K) \in \text{Rep}_{\mathbb{Q}_p}^{\text{cr,ab}} \text{Gal}_K$. For a finite extension $K \subset K'$, we have $(V, \rho) \in \text{Rep}_{\mathbb{Q}_p} T_K \subset \text{Rep}_{\mathbb{Q}_p} T_{K'}$ and the commutative diagram

$$\begin{array}{ccc} \text{Gal}_K & \xrightarrow{\chi_K} & T_K(\mathbb{Q}_p) \\ \text{Res}\uparrow & & N_{K'/K}\uparrow \\ \text{Gal}_{K'} & \xrightarrow{\chi_{K'}} & T_{K'}(\mathbb{Q}_p) \end{array}$$

induces another commutative diagram

$$\begin{array}{ccc} \operatorname{Rep}_{\mathbb{Q}_p} T_K & \longrightarrow & \operatorname{Rep}_{\mathbb{Q}_p}^{\mathrm{cr,ab}} \operatorname{Gal}_K \\ \downarrow & & \downarrow \\ \operatorname{Rep}_{\mathbb{Q}_p} T_{K'} & \longrightarrow & \operatorname{Rep}_{\mathbb{Q}_p}^{\mathrm{cr,ab}} \operatorname{Gal}_{K'} \end{array}$$

thus $\rho \circ \chi_K \in \mathrm{CrAb}(V)$ does not depend upon the chosen $K$, so $V$ is a well-defined $\otimes$-functor. It is plainly compatible with fiber functors and faithful. The diagram is commutative, and any object or arrow on the bottom right category comes from an object or arrow on the bottom left category for a sufficiently large $K$, so $V$ is also full and essentially surjective. □

## 5.3 Wintenberger's functor

### 5.3.1 Universal norms

Let $K_0^{\mathrm{a}}$ be the composite of $K_0$ and $\overline{\mathbb{Q}}_p$ in $\overline{K}_0$. First we need to prove that there exists a norm compatible system of uniformizers $\pi = (\pi_K) \in \varprojlim K^\times$, for $K \subset K_0^{\mathrm{a}}$. For each finite extension $K_0 \subset K \subset \overline{K}_0$, let $v_K : K^\times \to \mathbb{Z}$ be the normalized valuations of $K$. For $K_0 \subset K_1 \subset K_2 \subset \overline{K}_0$, the norm map yields a morphism of exact sequences

$$\begin{array}{ccccccccc} 1 & \to & \mathcal{O}_{K_2}^\times & \to & K_2^\times & \stackrel{v_{K_2}}{\to} & \mathbb{Z} & \to & 0 \\ & & \downarrow & & \downarrow & & \| & & \\ 1 & \to & \mathcal{O}_{K_1}^\times & \to & K_1^\times & \stackrel{v_{K_1}}{\to} & \mathbb{Z} & \to & 0 \end{array}$$

It is known by [36, 6.5.6 and 6.5.8] that the norm maps $N_{K_2/K_1} : \mathcal{O}_{K_2}^\times \to \mathcal{O}_{K_1}^\times$ are all surjective. Every extension of $K_0$ contained in $K_0^{\mathrm{a}}$ can be written as $K_0 E$, for $E \subset \overline{\mathbb{Q}}_0$ a finite extension of $\mathbb{Q}_p$. There are finitely many finite extensions $E \subset \overline{\mathbb{Q}}_p$ of $\mathbb{Q}_p$ of a given degree, so $K_0$ has only finitely many finite extensions of any given degree. Writing $K_n$ for the composite of all finite extensions of $K_0$ whose degree is less than $n$, we thus obtain a totally ordered and countable cofinal subset of $\{K \mid K_0 \subset K \subset K_0^{\mathrm{a}}\}$, so we can use the Mittag-Leffler condition on this system to get an exact sequence

$$1 \to \varprojlim_{K_0 \subset K \subset K_0^{\mathrm{a}}} \mathcal{O}_K^\times \to \varprojlim_{K_0 \subset K \subset K_0^{\mathrm{a}}} K^\times \stackrel{v}{\to} \mathbb{Z} \to 0..$$

We fix $\pi = (\pi_K) \in \varprojlim_{K_0 \subset K \subset K_0^{\mathrm{a}}} K^\times$ such that $v(\pi) = 1$, i.e. a norm compatible system of uniformizers. From now on, we also require that $\pi_{K_0} = p$.

### 5.3.2 The element $b \in T(K_0)$

There is a morphism

$$\varprojlim_{K_0 \subset K \subset K_0^{\mathrm{a}}} K^\times \to T(K_0)$$

which maps

$$\{x_K\} \in \varprojlim_{K_0 \subset K \subset K_0^{\mathrm{a}}} K^\times$$

to the element

$$\{y_E\} \in T(K_0) = \varprojlim_{\mathbb{Q}_p \subset E \subset \overline{\mathbb{Q}}_p} T_E(K_0) = \varprojlim_{\mathbb{Q}_p \subset E \subset \overline{\mathbb{Q}}_p} (E \otimes K_0)^\times$$

defined as follows. Let $E_0 = E \cap K_0 = E \cap \mathbb{Q}_p^{\mathrm{nr}} \subset E$ be the maximal unramified extension of $\mathbb{Q}_p$ in $E$, and let $K \simeq E \otimes_{E_0} K_0$ be the compositum of $E$ and $K_0$ in $K_0^{\mathrm{a}}$. Then

$$E \otimes K_0 = E \otimes_{E_0} (E_0 \otimes_{\mathbb{Q}_p} K_0) = \prod_{i=0}^{[E_0:\mathbb{Q}_p]-1} E \otimes_{E_0,\sigma^i} K_0.$$

For $i = 0$, we have $E \otimes_{E_0} K_0 \simeq K$ by $e \otimes b \mapsto eb$. This construction allows us to define an inclusion

$$K^\times \simeq (E \otimes_{E_0} K_0)^\times \subset (E \otimes K_0)^\times = T_E(K_0).$$

We set $y_E$ as the image of $x_K$ under this inclusion.

**Lemma 5.2.** *This construction yields a well-defined morphism* $\varprojlim_{K_0 \subset K \subset K_0^{\mathrm{a}}} K^\times \to T(K_0)$.

*Proof.* Fix $E \subset E'$, giving $K \subset K'$ and $E_0 \subset E_0'$. We have to show that the norm $N_{E'/E} : T_{E'}(K_0) \to T_E(K_0)$ maps $y_{E'}$ to $y_E$. Considering the decomposition $E \otimes K_0 = \prod_i L_i$, where $L_i = E \otimes_{E_0, \sigma^i} K_0$ for $0 \le i \le a-1$ with $a = [E_0 : \mathbb{Q}_p]$, so that

$$\begin{aligned} E' \otimes K_0 &= E' \otimes_E (E \otimes K_0) \\ &= \prod_i E' \otimes_E (E \otimes_{E_0, \sigma^i} K_0) \\ &= \prod_i E' \otimes_{E_0, \sigma^i} K_0 \\ &= \prod_i E' \otimes_{E_0'} (E_0' \otimes_{E_0, \sigma^i} K_0) \\ &= \prod_i E' \otimes_{E_0', \sigma^{aj+i}} K_0 \\ &= \prod_i \prod_j L_{i,j} \end{aligned}$$

with $L_{i,j} = E' \otimes_{E_0', \sigma^{aj+i}} K_0$ for $0 \le j \le b-1$ with $b = [E_0' : E]$. We obtain $N_{E'/E}(x_{i,j}) = (\prod_j N_{E'/E}(x_{i,j}))_i$. For $y_{E'}$, all $x_{i,j} = 1$ except $x_{0,0} = \pi_{K'}$, so $N_{E'/E}(x_{i,j})$ has all components 1, except at $i = 0$, where we get $\prod_j N_{E'/E}(x_{0,j}) = N_{E'/E}(x_{0,0}) = N_{E'/E}(\pi_{K'}) = \pi_K$. Thus $N_{E'/E}(y_{E'}) = y_E$. □

We define $b$ as the image of $\pi$ in $T(K_0)$ by this morphism. Let $b_E$ and $b_K$ be the images of $b$ in $T_E(K_0)$ and $T_K(K_0)$, respectively.

**Lemma 5.3.** *We have $b_K = N_{K/K_0}(\mu_K(\pi_K))$ in $T_K(K_0)$, for $K \subset K_0^{\mathrm{a}}$.*

*Proof.* We look at the images of both elements in $T_E(K_0)$ for a $E$ large enough such that $K = E \cdot K_0$, using

$$\begin{aligned} E \otimes K_0 &= \prod_{i=0}^{[E_0:\mathbb{Q}_p]-1} E \otimes_{E_0, \sigma^i} K_0 = \prod_{i=0}^{[E_0:\mathbb{Q}_p]-1} L_i \\ E \otimes K &= \prod_{i=0}^{[E_0:\mathbb{Q}_p]-1} L_i \otimes_{K_0} K \end{aligned}$$

Then, $b$ maps to $b_E = (\pi_K, 1, \ldots, 1)$ in $T_E(K_0) = \prod L_i^\times$, $\mu_K(\pi_K)$ maps to $(x, 1, \ldots, 1)$ in $T_E(K) = \prod (L_i \otimes_{K_0} K)^\times$ and $N_{K/K_0}(\mu_K(\pi_K))$ maps to $(N_{K/K_0}(x), 1, \ldots, 1)$ in $\prod L_i^\times$, where $x = (\pi_K, 1)$ in $L_0 \otimes K = K \times K'$, thus $N_{K/K_0}(x) = \pi_K$ in $L_0 = K$, i.e. $b_K$ and $N_{K/K_0}(\mu_K(\pi_K))$ have the same image in $T_E(K_0)$ for every $E \subset K$ and the lemma holds. □

### 5.3.3 Wintenberger's functor

For $K \subset \overline{K}_0$, there is a strictly commutative diagram

$$\begin{array}{ccc} \mathrm{Rep}_{\mathbb{Q}_p} T_K & \hookrightarrow & \mathrm{Rep}_{\mathbb{Q}_p} T \\ \downarrow{D_{\pi_K}} & & \downarrow{D_\pi} \\ {}^{\mathrm{wa}}\mathrm{MF}_K^\sigma & \hookrightarrow & {}^{\mathrm{wa}}\mathrm{MF}_{\overline{K}_0}^\sigma \end{array}$$

of $\otimes$-categories. The two horizontal arrows are the obvious $\otimes$-functors. The left hand side of the diagram is given by Wintenberger's $\otimes$-functor

$$D_{\pi_K} \ : \ \mathrm{Rep}_{\mathbb{Q}_p} T_K \to{}^{\mathrm{wa}} \mathrm{MF}_K^\sigma$$

defined by taking, for a representation $\tau = (V, \rho) \in \mathrm{Rep}_{\mathbb{Q}_p} T_K$, the filtered isocrystal $D_{\pi_K}(\tau)$ given by:

- The underlying $K_0$-vector space $D_{\pi_K}(V) = V \otimes_{\mathbb{Q}_p} K_0$,
- The Frobenius morphism $\sigma_{D_{\pi_K}(V)} = \rho(b_K) \otimes \sigma$, for $b_K \in T_K(K_0)$ defined as above.
- The Hodge filtration given by $\mathcal{F}_{\mathrm{H}}^i(D_{\pi_K}(V)_K) = \oplus_{j \ge i} V_{K,j}$ for every $i \in \mathbb{Z}$, where $V \otimes_{\mathbb{Q}_p} K = \oplus_{i \in \mathbb{Z}} V_{K,i}$ is the weight decomposition attached to $\rho \circ \mu_K : \mathbb{G}_{m,K} \to T_K \to \mathrm{GL}(V)_K$.

The right hand side of the diagram is given by the $\otimes$-functor

$$D_\pi \ : \ \mathrm{Rep}_{\mathbb{Q}_p} T \to{}^{\mathrm{wa}} \mathrm{MF}_{\overline{K}}^\sigma$$

defined as follows: for every $(V, \rho) \in \mathrm{Rep}_{\mathbb{Q}_p} T$,

- The underlying module $D_\pi(V) = V \otimes_{\mathbb{Q}_p} K_0$,
- The Frobenius on $D_\pi(V)$ given by $\sigma_{D_\pi(V)} = \rho(b) \otimes \sigma$, for $b \in T(K_0)$ defined as above.

- The Hodge filtration on $D_\pi(V)_{\overline{K}_0}$ defined by taking $\mathcal{F}_H^i(D_\pi(V)_{\overline{K}_0}) = \oplus_{j \geq i} V_{\overline{K}_0, j}$ for every $i \in \mathbb{Z}$, where $V \otimes \overline{K}_0 = \oplus_{i \in \mathbb{Z}} V_{\overline{K}_0, i}$ is the weight decomposition attached to $\rho \circ \mu : \mathbb{G}_{m, \overline{K}_0} \to T_{\overline{K}_0} \to \mathrm{GL}(V)_{\overline{K}_0}$.

We may view $D_\pi$ as a trivial filtered isocrystal with $T$-structure.

**Proposition 5.4.** *The Hodge filtration $\mathcal{F}_H \circ D_\pi$ is split by $\mu : \mathbb{G}_{m, \overline{K}_0} \to T_{\overline{K}_0}$. The Newton graduation $\mathcal{G}_N \circ D_\pi$ is given by the cocharacter $\nu : \mathbb{D}(\mathbb{Q}) \to T$.*

*Proof.* The statement about the Hodge filtration is obvious from the definition. For the newton graduation, we need to prove that for every $\tau = (V, \rho) \in \mathrm{Rep}_{\mathbb{Q}_p} T$, the weight decomposition $V = \oplus_{\lambda \in \mathbb{Q}} V_\lambda$ given by $\rho \circ \nu$ is the Newton graduation of $D_\pi(\tau)$, for every $\lambda \in \mathbb{Q}$.

We may assume that $\tau$ factors through $T_E$ for a finite extension $E$. By definition, $\nu_E$ is the average for the Galois orbits of $\mu_E$. On the other hand, Kottwitz proves in [26, 2.8.1] that the average map has a factorization
$$X_*(T_E) \twoheadrightarrow X_*(T_E)_\Gamma \simeq B(T_E) \hookrightarrow (X_*(T_E) \otimes \mathbb{Q})^\Gamma$$
where $\Gamma = \mathrm{Gal}_{\mathbb{Q}_p}$, $X_*(T_E)^\Gamma$ and $X_*(T_E)_\Gamma$ are, respectively, the $\Gamma$-invariants and $\Gamma$-coinvariants of $X_*(T_E)$, and $B(T_E)$ is the set of $\sigma$-conjugacy classes in $T_E(K_0)$ and the last map is the Newton map (with values in $(X_*(T_E) \otimes \mathbb{Q})^\Gamma = \mathbf{C}^\mathbb{Q}(T_E)^\Gamma$). So it suffices to prove that the image of $\mu_E$ in $B(T_E)$ under the first map is $\sigma$-conjugated to $b_E$. In [26, 2.5], we see that the image of $\mu_E$ in $B(T_E)$ is given by the $\sigma$-conjugation class of $b'_E = N_{E/E_0}(\mu_E(\pi'_E))$ for $\pi'_E$ an uniformizer of $E$. Now, the image of $b'_E$ in $T_E(K_0) = \prod_{i=0}^{s-1} L_i^\times$, where $L_i = E \otimes_{E_0, \sigma^i} K_0$, is $(\pi'_E, 1, \ldots, 1)$ while the image of $b_E$ is $(\pi_K, 1, \ldots, 1)$.

Let $s = [E : \mathbb{Q}_p]$. It remains to prove that there exists $x = (x_0, \ldots, x_{s-1}) \in T_E(K_0)$ such that $(\sigma - 1)(x)(b'_E) = b_E$. Since we have
$$(\sigma - 1)(x_0, \ldots, x_{s-1}) = \left( \frac{\sigma x_{s-1}}{x_0}, \frac{\sigma x_0}{x_1}, \ldots, \frac{\sigma x_{s-2}}{x_{s-1}} \right) \quad \text{in} \quad T_E(K_0),$$
then it suffices to find $x_0 \in K_0$ such that $(\sigma^s - 1)x_0 = u$ where $u \in \mathcal{O}_{K_0}^\times$ such that $u\pi'_E = \pi_E$. Since $(\sigma^s - 1)$ is surjective on $\mathcal{O}_{K_0}^\times$, this element $x_0 \in K_0$ exist and we can take
$$x = (x_0, \sigma x_0, \ldots, \sigma^{s-1} x_0) \in T_E(K_0)$$
to obtain the $\sigma$-conjugation between $b'_E$ and $b_E$. $\square$

The next proposition gives us the relation between Wintenberger's functor $D_\pi$ and Fontaine's functor $D_{\mathrm{cris}} \circ V^u$.

**Proposition 5.5.** *For any finitely generated tensor subcategory $\mathcal{V}$ of $\mathrm{Rep}_{\mathbb{Q}_p} T$, there is an isomorphism of $\otimes$-functors $(D_{\mathrm{cris}} \circ V^u)_{|\mathcal{V}} \simeq D_{\pi | \mathcal{V}}$.*

*Proof.* This immediately follows from the main result of [43]: any such $\mathcal{V}$ is contained in $\mathrm{Rep}_{\mathbb{Q}_p} T_K$ for a sufficiently large finite extension $K \subset K_0^a$ of $K_0$. $\square$

As a consequence, we get the following result about the Fargues filtration of $V^u$.

**Proposition 5.6.** *The filtration $\mathcal{F}_F(V^u)$ is split by $\nu^\iota = \nu^{-1}$ and $\mathcal{F}_F(D_{\mathrm{cris}} \circ V^u) = \mathcal{F}_N^\iota(D_{\mathrm{cris}} \circ V^u)$.*

*Proof.* Write $V = V^u$. Fix $\tau \in \mathrm{Rep}_{\mathbb{Q}_p} T$. Since $T$ is commutative, we have a decomposition $\tau = \oplus_{\lambda \in \mathbb{Q}} \tau_\lambda$ according to $\nu^\iota$-weights. We first want to show $\mathcal{F}_F^\gamma(V(\tau)) = \oplus_{\lambda \geq \gamma} V(\tau_\lambda)$. We obviously have
$$\mathcal{F}_F^\gamma(V(\tau)) = \oplus_{\lambda \in \mathbb{Q}} \mathcal{F}_F^\gamma(V(\tau_\lambda)),$$
so we need to show $\mathcal{F}_F^\gamma(V(\tau_\lambda)) = V(\tau_\lambda)$ if $\lambda \geq \gamma$ and 0 otherwise, i.e. we want to show that $V(\tau_\lambda)$ is semi-stable of slope $\lambda$. This is equivalent to prove that the weakly admissible filtered isocrystal $D_{\mathrm{cris}}(V(\tau_\lambda))$ is semi-stable of slope $\lambda$. Now, we have $D_{\mathrm{cris}}(V(\tau_\lambda)) \simeq D_\pi(\tau_\lambda)$ which is Newton isoclinic of slope $-\lambda$, by Proposition 5.4thm.5.4. Since $\mathbf{t}_F(D_{\mathrm{cris}}(V(\tau_\lambda))) \leq \mathbf{t}_N^\iota(D_{\mathrm{cris}}(V(\tau_\lambda)))$, the equality must hold and $D_{\mathrm{cris}}(V(\tau_\lambda))$ is indeed semi-stable of slope $\lambda$.

We next want to show that $\mathcal{F}_F(D_{\mathrm{cris}} \circ V(\tau)) = \mathcal{F}_N^\iota(D_{\mathrm{cris}} \circ V(\tau))$ and we know that $\mathcal{F}_F \circ D_{\mathrm{cris}} = D_{\mathrm{cris}} \circ \mathcal{F}_F$, so it suffices to show that $\mathcal{F}_N^{\iota, \gamma}(D_{\mathrm{cris}} \circ V(\tau)) = \oplus_{\lambda \geq \gamma} D_{\mathrm{cris}} \circ V(\tau_\lambda)$. Let $\mathcal{V}$ be the $\otimes$-category generated by $\tau$ and $\tau_\lambda$ for each $\lambda$ appearing in the decomposition of $\tau$ (since there is a finite number of such $\tau_\lambda$'s, $\mathcal{V}$ is finitely generated). By last proposition, we can replace $D_{\mathrm{cris}} \circ V$ by $D_\pi$ and it suffices to show that $\mathcal{F}_N^{\iota, \gamma}(D_\pi(\tau)) = \oplus_{\lambda \geq \gamma} D_\pi(\tau_\lambda)$, which is obvious since $\mathcal{F}_N^\iota \circ D_\pi$ is split by $\nu^\iota$. $\square$

## 5.4 Germs of abelian crystalline representations with $G$-structure

We want to add a $G$-structure to the germs of abelian crystalline representation, as we did in the sections before for general germs of crystalline representations. Let $G$ be a reductive group over $\mathbb{Z}_p$. A germ of abelian crystalline representation with $G$-structure is an exact and faithful $\otimes$-functor

$$V \;:\; \operatorname{Rep}_{\mathbb{Q}_p} G \to \operatorname{Rep}_{\mathbb{Q}_p}^{\mathrm{cr,ab}}\{\mathrm{Gal}_{K_0}\}.$$

We say that $V$ is trivial when $\omega(V) = \omega_{G,\mathbb{Q}_p}$.

**Proposition 5.7.** *The following data are equivalent:*

1. *A morphism $x : T \to G$,*

2. *A trivial germ of abelian crystalline representation with $G$-structure $V_x$,*

3. *A cocharacter $\mu_x \;:\; \mathbb{G}_{m,\overline{\mathbb{Q}}_p} \to G_{\overline{\mathbb{Q}}_p}$ such that the Mumford-Tate group $\mathrm{MT}(\mu_x)$ is a torus.*

*Moreover, $K$ is a field definition of $V_x$ if and only if $x$ factors tough $T_K$, if and only if $\mu_x$ is defined over $K \cap \overline{\mathbb{Q}}_p$.*

*Proof.* $(1) \Leftrightarrow (2)$: From a morphism $x : T \to G$ we get

$$V_x = V^u \circ x \;:\; \operatorname{Rep}_{\mathbb{Q}_p} G \to \operatorname{Rep}_{\mathbb{Q}_p} T \simeq \operatorname{Rep}_{\mathbb{Q}_p}^{\mathrm{cr,ab}}\{\mathrm{Gal}_{K_0}\}.$$

Given $V_x$, we obtain a map $\operatorname{Aut}^{\otimes}(\omega(V_x)) \to G$. On the other hand, Fontaine gives an isomorphism $\operatorname{Aut}^{\otimes}(\omega(V_x)) \to T$. Composing the inverse of the latter with the first map, we obtain $x : T \to G$.

$(1) \Leftrightarrow (3)$: Now, given $x$, we obtain $\mu_x = x \circ \mu$. Then, $\mathrm{MT}(\mu_x) \subset T$ since $\mathrm{MT}(\mu_x)$ is defined as the smallest algebraic subgroup of $G$, defined over $\mathbb{Q}_p$ such that $\mu_x$ factors through it. Moreover, $\mathbb{G}_{m,\overline{\mathbb{Q}}_p}$ is a connected group, so we have a factorization

$$\mathbb{G}_{m,\overline{\mathbb{Q}}_p} \xrightarrow{\mu_x} \mathrm{MT}(\mu_x)^{\circ}_{\overline{\mathbb{Q}}_p} \subset \mathrm{MT}(\mu_x)_{\overline{\mathbb{Q}}_p}.$$

Then $\mathrm{MT}(\mu_x)^{\circ}_{\overline{\mathbb{Q}}_p} = \mathrm{MT}(\mu_x)_{\overline{\mathbb{Q}}_p}$ and $\mathrm{MT}(\mu_x)$ is connected. Over $\mathbb{Q}_p$, every connected subgroup of a torus is a torus, so $\mathrm{MT}(\mu_x)$ is a torus. Conversely, starting with $\mu_x$ defined over $E$, pick a finite extension $E \subset \overline{\mathbb{Q}}_p$ of $\mathbb{Q}_p$ such that $\mu_x$ is defined over $E$. Then $\mu_x$ yields an element of

$$\begin{aligned}
\operatorname{Hom}_{\mathbb{Z}[\mathrm{Gal}_E]}(X^*(\mathrm{MT}(\mu_x)), \mathbb{Z}) &= \operatorname{Hom}_{\mathbb{Z}[\mathrm{Gal}_{\mathbb{Q}_p}]}(X^*(\mathrm{MT}(\mu_x)), \mathcal{C}(\operatorname{Hom}(E, \overline{\mathbb{Q}}_p), \mathbb{Z})) \\
&= \operatorname{Hom}(T_E, \mathrm{MT}(\mu_x))
\end{aligned}$$

hence $x : T \to T_E \to \mathrm{MT}(\mu_x) \subset G$. In the display, the first map is $f \mapsto \tilde{f}$, $\tilde{f}(a)(\iota) = f(\sigma a)$ if $\sigma \iota = \iota_0$.

The statement about the field of definition is obvious by construction of the equivalences. □

The next proposition gives us the Hodge and Newton types for a trivial germ of abelian crystalline representations with $G$-structure.

**Proposition 5.8.** *Let $V_x$ be a trivial germ of abelian crystralline representation with $G$-structure, with associated morphism $x : T \to G$. Then, we obtain the Hodge, Newton and Fargues types by*

$$\begin{aligned}
\mathbf{t}_{\mathrm{H}}(V_x) &= [x \circ \mu] &\text{in}&\quad \mathbf{C}^{\mathbb{Z}}(G_{\overline{\mathbb{Q}}_p}) \\
\mathbf{t}_{\mathrm{N}}(V_x) &= [x \circ \nu] &\text{in}&\quad \mathbf{C}^{\mathbb{Q}}(G_{\mathbb{Q}_p}). \\
\mathbf{t}_{\mathrm{F}}(V_x) &= [x \circ \nu^{\iota}] &\text{in}&\quad \mathbf{C}^{\mathbb{Q}}(G_{\mathbb{Q}_p})
\end{aligned}$$

*Therefore, $V_x$ is ordinary when $[x \circ \nu] = [x \circ \mu]^{\#}$.*

*Proof.* As an abuse of notation, let $x : \operatorname{Rep}_{\mathbb{Q}_p} G \to \operatorname{Rep}_{\mathbb{Q}_p} T$ and let $\mathcal{V}_x$ be the essential image of this map. It is a finitely generated $\otimes$-subcategory of $\operatorname{Rep}_{\mathbb{Q}_p} T$, thus $D_{\mathrm{cris}} \circ V_x = D_{\mathrm{cris}} \circ V^u \circ x \simeq D_{\pi} \circ x$, by Proposition 5.5thm.5.5. The functor $D_{\pi} \circ x$ is trivial, and we know that $\mathcal{F}_{\mathrm{H}}$ is split by $x \circ \mu$ and $\mathcal{G}_{\mathrm{N}}$ is given by $x \circ \nu$, by Proposition 5.4thm.5.4 and that $\mathcal{F}_{\mathrm{F}}$ is split by $x \circ \nu^{\iota}$ by Proposition 5.6thm.5.6. □

Since $G$ is defined over $\mathbb{Z}_p$, it is quasi-split over $\mathbb{Z}_p$ and $\mathbb{Q}_p$: there are Borel pairs $(T_0, B)$ in $G$ and they are all conjugated. In particular, we obtain a conjugation class of distinguished maximal tori in $G$, formed by the ones contained in a Borel pair. Let $(T_0, B)$ be a Borel pair in $G$. A cocharacter $\mu \in X_*(T_0)$ is called $B$-dominant if for every positive root $\alpha$ of $T_0$ in $B$ we have $\langle \alpha, \mu \rangle \geq 0$, where $\langle \cdot, \cdot \rangle$ is the pairing between $X_*(T_0)$ and $X^*(T_0)$. We denote by $X_*(T_0)^{\mathrm{dom}}$ the set of $B$-dominant cocharacters of $T_0$.

**Proposition 5.9.** *If the image of $x$ is contained in a distinguished maximal torus, then $V_x$ is ordinary. In this case, the definition field is an unramified extension of $\mathbb{Q}_p$ and we have $x(b) = \mu_x(p)$.*

*Proof.* By assumption, there is a Borel pair $(T_0, B) \subset G$ defined over $\mathbb{Q}_p$ such that $\mu_x = x \circ \mu$ factors through a $B$-dominant cocharacter of $T_0$, i.e. $\mu_x \in X_*(T_0)^{\mathrm{dom}}$. The morphism $T_0 \hookrightarrow G$ induces a $\mathrm{Gal}_{\mathbb{Q}_p}$-equivariant map
$$X_*(T_0) \otimes \mathbb{Q} = \mathbf{C}^{\mathbb{Q}}(T_{0, \overline{\mathbb{Q}}_p}) \to \mathbf{C}^{\mathbb{Q}}(G_{\overline{\mathbb{Q}}_p})$$
whose restriction to the $B$-dominant cone is also compatible with the monoid structures: it is a $\mathrm{Gal}_{\mathbb{Q}_p}$-equivariant isomorphism of monoids $(X_*(T_0) \otimes \mathbb{Q})^{\mathrm{dom}} \to \mathbf{C}^{\mathbb{Q}}(G_{\overline{\mathbb{Q}}_p})$. Plainly, $(x \circ \mu)^\# = x \circ \nu$ on the left hand side, thus also $[x \circ \mu]^\# = [x \circ \nu]$ on the right hand side, i.e. $\mathbf{t}_{\mathrm{H}}(V_x) = \mathbf{t}_{\mathrm{N}}(V_x)$ by proposition 5.8thm.5.8 and $V_x$ is ordinary.

In this case, $\mu_x$ is a cocharacter of $T_0$, which is split over $K_0$, so $\mu_x$ is defined over an unramified extension $E$ of $\mathbb{Q}_p$. Thus $x$ factors through $T_{K_0}$ by construction and
$$x(b) = x_{K_0}(b_{K_0}) = x_{K_0} \circ \mu_{K_0}(\pi_{K_0}) = x_{K_0} \circ \mu_{K_0}(p) = (x \circ \mu)(p) = \mu_x(p),$$
where $x_{K_0} : T_{K_0} \to T_0 \hookrightarrow G$ is the factorization of $x$. $\square$

## 5.5 The reduction map

For the rest of the section, we fix a trivial germ of abelian crystalline representations $V_x$ with $G$-structure such that $\mathrm{MT}(\mu_x)$ is contained in $T_0$, for $(T_0, B)$ a Borel pair of $G_{\mathbb{Z}_p}$ with $\mu_x$ dominant. We set
$$\nu_x = x \circ \nu \; : \; \mathbb{D}(\mathbb{Q})_{|\mathbb{Q}_p} \to T_{|\mathbb{Q}_p} \to G_{|\mathbb{Q}_p}$$
and $M_x = Z_G(\nu_x)$. The cocharacter $\nu_x$ is defined over $\mathbb{Z}_p$, by [8, Proposition 3 and Lemma 4], so $M_x$ is also defined over $\mathbb{Z}_p$. By Proposition 5.8thm.5.8, $M_x$ is a Levi subgroup of the parabolic subgroup $P_x = P_{\mathcal{F}}$ stabilizing $\mathcal{F}_{\mathrm{F}}(V)$. Note also that $B \subset P_{\mathcal{F}}$ since $\nu_x$ is $B$-dominant, since $\mu_x$ and all its conjugates are.

As we have seen in last section, the map red has a factorization
$$\begin{array}{ccc} \mathcal{L}(V_x) & \xrightarrow{\mathrm{red}} & \mathcal{L}(D_x, \mathbf{t}_{\mathrm{H}}^\iota(D_x)) \\ & \searrow \quad \nearrow & \\ & U_{\mathcal{F}_F}(\mathbb{Q}_p) \backslash \mathcal{L}(V_x) & \end{array}$$
where $D_x = D_{\mathrm{cris}}(V_x)$, $\mathcal{F}_F$ is the Fargues filtration on the germ of $G$-crystalline representation $V_x$, and the map red is compatible with the morphism $\mathrm{Aut}^\otimes(V_x) \to \mathrm{Aut}^\otimes(D_x)$. The aim of the rest of the section is to give a more explicit description of the diagram above.

### 5.5.1 The source

Since $V_x$ is trivial, we have a canonical isomorphism
$$\mathcal{L}(V_x) \simeq G(\mathbb{Q}_p)/G(\mathbb{Z}_p),$$
since $\mathcal{L}(V_x) = \mathcal{L}(\omega_{G, \mathbb{Q}_p})$ and $\mathrm{Aut}^\otimes(\omega_{G, \mathbb{Q}_p}) = G(\mathbb{Q}_p)$ acts transitively on $\mathcal{L}(\omega_{G, \mathbb{Q}_p})$, with stabilizer $\mathrm{Aut}^\otimes(\omega_{G, \mathbb{Z}_p}) = G(\mathbb{Z}_p)$, for $\omega_{G, \mathbb{Z}_p} \in \mathcal{L}(V_x)$ the trivial lattice. Next proposition describes the automorphism group of $V_x$ as a subgroup of $\mathrm{Aut}^\otimes(\omega(V_x)) = G(\mathbb{Q}_p)$.

**Proposition 5.10.** *For $V_x$ as above, we have $\mathrm{Aut}^\otimes(V_x) = M_x(\mathbb{Q}_p)$.*

*Proof.* For $g \in G(\mathbb{Q}_p)$, we denote by $\mathrm{Int}(g)$ the inner automorphism defined by $g$. Then, an element $g \in \mathrm{Aut}^\otimes(\omega(V_x))$ defines an isomorphism $V_x \to V_{\mathrm{Int}(g) \cdot x}$, since

$$\begin{aligned}
\tau(g)(\sigma \cdot_x v) &= \tau(g)(\tau(x \circ \chi_K(\sigma)) \cdot v) \\
&= \tau(g \cdot x \circ \chi_K(\sigma))(v) \\
&= \tau(g \cdot x \circ \chi_K(\sigma) \cdot g^{-1} \cdot g)(v) \\
&= \tau(\mathrm{Int}(g)(x) \circ \chi_K(\sigma) \cdot g)(v) \\
&= \tau(\mathrm{Int}(g)(x) \circ \chi_K(\sigma))(\tau(g) \cdot v) \\
&= \sigma \cdot_{\mathrm{Int}(g)(x)} \tau(g)v
\end{aligned}$$

for $\tau \in \mathrm{Rep}_{\mathbb{Q}_p} G$, $v \ni \omega_G(\tau)$, $K$ a field of definition of $V_x$ and $\sigma \in \mathrm{Gal}_K$ (we may take $K = K_0$). By definition, $\mathrm{Aut}^\otimes(V_x)$ is the set of $g \in \mathrm{Aut}^\otimes(\omega(V_x))$ inducing an isomorphism of $V_x$, i.e. the set of $g \in G(\mathbb{Q}_p)$ such that $V_x = V_{\mathrm{Int}(g)(x)}$, and by the classification of trivial germs of abelian crystalline representations with $G$-structure given in 5.7thm.5.7, this is equivalent to the set of $g \in G(\mathbb{Q}_p)$ such that $x = \mathrm{Int}(g)(x)$ or, with respect to $\mu_x$, the set of $g \in G(\mathbb{Q}_p)$ such that $\mu_x = \mathrm{Int}(g)(\mu_x)$, i.e.

$$\mathrm{Aut}^\otimes(V_x) = \{ g \in G(\mathbb{Q}_p) \mid g \cdot x \cdot g^{-1} = x \} = \{ g \in G(\mathbb{Q}_p) \mid g \cdot \mu_x \cdot g^{-1} = \mu_x \}.$$

Now, the faithful $\otimes$-functor $D_\pi$ is compatible with fiber functors and induces a commutative diagram

$$\begin{array}{ccc}
\mathrm{Aut}^\otimes(V_x) & \longrightarrow & \mathrm{Aut}^\otimes(D_\pi(V_x)) \\
\downarrow & & \downarrow \\
\mathrm{Aut}^\otimes(\omega(V_x)) & \longrightarrow & \mathrm{Aut}^\otimes(\omega(D_\pi(V_x))) \\
\| & & \| \\
G(\mathbb{Q}_p) & \longrightarrow & G(K_0).
\end{array}$$

If $g \in \mathrm{Aut}^\otimes(V_x)$, $g \in G(K_0)$ preserves in particular the underlying isocrystal associated to $D_\pi(V_x)$, so it also preserves the Newton graduation, which is given by $\nu_x$, thus $g$ normalizes $\nu_x$, i.e.

$$g \in M_x(K_0) \cap G(\mathbb{Q}_p) = M_x(\mathbb{Q}_p).$$

On the other hand, we have a factorization (all groups over $\overline{\mathbb{Q}}_p$)

$$\begin{array}{ccccc}
\mathbb{G}_m & \xrightarrow{\mu_x} & T_0 & \longrightarrow & G \\
& \searrow & \downarrow & \nearrow & \\
& & Z_G(M_x) \hookrightarrow M_x & &
\end{array}$$

Indeed, $\nu_x$ factors through $T_0$, which is commutative, so $T_0 \subset M_x$, thus $\mu_x$ factors through $M_x$ since it factors through $T_0 \subset M_x$. Let $\mathcal{M} = \mathrm{Lie}(M_x)$. Let $\Phi(T_0, M_x)$ be the system of roots of $T_0$ in $M_x$, i.e. characters such that $\mathcal{M} = \mathcal{M}^0 \oplus \oplus_{\alpha \in \Phi(T_0, M_x)} \mathcal{M}^\alpha$, where $\mathcal{M}^0 = \mathrm{Lie}\, T_0$. Let $\Phi^+(T_0, M_x) = \Phi(T_0, M_x) \cap \Phi^+(T_0, G)$ be the roots of $T_0$ in $\mathcal{M}$ which are positive with respect to $B \cap M_x$, a Borel subgroup of $M_x$ with Levi $T_0$. For every $\alpha \in \Phi^+(T_0, M_x)$ and $\mu'$ in the $\mathrm{Gal}_{\mathbb{Q}_p}$-orbit $O(\mu_x)$ of $\mu_x$, we have $\langle \alpha, \mu' \rangle \geq 0$, but $\sum_{O(\mu_x)} \langle \alpha, \mu' \rangle = \langle \alpha, \sum_{O(\mu_x)} \mu' \rangle = \#O(\mu_x)\langle \alpha, \nu_x \rangle = 0$ since $\nu_x$ is central in $M_x$, thus $\langle \alpha, \mu' \rangle = 0$ for every $\mu' \in O(\mu_x)$. In particular $\langle \alpha, \mu_x \rangle = 0$ for every $\alpha \in \Phi^+(T_0, M_x)$, which implies that $\mu_x$ is central in $M_x$. Thus $M_x(\mathbb{Q}_p) \subset \mathrm{Aut}^\otimes(V_x)$ which proves that $M_x(\mathbb{Q}_p) = \mathrm{Aut}^\otimes(V_x)$. $\square$

The last result can be reformulated as follows: since $x$ factors through $T_0$, it factors through $M_x$ and we have $V_x : \mathrm{Rep}_{\mathbb{Q}_p} G \to \mathrm{Rep}_{\mathbb{Q}_p} M_x \xrightarrow{V_x^{M_x}} \mathrm{Rep}_{\mathbb{Q}_p}^{\mathrm{cr,ab}} \{\mathrm{Gal}_{K_0}\}$, which induces a diagram

$$\begin{array}{ccc}
\mathrm{Aut}^\otimes(V_x^{M_x}) & \xrightarrow{\simeq} & \mathrm{Aut}^\otimes(V_x) \\
\| & & \downarrow \\
\mathrm{Aut}^\otimes(\omega(V_x^{M_x})) & \hookrightarrow & \mathrm{Aut}^\otimes(\omega(V_x)) \\
\| & & \| \\
M_x(\mathbb{Q}_p) & \longrightarrow & G(\mathbb{Q}_p).
\end{array}$$

Note that $V_x^{M_x}$ is also ordinary by Proposition 5.9thm.5.9, since the corresponding morphism $x : T \to M_x$ factors through the distinguished torus $T_0$ of $M_x$. Moreover, we have an $M_x(\mathbb{Q}_p)$-equivariant diagram

$$\begin{array}{ccccc}
\mathcal{L}(V_x^{M_x}) & \longrightarrow & \mathcal{L}(V_x) & \longrightarrow & U_{\mathcal{F}_{\mathrm{F}}}(\mathbb{Q}_p)\backslash\mathcal{L}(V_x) \\
\| & & \| & & \| \\
M_x(\mathbb{Q}_p)/M_x(\mathbb{Z}_p) & \to & G(\mathbb{Q}_p)/G(\mathbb{Z}_p) & \to & U_{\mathcal{F}_{\mathrm{F}}}(\mathbb{Q}_p)\backslash G(\mathbb{Q}_p)/G(\mathbb{Z}_p).
\end{array}$$

The next proposition gives us a description of the quotient:

**Proposition 5.11.** *The top map $\mathcal{L}(V_x^{M_x}) \to U_{\mathcal{F}_{\mathrm{F}}}(\mathbb{Q}_p)\backslash\mathcal{L}(V_x)$ bijective.*

*Proof.* It is well known. The Iwasawa decomposition $G(\mathbb{Q}_p) = B(\mathbb{Q}_p)G(\mathbb{Z}_p)$ gives us

$$G(\mathbb{Q}_p) = P_{\mathcal{F}_{\mathrm{F}}}(\mathbb{Q}_p)G(\mathbb{Z}_p) = U_{\mathcal{F}_{\mathrm{F}}}(\mathbb{Q}_p)M_x(\mathbb{Q}_p)G(\mathbb{Z}_p)$$

as $B(\mathbb{Q}_p) \subset P_{\mathcal{F}_{\mathrm{F}}}(\mathbb{Q}_p) = U_{\mathcal{F}_{\mathrm{F}}}(\mathbb{Q}_p) \rtimes M_x(\mathbb{Q}_p)$. Thus we have the surjectivity of the bottom map. The injectivity is given by [8, 6.3.3]. □

### 5.5.2 The target

It remains to study the description of both the set $\mathcal{L}(D_x, \mathbf{t}_{\mathrm{H}}^{\iota}(D_x))$ and the action of the automorphism group $\mathrm{Aut}^{\otimes}(D_x)$. We may assume that $D_x$ is isomorphic to a trivial $G$-isocrystal, since we can fix an isomorphism $D_x = D_{\mathrm{cris}} \circ V_x \simeq D_\pi \circ x$ and use Wintenberger's results in [43]. Applying Fontaine's functor $D_{\mathrm{cris}}$ to

$$V_x \;:\; \mathrm{Rep}_{\mathbb{Q}_p} G \to \mathrm{Rep}_{\mathbb{Q}_p} M_x \xrightarrow{V_x^{M_x}} \mathrm{Rep}_{\mathbb{Q}_p}^{\mathrm{cr,ab}}\{\mathrm{Gal}_{K_0}\}$$

we obtain

$$D_x \;:\; \mathrm{Rep}_{\mathbb{Q}_p} G \to \mathrm{Rep}_{\mathbb{Q}_p} M_x \xrightarrow{D_x^{M_x}} \mathrm{Mod}_{K_0}^\sigma .$$

**Proposition 5.12.** *The functor $D_{\mathrm{cris}}$ induces a bijection $\mathrm{Aut}^{\otimes}(V_x) \simeq \mathrm{Aut}^{\otimes}(D_x)$.*

*Proof.* We may replace $D_x$ by $D_\pi \circ x$. As we have seen, the $\otimes$-automorphisms of $D_\pi \circ x$ preserve the Newton cocharacter, so their image is in $M_x(K_0)$ and we get

$$\mathrm{Aut}^{\otimes}(D_\pi \circ x) = \{m \in M_x(K_0) \mid mx(b) = x(b)\sigma(m)\}.$$

Now, by ordinarity, we know by Proposition 5.9thm.5.9 that $x(b) = x(\mu(p)) = \mu_x(p)$, which is central in $M_x(K_0)$, thus

$$\mathrm{Aut}^{\otimes}(D_\pi(V_x)) = \{m \in M_x(K_0) \mid m = \sigma(m)\} = M_x(\mathbb{Q}_p).$$

Our claim now follows from Proposition 5.10thm.5.10. □

The previous result and Proposition 5.10thm.5.10 together imply that $D_{\mathrm{cris}}$ and the embedding $M_x \hookrightarrow G$ yield group isomorphisms

$$\begin{array}{ccc}
\mathrm{Aut}^{\otimes}(V_x^{M_x}) & \xrightarrow{\simeq} & \mathrm{Aut}^{\otimes}(D_x^{M_x}) \\
\downarrow \simeq & & \downarrow \simeq \\
\mathrm{Aut}^{\otimes}(V_x) & \xrightarrow{\simeq} & \mathrm{Aut}^{\otimes}(D_x).
\end{array}$$

A similar result holds for lattices. We will use the script $M_x$ or $G$ to indicate if we are working with the vectorial distance defined on $\mathbf{B}^e(M_{xK})$ or $\mathbf{B}^e(G_K)$, respectively.

**Proposition 5.13.** *The functor $D_{\mathrm{cris}}$ and the embedding $M_x \hookrightarrow G$ yield a commutative diagram*

$$\begin{array}{ccc}
\mathcal{L}(V_x^{M_x}) & \xrightarrow{\simeq} & \mathcal{L}\left(D_x^{M_x}, \mathbf{t}_{\mathrm{H}}^{\iota, M_x}(D_x^{M_x})\right) \\
\cap & & \downarrow \simeq \\
\mathcal{L}(V_x) & \xrightarrow{\mathrm{red}} & \mathcal{L}\left(D_x, \mathbf{t}_{\mathrm{H}}^{\iota, G}(D_x)\right).
\end{array}$$

*Proof.* To compute the right hand side map, we may again replace $D_x$ by $D_\pi \circ x$, that we call $D$ to easy the notations. Also, we set $V = V_x$, $M = M_x$, $\mathbf{t}_M = \mathbf{t}_{\mathrm{H}}^{\iota, M_x}(D_x^{M_x})$ in $\mathbf{C}^{\mathbb{Z}}(M)$ and $\mathbf{t}_G = \mathbf{t}_{\mathrm{H}}^{\iota, G}(D_x)$ in $\mathbf{C}^{\mathbb{Z}}(G)$. We have to show that the embedding $\mathcal{L}(V^M) \subset \mathcal{L}(V)$ yields a bijection $\mathcal{L}(D^M, \mathbf{t}_M) \simeq \mathcal{L}(D, \mathbf{t}_G)$. This is well-known result, given by Kottwitz in [29, Theorem 1.1]. We propose here an alternative proof. We have a diagram

$$\begin{array}{ccc}
\mathcal{L}(D^M, \mathbf{t}_M) & \hookrightarrow & \mathcal{L}(D, \mathbf{t}_G) \\
\downarrow & & \downarrow \\
\mathbf{B}(\omega_M, K_0) & \hookrightarrow & \mathbf{B}(\omega_G, K_0) \\
\simeq \uparrow & & \simeq \uparrow \\
\mathbf{B}^e(M_{K_0}) & \hookrightarrow & \mathbf{B}^e(G_{K_0})
\end{array}$$

which is cartesian, i.e. $\mathcal{L}(D^M, \mathbf{t}_M) = \mathcal{L}(D, \mathbf{t}_G) \cap \mathbf{B}^e(G_{K_0})$. From [10, Theorem 7], we thus already obtain $\mathcal{L}(D, \mathbf{t}_G) \subset \mathcal{L}(D^M)$, therefore

$$\mathcal{L}(D, \mathbf{t}_G) = \coprod\nolimits_{\mathbf{t}'_M \mapsto \mathbf{t}_G} \mathcal{L}(D^M, \mathbf{t}'_M)$$

where $\mathbf{t}'_M$ runs through the fiber $\mathbf{C}^{\mathbb{Z}}(M) \twoheadrightarrow \mathbf{C}^{\mathbb{Z}}(G)$ over $\mathbf{t}_G$. So, it suffices to show that $\mathcal{L}(D^M, \mathbf{t}'_M) = \emptyset$ whenever $\mathbf{t}'_M \neq \mathbf{t}_M$. Note that $\mathcal{L}(D^M, \mathbf{t}_M)$ is indeed non empty, since the $M$-isocrystal $D^M$ is $\mathbf{t}_M$-ordinary by 4.26thm.4.26, it is therefore sufficient to establish that $\mathbf{t}_M$ is the only point in the intersection of the relevant fibers of the Kottwitz map $[-]_{M,\Gamma} : \mathbf{C}^{\mathbb{Z}}(M) \to \pi_1(M)_\Gamma$ and the projection $\mathbf{C}^{\mathbb{Z}}(M) \to \mathbf{C}^{\mathbb{Z}}(G)$.

Let $W_G(T_0) = N_G(T_0)/Z_G(T_0)$ be the Weyl group of $T_0$ in $G$, denote by $R_G^\vee = \Phi^\vee(T_0, G)$ the system of coroots of $T_0$ in $G$, $\Delta_G$ the simple coroots of $T_0$ in $G$ and $X_*(T_0)^{G-\mathrm{dom}}$ the cone of dominant cocharacters defined by $B$. We use analogous notation for the corresponding objects defined with respect to $M$ and $M \cap B$. Finally, denote by $X_\Gamma$ the coinvariants of $\Gamma = \mathrm{Gal}_{\mathbb{Q}_p}$ acting on a group $X$. We then have the following commutative diagram, with exact second and third rows (compare with [29]):

$$\begin{array}{ccccccccc}
& & \mathbf{C}^{\mathbb{Z}}(M) \simeq W_{M_x}(T_0) \backslash X_*(T_0) & \to & W_G(T_0) \backslash X_*(T_0) \simeq \mathbf{C}^{\mathbb{Z}}(G) & & \\
& & \downarrow & & \downarrow & & \\
0 & \longrightarrow & \bigoplus_{\delta \in \Delta_G \backslash \Delta_M} \mathbb{Z} \cdot \delta & \longrightarrow & X_*(T_0)/\mathbb{Z} \cdot R_M^\vee & \longrightarrow & X_*(T_0)/\mathbb{Z} \cdot R_G^\vee & \longrightarrow & 0 \\
& & \downarrow & & \downarrow & & \downarrow & & \\
0 & \to & \bigoplus_{\Gamma \cdot \delta \in \Gamma \backslash (\Delta_G \backslash \Delta_M)} \mathbb{Z} \cdot \Gamma \cdot \delta & \longrightarrow & (X_*(T_0)/\mathbb{Z} \cdot R_M^\vee)_\Gamma & \longrightarrow & (X_*(T_0)/\mathbb{Z} \cdot R_G^\vee)_\Gamma & \longrightarrow & 0
\end{array}$$

By construction $\mathbf{t}_M = \mu_x$ in $X_*(T_0)^{G-\mathrm{dom}} \subset X_*(T_0)^{M-\mathrm{dom}}$. Let $\mathbf{t}'_M \in X_*(T_0)^{M-\mathrm{dom}}$ be any other element in the fiber of $\mathbf{C}^{\mathbb{Z}}(M) \twoheadrightarrow \mathbf{C}^{\mathbb{Z}}(G)$ over $\mathbf{t}_G$. Then, $\mathbf{t}'_{M_x} = w \cdot \mathbf{t}_{M_x}$ for some $w \in W_G(T_0)$, thus $\mathbf{t}_M - \mathbf{t}'_M \in X_*(T_0)$ is a linear combination of elements of $\Delta_G$ with non-negative coefficients. In the second line of our diagram, we thus obtain

$$\mathbf{t}_{M_x} - \mathbf{t}'_{M_x} = \sum\nolimits_{\delta \in \Delta_G \backslash \Delta_{M_x}} n_\delta \cdot \delta \in \bigoplus_{\delta \in \Delta_G \backslash \Delta_{M_x}} \mathbb{Z} \cdot \delta, \quad \text{with } n_\delta \geq 0.$$

The map

$$\bigoplus_{\delta \in \Delta_G \backslash \Delta_{M_x}} \mathbb{Z} \cdot \delta \to \bigoplus_{\Gamma \cdot \delta \in \Gamma \backslash (\Delta_G \backslash \Delta_{M_x})} \mathbb{Z} \cdot \Gamma \cdot \delta$$

is given by

$$\sum_{\delta \in \Delta_G \backslash \Delta_{M_x}} n_\delta \cdot \delta \mapsto \sum_{\delta \in \Delta_G \backslash \Delta_{M_x}} \left( \sum_{\delta' \in \Gamma \cdot \delta} n_{\delta'} \right) \cdot \delta.$$

Thus if also $[\mathbf{t}_M]_{M,\Gamma} = [\mathbf{t}'_M]_{M,\Gamma}$ in $(X_*(T_0)/\mathbb{Z} \cdot R_M^\vee)_\Gamma$ (which is the case if $\mathcal{L}(D^M, \mathbf{t}'_M) \neq \emptyset$), then actually $\mathbf{t}'_M = \mathbf{t}_M$ as was to be shown, since all the $n_\delta$ must be 0.

We still have to show that $\mathcal{L}(V^M) \hookrightarrow \mathcal{L}(D^M)$ induces a bijection $\mathcal{L}(V^M) \simeq \mathcal{L}(D^M, \mathbf{t}_M)$. Inside $\mathcal{L}(D^M) \simeq M(K_0)/M(W(\mathbb{F}))$, we have

$$\mathcal{L}(D^M, \mathbf{t}_M) = \{m \in M(K_0)/M(W(\mathbb{F})) \mid m^{-1}\mu_x(p)\sigma(m) \in G(W(\mathbb{F}))\mu_x(p)G(W(\mathbb{F}))\}$$

which equals

$$\{m \in M(K_0)/M(W(\mathbb{F})) \mid m^{-1}\sigma(m) \in G(W(\mathbb{F}))\}$$

since $\mu_x(p)$ is central in $M$. In other words,

$$\mathcal{L}(D^M, \mathbf{t}_M) = \{y \in \mathcal{L}(D^M) \mid \sigma y = y\} = \mathcal{L}(V^x)$$

using, for instance, [10, 2.5.5] for the last equality. $\square$

Let $s$ be the order of the $\Gamma$-orbit of $\mu_x$ in $X_*(T_0)$. Then $s\nu_x \in X_*(T_0)$ and $\sigma^s(\mu_x) = \mu_x$, thus

$$\sigma^s \mathbf{t}_H^\iota(D_x) = \mathbf{t}_H^\iota(D_x) \quad \text{in} \quad \mathbf{C}^{\mathbb{Z}}(G).$$

It follows that the $s$-power of the Frobenius of $D_x$ acts on $\mathcal{L}(D_x, \mathbf{t}_H^\iota(D_x))$, giving an automorphism

$$F \; : \; \mathcal{L}(D_x, \mathbf{t}_H^\iota(D_x)) \to \mathcal{L}(D_x, \mathbf{t}_H^\iota(D_x)).$$

It is related to the operator $\Phi_{\mathrm{cris}}^s$ of section 7.3.5 as follows:

**Proposition 5.14.** *We have $F = \Phi_{\mathrm{cris}}^s$ on $\mathcal{L}(D_x, \mathbf{t}_H^\iota(D_x))$.*

*Proof.* We may as usual replace $D_x$ by $D_\pi \circ x$ and call it $D$ to ease the notations. Let also $V = V_x$ and $M = M_x$. Then $F$ acts on
$$\mathcal{L}(D) \simeq \mathcal{L}(\omega_{G,K_0}) \simeq G(K_0)/G(\mathbb{Z}_p)$$
by the element $(\mu_x(p), \sigma)^s = (s\nu_x(p), \sigma^s) \in G(K_0) \rtimes \langle \sigma \rangle^{\mathbb{Z}}$. It thus acts on
$$\mathcal{L}(D, \mathbf{t}_H^\iota(D)) \simeq \mathcal{L}(V^M) \simeq M(\mathbb{Q}_p)/M(\mathbb{Z}_p)$$
by $s\nu_x(p) \in M(\mathbb{Q}_p)$. For any $y \in \mathcal{L}(D, \mathbf{t}_H^\iota(D))$ and $\tau \in \mathrm{Rep}_{\mathbb{Z}_p} G$ with weight decomposition $\tau|M = \oplus_{a \in \mathbb{Z}} \tau_a$ with respect to the central cocharacter $s\nu_x : \mathbb{G}_{m,\mathbb{Z}_p} \to M$, we thus have $F(y)(\tau) = \oplus_{a \in \mathbb{Z}} p^a \cdot y(\tau_a)$. On the other hand, since $\nu_x^{-1}$ splits $\mathcal{F}_N^\iota(D)$,

$$\begin{array}{rcl}
\Phi_{\mathrm{cris}}^s(y)(\tau) & = & \sum_{i \in \mathbb{Z}} p^{-i} y(\tau) \cap (s\mathcal{F}_N^\iota(D))^i(\tau) \\
& = & \sum_{i \in \mathbb{Z}} \left( \oplus_{a \in \mathbb{Z}} p^{-i} y(\tau_a) \right) \cap \left( \oplus_{a \leq -i} y(\tau_a) \otimes K_0 \right) \\
& = & \oplus_{a \in \mathbb{Z}} p^a \cdot y(\tau_a).
\end{array}$$

Thus indeed $F = \Phi_{\mathrm{cris}}^s$ on $\mathcal{L}(D_x, \mathbf{t}_H^\iota(D_x))$. $\square$

*Remark* 18. This also proves that $\Phi_{\mathrm{cris}}^s$ is bijective on $\mathcal{L}(D_x, \mathbf{t}_H^\iota(D_x))$, as we claimed in the proof of 4.32thm.4.32.

### 5.5.3 The main result

Putting all the results in this section together, we have obtained a commutative diagram

$$\begin{array}{ccccc}
M_x(\mathbb{Q}_p)/M_x(\mathbb{Z}_p) & \xrightarrow{\simeq} & \mathcal{L}(V_x^{M_x}) & \xrightarrow{\simeq} & \mathcal{L}(D_x^{M_x}, \mathbf{t}_H^\iota(D_x^{M_x})) \\
\cap & & \cap & & \downarrow \simeq \\
G(\mathbb{Q}_p)/G(\mathbb{Z}_p) & \xrightarrow{\simeq} & \mathcal{L}(V_x) & \xrightarrow{\mathrm{red}} & \mathcal{L}(D_x, \mathbf{t}_H^\iota(D_x))
\end{array}$$

which is equivariant with respect to the actions of

$$\begin{array}{ccccccc}
M_x(\mathbb{Q}_p) & \xrightarrow{\simeq} & \mathrm{Aut}^\otimes(\omega(V_x^{M_x})) & \xrightarrow{\simeq} & \mathrm{Aut}^\otimes(V_x^{M_x}) & \xrightarrow{\simeq} & \mathrm{Aut}^\otimes(D_x^{M_x}) \\
\cap & & \cap & & \downarrow \simeq & & \downarrow \simeq \\
G(\mathbb{Q}_p) & \xrightarrow{\simeq} & \mathrm{Aut}^\otimes(\omega(V_x)) & \supset & \mathrm{Aut}^\otimes(V_x) & \xrightarrow{\simeq} & \mathrm{Aut}^\otimes(D_x).
\end{array}$$

In particular:

**Theorem 5.15.** *The reduction map*

$$\mathrm{red} \; : \; \mathcal{L}(V_x) \to \mathcal{L}(D_x, \mathbf{t}_H^\iota(D_x))$$

*factors through an $M_x(\mathbb{Q}_p)$-equivariant bijection*

$$U_{\mathcal{F}_F}(\mathbb{Q}_p) \backslash \mathcal{L}(V_x) \simeq \mathcal{L}(D_x, \mathbf{t}_H^\iota(D_x)).$$

# References

[1] *Schémas en groupes (SGA 3). Tome III. Structure des schémas en groupes réductifs*, Documents Mathématiques (Paris) [Mathematical Documents (Paris)], 8, Société Mathématique de France, Paris, 2011. Séminaire de Géométrie Algébrique du Bois Marie 1962–64. [Algebraic Geometry Seminar of Bois Marie 1962–64].

[2] Y. ANDRÉ, *Slope filtrations*, Confluentes Math., 1 (2009), pp. 1–85.

[3] M. AUSLANDER AND D. A. BUCHSBAUM, *Homological dimension in Noetherian rings*, Proc. Nat. Acad. Sci. U.S.A., 42 (1956), pp. 36–38.

[4] B. BHATT, M. MORROW, AND P. SCHOLZE, *Integral p-adic hodge theory*, (2016).

[5] O. BRINON AND B. CONRAD, *CMI summer school notes on p-adic Hodge theory*, 2009.

[6] M. BROSHI, *G-torsors over a Dedekind scheme*, J. Pure Appl. Algebra, 217 (2013), pp. 11–19.

[7] J.-L. COLLIOT-THÉLÈNE AND J.-J. SANSUC, *Fibrés quadratiques et composantes connexes réelles*, Math. Ann., 244 (1979), pp. 105–134.

[8] C. CORNUT, *Filtrations and buildings*, (2014).

[9] ———, *Mazur's inequality and Laffaille's theorem*, (2015).

[10] C. CORNUT AND M.-H. NICOLE, *Cristaux et immeubles*, Bull. Soc. Math. France, 144 (2016), pp. 125–143.

[11] J.-F. DAT, S. ORLIK, AND M. RAPOPORT, *Period domains over finite and p-adic fields*, vol. 183 of Cambridge Tracts in Mathematics, Cambridge University Press, Cambridge, 2010.

[12] P. DELIGNE, *Catégories tannakiennes*, in The Grothendieck Festschrift, Vol. II, vol. 87 of Progr. Math., Birkhäuser Boston, Boston, MA, 1990, pp. 111–195.

[13] P. DELIGNE, J. S. MILNE, A. OGUS, AND K.-Y. SHIH, *Hodge cycles, motives, and Shimura varieties*, vol. 900 of Lecture Notes in Mathematics, Springer-Verlag, Berlin-New York, 1982.

[14] C. W. ERICKSON AND B. LEVIN, *A Harder-Narasimhan theory for Kisin modules*, 2015.

[15] L. FARGUES, *Théorie de la réduction pour les groupes p-divisibles*. Preprint.

[16] L. FARGUES, *La filtration de Harder-Narasimhan des schémas en groupes finis et plats*, J. Reine Angew. Math., 645 (2010), pp. 1–39.

[17] J.-M. FONTAINE, *Représentations p-adiques des corps locaux. I*, in The Grothendieck Festschrift, Vol. II, vol. 87 of Progr. Math., Birkhäuser Boston, Boston, MA, 1990, pp. 249–309.

[18] Q. R. GASHI, *On a conjecture of Kottwitz and Rapoport*, Ann. Sci. Éc. Norm. Supér. (4), 43 (2010), pp. 1017–1038.

[19] A. GENESTIER AND V. LAFFORGUE, *Structures de Hodge-Pink pour les $\varphi/\mathfrak{S}$-modules de Breuil et Kisin*, Compos. Math., 148 (2012), pp. 751–789.

[20] N. M. KATZ, *Slope filtration of F-crystals*, in Journées de Géométrie Algébrique de Rennes (Rennes, 1978), Vol. I, vol. 63 of Astérisque, Soc. Math. France, Paris, 1979, pp. 113–163.

[21] M. KISIN, *Mod p points on shimura varieties of abelian type*.

[22] M. KISIN, *Crystalline representations and F-crystals*, in Algebraic geometry and number theory, vol. 253 of Progr. Math., Birkhäuser Boston, Boston, MA, 2006, pp. 459–496.

[23] ———, *Integral canonical models of Shimura varieties*, J. Théor. Nombres Bordeaux, 21 (2009), pp. 301–312.

[24] ———, *Integral models for Shimura varieties of abelian type*, J. Amer. Math. Soc., 23 (2010), pp. 967–1012.

[25] M. Kisin and G. Pappas, *Integral models of shimura varieties with parahoric level structure*, (2015).

[26] R. E. Kottwitz, *Isocrystals with additional structure*, Compositio Math., 56 (1985), pp. 201–220.

[27] ———, *Shimura varieties and λ-adic representations*, in Automorphic forms, Shimura varieties, and *L*-functions, Vol. I (Ann Arbor, MI, 1988), vol. 10 of Perspect. Math., Academic Press, Boston, MA, 1990, pp. 161–209.

[28] ———, *Isocrystals with additional structure. II*, Compositio Math., 109 (1997), pp. 255–339.

[29] ———, *On the Hodge-Newton decomposition for split groups*, Int. Math. Res. Not., (2003), pp. 1433–1447.

[30] R. P. Langlands, *Some contemporary problems with origins in the Jugendtraum*, in Mathematical developments arising from Hilbert problems (Proc. Sympos. Pure Math., Vol. XXVIII, Northern Illinois Univ., De Kalb, Ill., 1974), Amer. Math. Soc., Providence, R. I., 1976, pp. 401–418.

[31] R. P. Langlands and M. Rapoport, *Shimuravarietäten und Gerben*, J. Reine Angew. Math., 378 (1987), pp. 113–220.

[32] T. Liu, *Compatibility of kisin modules for different uniformizers*, (2013).

[33] B. Mazur, *Frobenius and the Hodge filtration*, Bull. Amer. Math. Soc., 78 (1972), pp. 653–667.

[34] J. Milne, *Points on some Shimura varieties over finite fields: the conjecture of langlands and rapoport*, (2007).

[35] B. Moonen, *Serre-Tate theory for moduli spaces of PEL type*, Ann. Sci. École Norm. Sup. (4), 37 (2004), pp. 223–269.

[36] J. Neukirch, A. Schmidt, and K. Wingberg, *Cohomology of number fields*, vol. 323 of Grundlehren der Mathematischen Wissenschaften [Fundamental Principles of Mathematical Sciences], Springer-Verlag, Berlin, 2000.

[37] M. Rapoport and M. Richartz, *On the classification and specialization of F-isocrystals with additional structure*, Compositio Math., 103 (1996), pp. 153–181.

[38] J.-P. Serre, *Abelian l-adic representations and elliptic curves*, McGill University lecture notes written with the collaboration of Willem Kuyk and John Labute, W. A. Benjamin, Inc., New York-Amsterdam, 1968.

[39] ———, *Groupes algébriques associés aux modules de Hodge-Tate*, in Journées de Géométrie Algébrique de Rennes. (Rennes, 1978), Vol. III, vol. 65 of Astérisque, Soc. Math. France, Paris, 1979, pp. 155–188.

[40] C. S. Seshadri, *Geometric reductivity over arbitrary base*, Advances in Math., 26 (1977), pp. 225–274.

[41] X. Shen, *On the Hodge-Newton filtration for p-divisible groups with additional structures*, Int. Math. Res. Not. IMRN, (2014), pp. 3582–3631.

[42] R. Steinberg, *Regular elements of semisimple algebraic groups*, Inst. Hautes Études Sci. Publ. Math., (1965), pp. 49–80.

[43] J.-P. Wintenberger, *Torseur entre cohomologie étale p-adique et cohomologie cristalline; le cas abélien*, Duke Math. J., 62 (1991), pp. 511–526.

Unité de mathématiques pures et appliquées de l'ENS Lyon, ENS Lyon Site Monod, 46 Allée d'Italie, 69364 Lyon Cedex 07, France
*E-mail address:* `macarena.peche@ens-lyon.fr`